\documentclass[10pt,letterpaper]{article}
\usepackage{color}
\usepackage{multicol}
\usepackage[letterpaper, margin=1in]{geometry}
\usepackage{subfig}
\usepackage[overload]{empheq}

\usepackage{xcolor}

\usepackage{listings}
\usepackage{amsmath, amsthm, amssymb, wasysym, verbatim, bbm, graphics}
\usepackage{enumerate}
\usepackage{url}
\usepackage{mathtools}
\usepackage[colorlinks,linkcolor=blue]{hyperref}
\usepackage{xfrac}
\usepackage{float}
\usepackage{amsmath,bm}
\usepackage[title]{appendix}
\usepackage{todonotes}
\usepackage{subcaption}
\usepackage{multirow}
\usepackage{caption}
\usepackage{cleveref}
\usepackage{enumitem}
\usepackage{booktabs}
\usepackage{graphicx}
\usepackage{algorithm}
\usepackage{algpseudocode}
\usepackage[utf8]{inputenc}      
\DeclareUnicodeCharacter{0300}{\`{}}  

\graphicspath{{figures/}}

\DeclareMathOperator*{\argmin}{arg\,min}

\renewcommand\midrule{\specialrule{1pt}{0pt}{0pt}}
\renewcommand\toprule{\specialrule{1pt}{0pt}{0pt}}
\renewcommand\bottomrule{\specialrule{1pt}{0pt}{0pt}}

\usepackage[backend=biber,maxnames=4]{biblatex}
\newtheorem{lemma}{Lemma}[section]

\newtheorem{example}[lemma]{Example}

\newtheorem*{maintheorem*}{Main Theorem}
\theoremstyle{definition}{}

\allowdisplaybreaks

\numberwithin{equation}{section}

\begin{document}

\title{What metric to optimize for suppressing instability in a Vlasov-Poisson system?}

\author{Martin Guerra\thanks{Department of Mathematics, University of Wisconsin-Madison, Madison, WI 53706 USA} \and Qin Li\footnotemark[1] \and Yukun Yue\footnotemark[1] \and Leonardo Zepeda-N\'u\~nez\thanks{Google Research, Mountain View, CA 94043 USA}}

\date{}
\renewcommand*\contentsname{Outline}

\maketitle


\begin{abstract}
Stabilizing plasma dynamics is a central challenge in magnetic confinement fusion. A common approach is to introduce external electric fields to suppress instabilities in the plasma distribution. However, efficiently identifying such stabilizing fields remains challenging, even for simplified kinetic models such as the Vlasov–Poisson (VP) system.

In this work we study plasma stabilization from the perspective of PDE-constrained optimization. Our goal is to understand how the choice of objective function and the underlying kinetic dynamics influence the optimization landscape. First, we analyze the dispersion relation of the VP system and show that it reveals the spectral structure of the dynamics; eliminating unstable modes provides parameter configurations that lie close to the global optimum and serve as effective initial guesses for optimization. Second, we investigate several objective functions for stabilization and compare their optimization landscapes through numerical experiments.

Our results show that while different objectives lead to similar stabilizing parameter configurations, objective functions incorporating time-integrated information exhibit more convex-like landscapes and are therefore more favorable for gradient-based optimization methods. These findings provide insight into the design of objective functions for optimization-based plasma control and suggest promising directions for future research on real-time stabilization of kinetic plasma models.
\end{abstract}

\section{Introduction}
Plasma, an ionized gas composed of electrons, ions, and neutral particles, is often referred to as the fourth state of matter~\cite{burm2012plasma,fitzpatrick2022plasma,frank2012plasma}. One of its most compelling applications is the generation of clean energy through nuclear fusion~\cite{ichimaru1993nuclear,miyamoto1980plasma}. A principal obstacle to realizing practical fusion energy is the inherent instability of plasma dynamics: in many configurations, small perturbations to an equilibrium state can grow exponentially, disrupting the system and rendering a fusion device inoperable. Consequently, developing robust strategies to suppress these instabilities is a critical area of research.

These instabilities are typically mitigated using control strategies, wherein an external actuator—such as a magnetic or electric field—is applied to the plasma to counteract the growth of disruptive modes.
Control strategies can be broadly classified by their temporal nature. \textit{Dynamic} control involves real-time feedback, where the actuator responds continuously to the evolving state of the plasma. This paradigm has been the subject of extensive research, employing both traditional optimal control techniques~\cite{albi2025instantaneous,albi2025robustfeedbackcontrolcollisional,Blum2016,Einkemmer2025} and modern reinforcement learning (RL) methods~\cite{nature2022,nature2024}. However, dynamic approaches face a fundamental limitation: a \textit{timescale mismatch} between the rapid growth of certain instabilities and the slower response time of physical actuators~\cite{Moret1988}. In response to this challenge, \textit{static} control has recently regained significant interest from the community. This paradigm involves designing a time-independent control to dampen a small initial perturbation, with the goal of stabilizing the system at an ulterior time $T$~\cite{Einkemmer2024,Knopf2018,knopf2020optimal}. This work focuses on the design of such static controls.

A common approach for designing a static control is to impose an external electric field, $H(x)$, to guide the plasma towards a desired stable distribution~\cite{Einkemmer2024}. In this work we adopt a one-dimensional periodic configuration in space. This task can be formulated as a PDE-constrained optimization problem:
\begin{equation}\label{eq:optimization_pb_1}
\begin{aligned}
\argmin_{H} & \quad \mathcal{J}(f[H]) \\
\text{s.t.} & \quad \left\{\begin{array}{l}
   \partial_{t}f + v\,\partial_{x}f - (E_{f}+H)\,\partial_{v}f = 0\,,\\
   E_{f} = \partial_{x}V_{f}\,, \\
   \partial_{xx} V_{f} = 1 - \rho_{f} = 1 - \displaystyle\int f\,\mathrm{d}v\,,
\end{array}\right.
\end{aligned}
\end{equation}
with periodic boundary conditions in $x$. Here, the constraint is the VP system in 1D phase space. The function $f:(t,x,v)\in\mathbb{R}_{+}\times[0,L_{x}]\times\mathbb{R}\rightarrow\mathbb{R}_{+}$ is the time-dependent particle distribution. The system's evolution is driven by the total electric field, which is a sum of the externally applied control field $H(x)$ and the self-generated electric field $E_{f}(t,x)$. This self-generated field is determined self-consistently: the charge density $\rho_{f}(t,x) = \int f(t,x,v)\,\mathrm{d}v$ determines the electric potential $V_f(t,x)$ via the Poisson equation, and $E_f$ is determined by its gradient. The objective function $\mathcal{J}$ is a functional designed to quantify the degree of instability in the plasma distribution $f$.

This PDE-constrained optimization framework is a cornerstone of modern control theory and has been employed in numerous studies of plasma stability, both theoretical~\cite{Einkemmer2025,Glass2003,Glass2012,Knopf2018,Knopf2019} and numerical~\cite{Einkemmer2024,albi2025instantaneous,kawata2019dynamic,Bartsch2024,Bialek2001,Blum2016,Strait2014}. Regardless of the specific physical system, the underlying optimization problem is invariably challenging. The nonlinear nature of the governing PDEs typically leads to highly non-convex objective landscapes, populated with numerous local minima. As such, the careful design of the objective function $\mathcal{J}$ is critical. Different choices for $\mathcal{J}$ emphasize different features of instability~\cite{holm1985nonlinear, mouhot2011Landau}; for instance, one metric might penalize the kinetic energy of unstable wave modes, while another might measure the deviation of the entire  distribution in phase space from its equilibrium state. This choice directly shapes the optimization landscape and determines which properties of the solution are promoted. This leads to the central question that this manuscript seeks to answer:
\begin{center}
    \emph{What is a good metric for suppressing instabilities in the Vlasov-Poisson system?}
\end{center}

In the context of fusion energy, the ultimate goal is to maintain the plasma's energy in kinetic form, which implies minimizing the energy stored in the self-generated electric field~\cite{Freidberg2007, miyamoto1980plasma}. A common approach in the literature is to target this goal indirectly by minimizing the $L^2$-misfit between the final state of the distribution, $f_T := f(T,\cdot,\cdot)$, and the desired equilibrium, $f_{\text{eq}}$~\cite{Knopf2018,knopf2020optimal}. The rationale behind this approach lies in the fact that if $f_T$ perfectly matches the spatially uniform equilibrium $f_\text{eq}$, the self-generated electric energy must be zero. However, as demonstrated numerically in~\cite{Einkemmer2024}, this $L^2$ objective function is highly non-convex and oscillatory, creating a landscape with a myriad of local minima that can trap optimization algorithms.

This motivates our search for alternative objective functions, which we explore along two main axes: the choice of discrepancy measure and its temporal structure. Along the first axis, we compare two types of measures:
\begin{itemize}
    \item To quantify the deviation from equilibrium, we consider the \textit{Kullback--Leibler (KL) divergence}. Since $f$ represents a probability distribution, the KL divergence provides a natural information-theoretic measure of disparity compared with, for example, the $L^2$ norm. Moreover, the KL divergence corresponds to relative entropy and therefore carries a clear physical interpretation.
    
    \item To more directly target the underlying physics, we also consider minimizing the \textit{self-generated electric energy}. This quantity is macroscopic and inexpensive to evaluate, making it convenient for numerical optimization.
\end{itemize}

Along the second axis, we investigate the temporal structure of the objective function. In particular, we compare objectives evaluated only at the final time $T$ with objectives that are \textit{integrated over the entire time interval}. This allows us to examine whether incorporating information from the full system dynamics can regularize the optimization problem and mitigate the non-convexity of the resulting landscape.

Our numerical experiments suggest that objective functions incorporating time-integrated information tend to produce more convex-like optimization landscapes. Nevertheless, when the search domain is large and the initial guess is far from the optimum, non-convexity remains unavoidable, making it difficult for gradient-based methods to reliably locate the global minimizer.

Fortunately, for the VP system, the underlying PDE structure provides additional guidance. By analyzing the spectrum of the linearized system around equilibrium, one can identify the fastest-growing unstable modes through Fourier-Laplace analysis~\cite{Einkemmer2025}. This analysis allows us to derive control fields that specifically suppress these unstable modes. Although this study is performed in the linearized regime, it provides an effective guideline for selecting control parameters that stabilize the nonlinear optimization problem~\eqref{eq:optimization_pb_1}.

In summary, this paper presents an integrated analytical--numerical framework for suppressing instabilities in the VP system. Our numerical experiments systematically explore different objective functions, optimization solvers, and initialization strategies, leading to three main findings:
\begin{itemize}
\item First, objective functions that incorporate time-integrated information produce smoother and more convex-like optimization landscapes near the global minimum than objectives evaluated only at a terminal time.
\item Second, macroscopic metrics such as the electric energy, while physically intuitive, may introduce spurious local minima that can trap optimization algorithms, particularly line-search-based methods~\cite{NocedalNumOpt}.
\item Third, linear dispersion-relation analysis provides a reliable strategy for constructing an effective initial guess for the control parameters, placing the optimization within the basin of attraction of the global minimum.
\end{itemize}
Together, these observations provide a practical framework for efficiently identifying external control fields that suppress plasma instabilities.

\subsection{Related work} \label{sec:related_work}
Due to its significance, plasma control has attracted considerable research interest, leading to the development of a wide range of control strategies. As one might expect, these strategies are fundamentally shaped by the mathematical models used to describe plasma behavior. These models span a spectrum from macroscopic fluid descriptions—such as magnetohydrodynamics (MHD)\cite{Bialek2001,Blum2016,Strait2014}—to microscopic kinetic models like the VP or Vlasov-Maxwell system. Owing to the high mathematical and computational complexity of kinetic models, many approaches rely on simplifications such as dimensionality reduction\cite{Einkemmer2024} or linearization~\cite{mouhot2011Landau,HanKwan2021,Einkemmer2025} to make the resulting control problems more tractable. These approximations allow the use of optimal control techniques within available computational resources.

Once a model is chosen, another key consideration is the nature of the control strategy itself. Broadly speaking, plasma control can be categorized into static and dynamic control. In static control—such as the approach pursued in this work—the system starts near equilibrium, with only a small perturbation that the control aims to dampen by a prescribed target time~\cite{Einkemmer2024,Knopf2018,knopf2020optimal}. By contrast, dynamic control involves real-time feedback mechanisms that respond continuously to the evolving state of the plasma. Theoretical developments in this direction can be found in~\cite{Einkemmer2025}, while recent experimental implementations have emerged under the optimal control paradigm~\cite{albi2025instantaneous,albi2025robustfeedbackcontrolcollisional,Blum2016} and the reinforcement learning (RL) paradigm~\cite{nature2022,nature2024}. These RL-based methods recast the control problem in a data-driven framework, training neural networks to predict instabilities~\cite{Seo2023} and develop reactive strategies. While RL offers the promise of real-time, adaptive control of nonlinear plasma dynamics, it also faces limitations—particularly the mismatch between actuation and measurement timescales (on the order of 25 ms~\cite{Moret1988}) and the much faster evolution of certain instabilities, which may occur on sub-millisecond or even nanosecond timescales.

Regardless of the control type, implementation in experimental or computational settings typically involves actuating either a magnetic field or an electric field. Magnetic field control is the dominant paradigm in magnetic confinement fusion experiments, particularly for managing MHD-scale behavior in tokamaks~\cite{Bialek2001}. It has also been extended to kinetic models, where theoretical and numerical studies have investigated magnetic control within the VP framework to confine plasma away from reactor walls and prevent material damage~\cite{albi2025instantaneous, albi2025robustfeedbackcontrolcollisional, Knopf2018, knopf2020optimal}. In contrast, electric field control—although less common in large-scale experimental devices—has proven powerful in theoretical and computational studies, particularly for kinetic models. Applications include beam shaping in particle accelerators~\cite{Einkemmer2025} and, most relevant to this work, the suppression of kinetic instabilities via PDE-constrained optimization~\cite{Einkemmer2024}. Consistent with experimental and theoretical practice, we model the external electric field \(H\) as a prescribed actuator input superposed on the plasma's self-consistent field. Explicit realizations that justify this modeling choice include laboratory implementations~\cite{Gilson2004}, physics models with a given focusing field~\cite{Mitchell2019}, and mathematical analyses of VP with an added external force~\cite{Glass2012}. These examples align with the assumption adopted in this work. Moreover, in one-dimensional configurations, electric fields constitute the only physically viable external control mechanism.

Finally, we should stress that different mathematical definitions of ``instability" (objective functions) show drastically different landscapes for a given shared physics problem is not a new phenomenon, and is especially pronounced for wave-type problems, see~\cite{Engquist2022,Symes1991,ChenDing2023}. In the context of plasma control, though the result is expected, we have not found literature that thoroughly studied this phenomenon.

\subsection{Organization}
The paper is organized as follows. In Section~\ref{sec:prelim} we properly define the PDE-constrained optimization problem that we aim to solve, along the different possible objective functions that we will test. We also present the forward solver that we will use to solve the forward problem and the two canonical examples that we will present in this work. In Section~\ref{sec:dispersion} we briefly present the dispersion relation and a linear stability analysis which provides important information on the frequency and scaling that our control $H$ should have. In Section~\ref{sec:landscapes} we thoroughly study the different landscapes for different objectives functions. In Section~\ref{sec:numerical_exp} we make substantial numerical simulations with different objectives, methods, and initializations of the parameters.

\section{Preliminaries}\label{sec:prelim}
In fusion applications, plasmas are typically confined in toroidal or axisymmetric devices (tokamaks or stellarators). These geometries are well approximated by a circular configuration, for which a one-dimensional periodic model is standard. Accordingly, and as used in \eqref{eq:optimization_pb_1}, we work on $x\in[0,L_{x}]$ equipped with periodic boundary conditions. The dynamics are governed by the one-dimensional VP system:

\begin{equation}\label{eq:vlasov-poisoon_system_ext_1d}
\left\{\begin{array}{l}
   \partial_{t}f + v\,\partial_{x}f - (E_{f}+H)\,\partial_{v}f = 0\,,\\
   E_{f} = \partial_{x}V_{f}\,, \\
   \partial_{xx} V_{f} = 1 - \rho_{f} = 1 - \displaystyle\int f\,\mathrm{d}v\,.
\end{array}\right.
\end{equation}

It is not hard to see that in the Vlasov-Poisson system defined in~\eqref{eq:optimization_pb_1} (and~\eqref{eq:vlasov-poisoon_system_ext_1d}), any initial data that is spatially independent is an equilibrium state. In other words, if $f(0,x,v) = f_{0}(x,v) = f_{\text{eq}}(v)$, then $f(t,x,v) = f_{\text{eq}}(v)$ is a solution and, $\partial_{x}f = 0$ and $E_{f} = \partial_{x}V_{f} = 0$.

We also formally define the electric energy in time $\mathcal{E}_{f}(t)$ as
\begin{equation}\label{eq:electric_energy_def}
\mathcal{E}_{f}(t) = \int_{0}^{L_{x}}[E_{f}(t,x)]^{2}\,\mathrm{d}x\,.
\end{equation}

It is well known that equilibrium states which do not satisfy the Penrose condition~\cite{penrose1960electrostatic} are classified as unstable. In such cases, even a small perturbation to the equilibrium state can lead to exponential growth of the self-generated electric field $E_{f}$. In this work, we study two canonical examples of such unstable equilibria that violate the Penrose condition \cite{chen1984introduction,nicholson1983introduction}: the Two Stream equilibrium (Section~\ref{sec:TS_example}) and the Bump-on-Tail equilibrium (Section~\ref{sec:BoT_example}).

\subsection{PDE-constrained optimization problem}\label{sec:PDE_opt}

In general, it is difficult to suppress the self-generated electric field for every point in time, so, we attempt to suppress the instability of the Vlasov-Poisson system at a given time $T$ by introducing an external electric field $H(x)$. For this, we define $f[H] = f[H](t,x,v)$ as the solution of \eqref{eq:vlasov-poisoon_system_ext_1d} for a given external electric field $H$ and, our self-generated electric field for that system is now defined as $E_{f[H]}$. Therefore, we want to solve the following PDE-constrained optimization problem.

\begin{equation}\label{eq:optimization_pb}
\begin{aligned}
    \argmin_{H} & \quad \mathcal{J}(f[H]) \\
    \text{s.t.} & \, f[H] \text{ is a solution of }\eqref{eq:vlasov-poisoon_system_ext_1d}
\end{aligned}
\end{equation}
where our objective functional will be either of the following four\footnote{Another straightforward definition of objective functional is to use the $L^2$ difference between the equilibrium and the solution $\mathcal{J}(f[H]) = \frac{1}{2}\|f_{T}[H] - f_\text{eq}\|^{2}_{L^{2}(x,v)} = \frac{1}{2}\int_{-L_{v}}^{L_{v}}\int_{0}^{L_{x}} (f[H](T,x,v) - f_{\text{eq}}(v))^{2}\,\mathrm{d}x\,\mathrm{d}v\,,$. Numerically we find this objective functional presents very similar landscape as the KL divergence, as will be discussed in~\ref{sec:L2_KL}.}
\begin{subequations}
\begin{equation}\tag{KL}\label{eq:KL_obj}
\mathcal{J}_{1}(f[H]) = \text{KL}(f_{T}[H]||f_{\text{eq}}) = \int_{-L_{v}}^{L_{v}}\int_{0}^{L_{x}} f[H](T,x,v)\log\left(\frac{f[H](T,x,v)}{f_{\text{eq}}(v)}\right)\,\mathrm{d}x\,\mathrm{d}v\,,
\end{equation}
\begin{equation}\tag{EE}\label{eq:EE_obj}
    \mathcal{J}_{2}(f[H]) = \mathcal{E}_{f[H]}(T)  = \int_{0}^{L_{x}} [E_{f[H]}(T,x)]^{2}\,\mathrm{d}x\,,
\end{equation}
\begin{equation}\tag{KLT}\label{eq:KLT_obj}
\mathcal{J}_{3}(f[H]) = \int_{0}^{T} \text{KL}(f_{t}[H]||f_{\text{eq}})\,\mathrm{d}t = \int_{0}^{T}\int_{-L_{v}}^{L_{v}}\int_{0}^{L_{x}} f[H](t,x,v)\log\left(\frac{f[H](t,x,v)}{f_{\text{eq}}(v)}\right)\,\mathrm{d}x\,\mathrm{d}v\,\mathrm{d}t\,,
\end{equation}
\begin{equation}\tag{EET}\label{eq:EET_obj}
    \mathcal{J}_{4}(f[H]) = \int_{0}^{T} \mathcal{E}_{f[H]}(t)\mathrm{d}t  =\int_{0}^{T}\int_{0}^{L_{x}} [E_{f[H]}(t,x)]^{2}\,\mathrm{d}x\,\mathrm{dt}\,,
\end{equation}
\end{subequations}
where $E_{f[H]}$ and $\mathcal{E}_{f[H]}$ are as defined in~\eqref{eq:vlasov-poisoon_system_ext_1d} and~\eqref{eq:electric_energy_def} respectively and we set our computational domain to be $(x,v)\in[0,L_x]\times[-L_v,L_v]$ with $L_v$ big enough and the tail cut does not affect the numerical results. We note that these four objective functionals can be grouped into two categories: one based on electric energy and the other based on KL divergence. The former is a natural, physics-inspired choice, while the latter is a valid approach that has been considered in a recent work on control for the Vlasov-Poisson system \cite{borzi2025optimal}.

To simplify the computations we define the external field \(H(x)\) as a linear combination of cosine and sine basis functions:
\begin{equation}\label{eq:external_field}
\begin{aligned}
H(x;\boldsymbol{a},\boldsymbol{b}) & = \sum_{k=1}^{N} a_k \cos\left(\frac{2\pi k x}{L_{x}}\right) + b_k \sin\left(\frac{2\pi k x}{L_{x}}\right)
\\
& =\sum_{k=1}^{N} a_k \cos\left(kk_0 x\right) + b_k \sin\left(kk_0 x\right)\,,
\end{aligned}
\end{equation}
where $k_{0} = 2\pi/L_{x}$ is the reference frequency to ensure the periodicity of the field and, \(\boldsymbol{a}=(a_{1},...,a_{N})\) and \(\boldsymbol{b} = (b_{1},...,b_{N})\) represent the parameters to be determined. To simplify notation we will drop the dependence on $\boldsymbol{a}$ and $\boldsymbol{b}$ for $H$, and our problem becomes,
\begin{equation}\label{eq:optimization_pb_simple}
\begin{aligned}\argmin_{(\boldsymbol{a},\boldsymbol{b})} & \quad  \quad \quad \mathcal{J}(f[H]) \\
    \text{s.t.} & \, f[H] \text{ is a solution of }\eqref{eq:vlasov-poisoon_system_ext_1d} \\
    & \, (\boldsymbol{a},\boldsymbol{b}) \in \mathbb{R}^{N}\times \mathbb{R}^{N}\,.
\end{aligned}
\end{equation}

\subsection{Solving the forward problem}\label{sec:numerical_scheme}

In order to solve the optimization problem~\eqref{eq:optimization_pb_simple}, we must repeatedly solve the VP system~\eqref{eq:vlasov-poisoon_system_ext_1d} in an efficient and stable manner. A variety of numerical methods have been developed for such kinetic equations, including Eulerian methods~\cite{filbet2003comparison}, Lagrangian (particle) methods~\cite{Verboncoeur2005}, and semi-Lagrangian methods~\cite{cheng1976integration,Rossmanith2011}. A general overview of these approaches can be found in~\cite{filbet2003numerical}.

In this work we employ a semi-Lagrangian discretization due to its favorable properties, including unconditional stability and the ability to take relatively large time steps~\cite{besse2003semi,cai2021high}. Since our focus is on the optimization framework rather than the development of new numerical schemes, we use standard off-the-shelf PDE solvers.

In a nutshell, semi-Lagrangian methods trace the characteristics exactly and perform an interpolation since the translated solution may not coincide with the grid. The literature proposes different ways to perform this interpolation, with spline based~\cite{cheng1976integration,sonnendrucker1999semi,filbet2003comparison} and Fourier-based methods~\cite{Klimas1994,Klimas2018} being the most famous. In our case, we deploy a simple linear interpolation.

We now describe the semi-Lagrangian method used in our implementation. The key idea is that each transport term in the Vlasov equation can be solved analytically along its characteristic curves, which allows the scheme to avoid the usual CFL stability restriction. 

To this end, we split the Vlasov equation into two subproblems,
\begin{subequations}
\begin{equation}\label{eq:x_splitting}
\partial_{t}f + v\partial_{x}f = 0,
\end{equation}
\begin{equation}\label{eq:v_splitting}
\partial_{t}f - (E_{f} + H)\partial_{v}f = 0,
\end{equation}
\end{subequations}
representing the spatial transport and force terms, respectively. 
Since the coefficients $v$ and $E_f+H$ are explicitly known within each stage, these equations can be solved analytically along characteristics~\cite{sonnendrucker1999semi,einkemmer2014convergence}.

For the numerical implementation, we discretize time with step size $\Delta t$ and denote $t^n = n\Delta t$. 
The phase space $(x,v)$ is discretized on an equispaced grid with $M_x$ and $M_v$ points in the spatial and velocity directions, respectively. 
For simplicity we assume $M_x = M_v = M$ and denote the grid by $\{(x_i,v_j)\}_{i,j=0}^{M-1}$. 
Let $\mathsf{f}$ denote the numerical approximation of the distribution function. 
Applying Strang splitting, one time step from $t^n$ to $t^{n+1}$ consists of the following three stages.

\begin{itemize}

\item \textbf{Step 1: half-step advection in space.}  
The spatial transport equation is solved using a semi-Lagrangian update
\[
\mathsf{f}^{\ast}(x_i,v_j) 
= \mathsf{f}\!\left(t^n,\,x_i - v_j\frac{\Delta t}{2},\,v_j\right),
\qquad i,j = 0,\dots,M-1 .
\]

\item \textbf{Step 2: velocity update.}  

First we compute the self-consistent electric field $E_{\mathsf{f}^{\ast}}$. 
The electric potential $V_{\mathsf{f}^{\ast}}$ satisfies
\begin{equation}\label{eq:self_gen_eef}
E_{\mathsf{f}^{\ast}} = \partial_x V_{\mathsf{f}^{\ast}}\,,\quad\text{with}\quad \partial_{xx} V_{\mathsf{f}^{\ast}} = 1 - \rho_{\mathsf{f}^{\ast}}\,,
\end{equation}
which is solved using a pseudo-spectral method. 
In Fourier space this yields
\[
\widehat{E}_{\mathsf{f}^{\ast}}(k) = \frac{i}{k}\widehat{\rho}_{\mathsf{f}^{\ast}}(k), 
\qquad k \neq 0,
\]
while $\widehat{E}_{\mathsf{f}^{\ast}}(k=0)=0$ to enforce charge neutrality. The Fourier transform is coded using FFT.

Once the electric field is obtained, the velocity transport equation is advanced for a full time step,
\[
\mathsf{f}^{\ast\ast}(x_i,v_j) =
\mathsf{f}^{\ast}\!\left(x_i,\, v_j + (E_{\mathsf{f}^{\ast}} + H)\Delta t\right),
\qquad i,j = 0,\dots,M-1 .
\]

\item \textbf{Step 3: half-step advection in space.}  
Finally, the remaining half-step of spatial transport is applied,
\[
\mathsf{f}(t^{n+1},x_i,v_j) =
\mathsf{f}^{\ast\ast}\!\left(x_i - v_j\frac{\Delta t}{2},\,v_j\right),
\qquad i,j = 0,\dots,M-1 .
\]

\end{itemize}

The complete semi-Lagrangian scheme is summarized in Algorithm~\ref{alg:semi_lagrangian}.\footnote{\url{https://github.com/maguerrap/Vlasov-Poisson}}.

\begin{algorithm}
\caption{Semi-Lagrangian with interpolation}\label{alg:semi_lagrangian}
\begin{algorithmic}[1]
\Require $\{x_{i},v_{j}\}_{i,\,j=0}^{M-1}$, $f_{\text{iv}}$, $\Delta t$, $N$, $f_{\text{eq}}$, $L_{x}$, $L_{v}$, $H$.
\Ensure $\boldsymbol{f}_{n}:= \{\mathsf{f}(t^{n},x_{i},v_{j})\}_{i,\,j=0}^{M-1}$ for $n=0,...,N$.
\For{$n=0,1,...,N$}
    \State Compute $\mathsf{f}^{\ast}(x_{i},v_{j}) = \mathsf{f}(t^{n},x_{i}-v_{j}\Delta t/2,v_{j})$, for $i,j=0,1,...,M-1$.
    \State  Compute $E_{\mathsf{f}^{\ast}}$ by solving~\eqref{eq:self_gen_eef}.
    \State Compute $\mathsf{f}^{\ast\ast}(x_{i},v_{j} + (E_{\mathsf{f}^{\ast}}+H)\Delta t)$, for $i,j=0,1,...,M-1$.
    \State Compute $\mathsf{f}(t^{n+1},x_{i},v_{j}) = \mathsf{f}^{\ast\ast}(x_{i}-v_{j}\Delta t/2,v_{j})$, for $i,j=0,1,...,M-1$.
\EndFor
\end{algorithmic}
\end{algorithm}

\subsection{Two canonical examples}
Plasma system, or the system modeled by the VP system, can have many different kinds of instabilities. We only study two most prominent ones: Two Stream instability and Bump-on-Tail instability. Both instabilities violate the Penrose condition~\cite{penrose1960electrostatic} and are regarded as two canonical examples.
\subsubsection{Two Stream example}\label{sec:TS_example}
As a first example, we consider the Two Stream equilibrium. This distribution arises when two streams (or beams) of charged particles moving relative to each other are interpenetrating, causing instability \cite{anderson2001}. Mathematically, the equilibrium distribution is given by
\begin{equation}\label{eq:equilibrium_two_stream}
f_{\text{eq}}(v) = \frac{\alpha \exp\left(-\frac{1}{2}(v - \mu)^2\right) + (1 - \alpha) \exp\left(-\frac{1}{2}(v + \mu)^2\right)}{\sqrt{2\pi}}\,.
\end{equation}
It is a distribution composed of two Gaussians having opposite centers. In our simulation we set $\alpha = 0.5$ and $\mu = 2.4$.

The initial data is set to be a small perturbation from this equilibrium:
\begin{equation}\label{eq:perturbation:example1}
f_{\text{iv}}(x, v) = \left(1 + \varepsilon \cos(k_0 x)\right) f_{\text{eq}}(v) \doteq f_{\text{eq}}(v)+f_{\text{p}}(t=0,x,v)\,,
\end{equation}
with the initial perturbation $f_{\text{p}}(t=0,x,v)$ being of order $\varepsilon=0.001$.

In Figure~\ref{fig:TS_no_H} we plot our simulation~\eqref{eq:vlasov-poisoon_system_ext_1d} using Algorithm~\ref{alg:semi_lagrangian} without an external electric field (i.e., $H\equiv 0$). We simulate until $T=30$ -- a time one can observe strong instability starting to appear and electric field has shown exponential growth. The phase-space domain for this simulation is the $L_{x}\times L_{v} = 10\pi \times 6.0$ box.

\begin{figure}[ht]
    \centering
    \includegraphics[width=1.0\linewidth]{TS_noH.png}
    \caption{Simulation of~\eqref{eq:vlasov-poisoon_system_ext_1d} with $H\equiv 0$ for the Two Stream equilibrium. From left to right we have $f_{\text{eq}}(v)$, $f(T=30,x,v)$, $E_{f}(t,x)$ and, $\mathcal{E}_{f}(t)$ and KL over time.}
    \label{fig:TS_no_H}
\end{figure}

\subsubsection{Bump-on-Tail example}\label{sec:BoT_example}
The Bump-on-Tail equilibrium is another steady state that is heavily investigated. This distribution appears when a small group of fast-moving electrons (a ``bump'') is superimposed on the tail of the background electron velocity distribution and this can cause instability \cite{o1971nonlinear}. This can happen due to particle injection or wave-particle interactions. The equilibrium is defined as
\begin{equation}\label{eq:equilibrium_bump}
f_{\text{eq}}(v) = \frac{9}{10\sqrt{2\pi}} \exp\left(-\frac{1}{2}(v - \bar{v}_1)^2\right) + \frac{\sqrt{2}}{10\sqrt{\pi}} \exp\left(-2(v - \bar{v}_2)^2\right)\,.
\end{equation}
with the equilibrium parameters set to
\begin{align*}
\bar{v}_1 &= -3.0, \quad \bar{v}_2 = 4.5.
\end{align*}

This equilibrium is also not stable. We add an initial perturbation to it as the initial data:
\begin{equation}\label{eq:perturbation_bump}
f_{\text{iv}}(x, v) = f_{\text{eq}}(v) + \underbrace{\frac{\sqrt{2}\,\varepsilon}{10\sqrt{\pi}} \exp\left(-2(v - \bar{v}_2)^2\right) \cos(0.1 x)}_{\doteq f_{\text{p}}(t=0,x,v)}\,,
\end{equation}
with the initial perturbation $f_{\text{p}}(t=0,x,v)$ being of order $\varepsilon = 0.001$, the electric energy will grow exponentially over time. In Figure~\ref{fig:BoT_no_H} we simulate~\eqref{eq:vlasov-poisoon_system_ext_1d} using Algorithm~\ref{alg:semi_lagrangian} without an external electric field (i.e., $H\equiv 0$). The solution demonstrates strong electric energy growth at about $T=40$ when we terminate the simulation. The simulation is conducted over a $L_x\times L_v= 20\pi\times 9.0$ box.

\begin{figure}[ht]
    \centering
    \includegraphics[width=1.0\linewidth]{BoT_noH.png}
    \caption{Simulation of~\eqref{eq:vlasov-poisoon_system_ext_1d} with $H\equiv 0$ for the Bump-on-Tail equilibrium. From left to right we have $f_{\text{eq}}(v)$, $f(T=40,x,v)$, $E_{f}(t,x)$ and, $\mathcal{E}_{f}(t)$ and KL over time.}
    \label{fig:BoT_no_H}
\end{figure}

\section{Dispersion relation and linear stability analysis}\label{sec:dispersion}
Examining the dispersion relation is a classical approach to investigating the stability properties of dynamical systems \cite{mouhot2011Landau}. Mathematically, this involves linearizing the system around a desired equilibrium state and analyzing the spectrum of the resulting linear PDE operator. If the operator has spectrum with positive real parts, the system exhibits exponential divergence from equilibrium. To counteract this, we introduce a control term $H$ aimed at modifying the spectrum—specifically, to eliminate its unstable (positive) components and thereby keep the system close to the desired state.

This strategy was proposed and studied in depth in~\cite{Einkemmer2025}, where the dispersion relation was derived both with and without an external field, and a general procedure was outlined to design a field capable of suppressing growth. However, the method targets complete suppression by eliminating all unstable roots, which is a stringent requirement. Achieving this demands highly flexible external fields—indeed, the designs in~\cite{Einkemmer2025} result in fields that oscillate significantly in time. In practice, actuator bandwidth and shot-duration constraints limit how fast fields can be driven: for example, on TCV the plasma position responds to poloidal-field coil currents on a $\sim\!10\,\mathrm{ms}$ timescale~\cite[Section 5.2]{hommen2014real}, and during ITER first-plasma operations feedback is intentionally limited. Short pulse durations and loop latencies constrain the achievable control bandwidth, so only relatively slow corrections are feasible~\cite[Section 2.3]{snipes2021iter}. Accordingly, strongly time-dependent fields are operationally impractical on present devices. This motivates our use of a time-independent external field as specified in~\eqref{eq:external_field}.

In contrast, we seek to impose a far less flexible external field: one that contains only finitely many Fourier modes and is time-independent. This constraint necessitates a relaxation of the original goal. Rather than eliminating all unstable roots, we shift our focus to suppressing the modes that grow the fastest. Our approach begins by identifying the specific modes that would exhibit the strongest exponential growth in the absence of control, then targeting their associated roots for elimination. This analysis builds on our previous derivation of the dispersion relation, now adapted to account for time-independent control.

Ultimately, the system remains nonlinear, so previously unexcited modes may be triggered over time. As a result, the present strategy is effective only in the short term. For long-term suppression, a return to a full PDE-constrained optimization framework is necessary. Nonetheless, our analysis offers a good initial guess for such optimization-based methods.

In the following section, we begin by reviewing the dispersion relation derived in~\cite{Einkemmer2025}, as presented in subsection~\ref{subsec:dispersion_Htx}. We then extend the analysis to the case of a time-independent external field, and propose our strategy in subsection~\ref{subsec:dispersion_Hx}. It provides a good initial guess for further optimization.

\subsection{Dispersion relation and control strategy when \texorpdfstring{$H=H(t,x)$}{H=H(t,x)}}\label{subsec:dispersion_Htx}
We first briefly review the dispersion relation for the Vlasov-Poisson system subject to an external field \(H(t,x)\) that varies in both time and space. The original derivation is complete in~\cite{Einkemmer2025}.

To do so, we first linearize the system around the desired equilibrium state, writing $f=f_{\text{eq}}+f_{\text{p}}$, then omitting the higher order terms in the expansion, the leading order solution to $f_{\text{p}}$ is:
\begin{equation}\label{eq:VP_system_pert_linear}
\left\{
    \begin{aligned}
        &\partial_t f_{\text{p}}(t,x,v) + v \partial_x f_{\text{p}}(t,x,v) - \left[E_{f_{\text{p}}}(t,x) + H(t,x)\right]\partial_v f_{\text{eq}}(v) = 0\,, \\
        &E_{f_{\text{p}}}(t,x) = \partial_{x} V_{f_{\text{p}}}(t,x)\,,\\
       &\partial_{xx}V_{f_{\text{p}}}(t,x) = -\rho_{\text{p}}(t,x) = -\int f_{\text{p}}(t,x,v)\,\mathrm{d}v\,.
    \end{aligned}
\right.
\end{equation}
Three crucial transforms will be heavily leaned on:
\begin{itemize}
    \item[--] Fourier transform in space, e.g.
   \begin{equation}\label{eq:Fourier-transform-space}
         \hat{F}(t,k) = \int_{-\infty}^{\infty} F(t,x)e^{-ikx}\,\mathrm{d}x\,.
    \end{equation}
    \item[--] Fourier transform in space and velocity, e.g:
    \begin{equation}\label{eq:Fourier-transform-space-velocity}
        \hat{F}(t,k,m) = \int_{-\infty}^{\infty} \int_{-\infty}^{\infty}F(t,x,v)e^{-ikx}e^{-imv}\,\mathrm{d}x\,\mathrm{d}v\,.
    \end{equation}
    \item[--] Laplace transform in time:
    \begin{equation}
    \label{eq:Laplace_transform_def}
    L[F](s) = \int_0^\infty e^{-st} F(t)\,\mathrm{d}t\,.
\end{equation}
\end{itemize}
Repeatedly deploying these transforms, one can deduce an important identity:
\begin{equation}
     \label{eq:external_Laplace}
     \left(L[\hat{\rho}_{\text{p}}(\cdot,k)](s)-L[\hat{S}(\cdot,k)](s)\right)\left( 1+L[\hat{U}(\cdot,k)](s) \right) = L[\hat{U}(\cdot,k)](s)\left(-L[\hat{S}(\cdot,k)](s) + ik L[\hat{H}(\cdot,k)](s)  \right)\,,
\end{equation}
where
\begin{itemize}
\item \( \hat{U}(t,k) \) is defined as:
\begin{equation}\label{eq:hatU_def}
 \hat{U}(t,k) :=   t\hat{f}_{\text{eq}}(kt)\,.
\end{equation}
It is a function of \( t \) and \( k \), and is uniquely determined by the equilibrium \(f_{\text{eq}}\).
\item \( \hat{S} \) is defined as
\begin{equation}\label{eq:hatS_def}
    \hat{S}(t,k) := \hat{f}_{\text{p}}(0,k, kt) \,.
\end{equation} 
Or equivalently $S(t,x) = \int f_{\text{p}}(t,x-vt,v)\,\mathrm{d}v$. It is the free-streaming solution.
\end{itemize}
According to the inversion Laplace transform formula, $\rho_{\text{p}}$ will have exponential growth if $1+L[\hat{U}(\cdot,k)](s)=0$ for some $s$ that has positive real part (i.e. $\Re(s)>0$)~\cite{Einkemmer2025}, when the source term on the right hand side of~\eqref{eq:external_Laplace} does not vanish at the same $s$. As a consequence, to eliminate the exponential growth, one should simply force the right hand side to be zero at the roots of $1+L[\hat{U}(\cdot,k)](s)$. In light of this discussion, in~\cite{Einkemmer2025}, the authors propose to set the control $H$ according to:
\begin{equation}
    \label{eq:external_field_choice}
    -L[\hat{S}(\cdot,k)](s) + ik L[\hat{H}(\cdot,k)](s) = h\left( 1 + L[\hat{U}(\cdot,k)](s) \right),
\end{equation}
where $h$ is a factor so that:
\begin{equation}\label{eqn:h_cond}
    h : \mathbb{C} \to \mathbb{C}\quad\text{s.t.}\quad h(0) = 0\quad \text{and}\quad \lim\limits_{|x| \to 0} \frac{|h(\cdot)|}{|\cdot|} = c \ \ \text{for some} \ c\in\mathbb{C}\,.
\end{equation}
This design effectively eliminates the positive roots, and thus the PDE's instability is suppressed. Further analytical derivations and comprehensive numerical validations are provided in~\cite{Einkemmer2025}.

\subsection{Dispersion relation and control strategy when \texorpdfstring{$H=H(x)$}{H=H(x)}}\label{subsec:dispersion_Hx}
It is evident that the control $H$ designed above has both time and space dependence. However, the $H$ we are seeking for only has spatial dependence, rendering these analysis results not immediately useful. In this section, we modify the analysis above and re-derive the dispersion relation with $H$ assumed to be time-independent.

Similar to the derivation above, we apply the Fourier-Laplace transform and find:

\begin{equation}\label{eq:external_Laplace_Hx}
\bigl(L[\hat{\rho}_{\mathrm{p}}(\cdot,k)](s) - L[\hat{S}(\cdot,k)](s)\bigr)\,
\bigl(1 + L[\hat{U}(\cdot,k)](s)\bigr)
= L[\hat{U}(\cdot,k)](s)\,
\bigl(-L[\hat{S}(\cdot,k)](s) + \tfrac{i k\,\hat{H}(k)}{s}\bigr).
\end{equation}
In the derivation, we used \(L[1](s)=\frac{1}{s}\) and the fact that $H$ does not have time dependence.

Following the strategy as in~\eqref{eq:external_field_choice}, we should set
\[
-\,L[\hat S(\cdot,k)](s) \;+\;\frac{i k\,\hat H(k)}{s} \;=\; 0
\]
for all complex frequencies \(s\). The equation is pointwise in $(k,s)$ while the function to be sought for only varies in $k$, rendering no solution that satisfies the entire equation. However, as discussed, in this regime when we have a far less flexible external field, we ought to relax our original goal. In this particular setting, we are to find the fastest growing mode and aim at suppressing that mode only.

Mathematically, to find the fastest growing mode, we set 
\begin{equation}\label{eqn:s0}
s_0(k) = \text{argmax}\{\Re(s)>0: \; 1 + L\bigl[\hat U(\cdot,k)\bigr](s) = 0\}
\end{equation}
for every fixed $k$. This is to identify the Laplace index who are roots for the operator and has the largest real component. We then require $\hat{H}$ only to eliminate that specific component: $-\,L\bigl[\hat S(\cdot,k)\bigr](s_0(k))\;+\;\frac{i k\,\hat H(k)}{s_0}=0$. This yields
\begin{equation}\label{eq:Hx_choice}
\hat H(k)
\;=\;
\frac{s_0(k)\,L\bigl[\hat S(\cdot,k)\bigr](s_0(k))}{i k}\,.
\end{equation}

We summarize our strategy as follows. For a given equilibrium distribution \(U\) and the initial perturbation \(S\), we:
\begin{enumerate}
  \item compute~\eqref{eqn:s0} for $s_0(k)$ for all \(k\) seen in the initial perturbation;
  \item compute the spectral amplitude $\hat{H}(k)$ using \eqref{eq:Hx_choice}.
\end{enumerate}
Since $\hat{S}$ is determined, according to~\eqref{eq:hatS_def}, by the initial perturbation. By definition~\eqref{eq:perturbation:example1} and~\eqref{eq:perturbation_bump}, it is of order $\mathcal{O}(\varepsilon)$, hence the guessed control $H$ should also be of $\mathcal{O}(\varepsilon)$.

\begin{figure}[ht]
  \centering
    \includegraphics[width=0.49\linewidth]{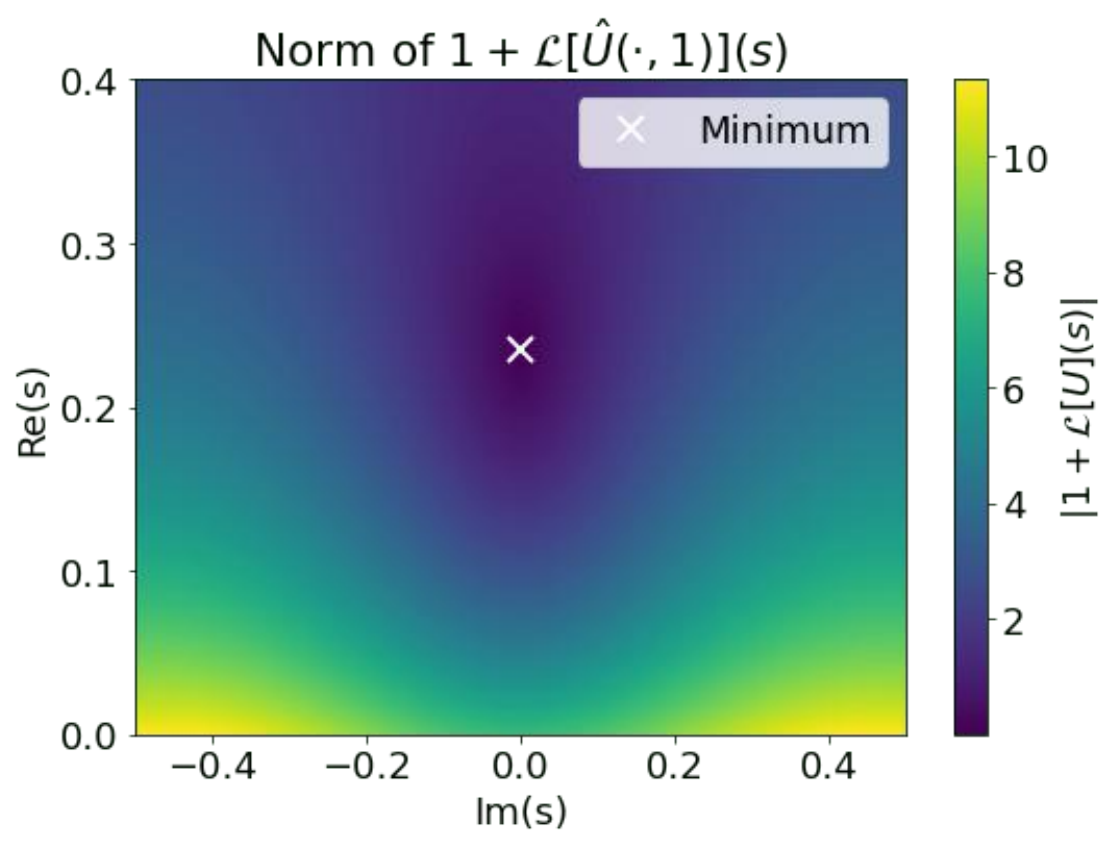}
    \includegraphics[width=0.49\linewidth]{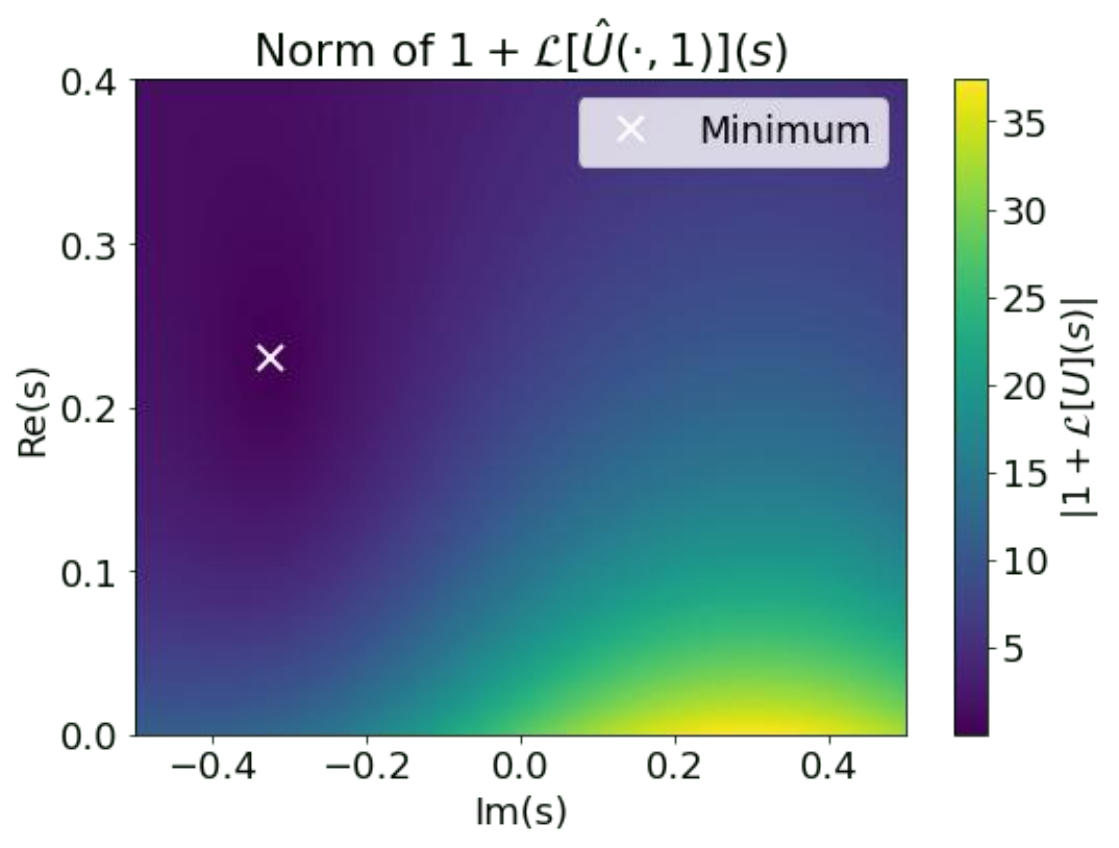}
  \caption{
Norm of $\|1 + \mathcal{L}[U](\cdot,1)(s)\|$: the Two Stream (left) and the Bump-on-Tail (right) examples. In Two Stream, the minimum occurs at $s_0 = 0.236 + 0i$ with $\|1 + \mathcal{L}[U](\cdot,1)(s_0)\| = 9.564044\times10^{-4}$. In Bump-on-Tail, the minimum occurs at $s_0 = 0.230 - 0.324i$ with $\|1 + \mathcal{L}[U](\cdot,1)(s_0)\| = 3.756278\times10^{-3}$. 
} 
  \label{fig:(1+LU)_laplace}
\end{figure}

\begin{example}\label{ex:initial_guess}
    For the two benchmark cases (Two-Stream and Bump-on-Tail), the initial perturbations excite only the modes $k=\pm1$. We therefore determine $s_0(\pm1)$ as the minimizers of $\big|1+\mathcal{L}[U](\cdot,1)(s)\big)$ (see Fig.~\ref{fig:(1+LU)_laplace}) and obtain $\widehat{H}(\pm1)$ from~\eqref{eq:Hx_choice}. The resulting fields $H(x)$ serve as the initial guesses for our optimization. The corresponding values of $s_0$ are reported alongside the plots in Fig.~\ref{fig:(1+LU)_laplace}.
\end{example}

These constructions yield a strong initial guess for the external field $H$. When applied to equation~\eqref{eq:vlasov-poisoon_system_ext_1d}, we observe significant suppression effects, as illustrated in Figure~\ref{fig:control_vs_nocontrol}. However, it is important to note that the electric energy eventually exhibits exponential growth over long times. This reflects the onset of nonlinear dynamics, which cannot be fully addressed by our initial design. As such, complete suppression requires the use of a PDE-constrained optimization framework. The field $H$ constructed here serves as an effective initial guess for such optimization, which we will discuss in detail in Section~\ref{sec:numerical_exp}. 

\begin{figure}[ht]
    \centering
    \includegraphics[width=1.0\linewidth]{TS_init_guess_30.png}
    \includegraphics[width=1.0\linewidth]{BoT_init_guess_40.png}
    \caption{Simulation of~\eqref{eq:vlasov-poisoon_system_ext_1d} for the Two Stream equilibrium up to $T=30$ (top) and Bump-on-Tail equilibrium up to $T=40$ (bottom). In each row, from left to right: $f_{\text{eq}}(v)$, $f(T,x,v)$ (for $H\equiv 0$), $f(T,x,v)$ (for $H$ obtained from Example~\ref{ex:initial_guess}) and $\mathcal{E}_{f}(t)$ (for $H\equiv 0$ and $H$ obtained from Example~\ref{ex:initial_guess}).}
    \label{fig:control_vs_nocontrol}
\end{figure}

\section{Landscape analysis of the objective function}\label{sec:landscapes}

In this section, we examine the landscape of the objective function. This is to evaluate the dependence of the four objective function (defined in~\eqref{eq:KL_obj},~\eqref{eq:EE_obj},~\eqref{eq:KLT_obj} and~\eqref{eq:EET_obj}) on parameters of the control (defined in~\eqref{eq:external_field}). Landscape analysis requires a sweep of PDE solves over all parameter choices in a pre-defined domain, and thus is computationally demanding.

In Section~\ref{sec:TS_landscape} and Section~\ref{sec:BoT_landscape} respectively, we provide the landscape for choices of two different parameter settings for the two different canonical examples (Sections \ref{sec:TS_example} and \ref{sec:BoT_example}) respectively.

\subsection{Two Stream example}\label{sec:TS_landscape}
We report results for the Two Stream instability problem in this section. Here we have the following two choices of parameter settings:

\begin{itemize}
    \item \textbf{Choice A:} 
    \begin{equation}\label{eq:H_sin}
        H(x) = b_{1}\sin(k_{0}x)\,,
    \end{equation}
    \item \textbf{Choice B:} 
    \begin{equation}\label{eq:H_sins}
        H(x) = b_{1}\sin (k_0 x) + b_{2}\sin(2k_0 x)\,.
    \end{equation}
\end{itemize}

\begin{figure}[ht]
    \centering
    \includegraphics[width=1.0\linewidth]{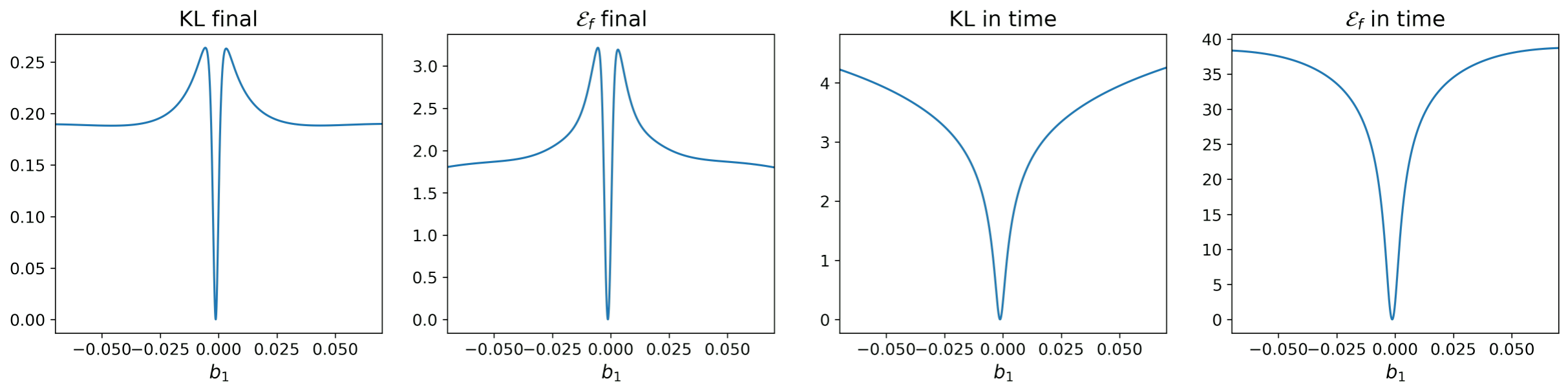}
    \caption{Landscape of the Two Stream instability on the domain \([-0.07,0.07]\) for \(b_{1}\) for~\eqref{eq:H_sin}. The objective functions are~\eqref{eq:KL_obj}(left),~\eqref{eq:EE_obj}(center-left),~\eqref{eq:KLT_obj}(center-right) and~\eqref{eq:EET_obj}(right).}
    \label{fig:landscape_two_stream_varepsilon_0.001_1D}
\end{figure}

For \textbf{Choice A}, with one parameter, we conduct a sweeping of PDE solutions by setting $b_{1}$ to be in the pre-defined range of $[-0.1,0.1]$, and compute the evaluation of the four objective functions (\eqref{eq:KL_obj},~\eqref{eq:EE_obj},~\eqref{eq:KLT_obj} and~\eqref{eq:EET_obj}). The result is plotted in Figure~\ref{fig:landscape_two_stream_varepsilon_0.001_1D}. Evidently, in this 1D situation, both~\eqref{eq:EE_obj} and~\eqref{eq:KL_obj} that only see the last frame of the solution give non-convex behavior in the neighborhood of the minimum, whereas~\eqref{eq:EET_obj} and~\eqref{eq:KLT_obj} that see the whole solution in time demonstrates a convex profile. This strongly suggests the information in-time is crucial in determining the convexity of the landscape of the objective functions.

\begin{figure}[ht]
    \centering
    \includegraphics[width=1.0\linewidth]{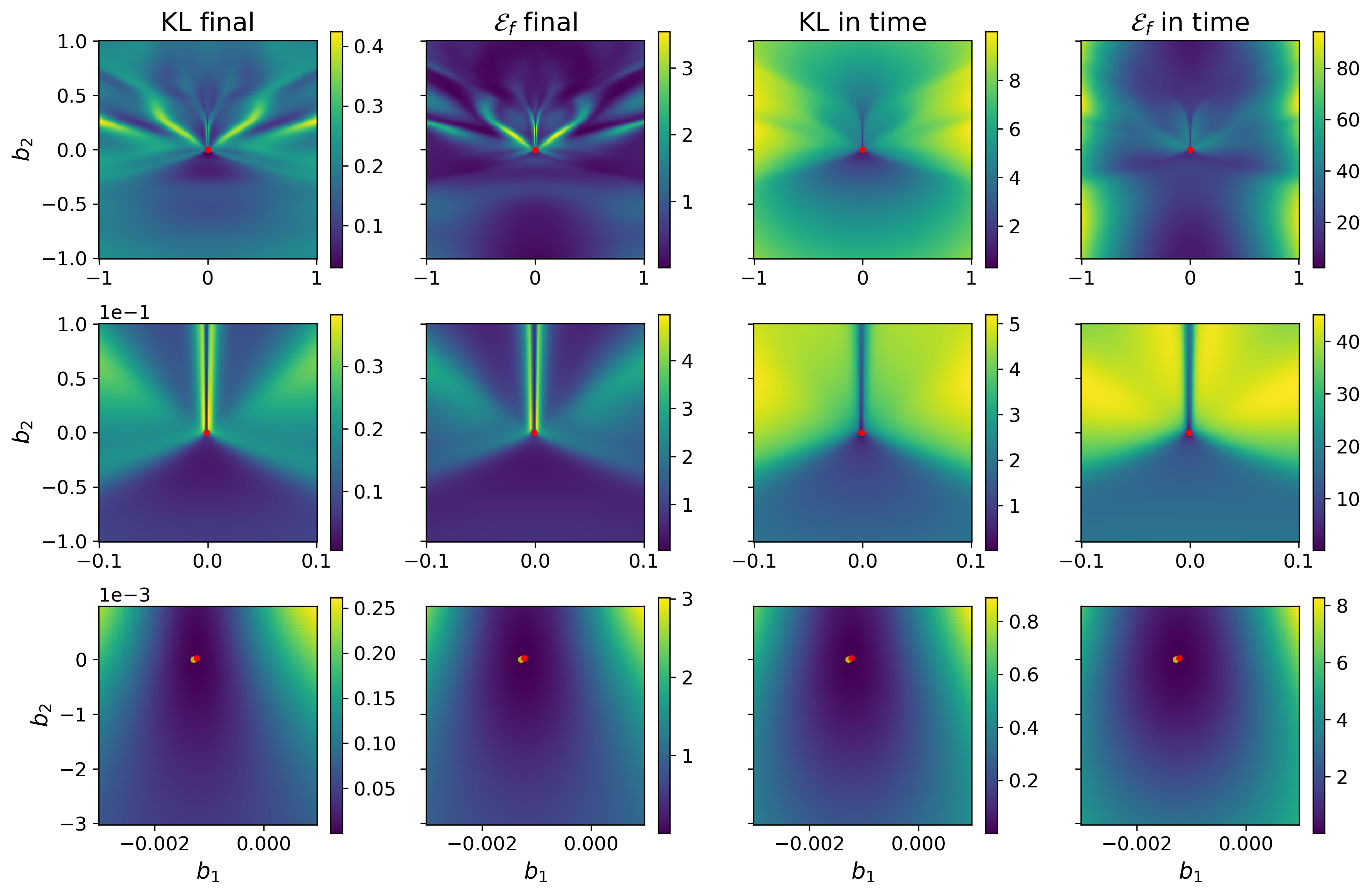}
    \caption{Landscape of the Two Stream instability with applied control over the domain \([-1,1]^{2}\)(top), \([-0.1,0.1]^{2}\)(center) and, \([-0.003,0.001]^{2}\)(bottom) for \((b_{1}, b_{2})\) for~\eqref{eq:H_sins}. The objective functions are~\eqref{eq:KL_obj}(left),~\eqref{eq:EE_obj}(center-left),~\eqref{eq:KLT_obj}(center-right) and~\eqref{eq:EET_obj}(right). Yellow dot in last row represents good initial guess from Example~\ref{ex:initial_guess} and red dot represents approximate global minimum.}
    \label{fig:landscape_two_stream_2D_varepsilon_0.001}
\end{figure}

Then, for \textbf{Choice B} with two parameters, we again perform a sweeping of PDE solutions by setting $(b_1,b_2)$ in three pre-defined domains:
\begin{itemize}
    \item[--]Far field: $[-1,1]\times[-1,1]$;
    \item[--]Mid-range: $[-0.1,0.1]\times[-0.1,0.1]$;
    \item[--]Near field: $[-0.003,0.001]\times[-0.003,0.001]$.
\end{itemize}
The wordings (``Far field, Mid-range, Near field") are selected based on the conclusions drawn from Section~\ref{sec:dispersion} of dispersion relation, where it was suggested the optimal $H$ should be in the range of $O(\varepsilon)$. In this case, considering~\eqref{eq:perturbation:example1}, $\varepsilon = 0.001$. In Figure~\ref{fig:landscape_two_stream_2D_varepsilon_0.001}, we plot the values for the four different objective functions in these three domains. The following conclusions are drawn:
\begin{itemize}
\item The landscape close to the global minimum (near-field) is convex for all four objective functions. To justify the convexity, we also compute Hessian at each point, and calculate the minimum eigenvalues. The values are shown in~\ref{sec:Hess_computations} Figure~\ref{fig:landscape_two_stream_2D_hess} and are evidently positive in the neighborhood of the global minimum.
\item In the mid-range set, both plots that only use last-frame information (\eqref{eq:KL_obj}-\eqref{eq:EE_obj}) present a pronounced bulge around the minimum that highlights a clear non-convex behavior along the first mode direction ($b_{1}$ variable), while the plots that integrate in time (\eqref{eq:KLT_obj}-\eqref{eq:EET_obj}) eliminate these bulge and present a convex objective landscape.
\item In the far field, both~\eqref{eq:KL_obj} and~\eqref{eq:EE_obj} are highly oscillatory, making the global minima hard to find. In comparison, the objective functions that integrates in time are much smoother. Comparing KL and electric energy, we also see that, in both final-frame and integrated-in-time case, objective functions defined by electric energy present minima for large values of $b_{2}$, suggesting large external electric field can also suppress instability. However, we should point out that this is a non-physical setting: The imposed electric field $H$ is large, and completely dominates the dynamics.
\end{itemize}

\subsection{Bump-on-Tail example}\label{sec:BoT_landscape}

In this section, we present results for the Bump-on-Tail instability problem with the following two choices of parameter settings:

\begin{itemize}
    \item \textbf{Choice C:} 
    \begin{equation}\label{eq:H_cos}
        H(x) = a_{1}\cos(k_{0}x)\,,
    \end{equation}
    \item \textbf{Choice D:} 
    \begin{equation}\label{eq:H_cos_sin}
        H(x) = a_{1}\cos (k_0 x) + b_{1}\sin(k_0 x)\,.
    \end{equation}
\end{itemize}

\begin{figure}[ht]
    \centering
    \includegraphics[width=1.0\linewidth]{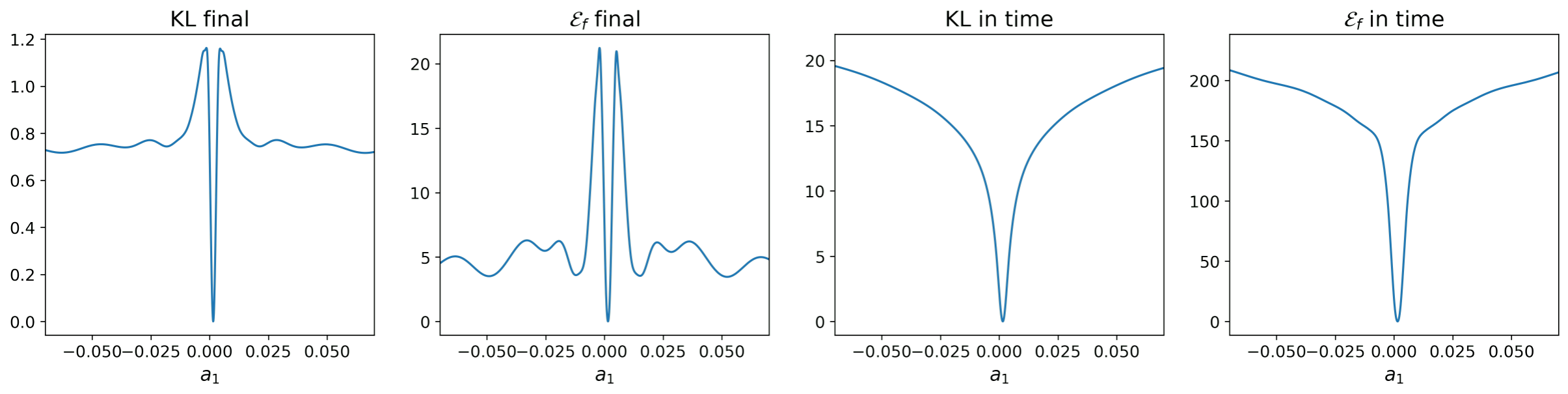}
    \caption{Landscape of the Bump-on-Tail instability on the domain \([-0.07,0.07]\) for \(a_{1}\) for~\eqref{eq:H_cos}. The objective functions are~\eqref{eq:KL_obj}(left),~\eqref{eq:EE_obj}(center-left),~\eqref{eq:KLT_obj}(center-right) and~\eqref{eq:EET_obj}(right).}
 \label{fig:landscape_bump-on-tail_varepsilon_0.001_1D}
\end{figure}

For \textbf{Choice C}, we perform a parameter sweep of the PDE solutions by varying $a_{1}$ over the pre-defined interval $[-0.1, 0.1]$. We then evaluate four different objective functions (\eqref{eq:KL_obj},~\eqref{eq:EE_obj},~\eqref{eq:KLT_obj} and~\eqref{eq:EET_obj}), in analogous to the Two Stream example. The results are displayed in Figure~\ref{fig:landscape_bump-on-tail_varepsilon_0.001_1D}. Notably, in this one-dimensional setting, both \eqref{eq:KL_obj} and~\eqref{eq:EE_obj} exhibit non-convex landscape, whereas~\eqref{eq:KLT_obj} and~\eqref{eq:EET_obj} display a convex profile, mirroring the behavior observed in the Two Stream example in Figure~\ref{fig:landscape_two_stream_varepsilon_0.001_1D}. This observation resonates with the situation for the Two Stream instability, and suggests information in-time is crucial in convexifying the objective functions.

For \textbf{Choice D}, we extend the analysis to a two-parameter sweep by varying the pair $(a_1,b_1)$ over three pre-defined domains ($[-1,1]^2$, $[-0.1,0.1]^2$, and $[-0.001,0.003]^2$). The near-field region is shifted due to the new location of the global minimum in this case.

\begin{figure}[ht]
    \centering
    \includegraphics[width=1.0\linewidth]{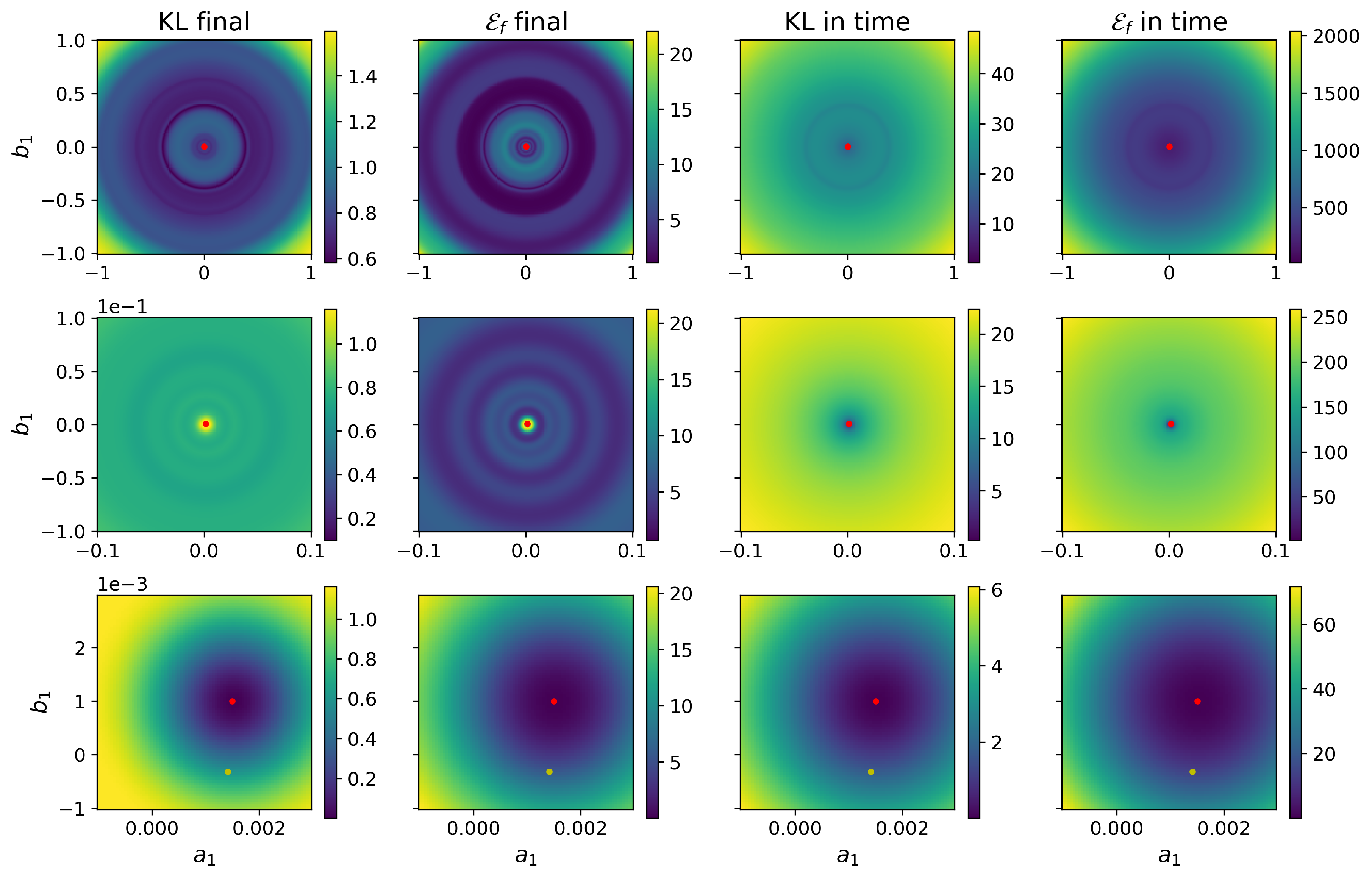}
    \caption{Landscape of the Bump-on-Tail instability with applied control over the domain \([-1,1]^{2}\)(top), \([-0.1,0.1]^{2}\)(center) and, \([-0.001,0.003]^{2}\)(bottom) for \((a_{1}, b_{1})\) for~\eqref{eq:H_cos_sin}. The objective functions are~\eqref{eq:KL_obj}(left),~\eqref{eq:EE_obj}(center-left),~\eqref{eq:KLT_obj}(center-right) and~\eqref{eq:EET_obj}(right). Yellow dot from last row represents good initial guess from Example~\ref{ex:initial_guess} and red dot represents approximate global minimum.}
    \label{fig:landscape_bump-on-tail_2D_varepsilon_0.001}
\end{figure}

Figure~\ref{fig:landscape_bump-on-tail_2D_varepsilon_0.001} displays the values of the four objective functions (\eqref{eq:KL_obj},~\eqref{eq:EE_obj},~\eqref{eq:EET_obj} and~\eqref{eq:KLT_obj}) across the three domains. The conclusions drawn from the Two Stream example continue to hold in the near-field and mid-range cases. In particular, to demonstrate the local convexity close to the global minimum, we plot in~\ref{sec:Hess_computations} Figure~\ref{fig:landscape_bump-on-tail_2D_hess} the minimum eigenvalue of the Hessians at each point and they are consistently above zero. In the far-field region,~\eqref{eq:EET_obj} and~\eqref{eq:KLT_obj} continue to be more convex than those of~\eqref{eq:KL_obj} and~\eqref{eq:EE_obj}. Moreover, they do not exhibit global minima for large values of $a_1$ and $b_1$ as seen in Figure~\ref{fig:landscape_two_stream_2D_varepsilon_0.001} for far-field region. This observation indicates that an excessively strong external electric field does not always produce a correspondingly small self-generated electric field.

\section{Solving the PDE-constrained optimization problem}\label{sec:numerical_exp}

In this section, we present the numerical results obtained from solving~\eqref{eq:optimization_pb_simple} using Algorithm~\ref{alg:semi_lagrangian} to address the forward system described by~\eqref{eq:vlasov-poisoon_system_ext_1d}. The computation was performed on an NVIDIA GeForce RTX 4090 GPU with 24 GB of memory. For the forward simulation, the phase space was discretized with $M = 256$ mesh points, and the time domain was discretized using a time step of $\Delta t = 0.1$. Each PDE solve can be completed within 5 seconds. Optimization was carried out using automatic differentiation, facilitated by the \textit{Python} library JAX~\cite{jax2018github}, to compute the required gradients. 

The plots shown in Section~\ref{sec:landscapes} suggest that the landscape for the objectives~\eqref{eq:KLT_obj} or~\eqref{eq:EET_obj} are rather similar. Therefore we expect the performance of optimization solvers with these objectives would be similar as well. For this reason, in the rest of this section, we perform numerical studies only on~\eqref{eq:KL_obj},~\eqref{eq:EE_obj} and~\eqref{eq:EET_obj}.

The landscape plots only provide some suggestions on the performance of optimization solvers. Indeed, for local solvers such as gradient descent, that can only see local behavior of objective functions, the solutions are trapped in local minima, and thus placing the initial guess in the global basin is crucial. On the other hand, solvers have their non-convex extensions. One such possibility is to integrate the line-search technique, which allows iterations to jump out of local minima. This can potentially relax the requirement for a good initial guess.

Another strategy typically deployed in optimization is ``over-parameterization''. This is to increase the number of unknowns and adjustable parameters, and thus formulating the problem in a higher-dimensional space. In this setting, it is very likely that the dimension of the global optima manifold increases, and so the chance of obtaining optimal points.

For these reasons, we numerically investigate the solvers' performance on two physics scenarios, Two Stream instability and Bump-on-Tail instability, respectively, in three distinct dimensions of variations:
\begin{itemize}
    \item Solvers can be either local or adaptive/global;
    \begin{itemize}
        \item Local solver. There are many kinds of local solvers. We adopt the simplest gradient descent with constant stepsize. The stepsize is small enough ($1e-8$, $1e-9$) that one can essentially view the optimization iteration as a gradient flow along the landscape. 
        \item Adaptive/global solver. There are also many choices of adaptive solvers. To be more comparable with the choice of local solver (GD), we select GD with line-search correction as our adaptive solver. The line-search will satisfy the strong Wolfe condition that we briefly describe below.
        
        Denoting the GD with stepsize $\alpha_n$ at $n$-th iteration for minimizing a function $F$:
        \[x_{n+1} = x_{n} - \alpha_{n}\nabla F(x_{n})\,.\]
        Line-search is imposed to ensure that the choice of $\alpha_{n}$ provides sufficient decrease. The ``sufficiency" is quantified by the following:
        \begin{align*}
            F(x_{n} - \alpha_{n}\nabla F(x_{n})) & \leq F(x_{n}) - c_{1}\alpha_{n}\|\nabla F(x_{n})\|^{2} \\
             |[\nabla F(x_{n} - \alpha_{n}\nabla F(x_{n}))]^{\top}\nabla F(x_{n})| & \leq c_{2}\|\nabla F(x_{n})\|^{2}
        \end{align*}
        for $0<c_{1}<c_{2}<1$. Namely, either the value of the new update is sufficiently smaller than the current value, or the new update has a gradient almost parallel to that of the old one. In particular, we use $c_{1} = 10^{-4}$ and $c_{2} = 0.9$.
    \end{itemize}
    \item Initial guess can be close, or far away from the optimum points. More specifically, the three kinds of initializations are
    \begin{itemize}
    \item Far initialization: $(\boldsymbol{a},\boldsymbol{b}) \in [-1,1]^{2N}$,
    \item Mid-range initialization: $(\boldsymbol{a},\boldsymbol{b}) \in [-0.05,0.05]^{2N}$,
    \item Near initialization:
    $(\boldsymbol{a},\boldsymbol{b}) \in [-0.003,0.001]^{2N}$ for Two Stream and $(\boldsymbol{a},\boldsymbol{b}) \in [-0.001,0.003]^{2N}$ for Bump-on-Tail.
\end{itemize}
The scalings are termed ``far", ``mid-range" and ``near" are closely tied to the landscapes we found on Section~\ref{sec:landscapes}.
    \item Optimization conducted with either $2$ (under or properly parameterized) or $14$ (over-parameterized) variables.

    According to the definition of the control term~\eqref{eq:external_field}, $H$ is composed of Fourier modes. The terminology ``under/properly" or ``over" have blurred scientific meanings. Roughly speaking, a system is ``over-parameterized" if one can find a configuration that already performs optimization well with a smaller set of parameters. We perform the study with
    \begin{itemize}
        \item Under/properly parameterized:
        \begin{equation}\label{eq:H_sins_under}
            H(x;b_{1},b_{2}) = b_{1}\sin (k_0 x) + b_{2}\sin(2k_0 x)\,, \quad \text{(Two Stream)}
        \end{equation}
        \begin{equation}\label{eq:H_cos_sin_under}
            H(x;a_{1},b_{1}) = a_{1}\cos (k_0 x) + b_{1}\sin(k_0 x)\,, \quad \text{(Bump-on-Tail)}
        \end{equation}
        \item Over-parameterized:
        \begin{equation}\label{eq:H_cos_sin_over}
            H(x;\boldsymbol{a},\boldsymbol{b}) = \sum_{k=1}^{14} \big(a_k \cos\left(kk_0 x\right) + b_k \sin\left(kk_0 x\right)\big)
        \end{equation}
        where $(\boldsymbol{a},\boldsymbol{b}) = ((a_{1},...,a_{14}),(b_{1},...,b_{14}))$.
    \end{itemize}
\end{itemize}

\subsection{Two Stream example}
We summarize results for controlling Two Stream instability, using the equilibrium state and initial condition presented in~\eqref{eq:equilibrium_two_stream}-\eqref{eq:perturbation:example1}. All numerical results are presented in Table~\ref{tab:TS_results} that shows the final distribution function $f_T[H]$ with $H$ found by the optimization solver. A more detailed summary of each experiment can be found in~\ref{sec:TS_summary}. The stepsize used for the constant stepsize gradient descent is $10^{-8}$.
\begin{table}[!ht]
    \centering
    \begin{tabular}{c|c|c|c|c|c|c|c}
    \toprule
    & & \multicolumn{3}{|c|}{Under-parametrized~\eqref{eq:H_sins_under}} & \multicolumn{3}{|c}{Over-parametrized~\eqref{eq:H_cos_sin_over}} \\
    \midrule
    Init. & Step & \eqref{eq:KL_obj} & \eqref{eq:EE_obj} & \eqref{eq:EET_obj} & \eqref{eq:KL_obj} & \eqref{eq:EE_obj} & \eqref{eq:EET_obj}  \\
    \midrule
    \multirow{4}{*}[0pt]{Far} & \multirow{2}{*}[22pt]{Adaptive} & \includegraphics[width=0.11\linewidth]{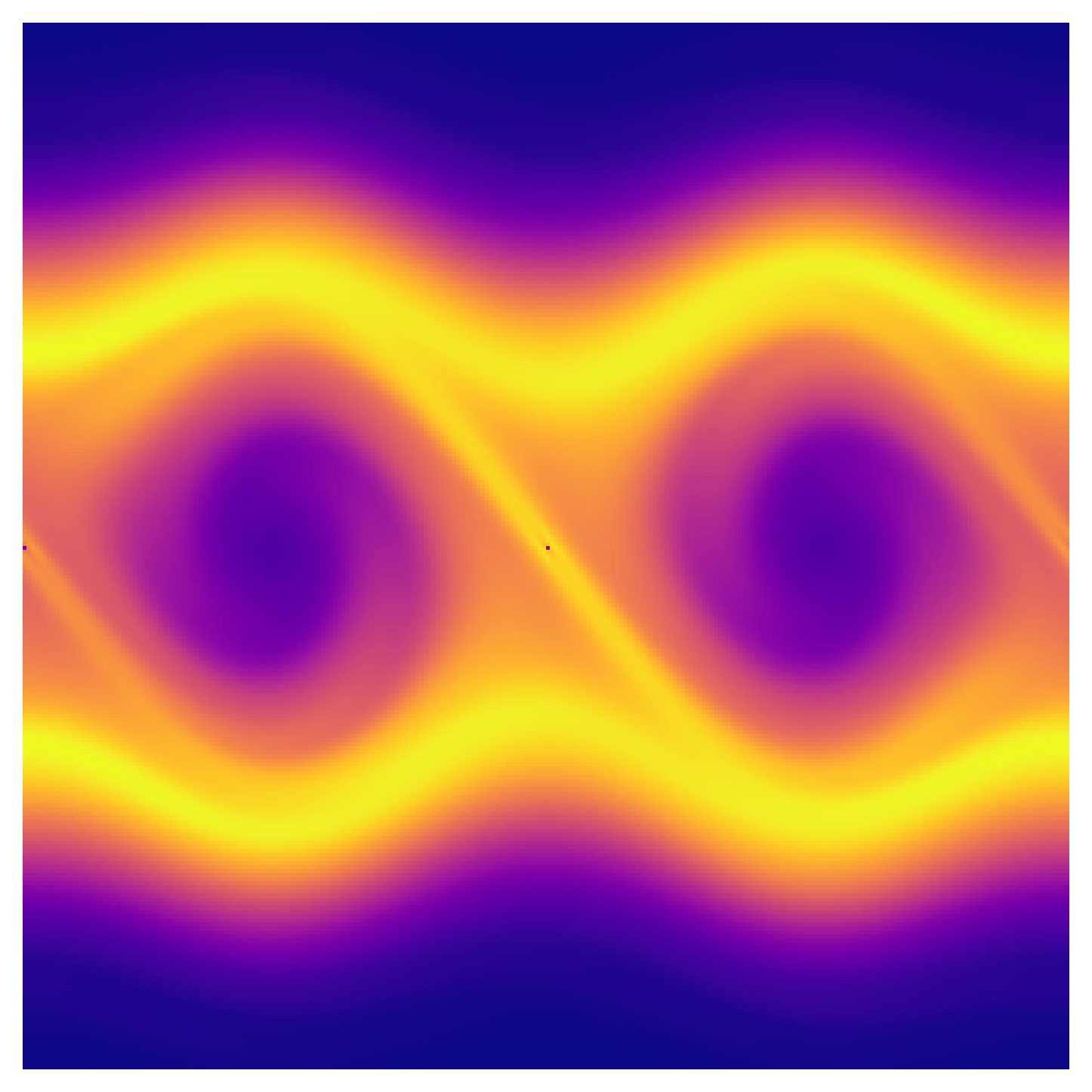} & \includegraphics[width=0.11\linewidth]{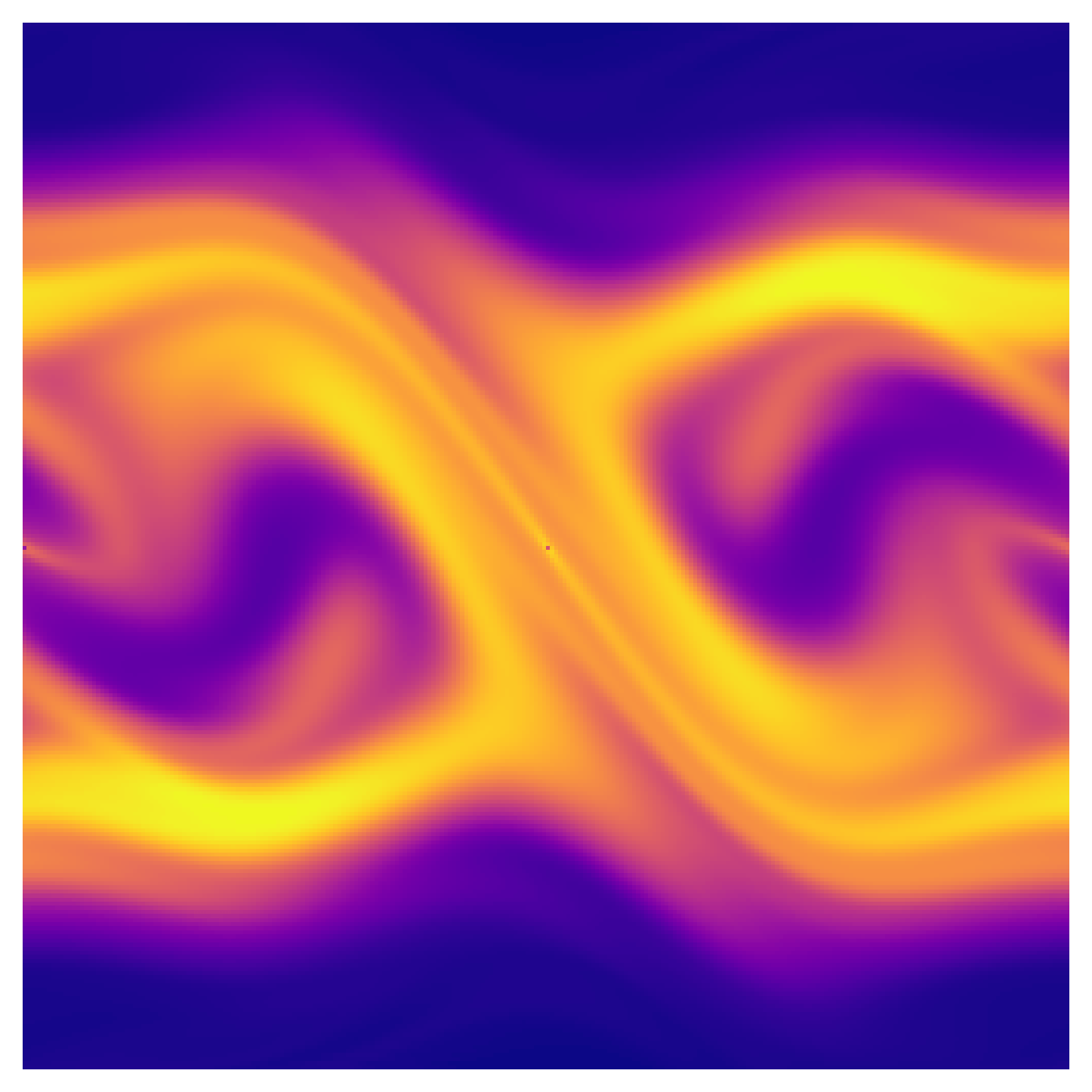} & \includegraphics[width=0.11\linewidth]{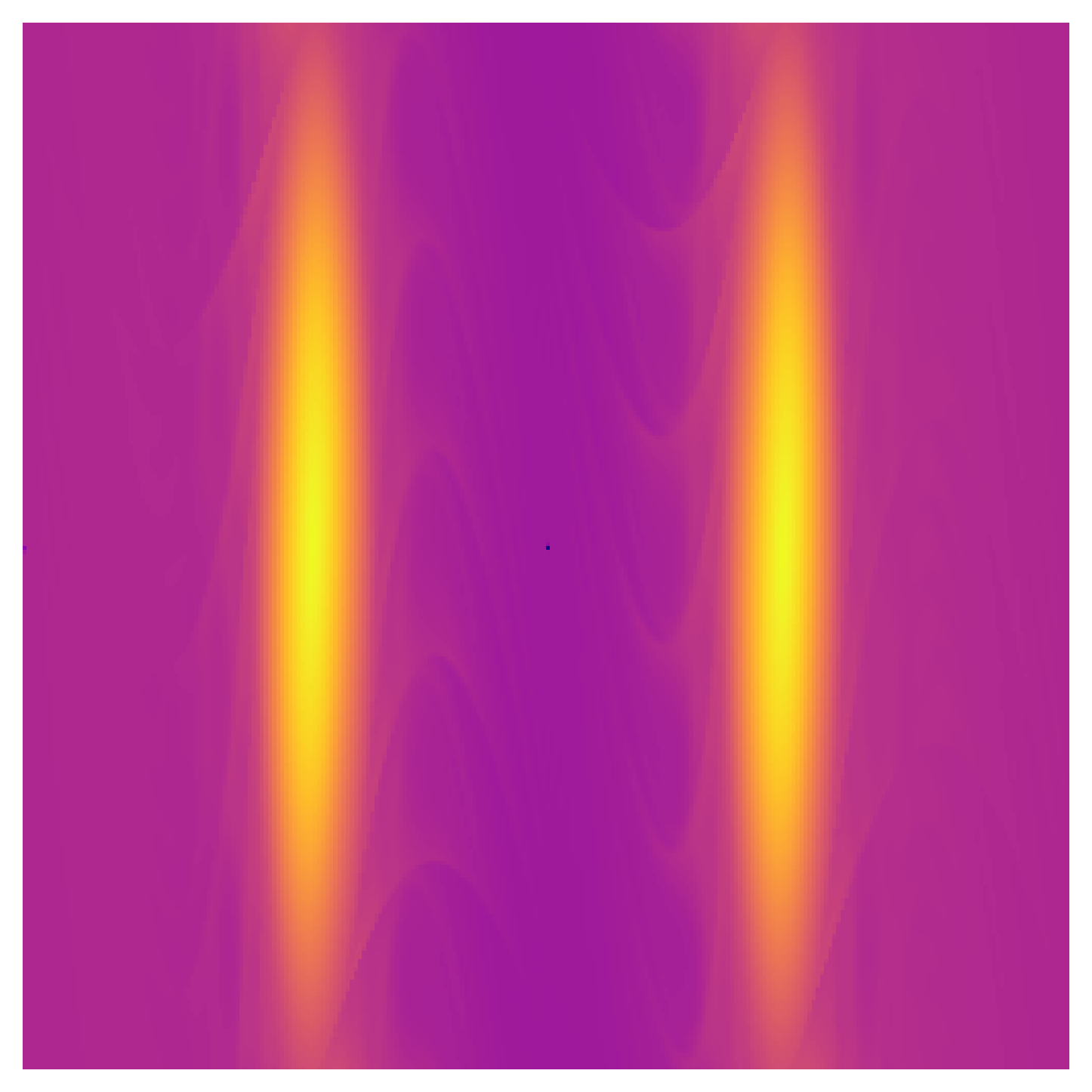} & \includegraphics[width=0.11\linewidth]{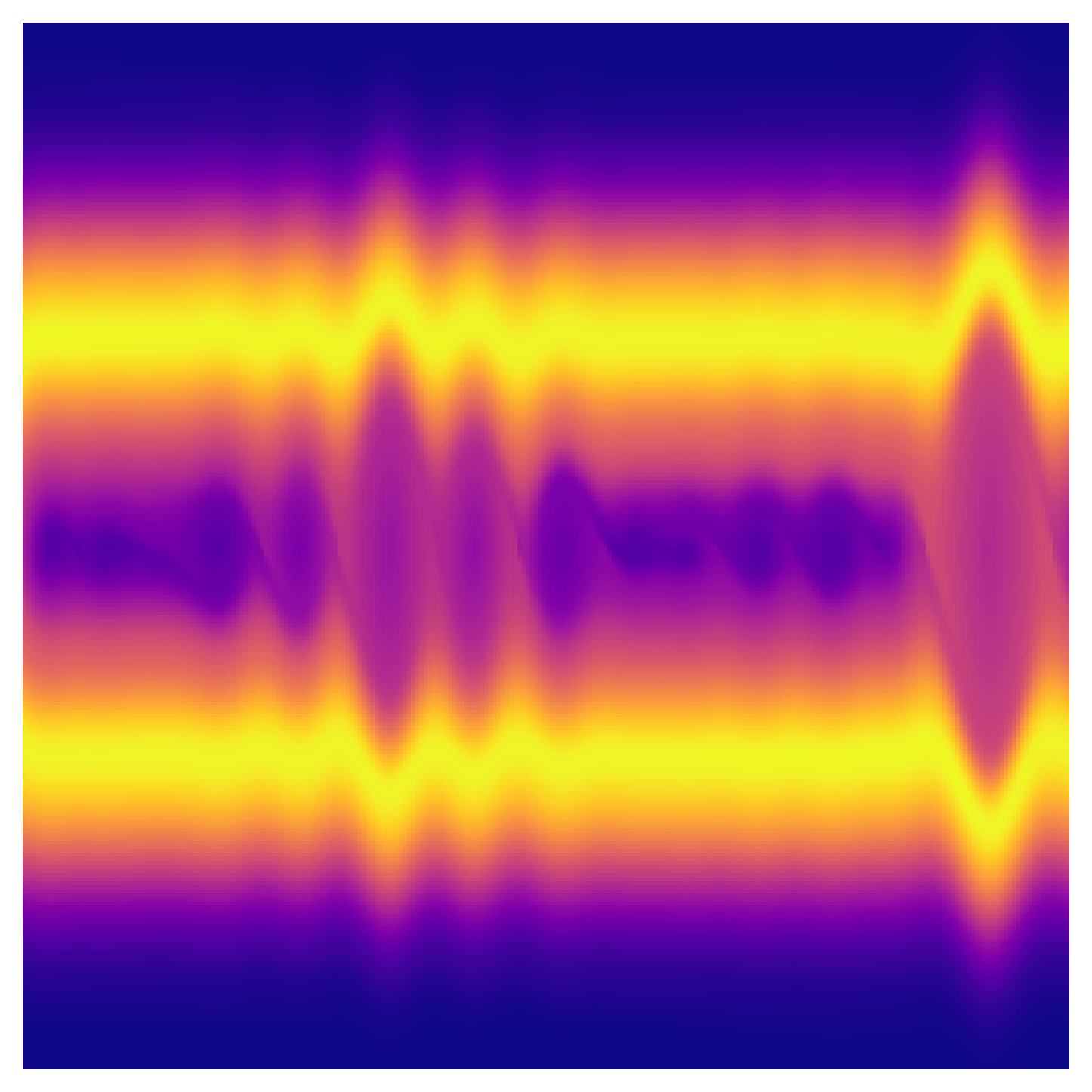} & \includegraphics[width=0.11\linewidth]{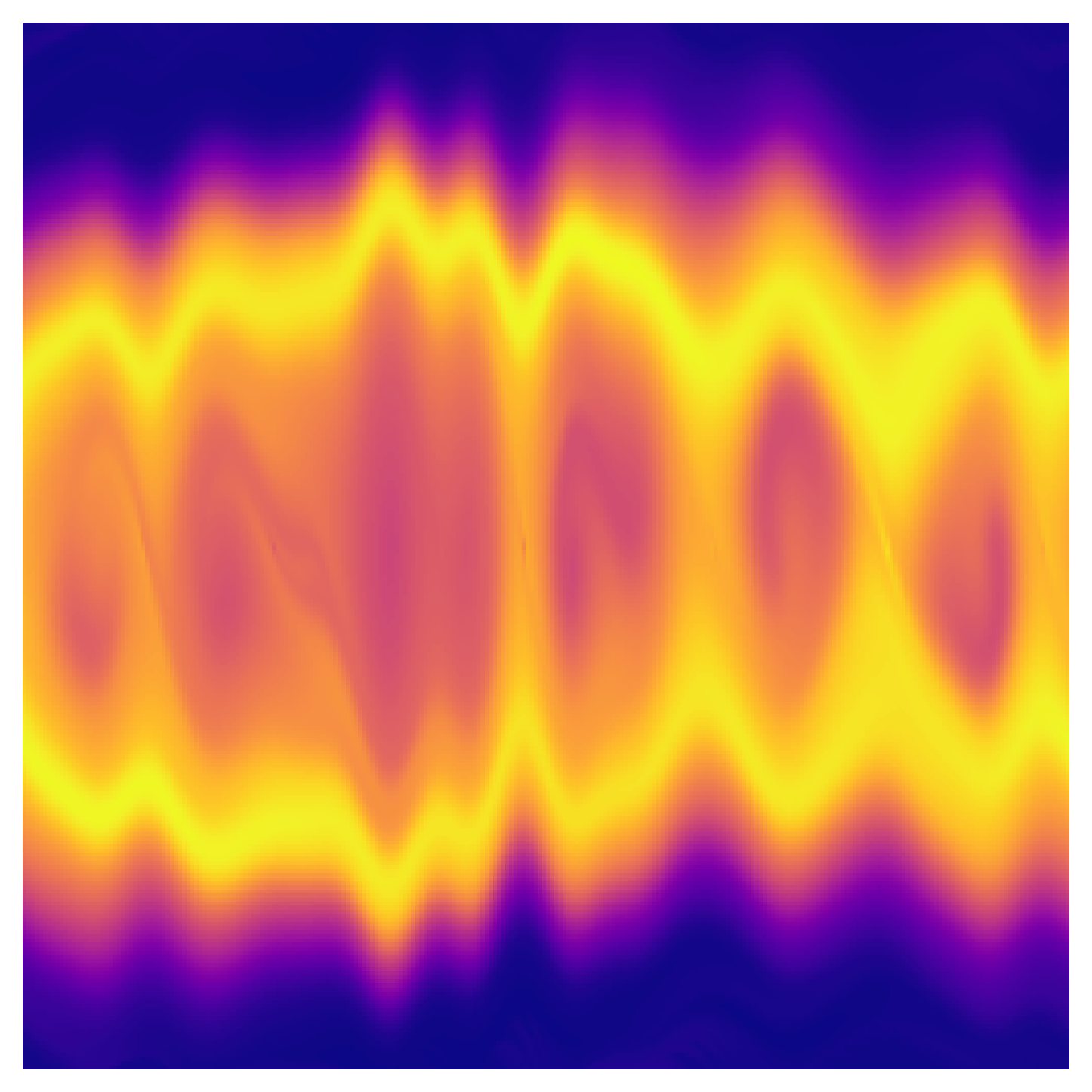} & \includegraphics[width=0.11\linewidth]{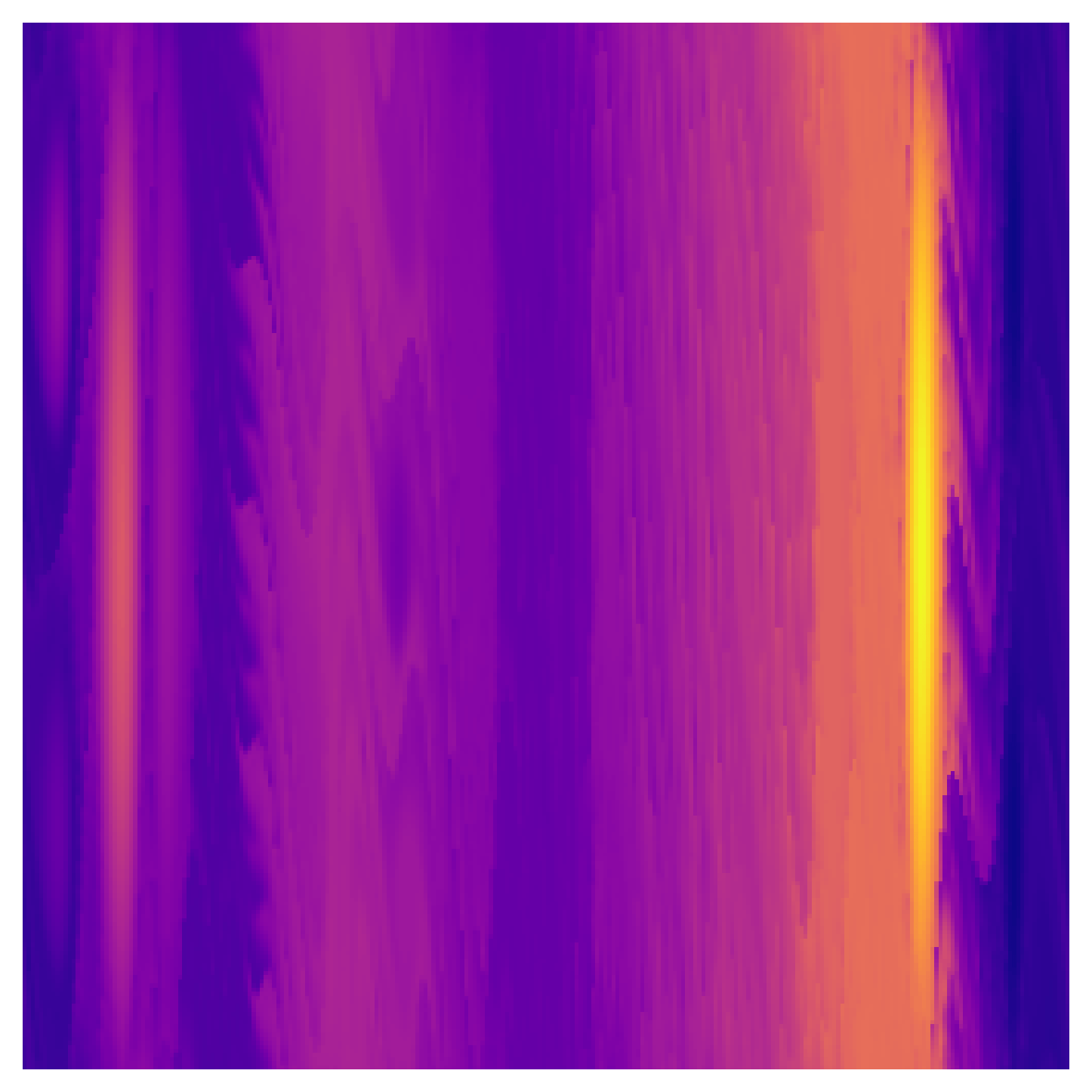} \\ 
    & & (Fig.~\ref{fig:TS_KL_GDL_far_under}) & (Fig.~\ref{fig:TS_ee_lf_GDL_far_under}) & (Fig.~\ref{fig:TS_ee_GDL_far_under}) &(Fig.~\ref{fig:TS_KL_GDL_far_over}) & (Fig.~\ref{fig:TS_ee_lf_GDL_far_over}) & (Fig.~\ref{fig:TS_ee_GDL_far_over}) \\ \cline{2-8}
    & \multirow{2}{*}[22pt]{Local} & \includegraphics[width=0.11\linewidth]{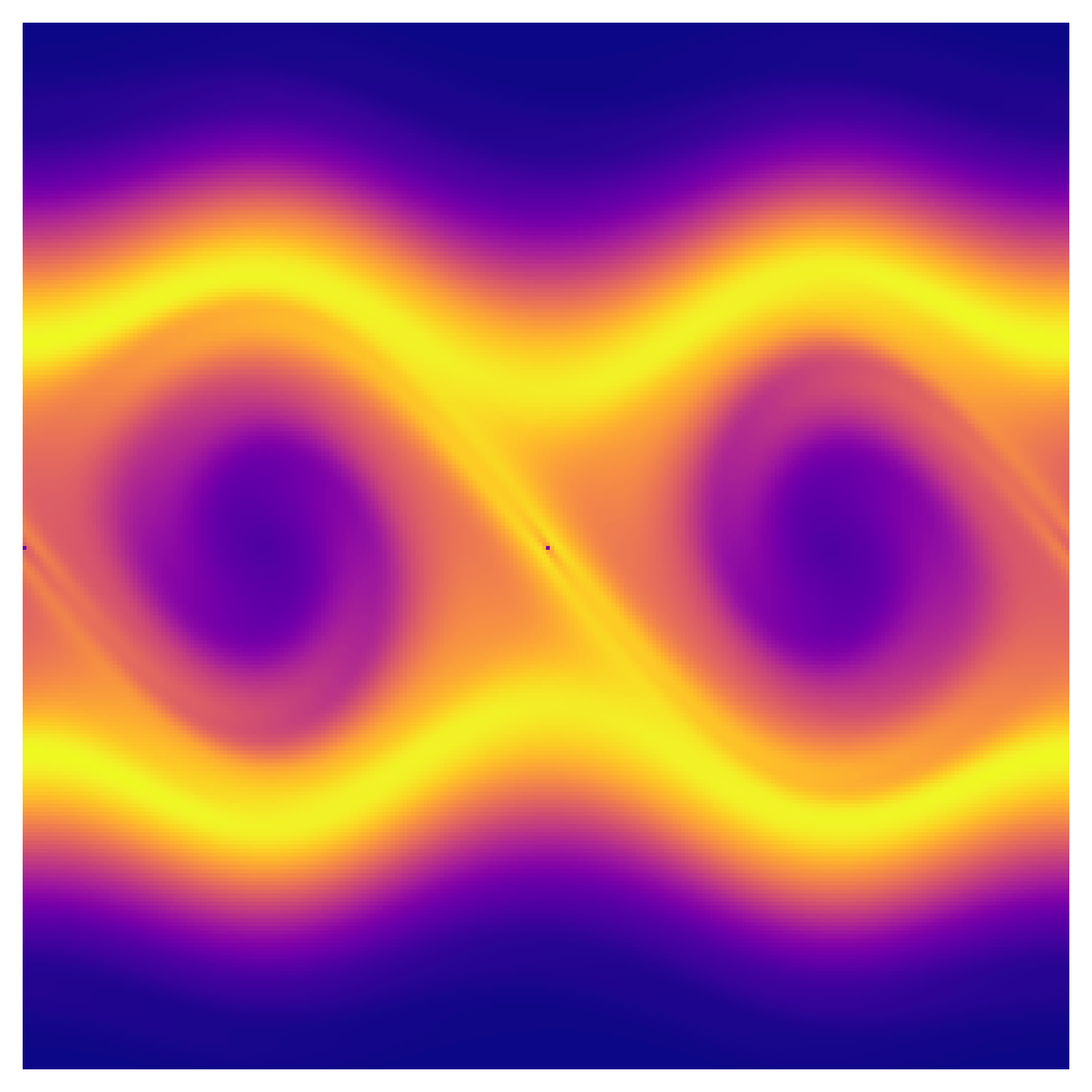} & \includegraphics[width=0.11\linewidth]{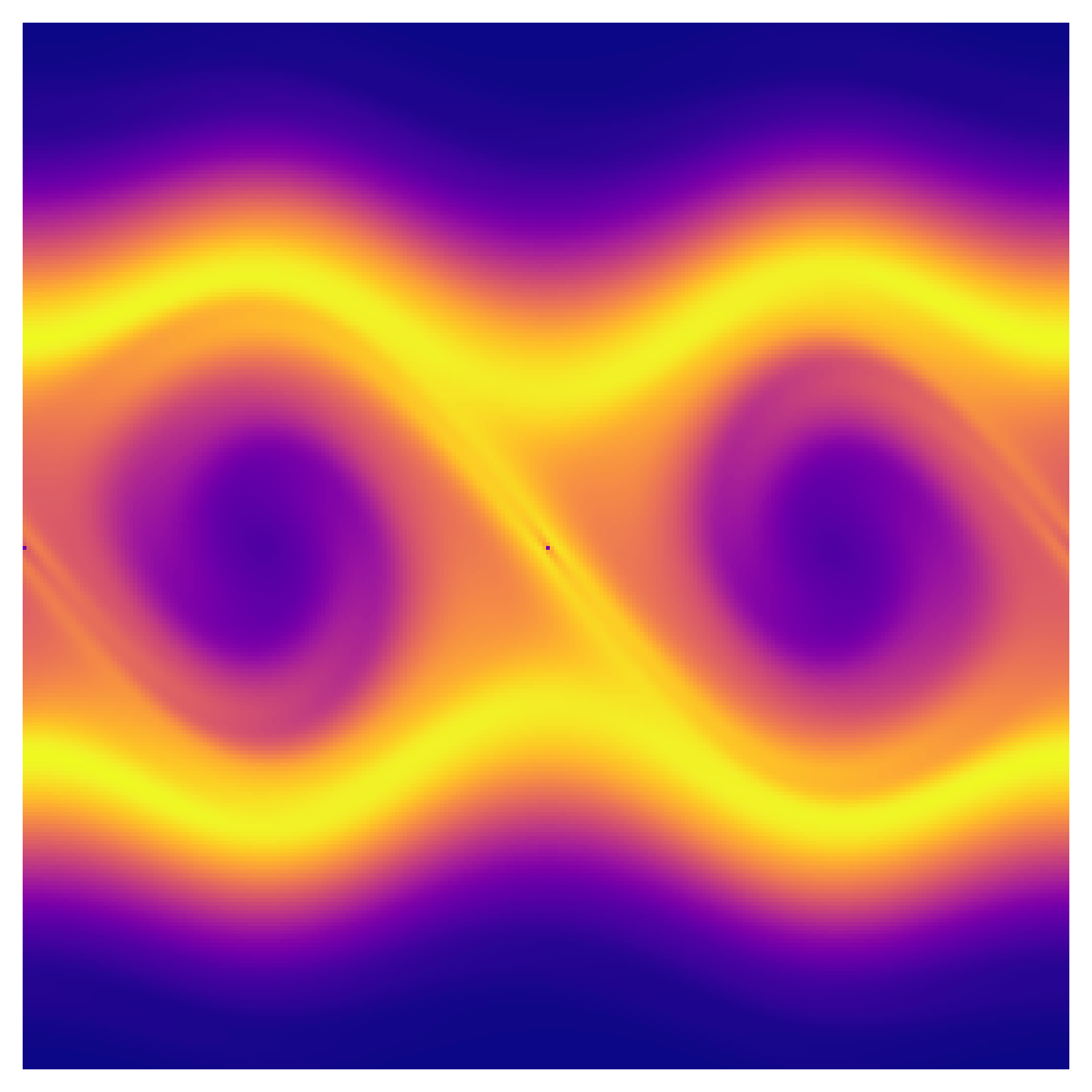} & \includegraphics[width=0.11\linewidth]{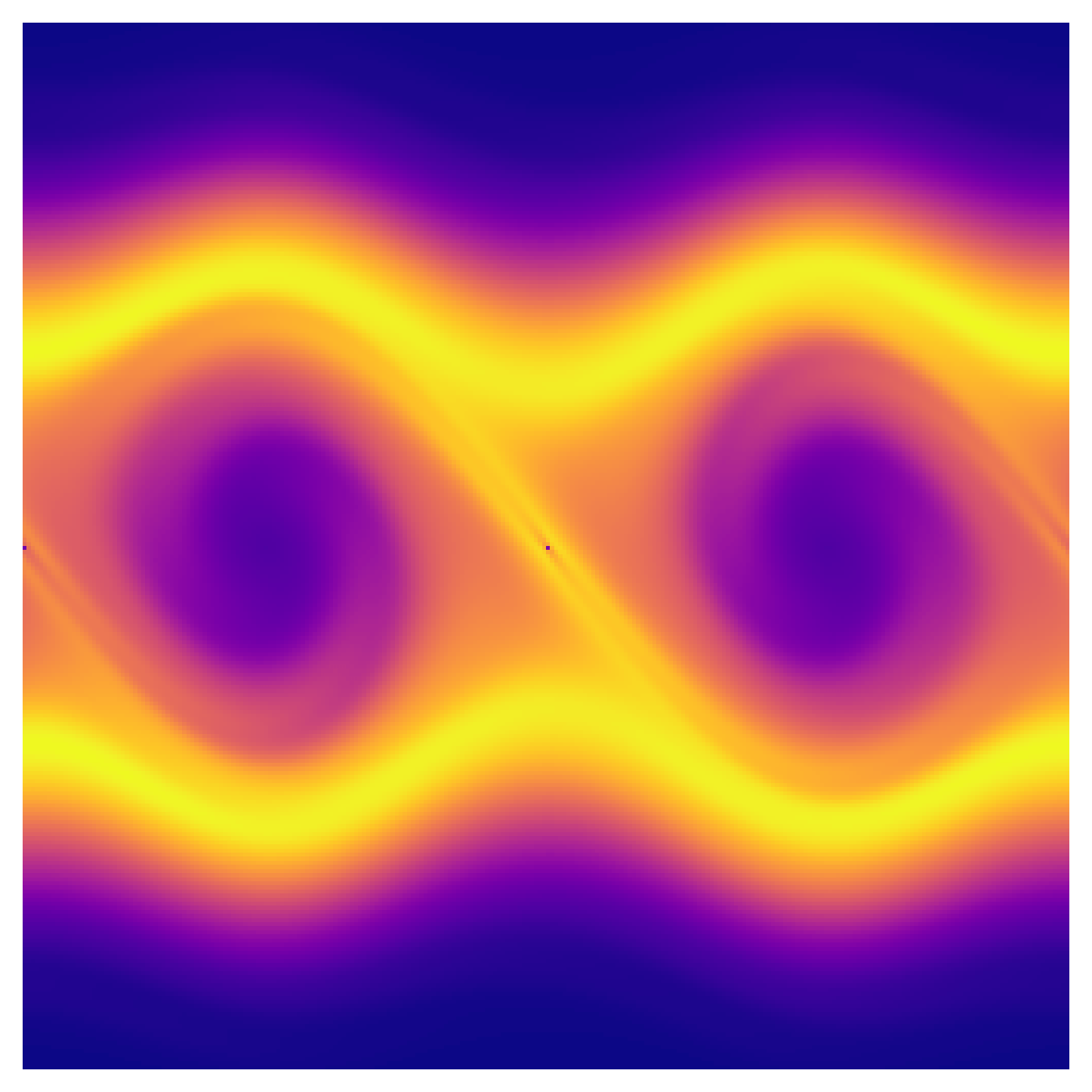} & \includegraphics[width=0.11\linewidth]{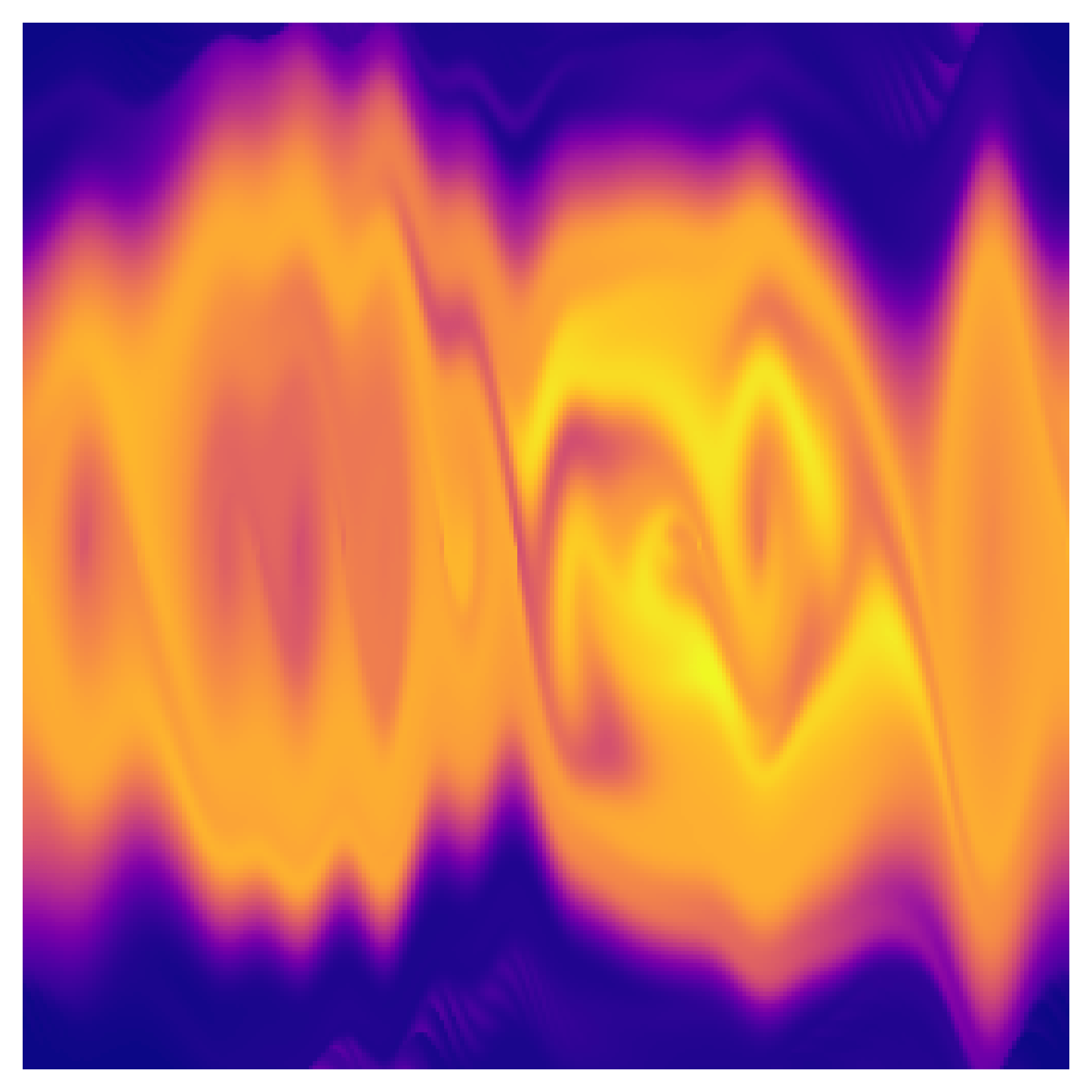} & \includegraphics[width=0.11\linewidth]{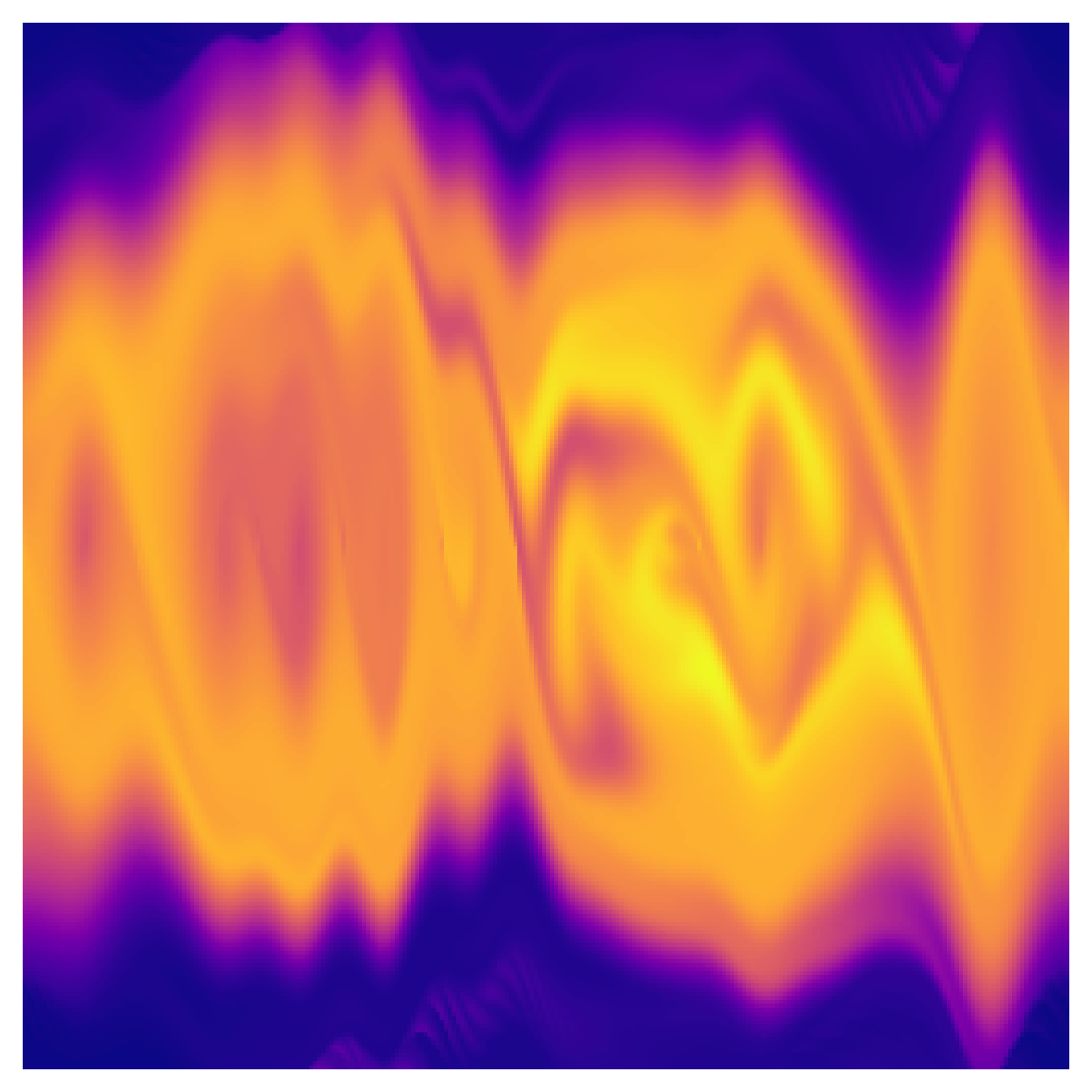} & \includegraphics[width=0.11\linewidth]{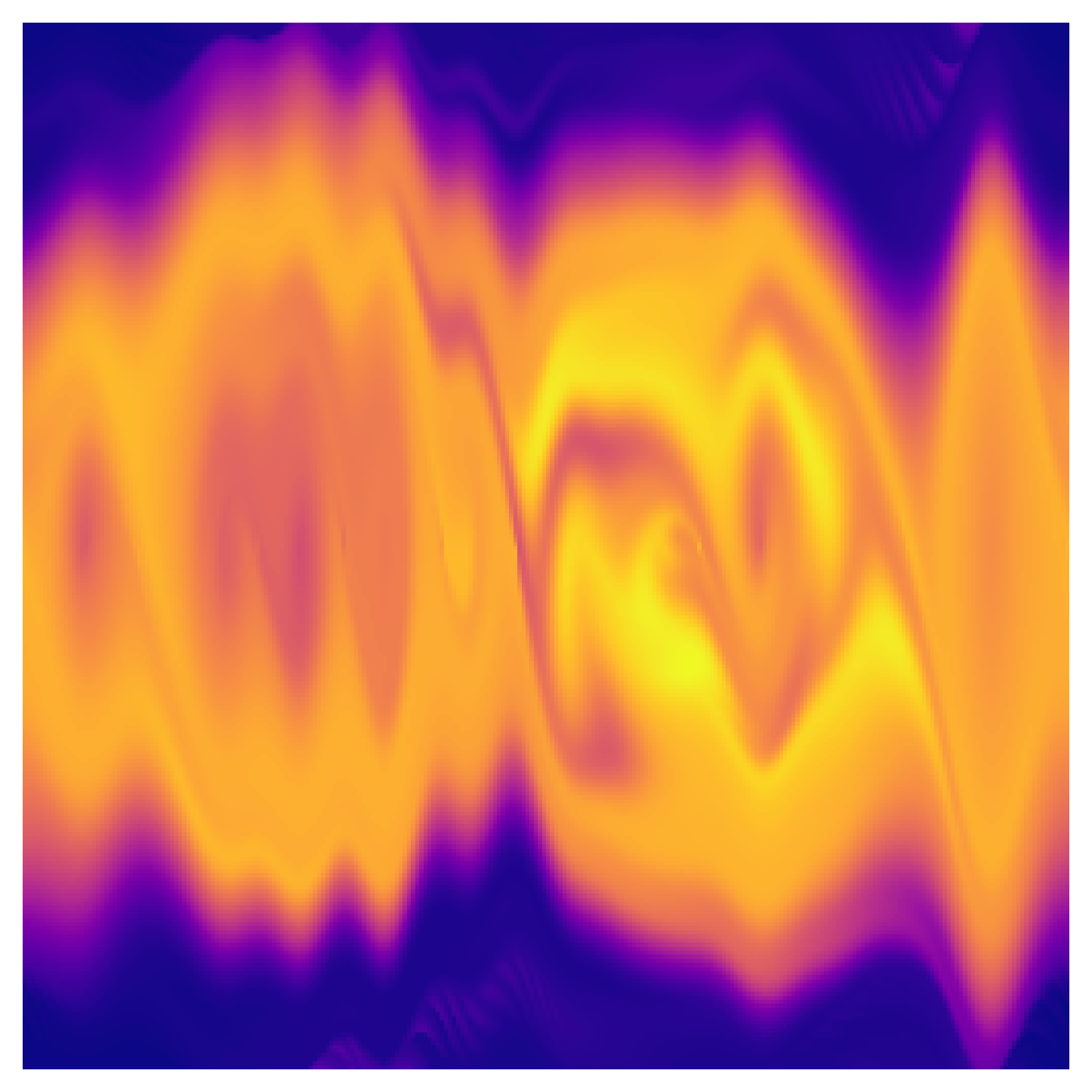} \\
    & & (Fig.~\ref{fig:TS_KL_GD_far_under}) & (Fig.~\ref{fig:TS_ee_lf_GD_far_under}) & (Fig.~\ref{fig:TS_ee_GD_far_under}) & (Fig.~\ref{fig:TS_KL_GD_far_over}) & (Fig.~\ref{fig:TS_ee_lf_GD_far_over}) & (Fig.~\ref{fig:TS_ee_GD_far_over}) \\
    \midrule
    \multirow{4}{*}[0pt]{Mid} & \multirow{2}{*}[22pt]{Adaptive} & \includegraphics[width=0.11\linewidth]{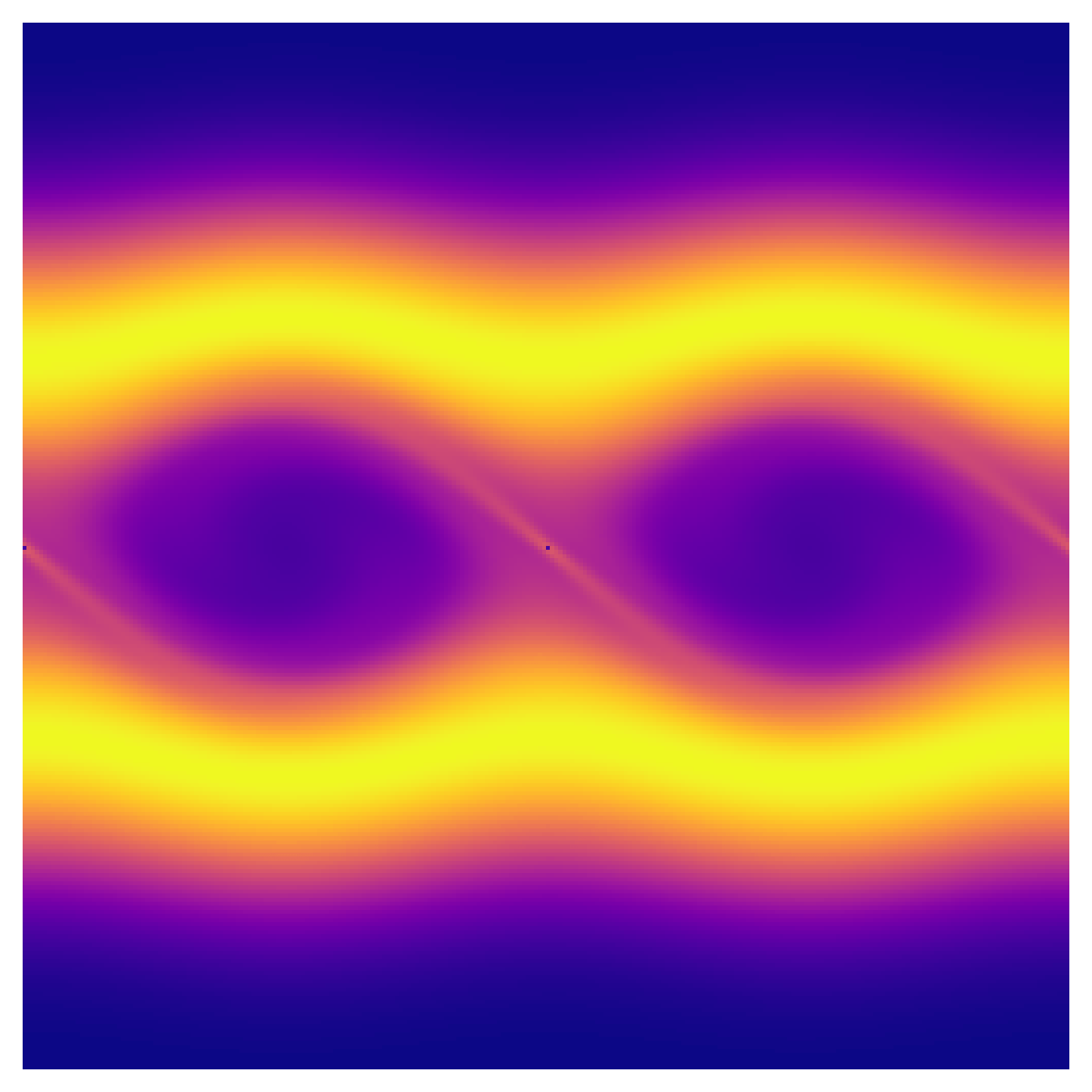} & \includegraphics[width=0.11\linewidth]{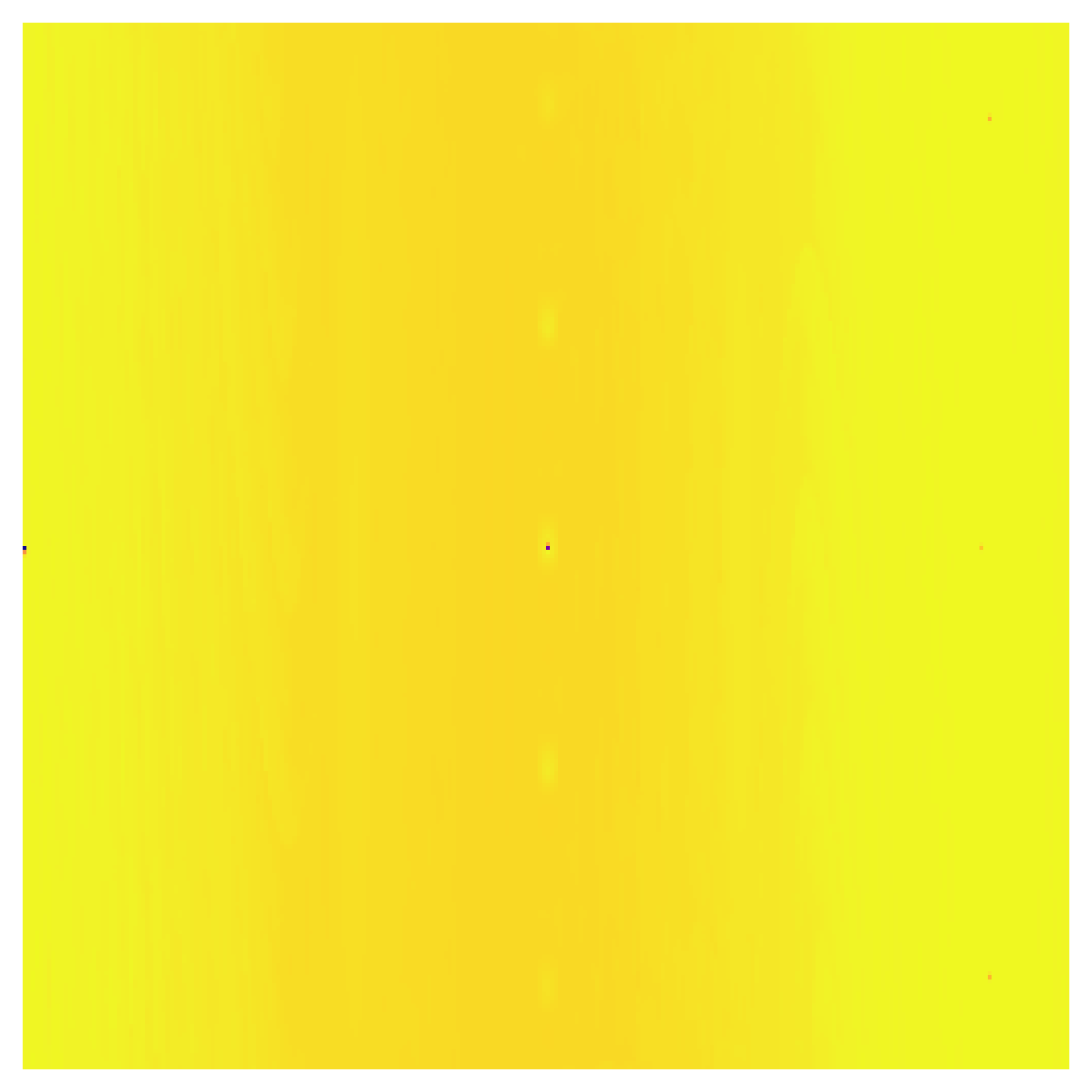} &
    \includegraphics[width=0.11\linewidth]{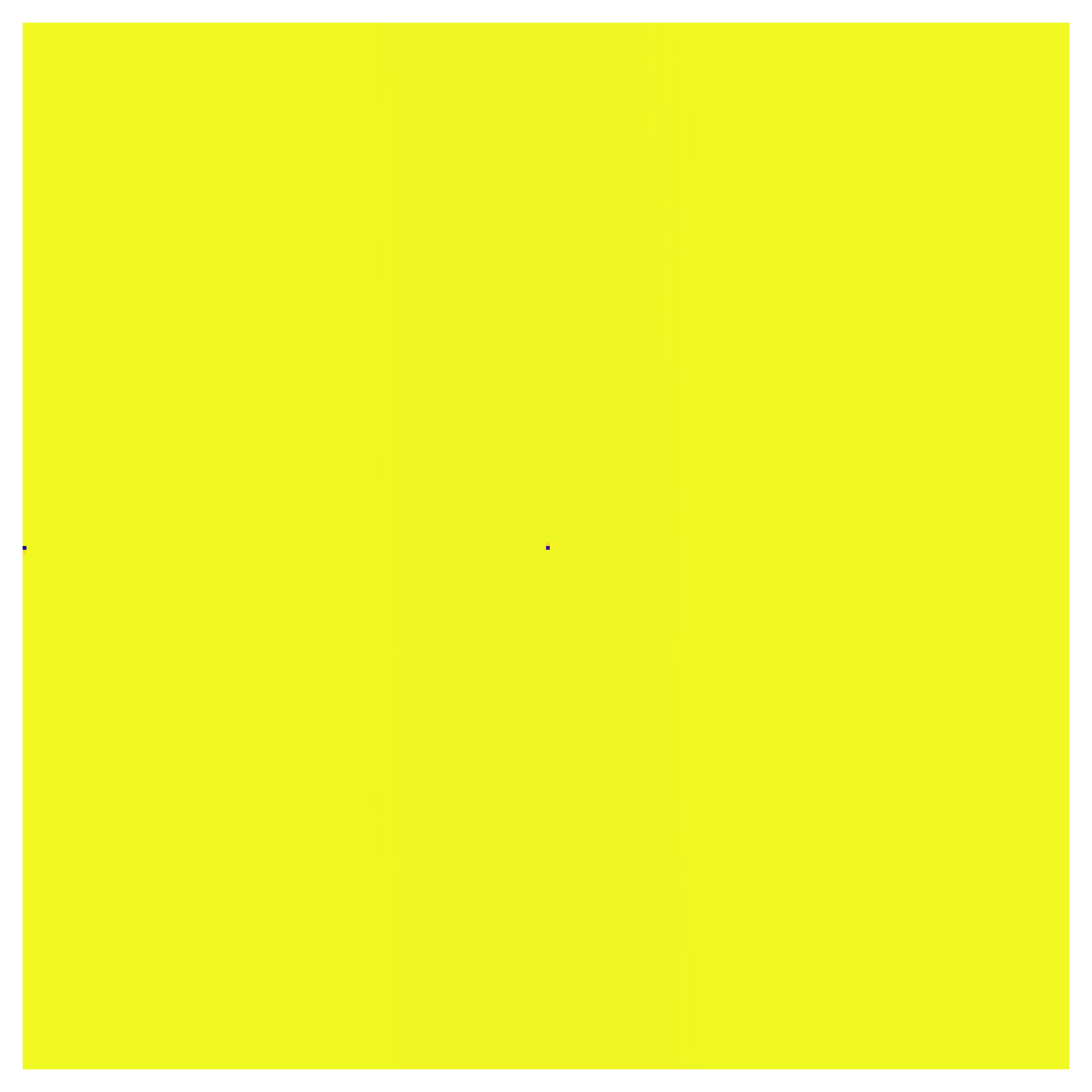} & \includegraphics[width=0.11\linewidth]{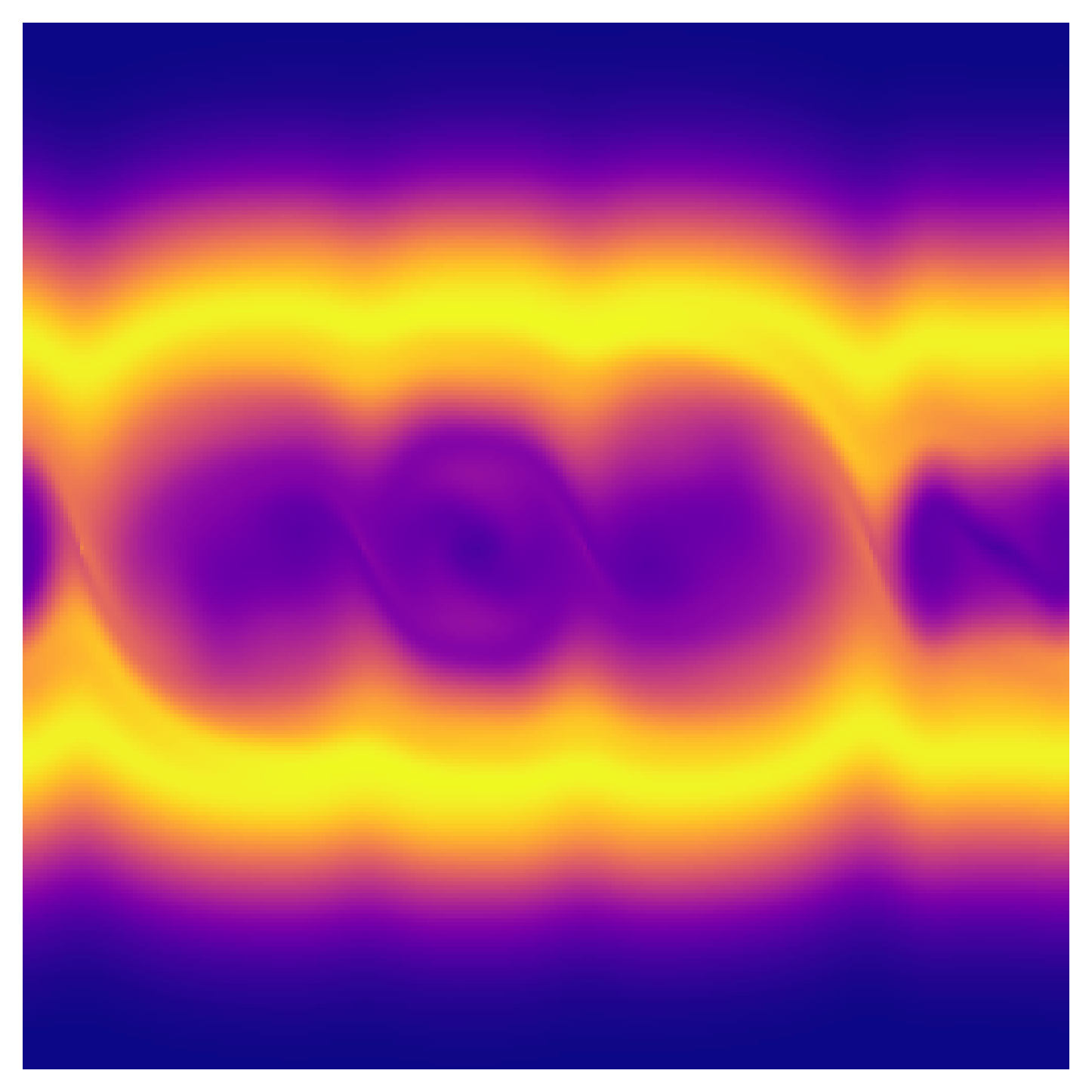} & \includegraphics[width=0.11\linewidth]{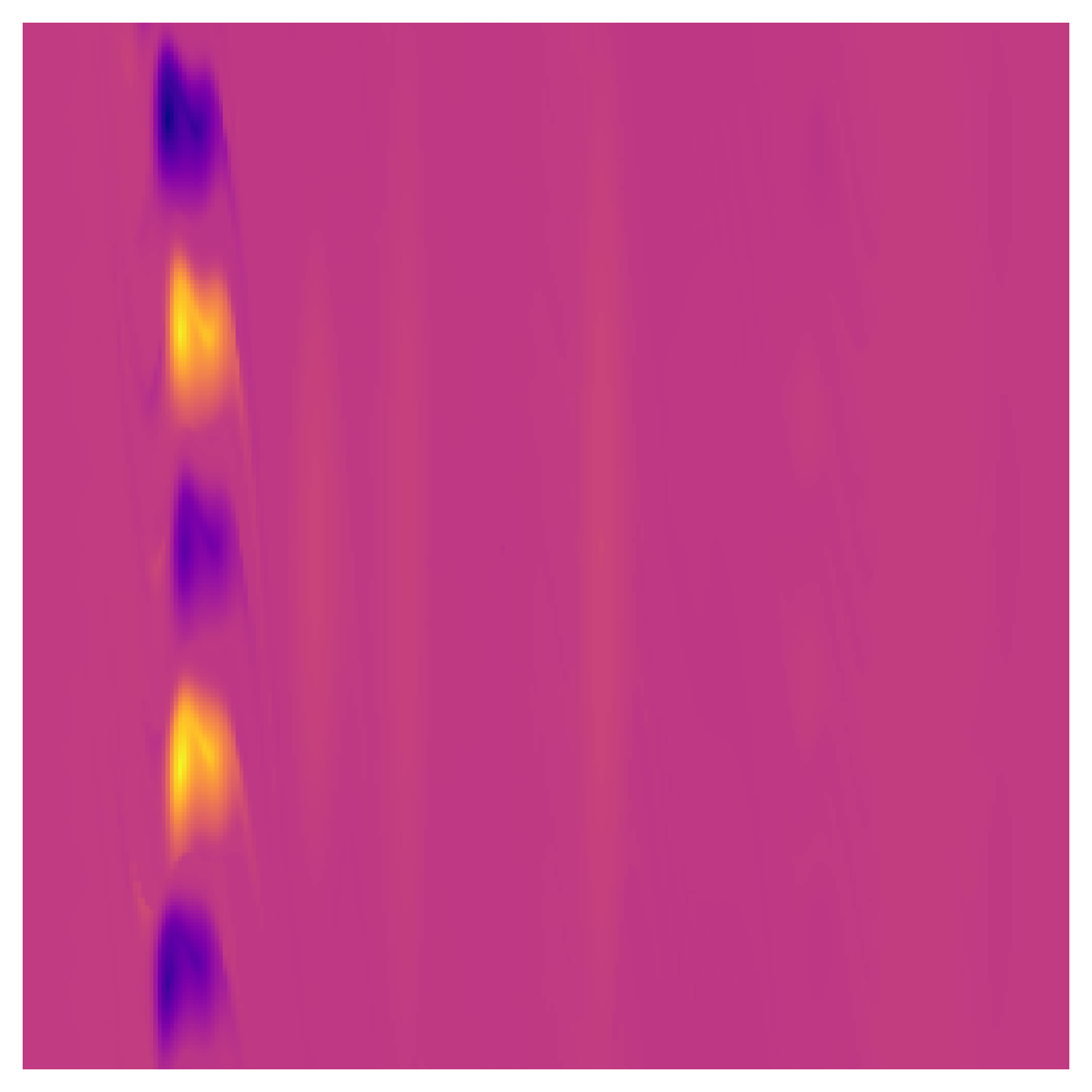} & \includegraphics[width=0.11\linewidth]{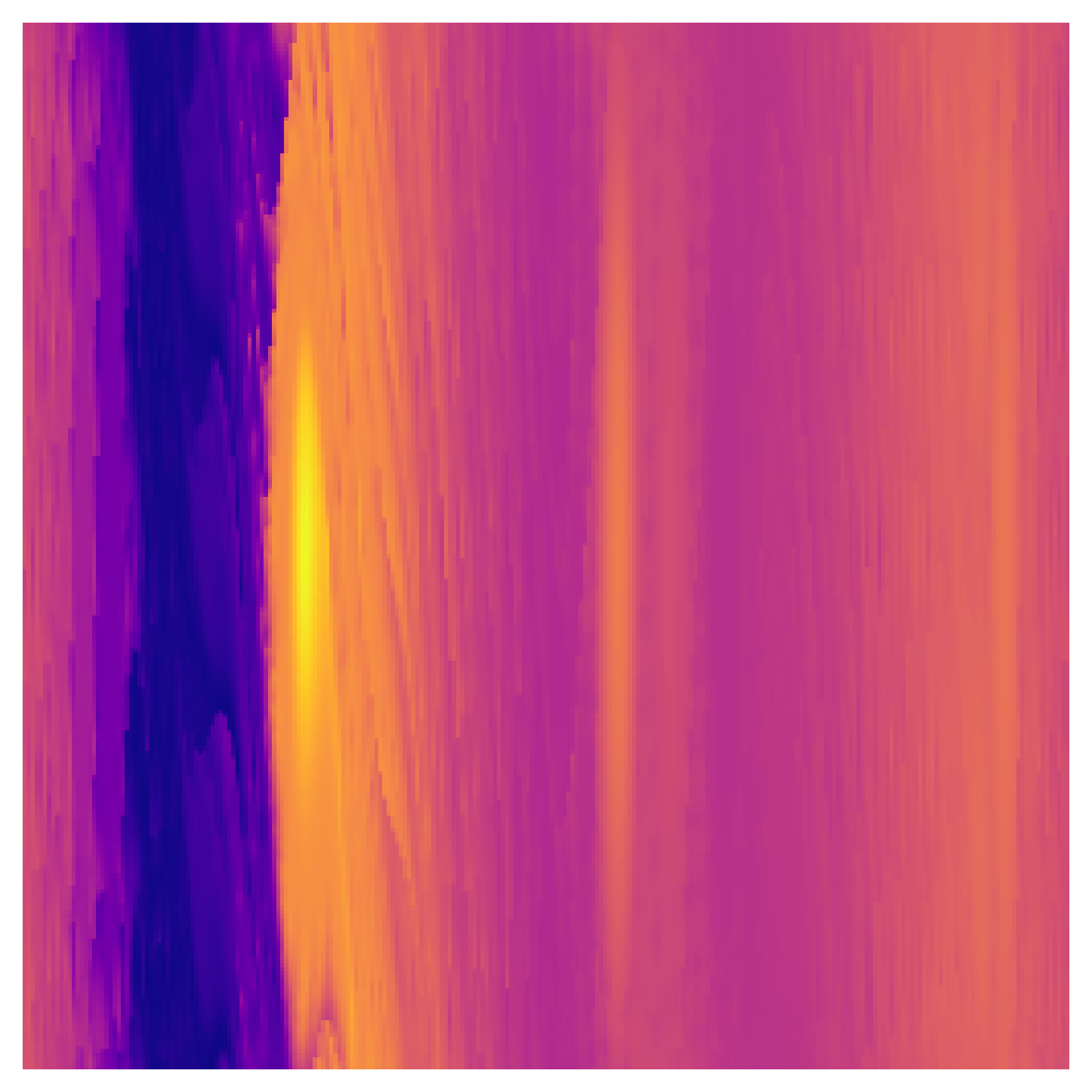} \\
    & & (Fig.~\ref{fig:TS_KL_GDL_near_under}) & (Fig.~\ref{fig:TS_ee_lf_GDL_near_under}) & (Fig.~\ref{fig:TS_ee_GDL_near_under}) & (Fig.~\ref{fig:TS_KL_GDL_near_over}) & (Fig.~\ref{fig:TS_ee_lf_GDL_near_over}) & (Fig.~\ref{fig:TS_ee_GDL_near_over}) \\ \cline{2-8}
    & \multirow{2}{*}[22pt]{Local} & \includegraphics[width=0.11\linewidth]{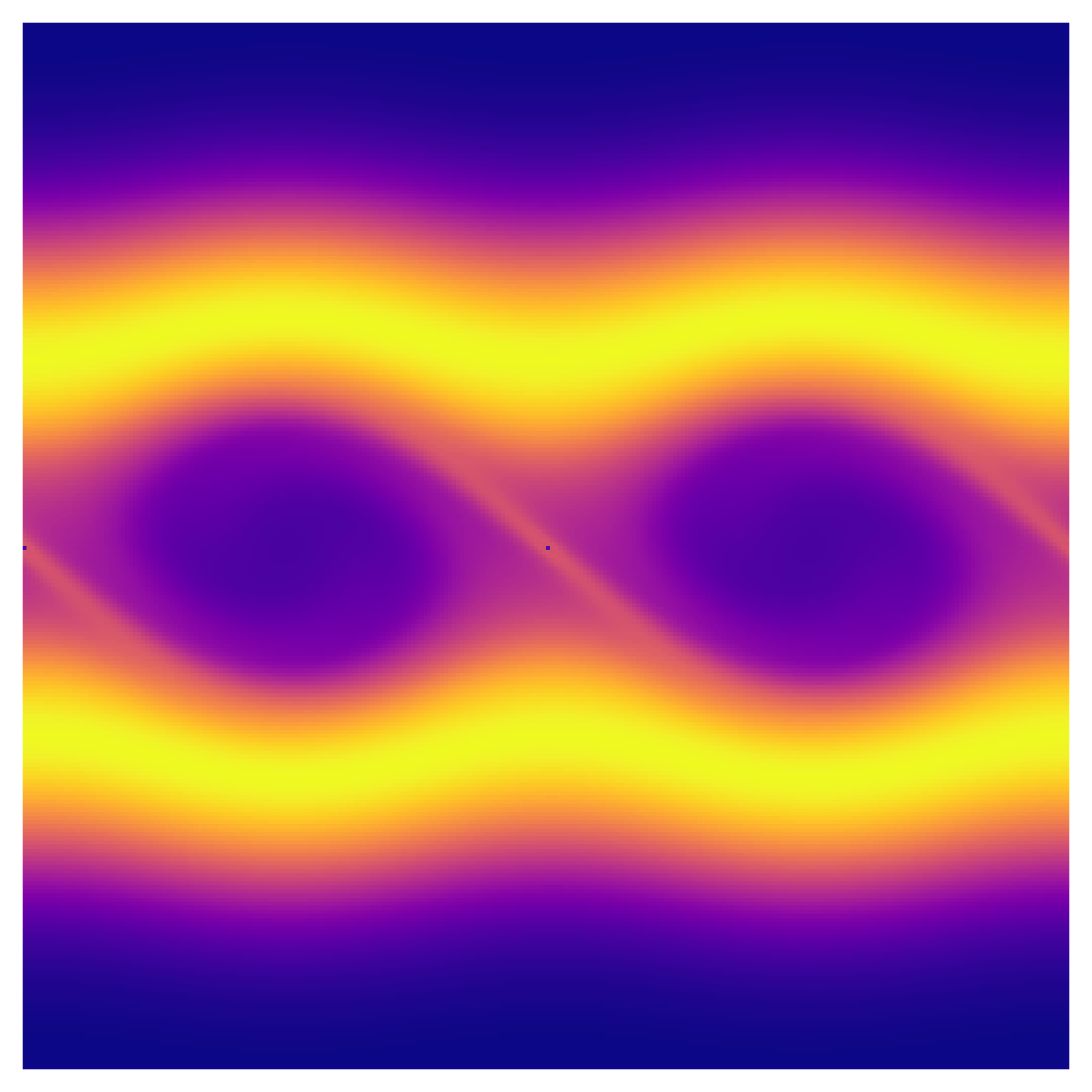} & \includegraphics[width=0.11\linewidth]{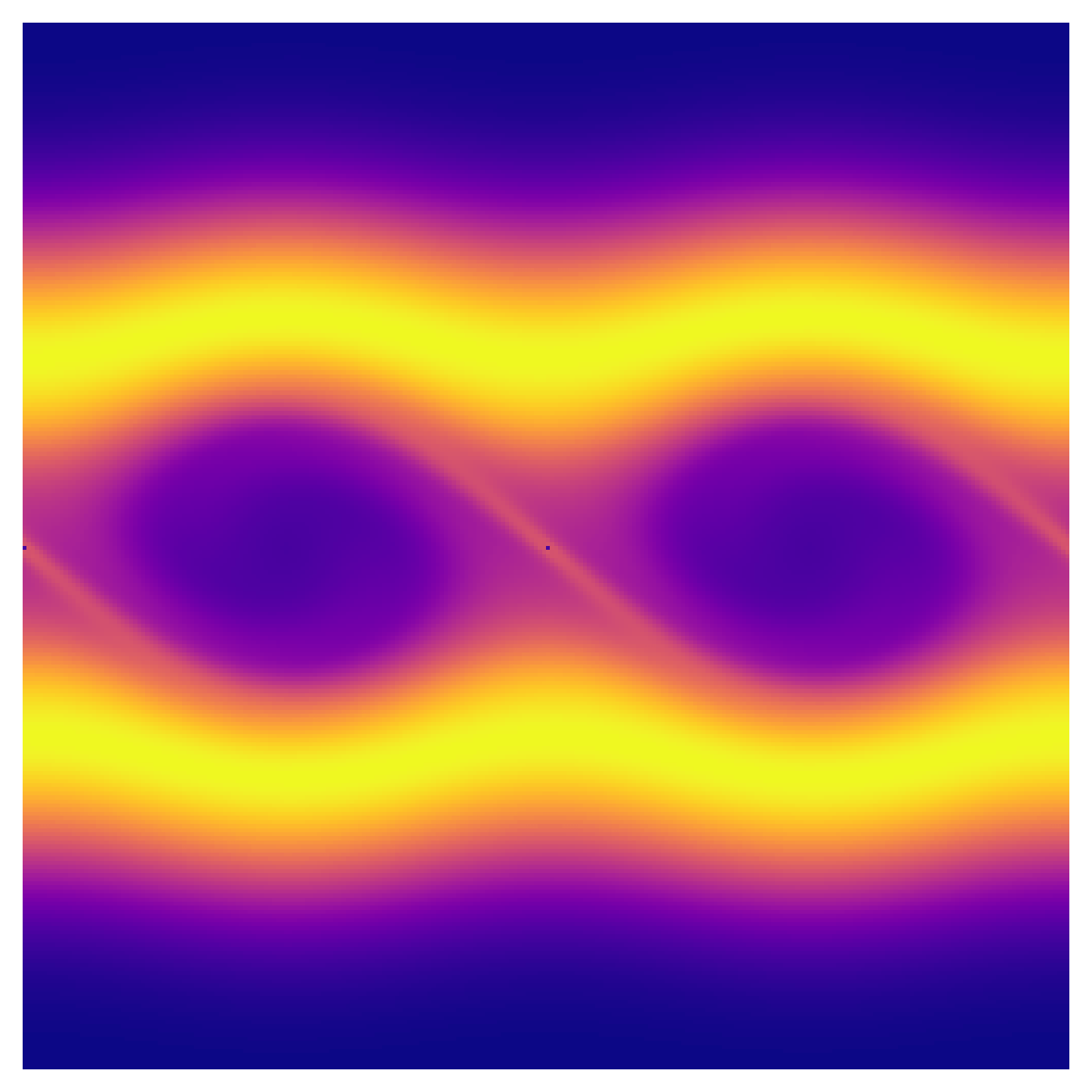} & \includegraphics[width=0.11\linewidth]{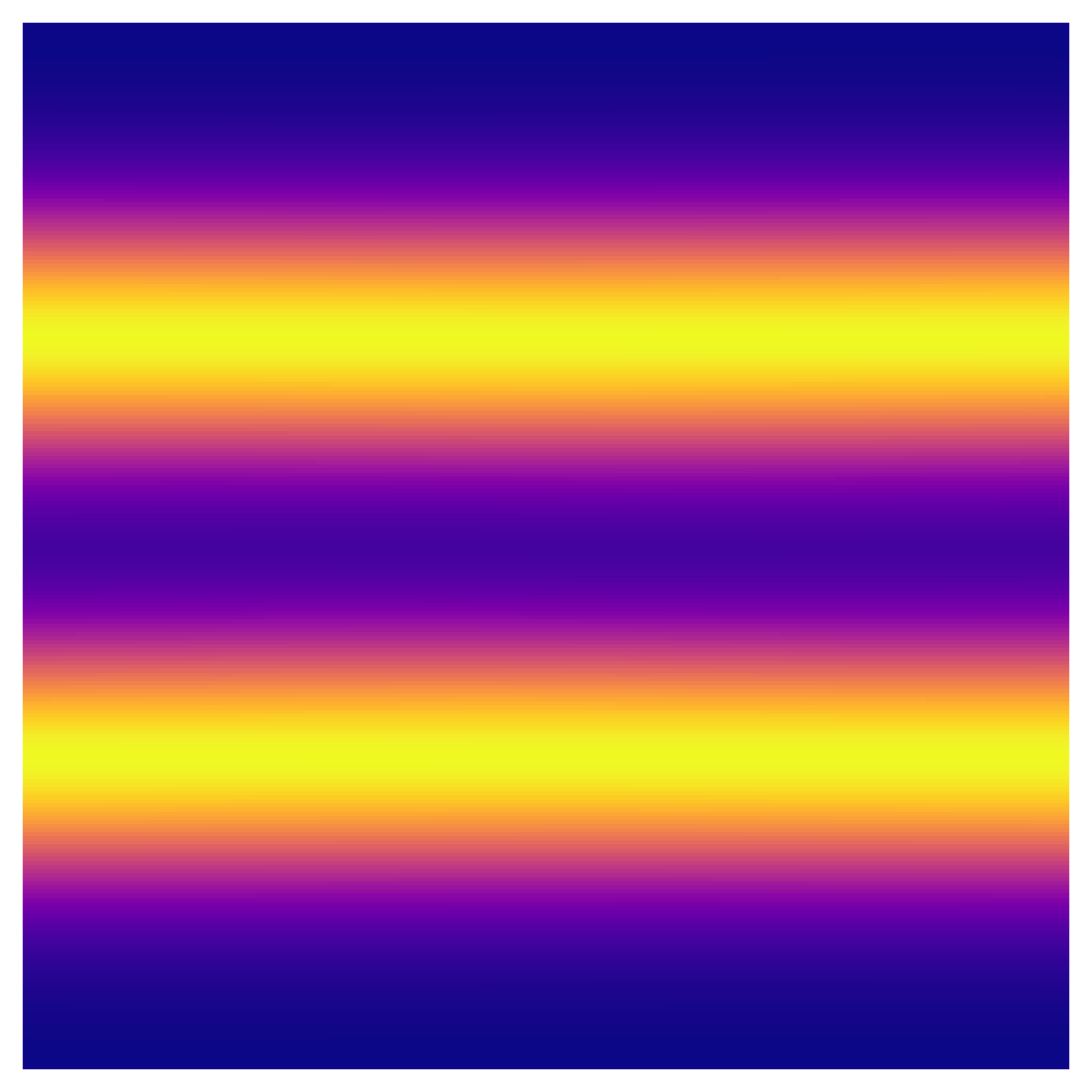} & \includegraphics[width=0.11\linewidth]{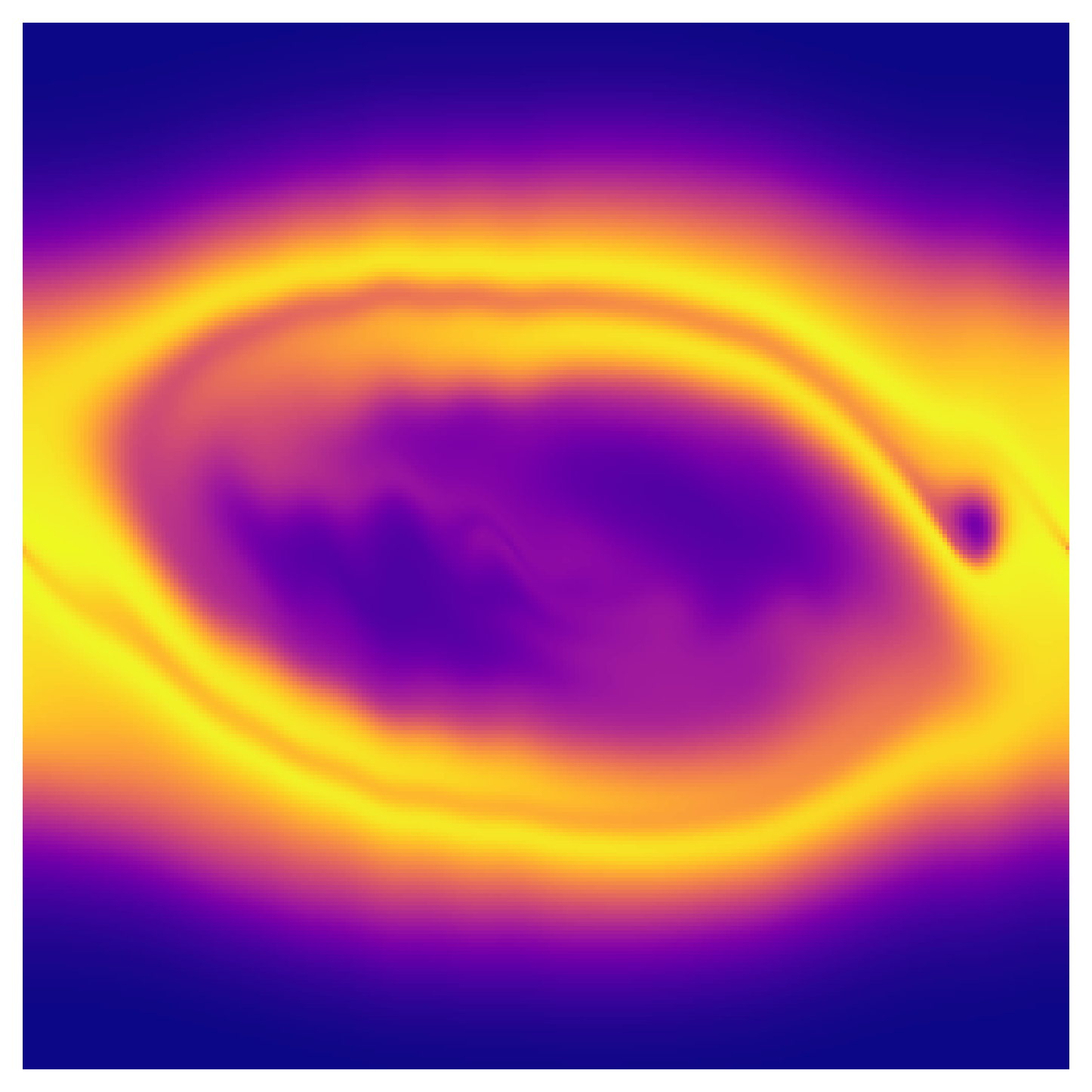} & \includegraphics[width=0.11\linewidth]{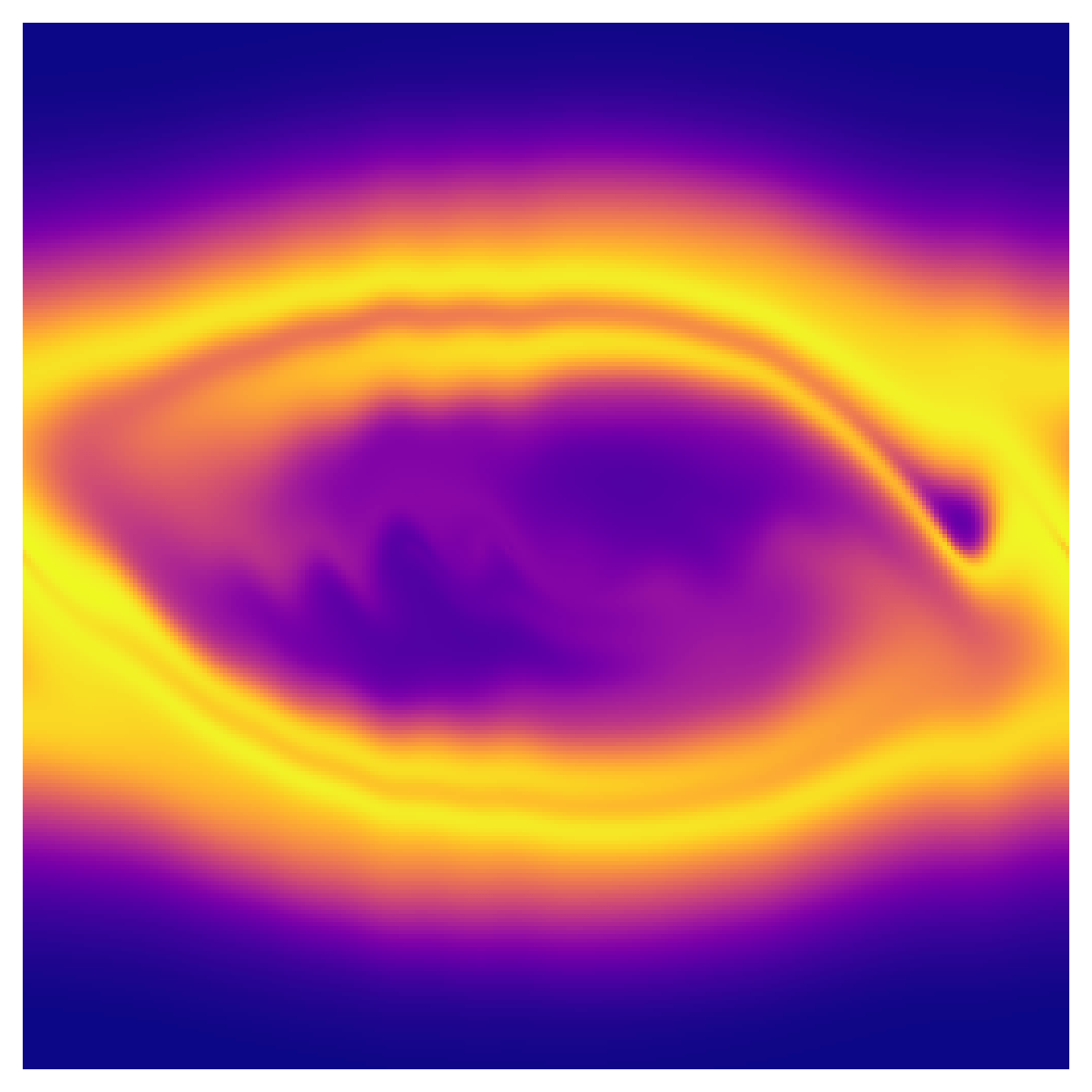} & \includegraphics[width=0.11\linewidth]{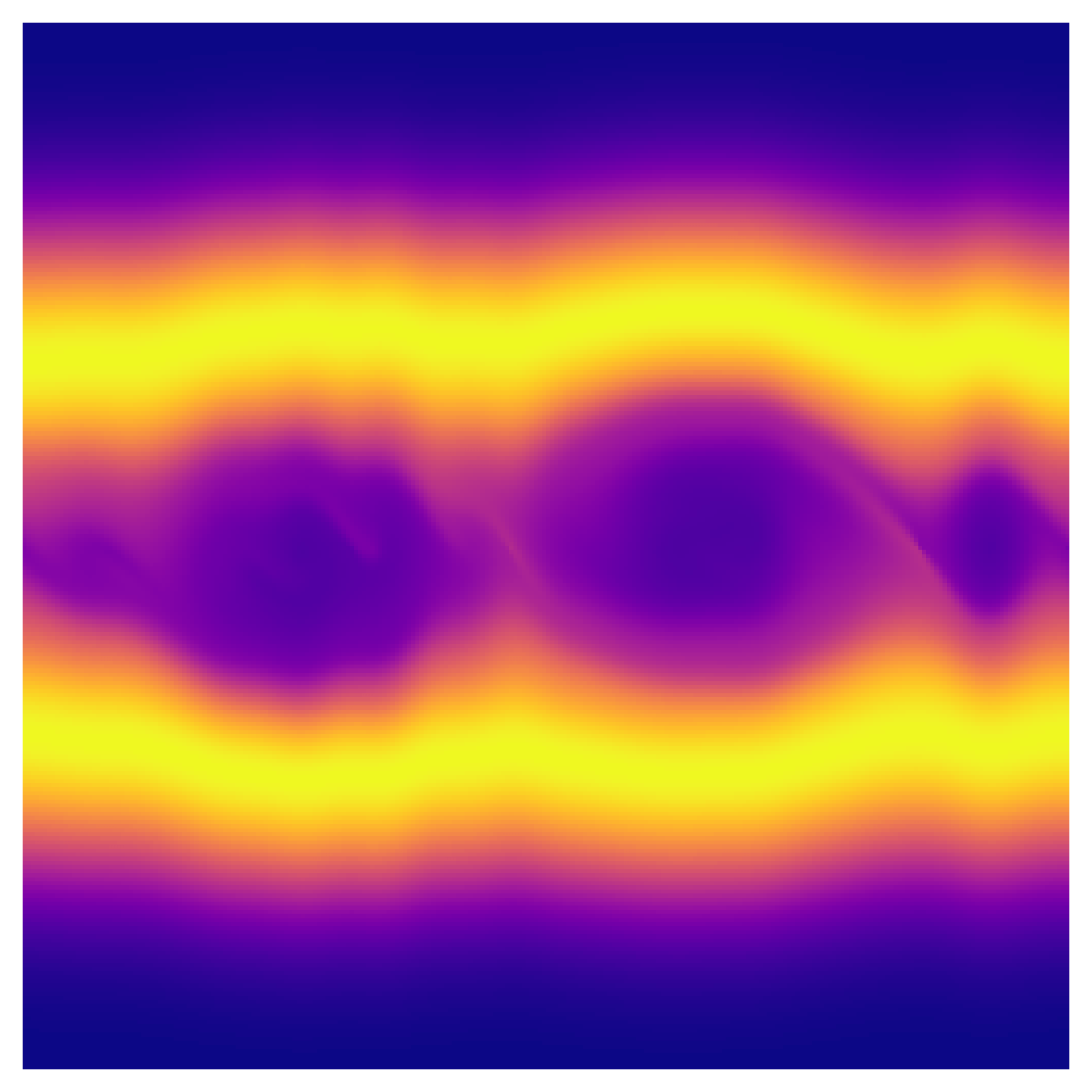} \\
    & & (Fig.~\ref{fig:TS_KL_GD_near_under}) & (Fig.~\ref{fig:TS_ee_lf_GD_near_under}) & 
    (Fig.~\ref{fig:TS_ee_GD_near_under}) & (Fig.~\ref{fig:TS_KL_GD_near_over}) & (Fig.~\ref{fig:TS_ee_lf_GD_near_over}) & (Fig.~\ref{fig:TS_ee_GD_near_over}) \\
    \midrule
    \multirow{4}{*}[0pt]{Near} & \multirow{2}{*}[22pt]{Adaptive} & \includegraphics[width=0.11\linewidth]{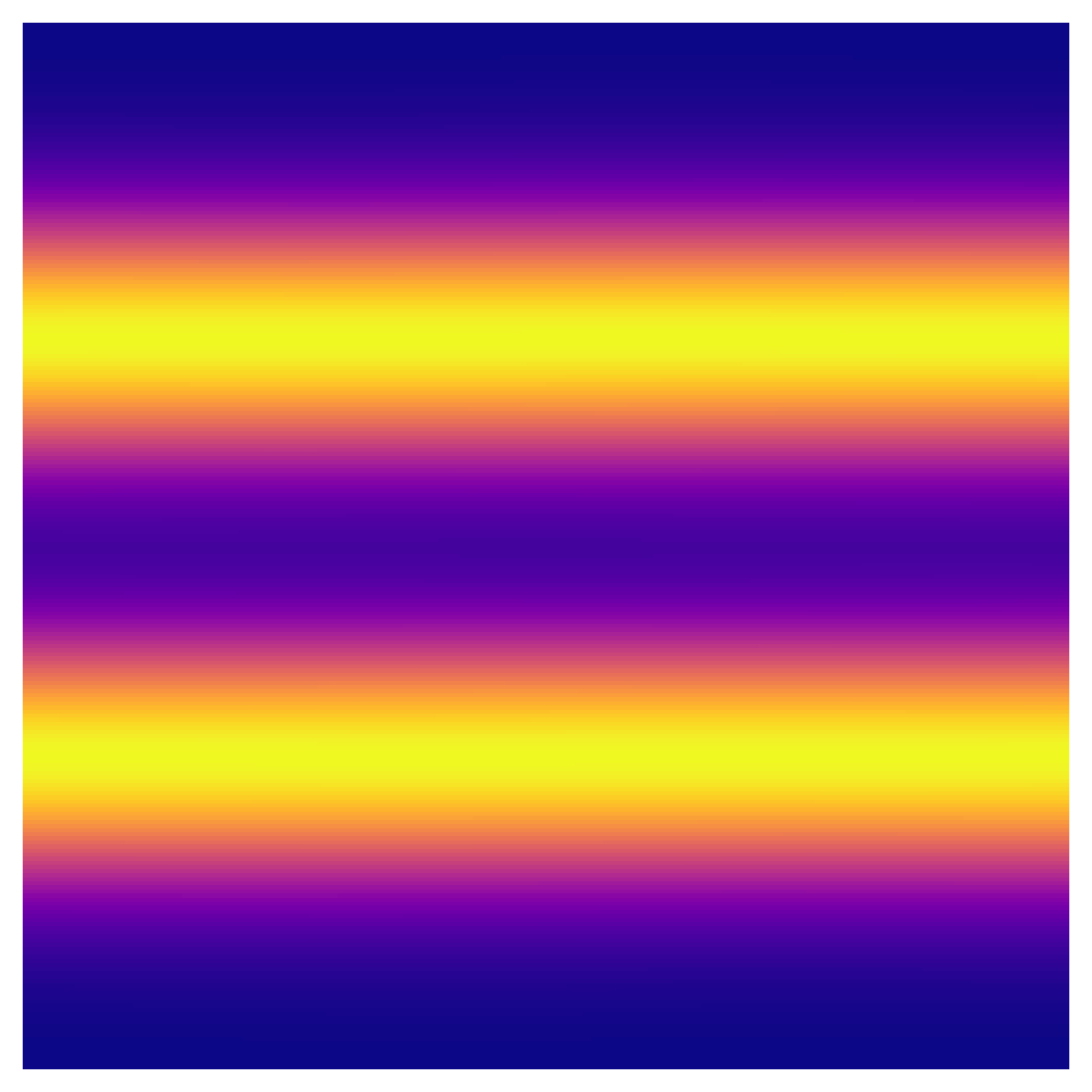}  &  \includegraphics[width=0.11\linewidth]{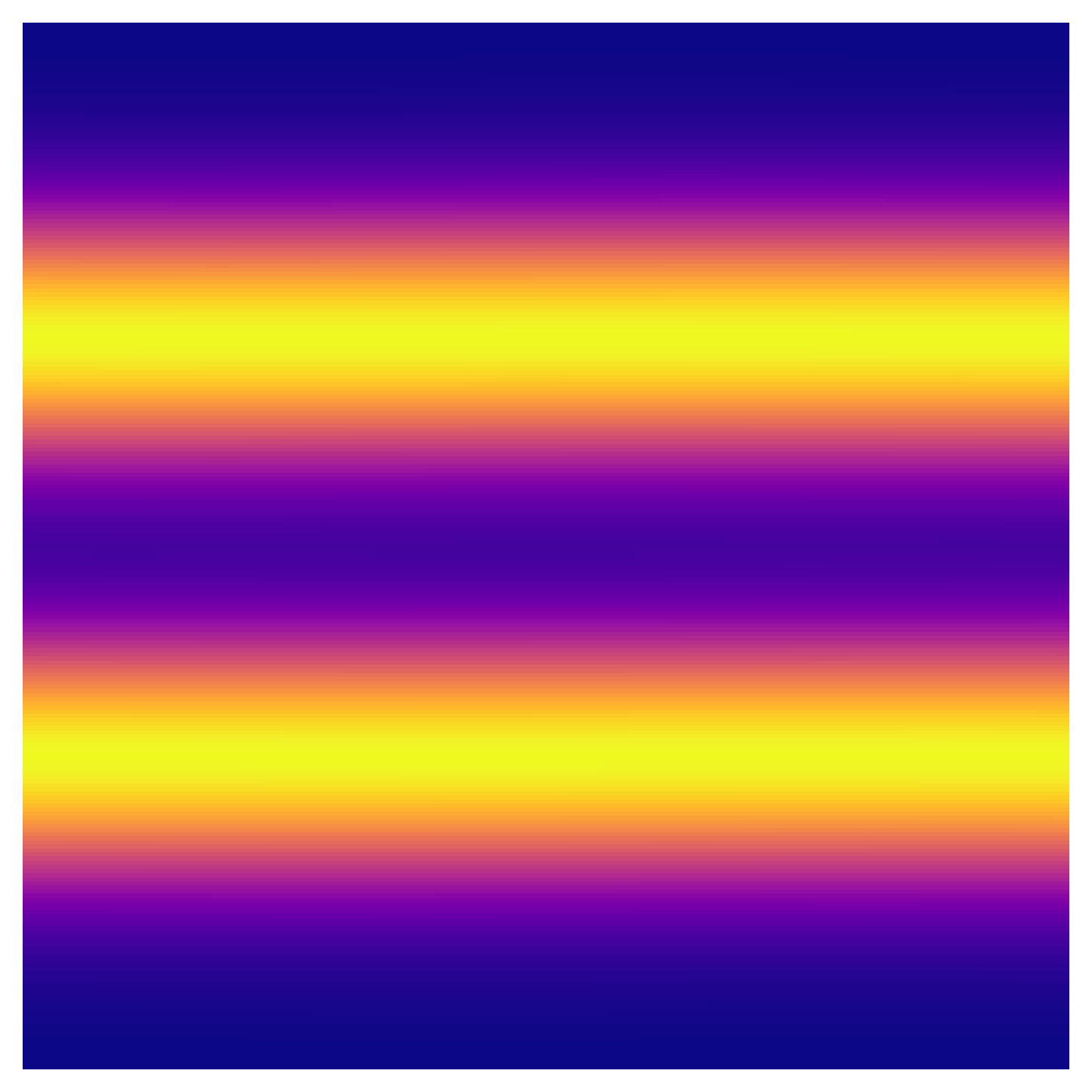}  &
    \includegraphics[width=0.11\linewidth]{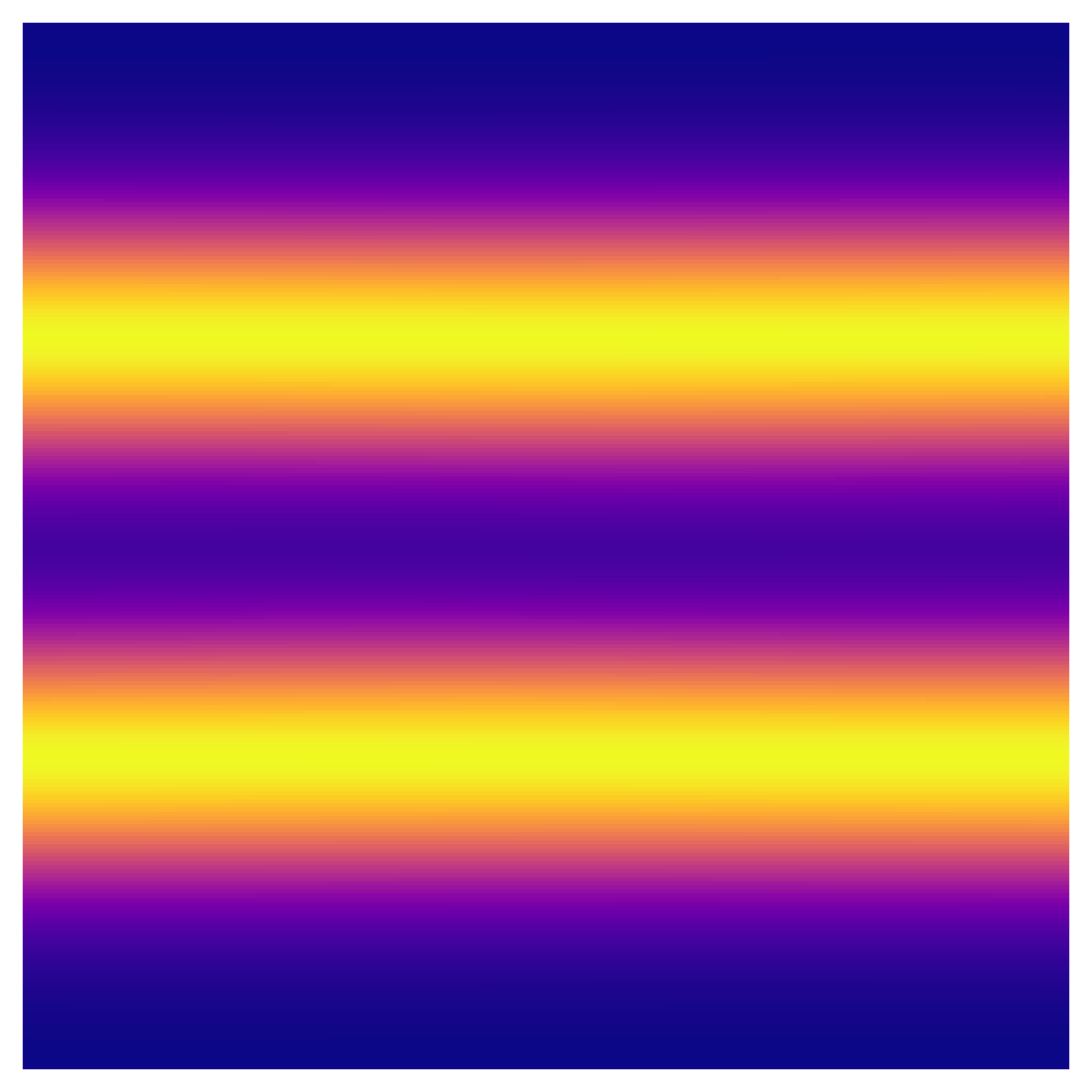}  &
    \includegraphics[width=0.11\linewidth]{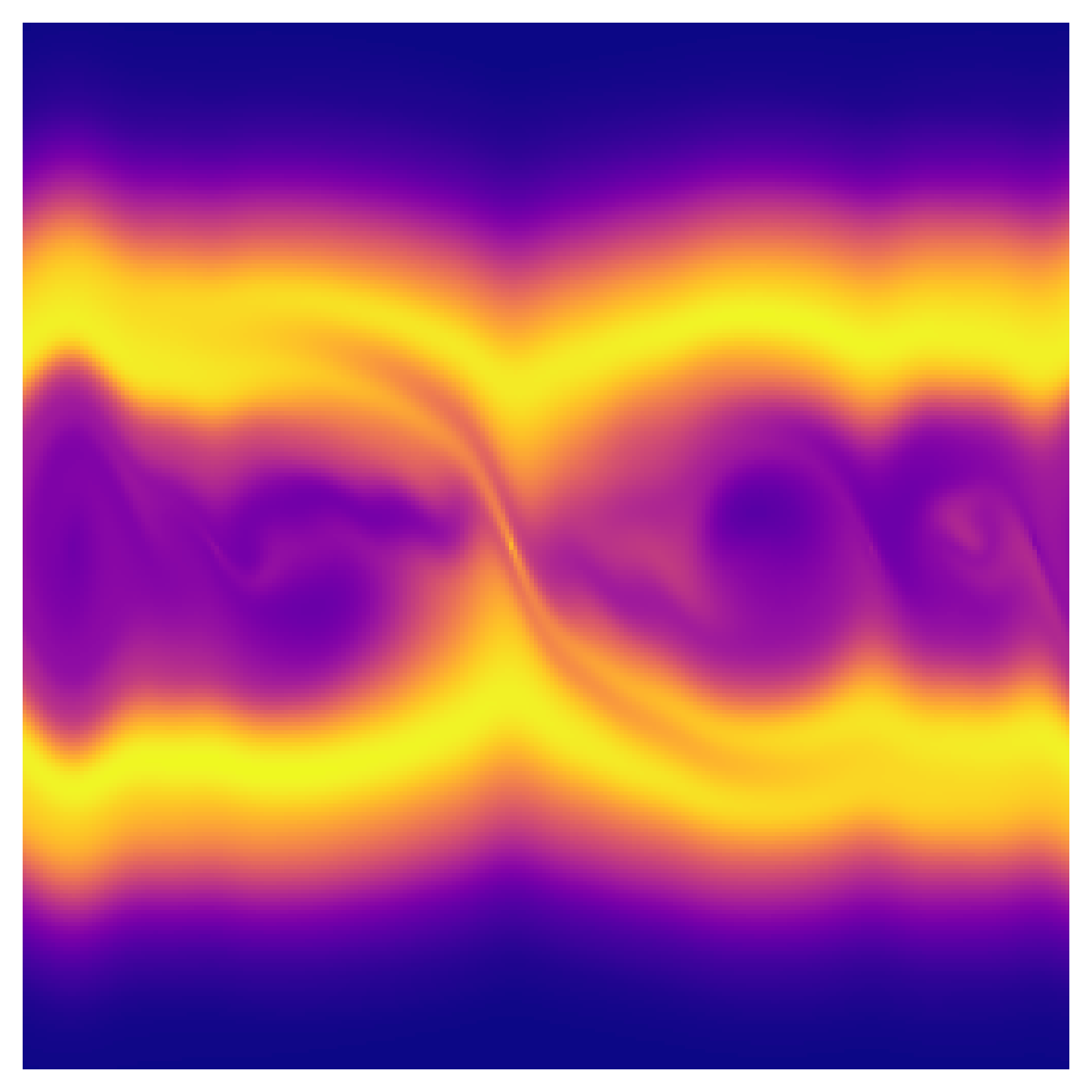} & \includegraphics[width=0.11\linewidth]{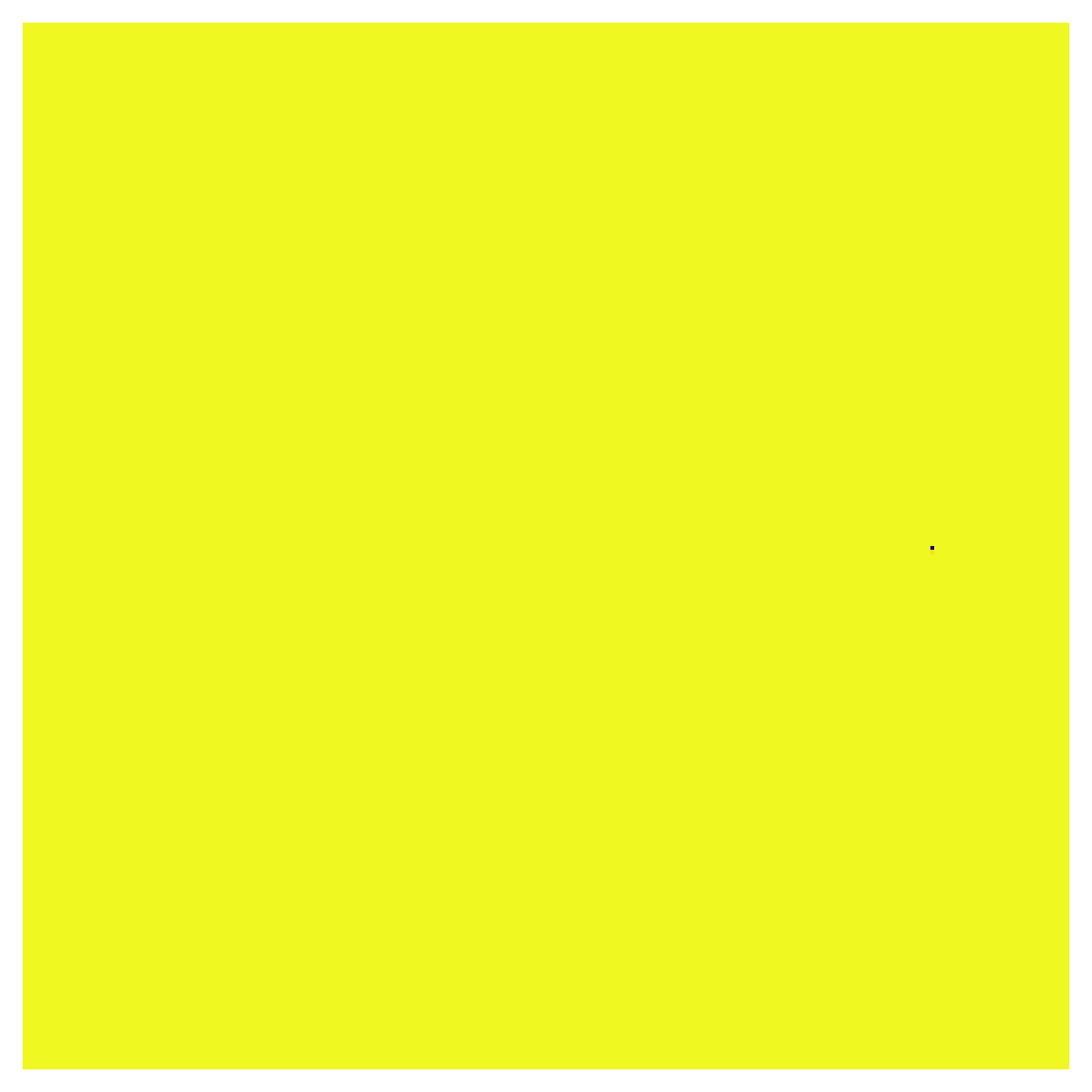} & \includegraphics[width=0.11\linewidth]{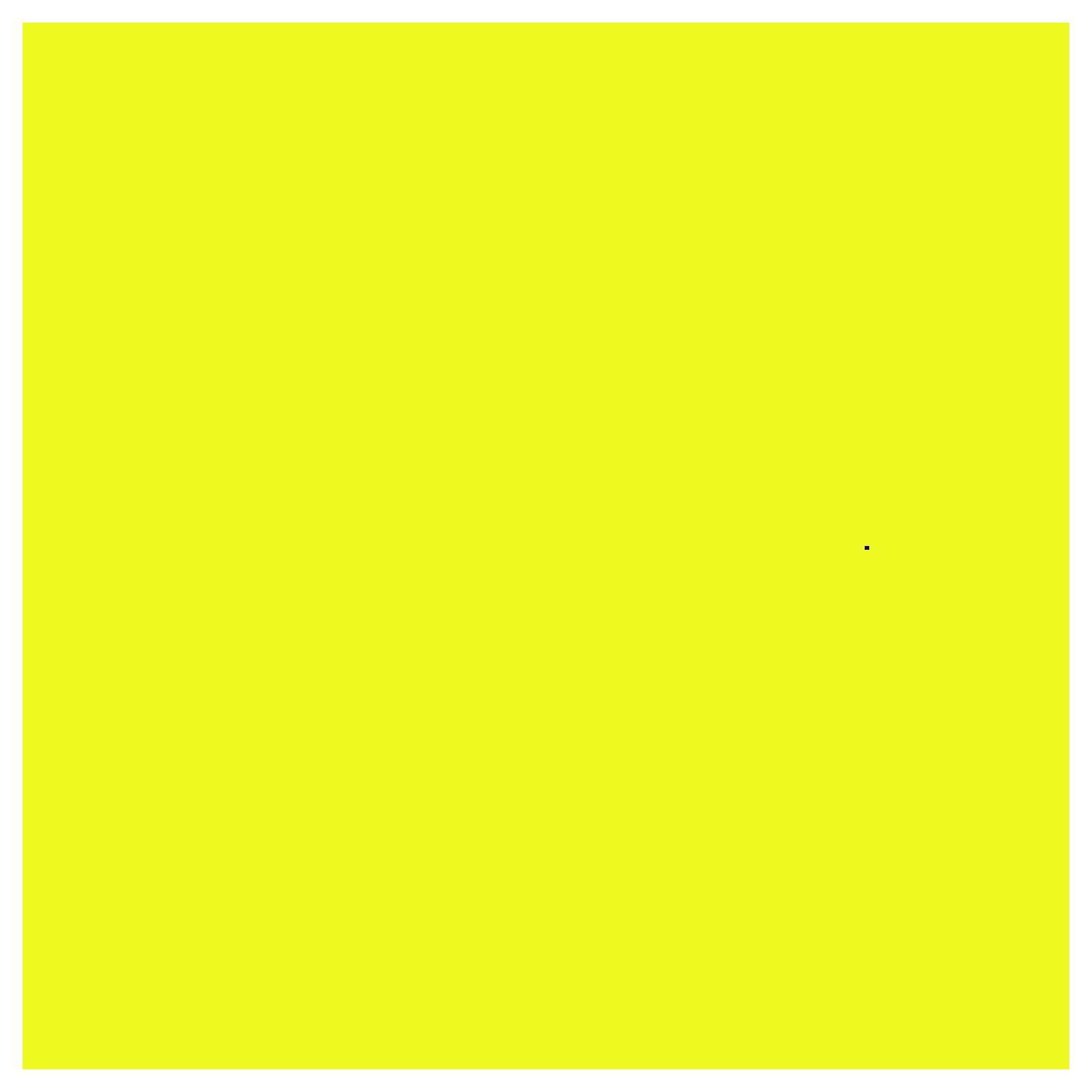} \\
    & & (Fig.~\ref{fig:TS_KL_GDL_local_under}) & (Fig.~\ref{fig:TS_ee_lf_GDL_local_under}) & 
    (Fig.~\ref{fig:TS_ee_GDL_local_under}) & (Fig.~\ref{fig:TS_KL_GDL_local_over}) & (Fig.~\ref{fig:TS_ee_lf_GDL_local_over}) & (Fig.~\ref{fig:TS_ee_GDL_local_over}) \\ \cline{2-8}
    & \multirow{2}{*}[22pt]{Local} & \includegraphics[width=0.11\linewidth]{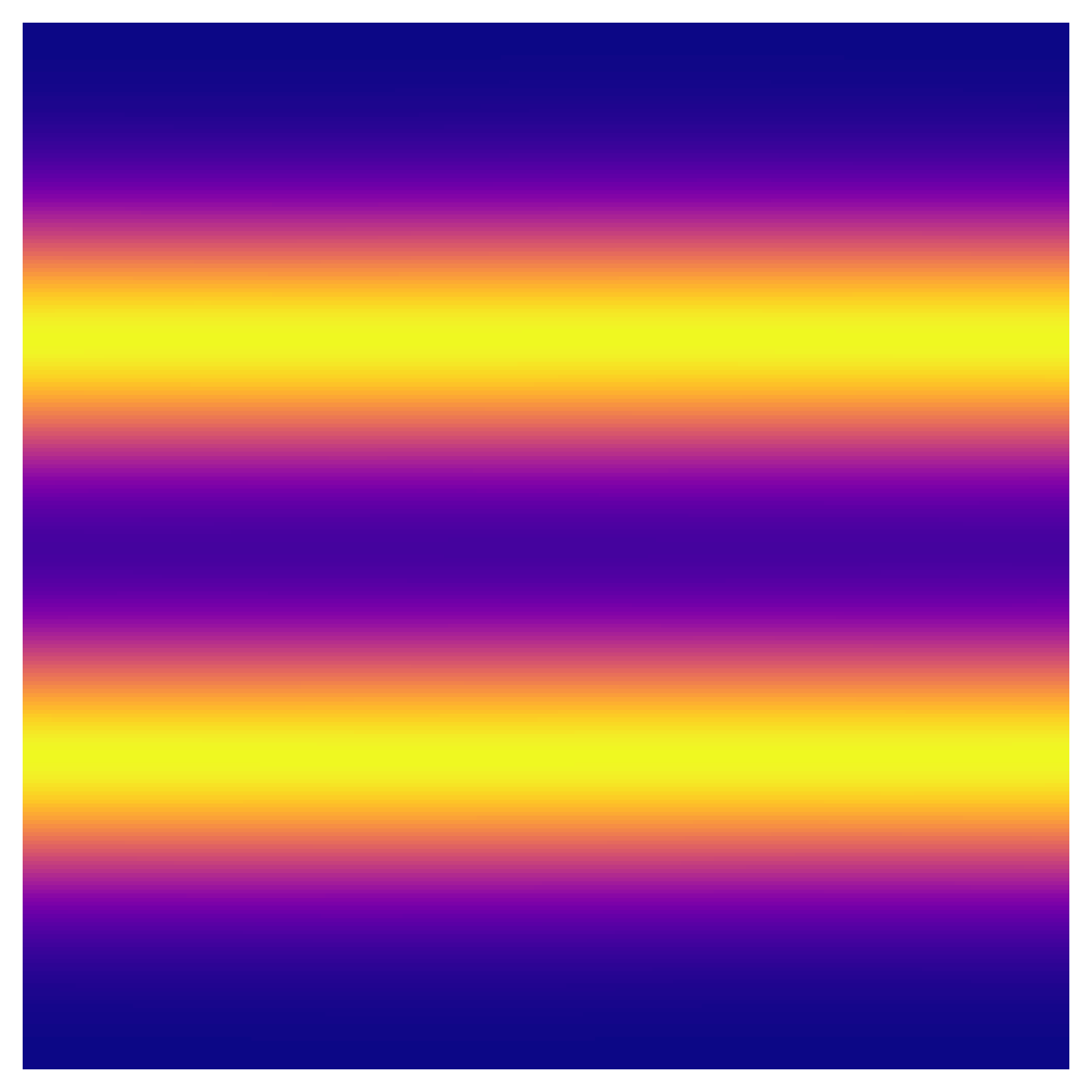} & \includegraphics[width=0.11\linewidth]{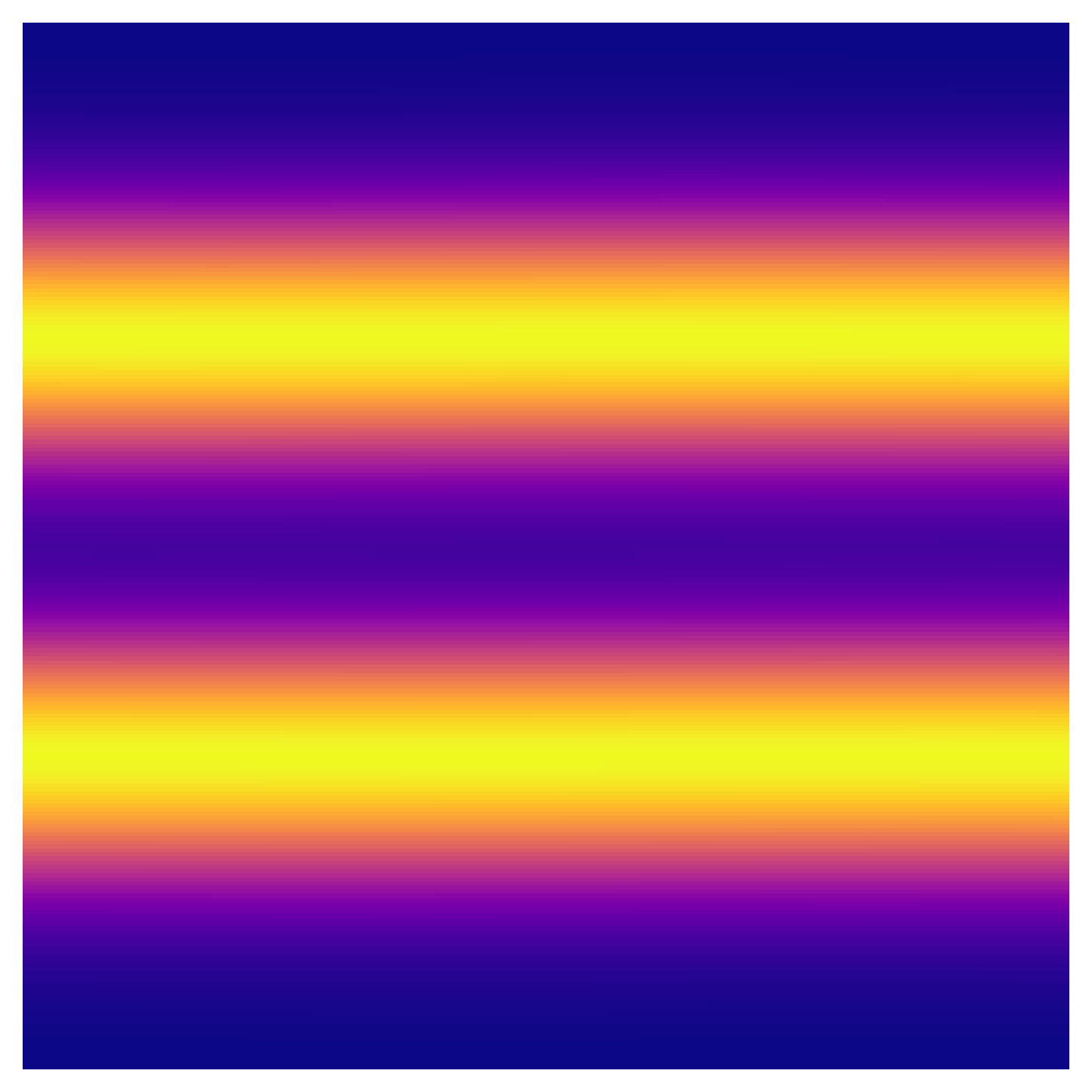} &
    \includegraphics[width=0.11\linewidth]{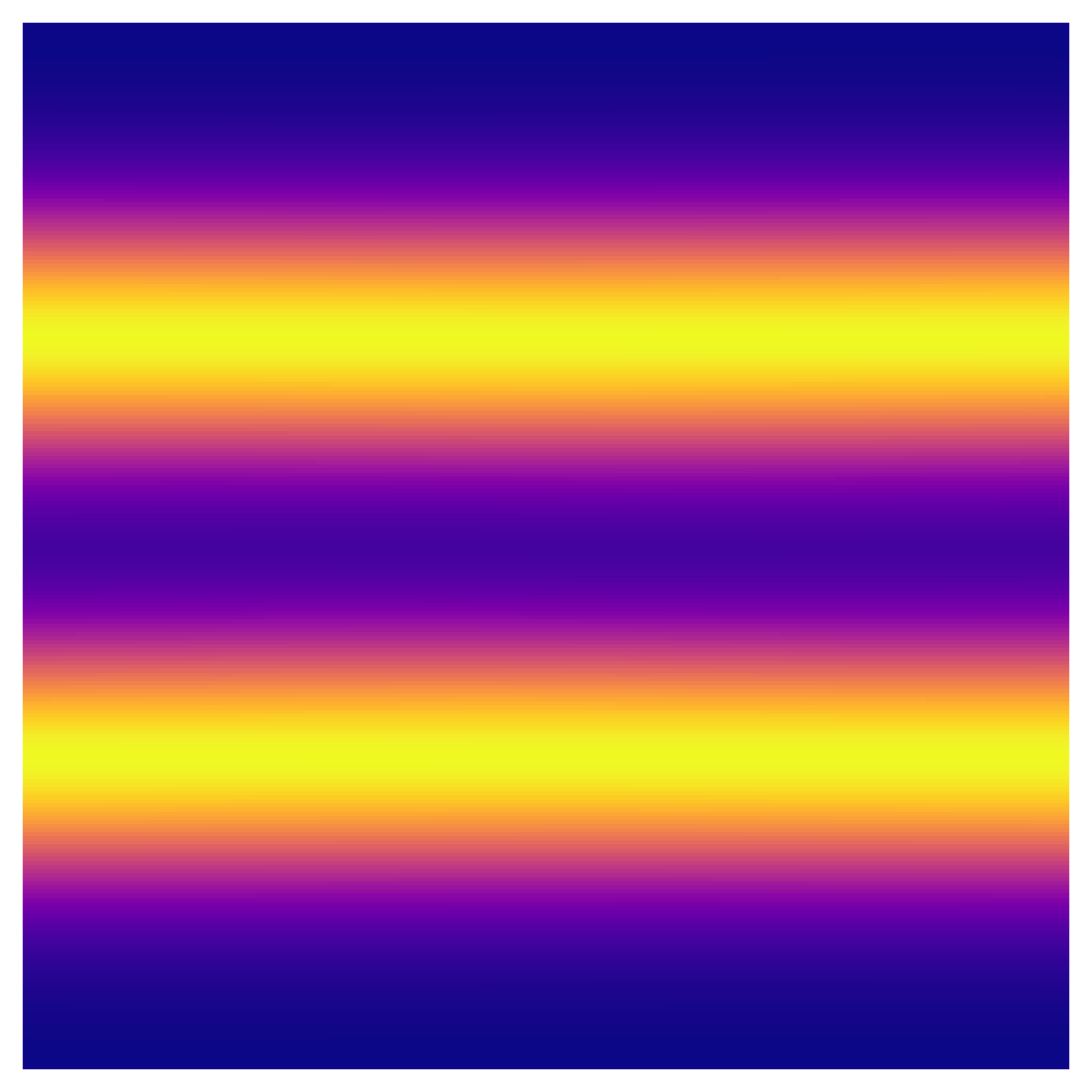} & \includegraphics[width=0.11\linewidth]{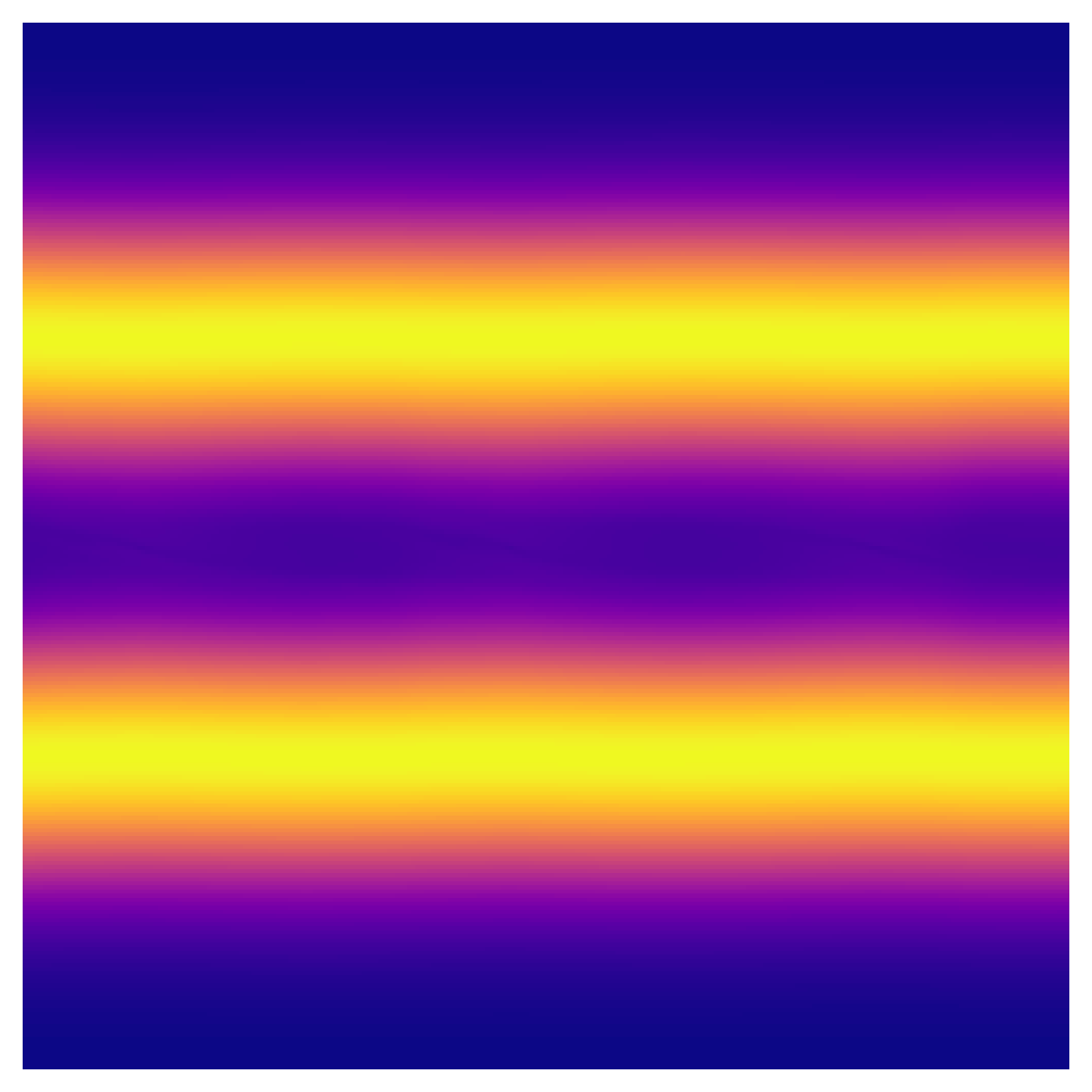} & \includegraphics[width=0.11\linewidth]{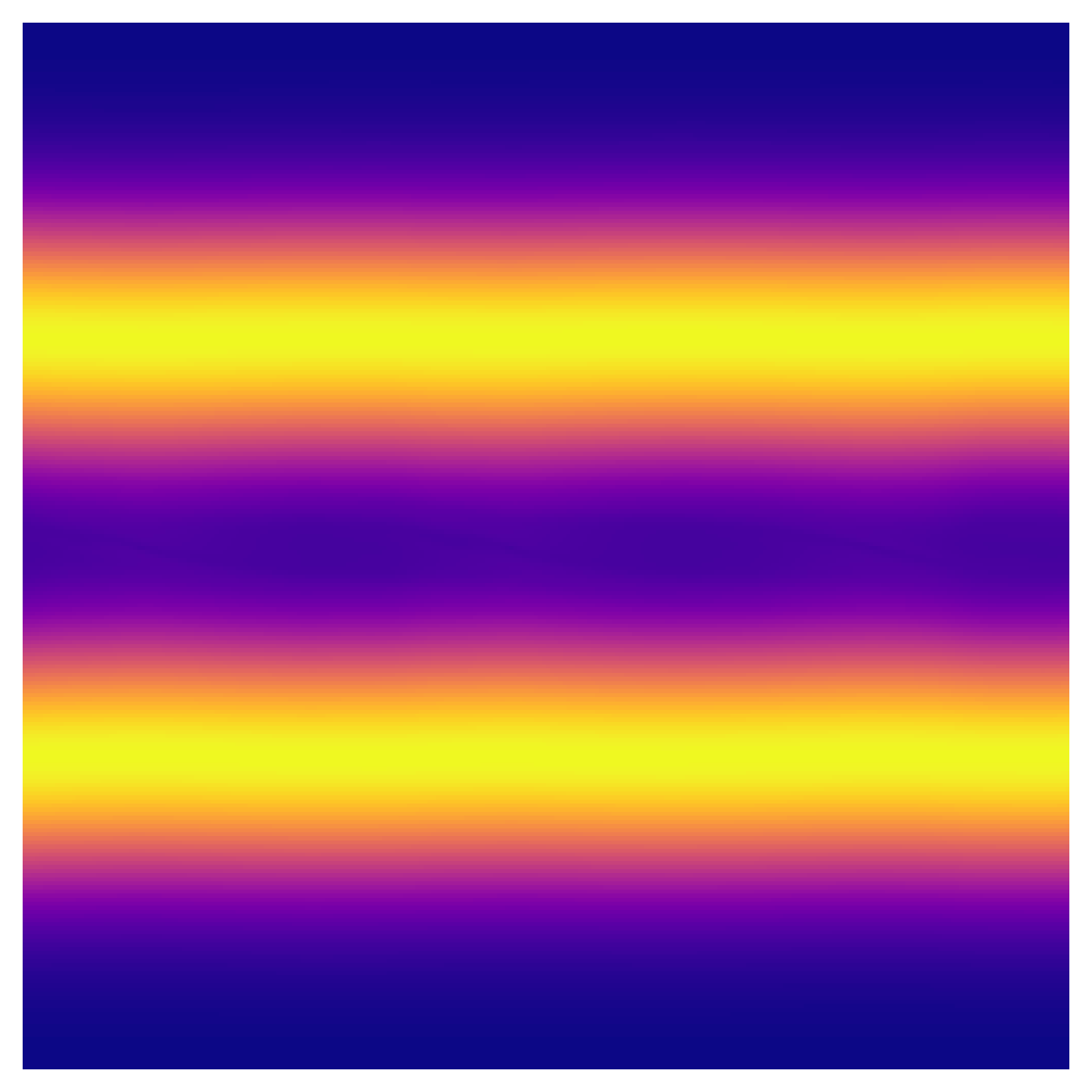} & \includegraphics[width=0.11\linewidth]{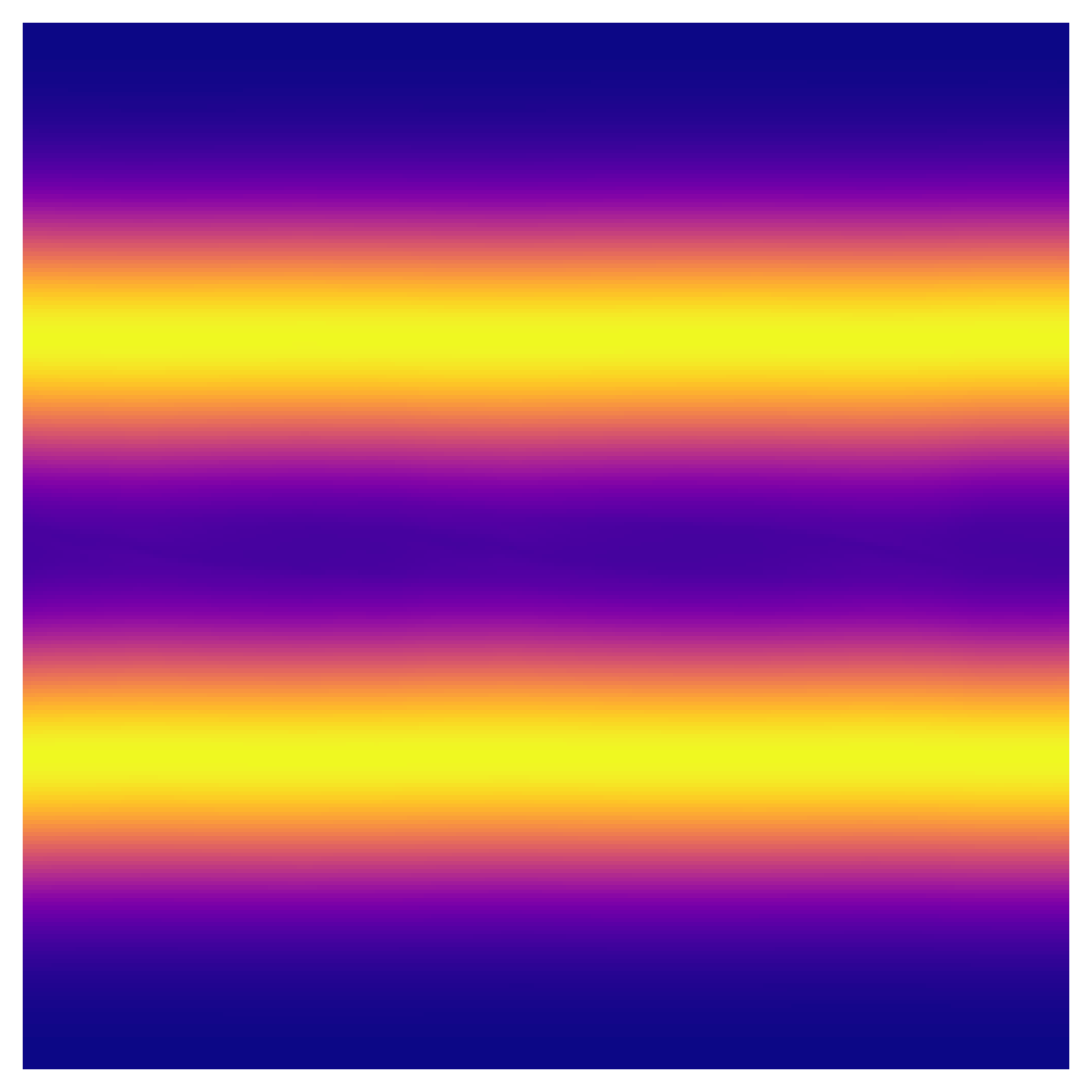} \\
    & & (Fig.~\ref{fig:TS_KL_GD_local_under}) & (Fig.~\ref{fig:TS_ee_lf_GD_local_under}) & (Fig.~\ref{fig:TS_ee_GD_local_under}) & (Fig.~\ref{fig:TS_KL_GD_local_over}) & (Fig.~\ref{fig:TS_ee_lf_GD_local_over}) & (Fig.~\ref{fig:TS_ee_GD_local_over}) \\
    \bottomrule
    \end{tabular}
    \caption{Summary of results for the Two Stream equilibrium.}
    \label{tab:TS_results}
\end{table}

Multiple conclusions can be drawn.
\begin{itemize}
    \item Local solver with near-optimal initialization (good initial guess) perform universally well, and the stability is achieved in both under-parameter regime, and over-parameter regime, as seen in Figures \ref{fig:TS_KL_GD_local_under}, \ref{fig:TS_ee_lf_GD_local_under}, \ref{fig:TS_ee_GD_local_under}, \ref{fig:TS_KL_GD_local_over}, \ref{fig:TS_ee_lf_GD_local_over} and \ref{fig:TS_ee_GD_local_over}. This signifies the importance of making a good initial guess. Making a good initial guess is usually not a feasible task for a classical PDE-constrained optimization, but is possible in this context by deploying the study of the dispersion relation, see discussion in Section~\ref{sec:dispersion}.
    \item Adaptive solver is not compatible with the objective functions~\eqref{eq:EE_obj} and~\eqref{eq:EET_obj}, as seen in Figures \ref{fig:TS_ee_GDL_far_under}, \ref{fig:TS_ee_GDL_far_over}, \ref{fig:TS_ee_lf_GDL_near_under}, \ref{fig:TS_ee_GDL_near_under}, \ref{fig:TS_ee_lf_GDL_near_over}, \ref{fig:TS_ee_GDL_near_over}, \ref{fig:TS_ee_lf_GDL_local_over} and \ref{fig:TS_ee_GDL_local_over}. This is already suggested by Figure~\ref{fig:landscape_two_stream_2D_varepsilon_0.001},~\eqref{eq:EET_obj} where non-physical local minima can be found when parameters take on very large values. When adaptive solvers such as GD with line-search are deployed, the searching can take the iteration very far into these large parameter regions. The found external field is extremely strong and it completely dominates the whole plasma behavior. When it happens, the self-generated electric field is indeed very small, but the solution renders meaningless physics. The only exception in this case is to set the initial data from a near-optimal initialization. The basin is convex enough to control the line-search's performance, see Figures~\ref{fig:TS_ee_GDL_local_under} and \ref{fig:TS_ee_lf_GDL_local_under}.
    \item The KL divergence does not seem to be a good objective functional for suppressing instability. The solution with good suppression is only observed when initialization is near-optimal already (see Figures~\ref{fig:TS_KL_GDL_local_under},~\ref{fig:TS_KL_GD_local_under} and~\ref{fig:TS_KL_GD_local_over}). With neither adaptive or local solver, in neither under-parameterized or over-parameterized regime, KL divergence returns a satisfactory solution as long as initialization is not near-optimal.
    \item Counter-intuitively, over-parametrization in this setting does not seem to bring too much benefit (comparing the last three columns with the first three columns), with the only exception being the use of an adaptive solver for minimizing KL with a far-away initialization (comparing Figure~\ref{fig:TS_KL_GDL_far_over} with Figure~\ref{fig:TS_KL_GDL_far_under}).
    \item Without good initial guess, the only good solution is found by optimizing~\eqref{eq:EET_obj} with simple GD solver in the properly parametrized setting when initialization is prepared in the mid-range (see Figure~\ref{fig:TS_ee_GD_near_under}).
\end{itemize}

\subsection{Bump-on-Tail example}
For the Bump-on-Tail example, we conduct the same experiments and summarize them in Table~\ref{tab:BoT_results} and~\ref{sec:BoT_summary}. Most observations drawn from the previous example still hold true. Yet we should identify the distinct features. In particular, our simulation seems to suggest that Bump-on-Tail instability is easier to be controlled, though the perfect recovery is more rare. More specifically, for Two Stream instability, when we have good initial guess, we obtained exact recovery when we use adaptive solvers in the under-parameterized regime and local solvers in the over-parameterized regime. We fail to do so in the Bump-on-Tail situation except for~\eqref{eq:EET_obj}. On the other hand, we find many more situations in which, though we are not being able to find exact recovery, the solutions are still relatively stable. This strongly suggests the local minima are good controls nevertheless, see trajectory data in Figures \ref{fig:BoT_KL_GDL_far_under}, \ref{fig:BoT_KL_GD_far_under}, \ref{fig:BoT_ee_GD_far_under}, \ref{fig:BoT_ee_lf_GD_far_under}, \ref{fig:BoT_KL_GD_local_under}, \ref{fig:BoT_ee_GD_local_under} and \ref{fig:BoT_ee_lf_GD_local_under}.

Finding good initialization is still instrumental and it brings out-performance across solvers and regimes. Adaptive solvers are still incompatible with~\eqref{eq:EE_obj} and~\eqref{eq:EET_obj}, and the performance is universally bad. Over-parameterization still does not necessarily bring benefits. In particular, with good initial guess and local solver, looking for best control in the over-parameterized regime brings even worse reconstruction, see Figures~\ref{fig:BoT_KL_GD_local_over} and \ref{fig:BoT_ee_GD_local_over} vs. Figures~\ref{fig:BoT_KL_GD_local_under} and \ref{fig:BoT_ee_GD_local_under}.

However, in comparison to the Two Stream case, KL divergence seems to be a reasonable objective function, by minimizing which, acceptable controls for Bump-on-Tail can still be found. This resonates our earlier discussion on local minima serving as good controls for Bump-on-Tail instability.

\begin{table}[!ht]
    \centering
    \begin{tabular}{c|c|c|c|c|c|c|c}
    \toprule
    & & \multicolumn{3}{|c|}{Under-parametrized~\eqref{eq:H_cos_sin_under}} & \multicolumn{3}{|c}{Over-parametrized~\eqref{eq:H_cos_sin_over}} \\
    \midrule
    Init. & Step & \eqref{eq:KL_obj} & \eqref{eq:EE_obj} & \eqref{eq:EET_obj} & \eqref{eq:KL_obj} & \eqref{eq:EE_obj} & \eqref{eq:EET_obj}  \\
    \midrule
    \multirow{4}{*}[0pt]{Far} & \multirow{2}{*}[22pt]{Adaptive} & \includegraphics[width=0.11\linewidth]{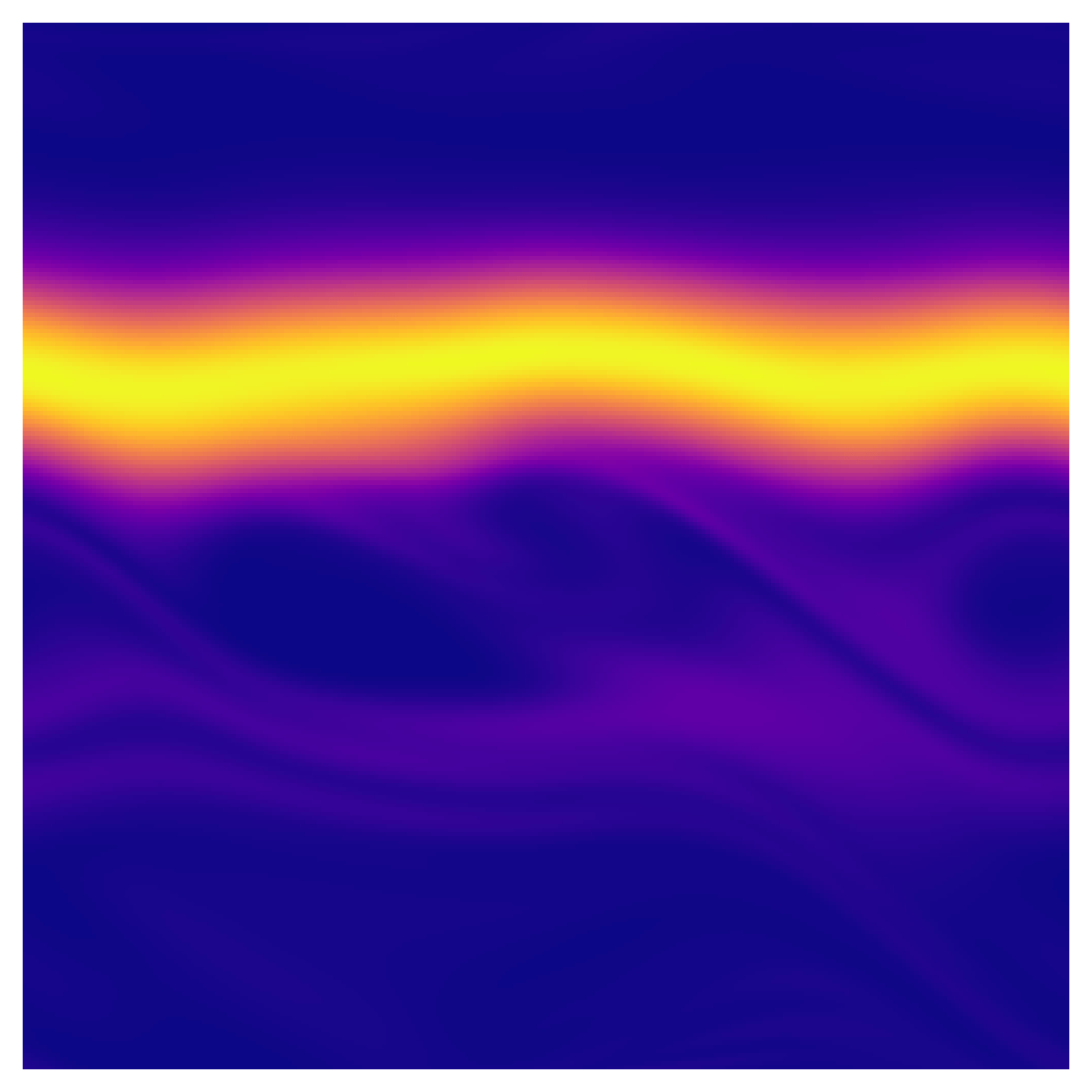} & \includegraphics[width=0.11\linewidth]{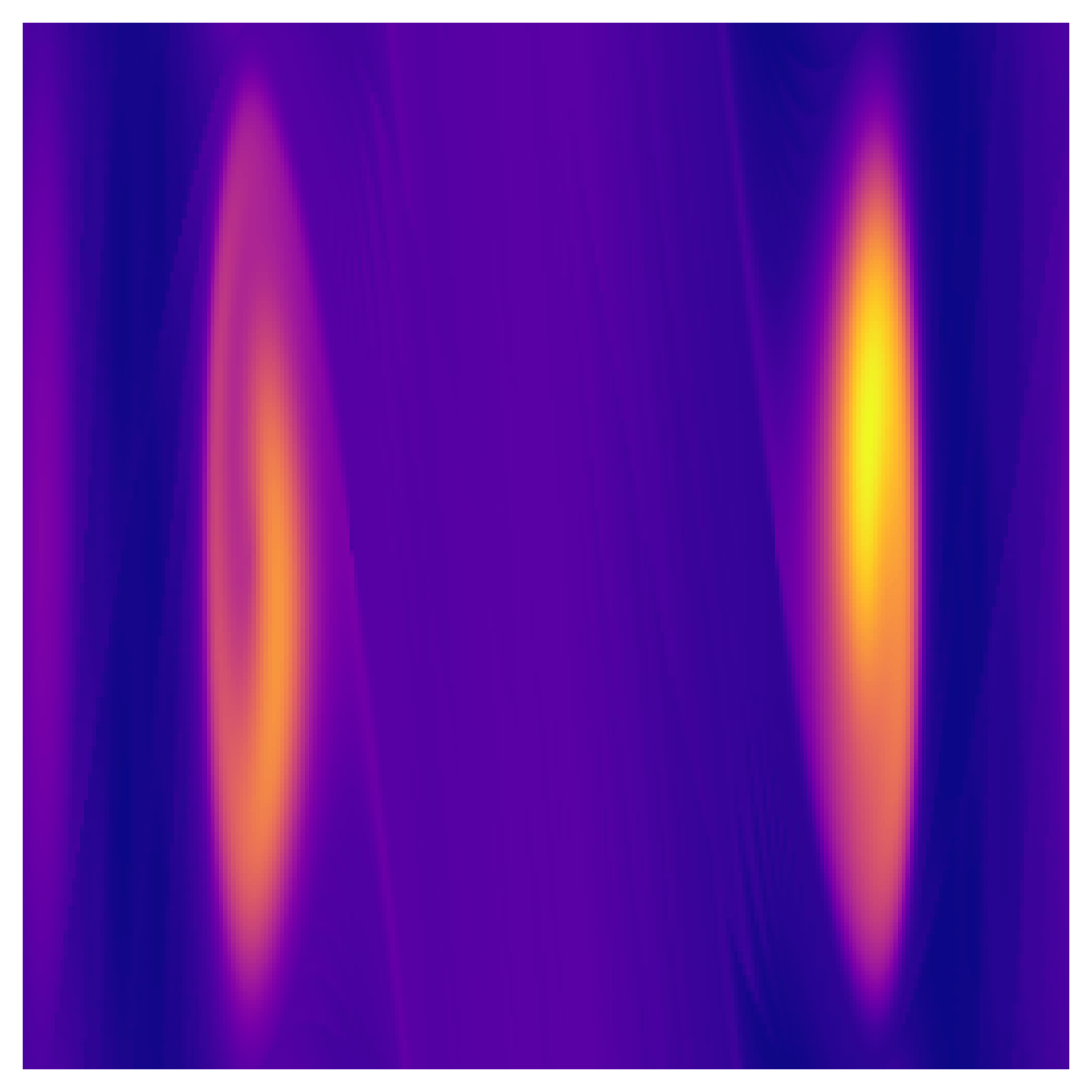} & \includegraphics[width=0.11\linewidth]{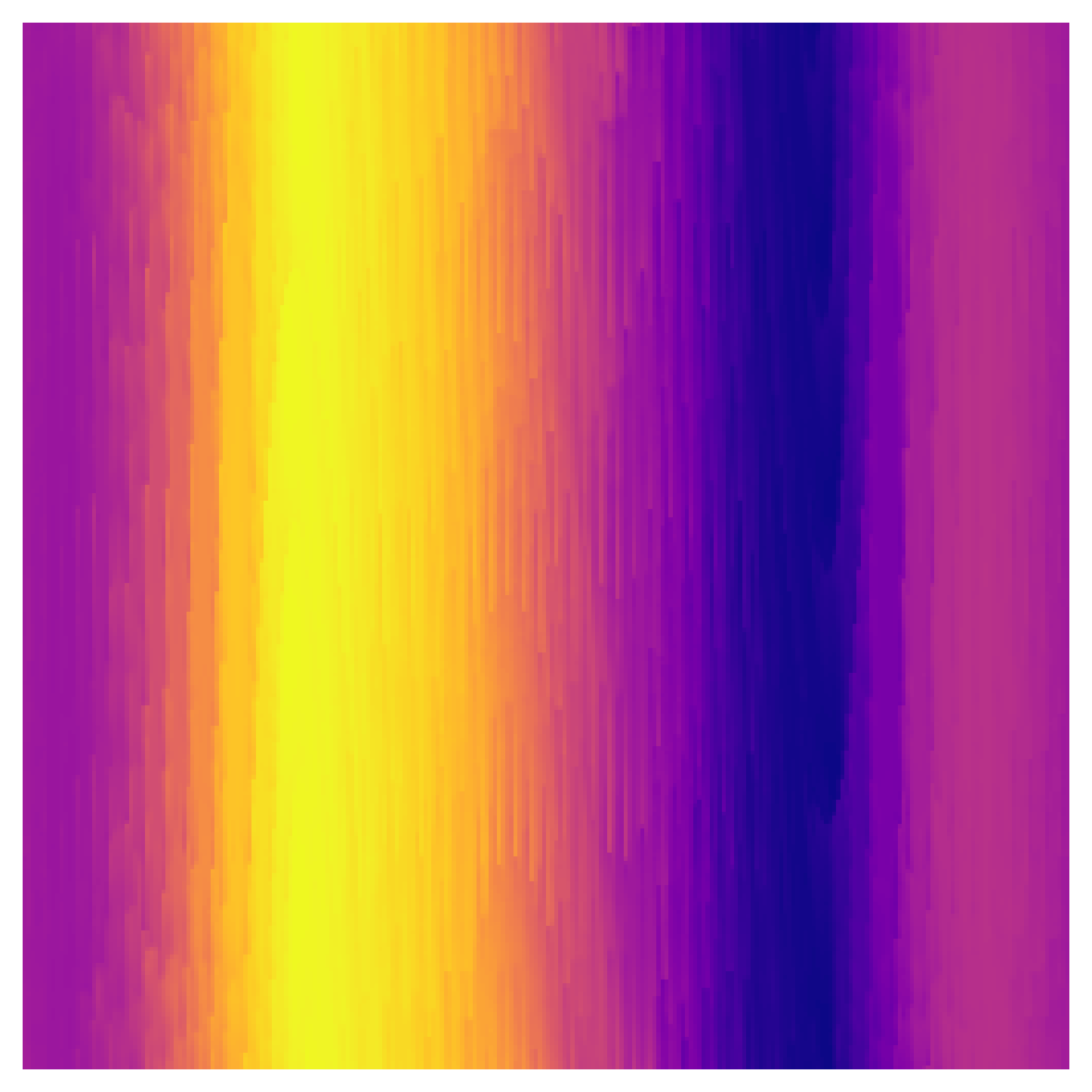} & \includegraphics[width=0.11\linewidth]{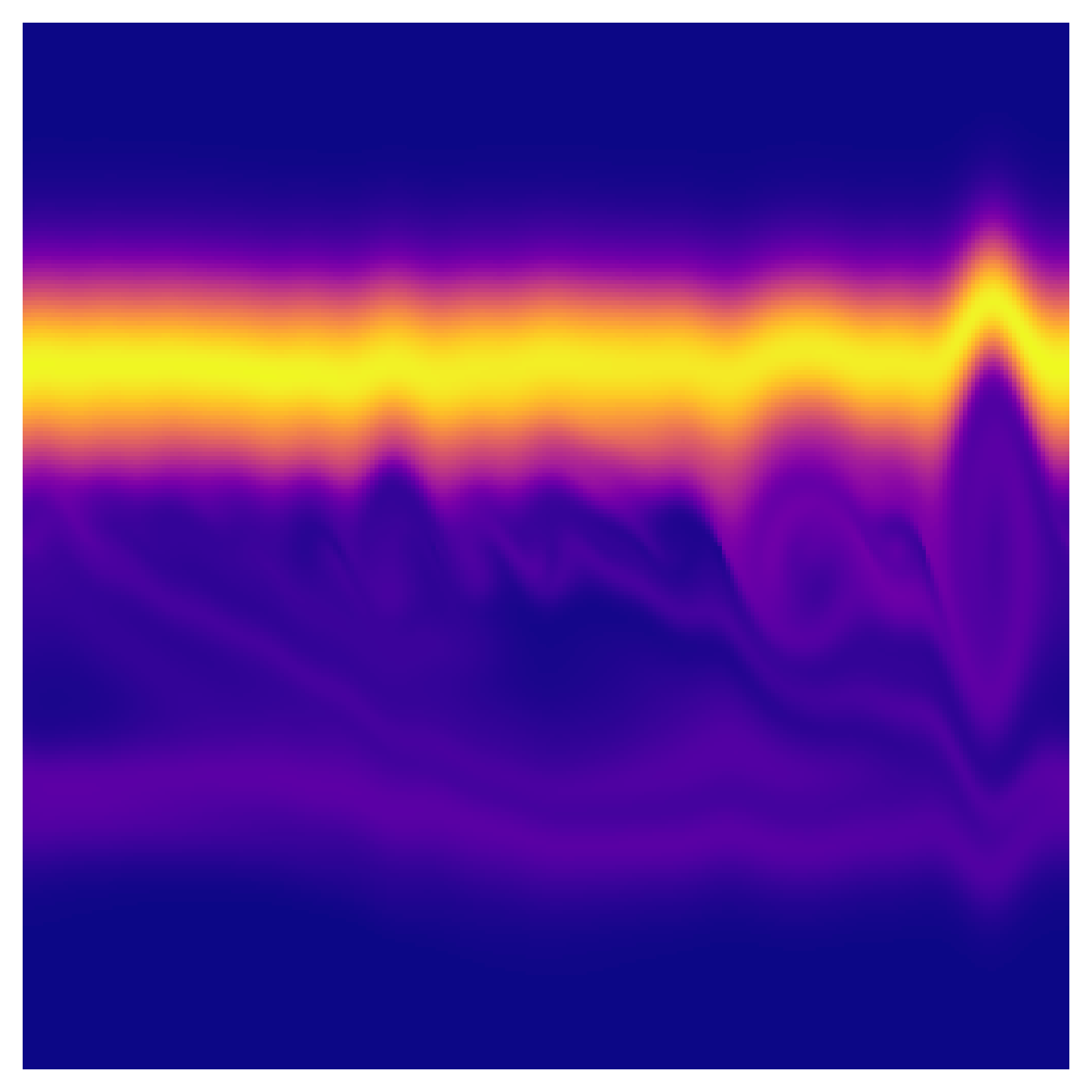} & \includegraphics[width=0.11\linewidth]{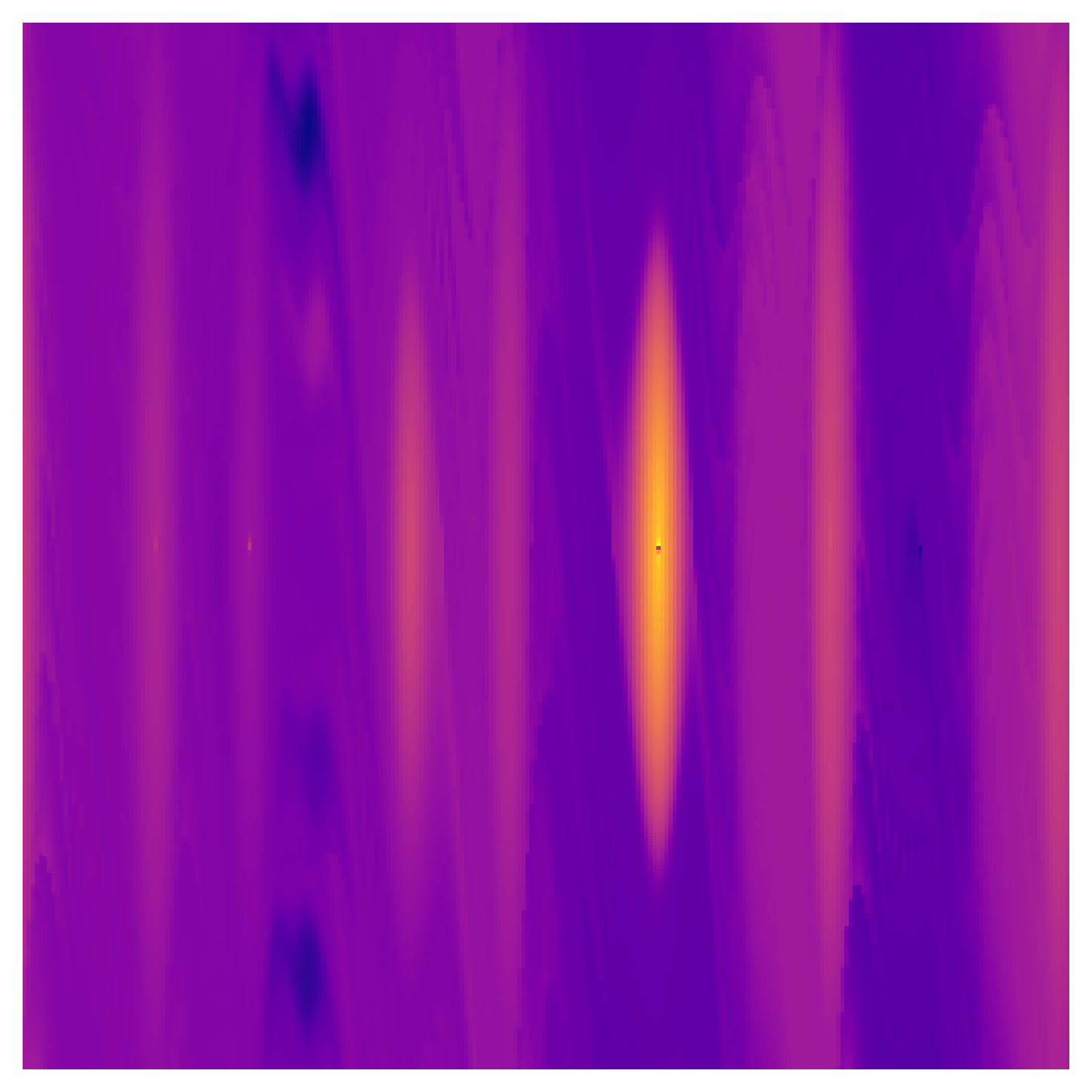} & \includegraphics[width=0.11\linewidth]{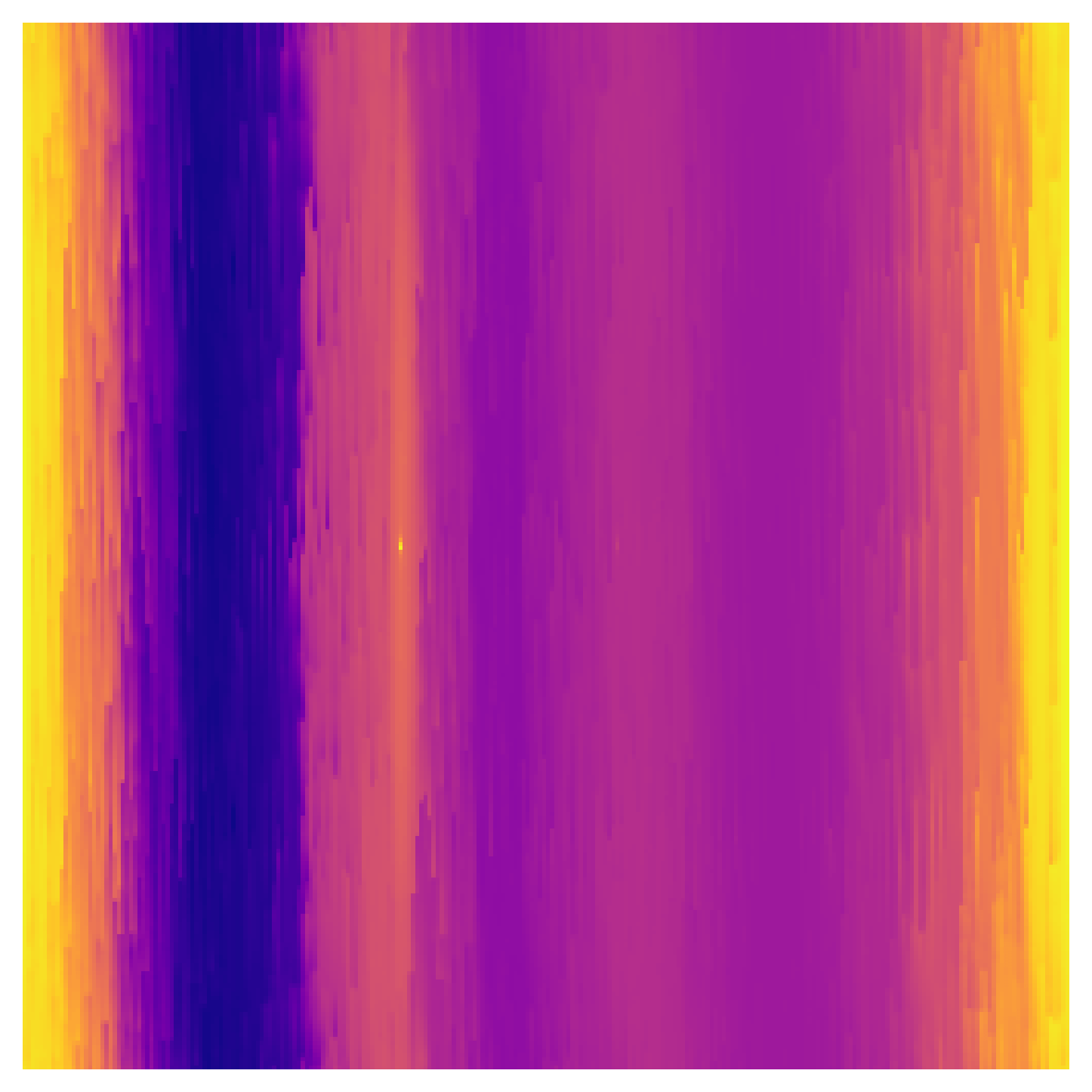} \\ 
    & & (Fig.~\ref{fig:BoT_KL_GDL_far_under}) & (Fig.~\ref{fig:BoT_ee_lf_GDL_far_under}) & (Fig.~\ref{fig:BoT_ee_GDL_far_under}) &(Fig.~\ref{fig:BoT_KL_GDL_far_over}) & (Fig.~\ref{fig:BoT_ee_lf_GDL_far_over}) & (Fig.~\ref{fig:BoT_ee_GDL_far_over}) \\ \cline{2-8}
    & \multirow{2}{*}[22pt]{Local} & \includegraphics[width=0.11\linewidth]{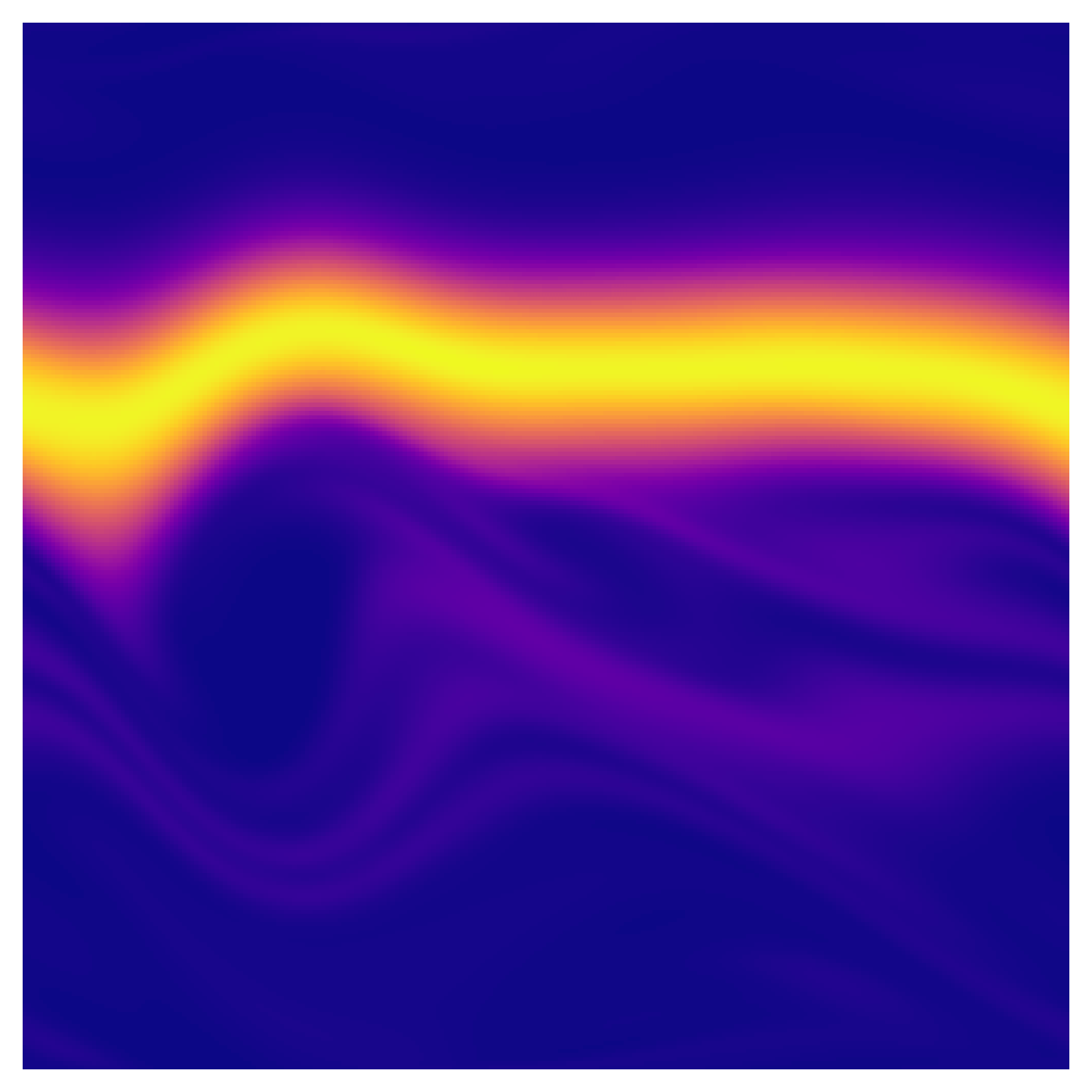} & \includegraphics[width=0.11\linewidth]{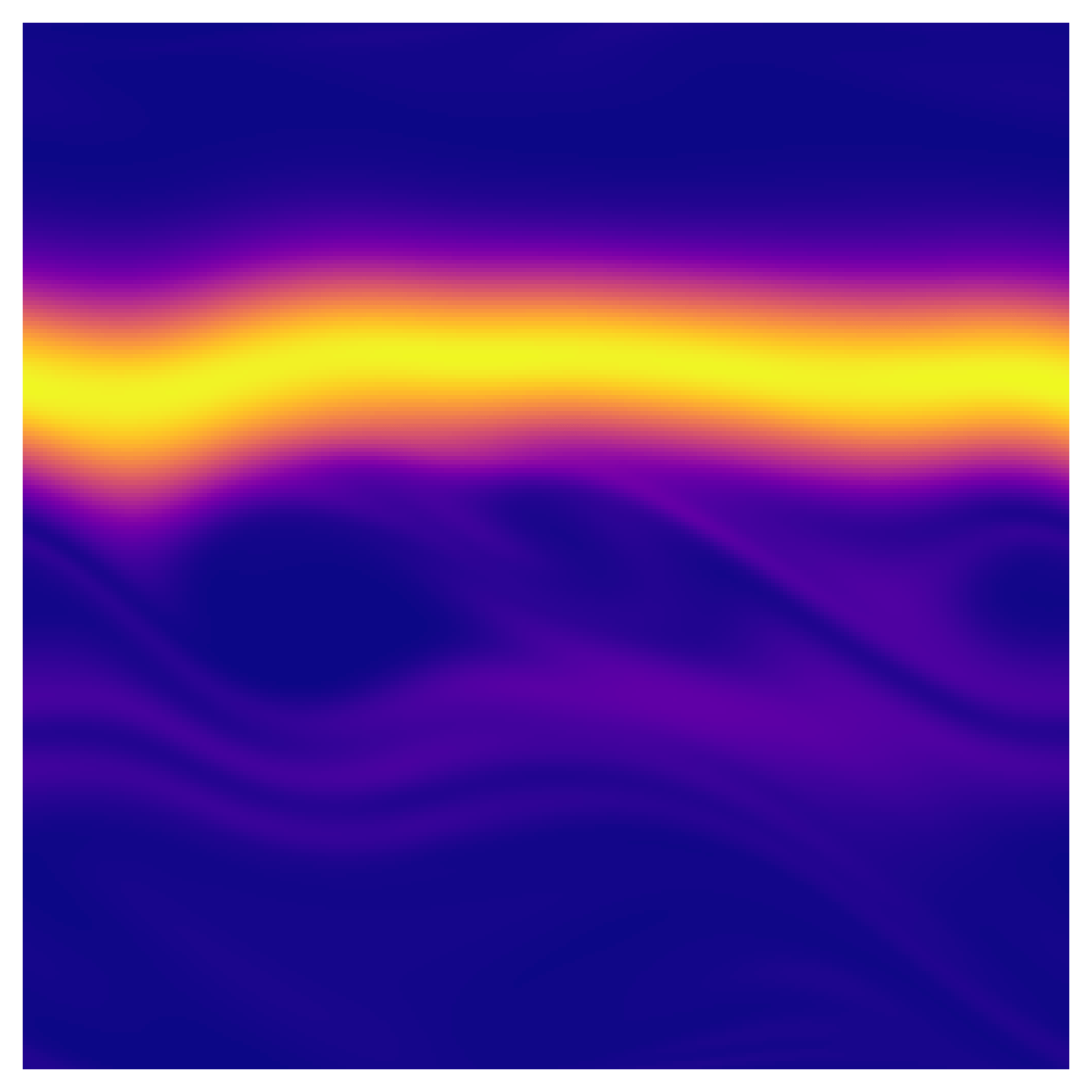} & \includegraphics[width=0.11\linewidth]{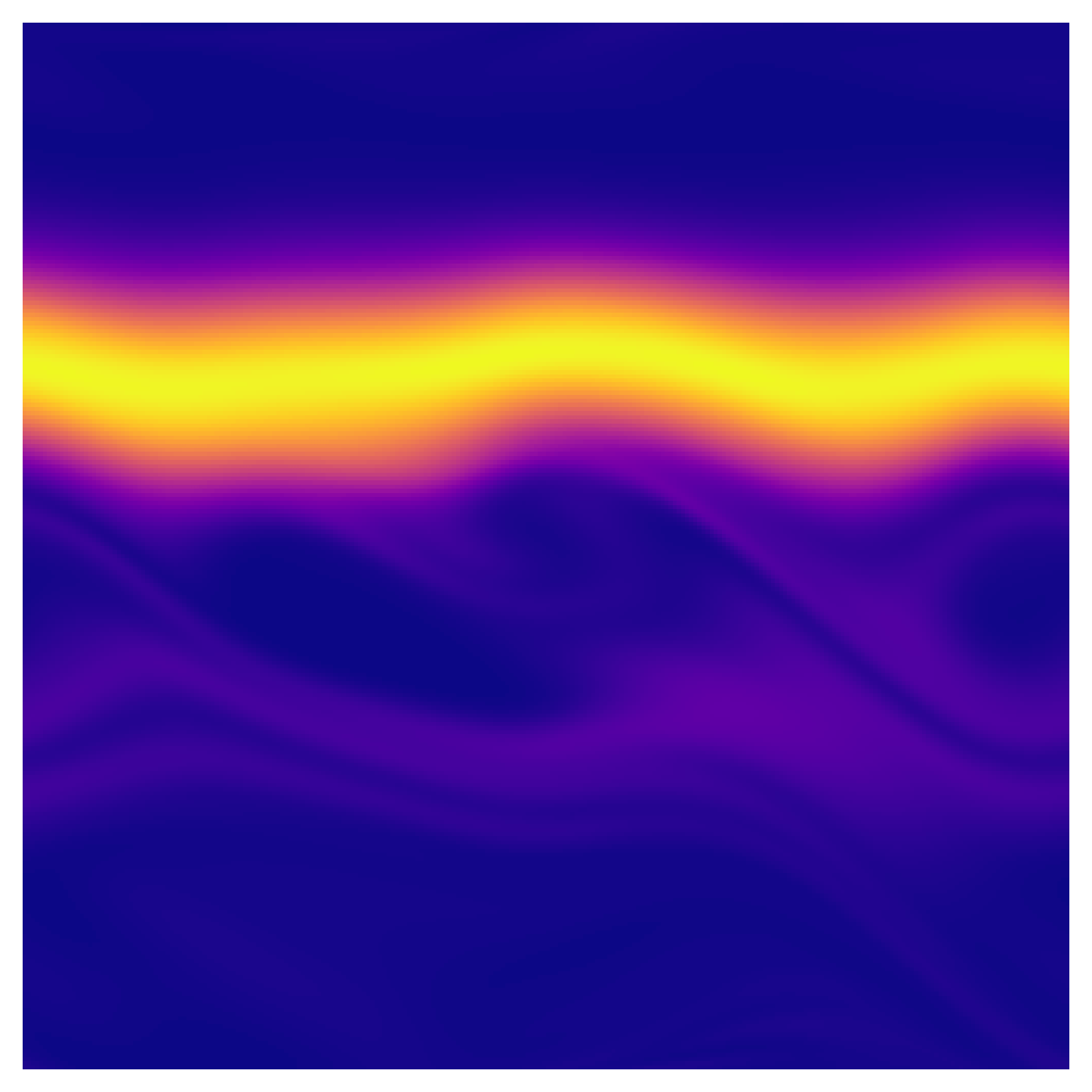} & \includegraphics[width=0.11\linewidth]{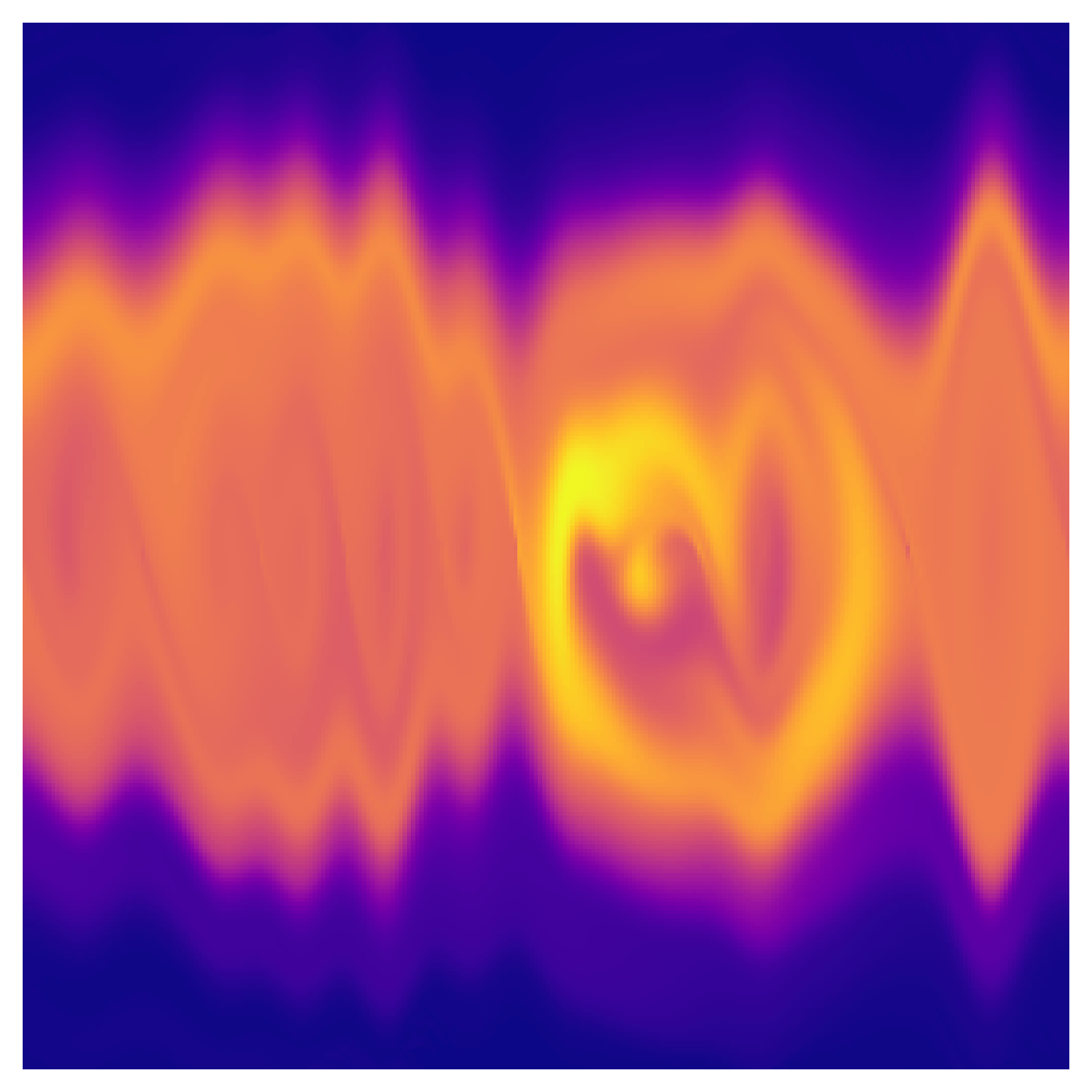} & \includegraphics[width=0.11\linewidth]{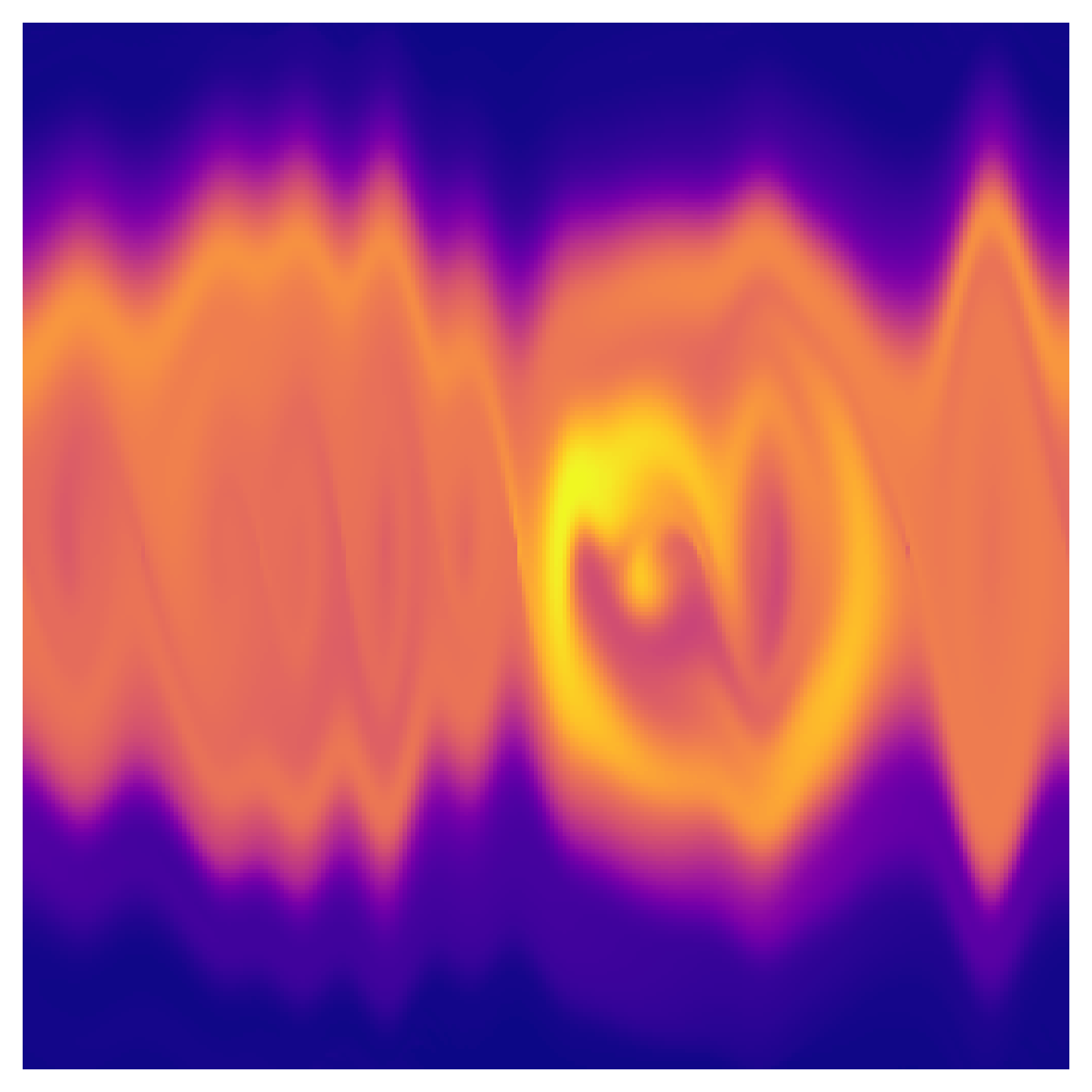} & \includegraphics[width=0.11\linewidth]{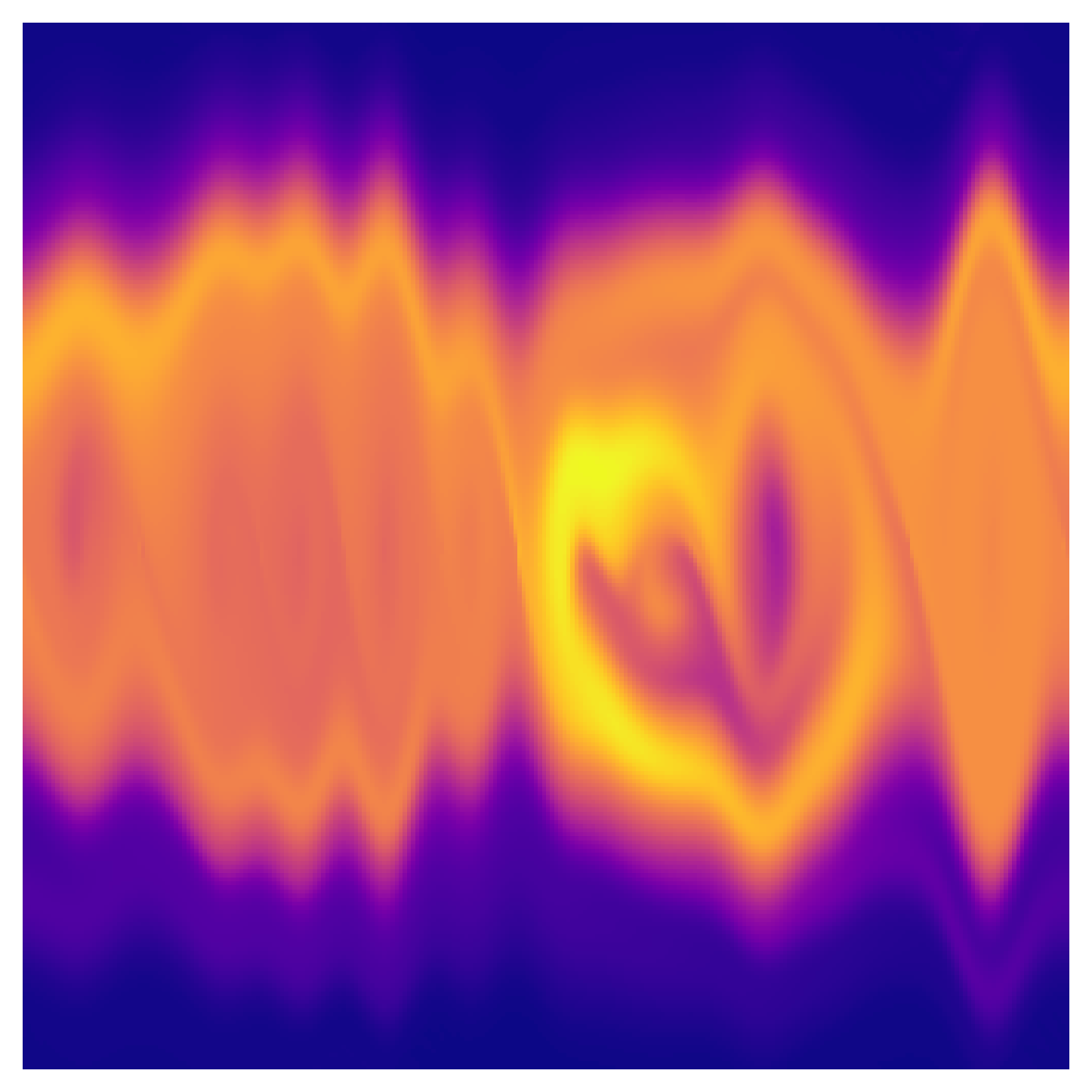} \\
    & & (Fig.~\ref{fig:BoT_KL_GD_far_under}) & (Fig.~\ref{fig:BoT_ee_lf_GD_far_under}) & (Fig.~\ref{fig:BoT_ee_GD_far_under}) & (Fig.~\ref{fig:BoT_KL_GD_far_over}) & (Fig.~\ref{fig:BoT_ee_lf_GD_far_over}) & (Fig.~\ref{fig:BoT_ee_GD_far_over}) \\
    \midrule
    \multirow{4}{*}[0pt]{Mid} & \multirow{2}{*}[22pt]{Adaptive} & \includegraphics[width=0.11\linewidth]{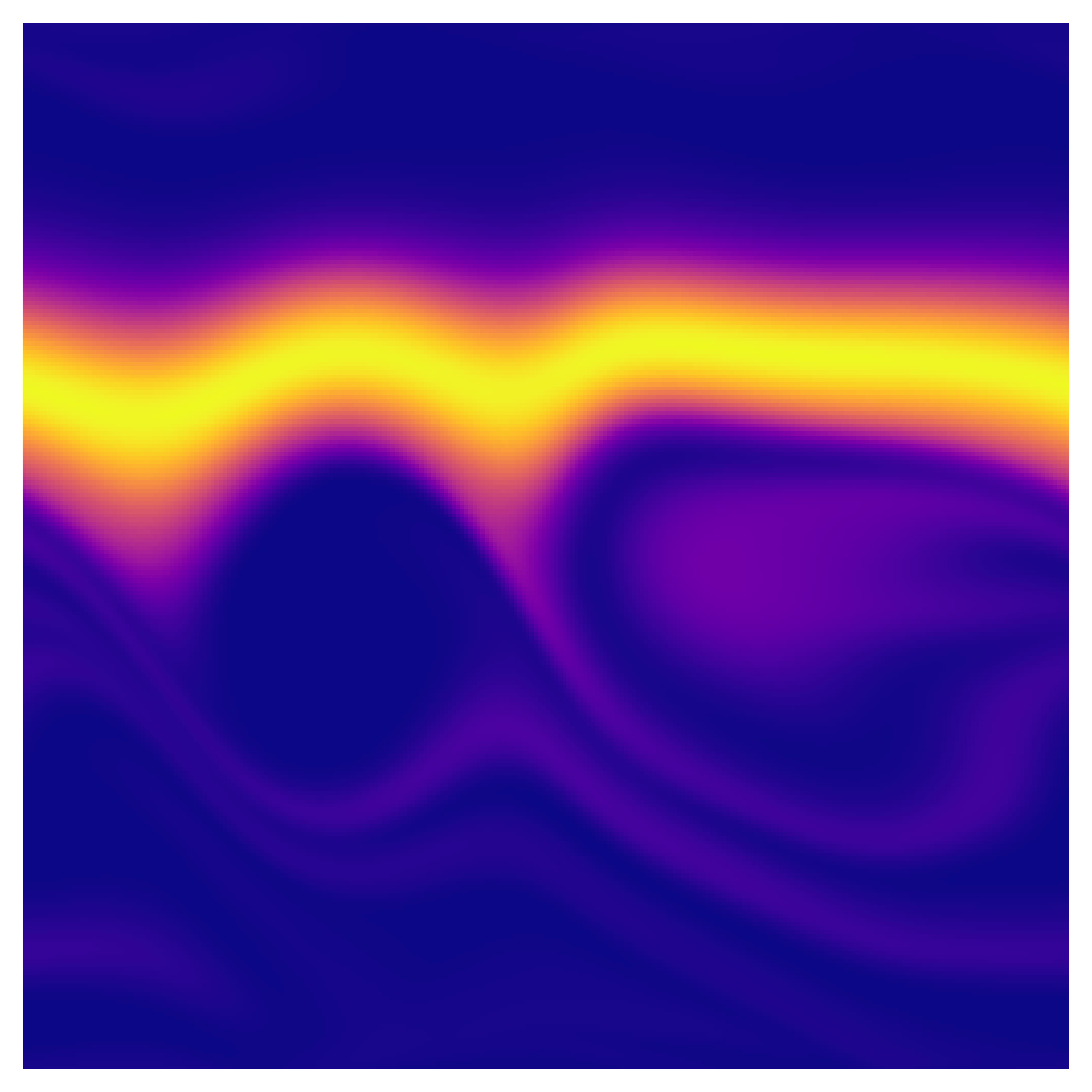} & \includegraphics[width=0.11\linewidth]{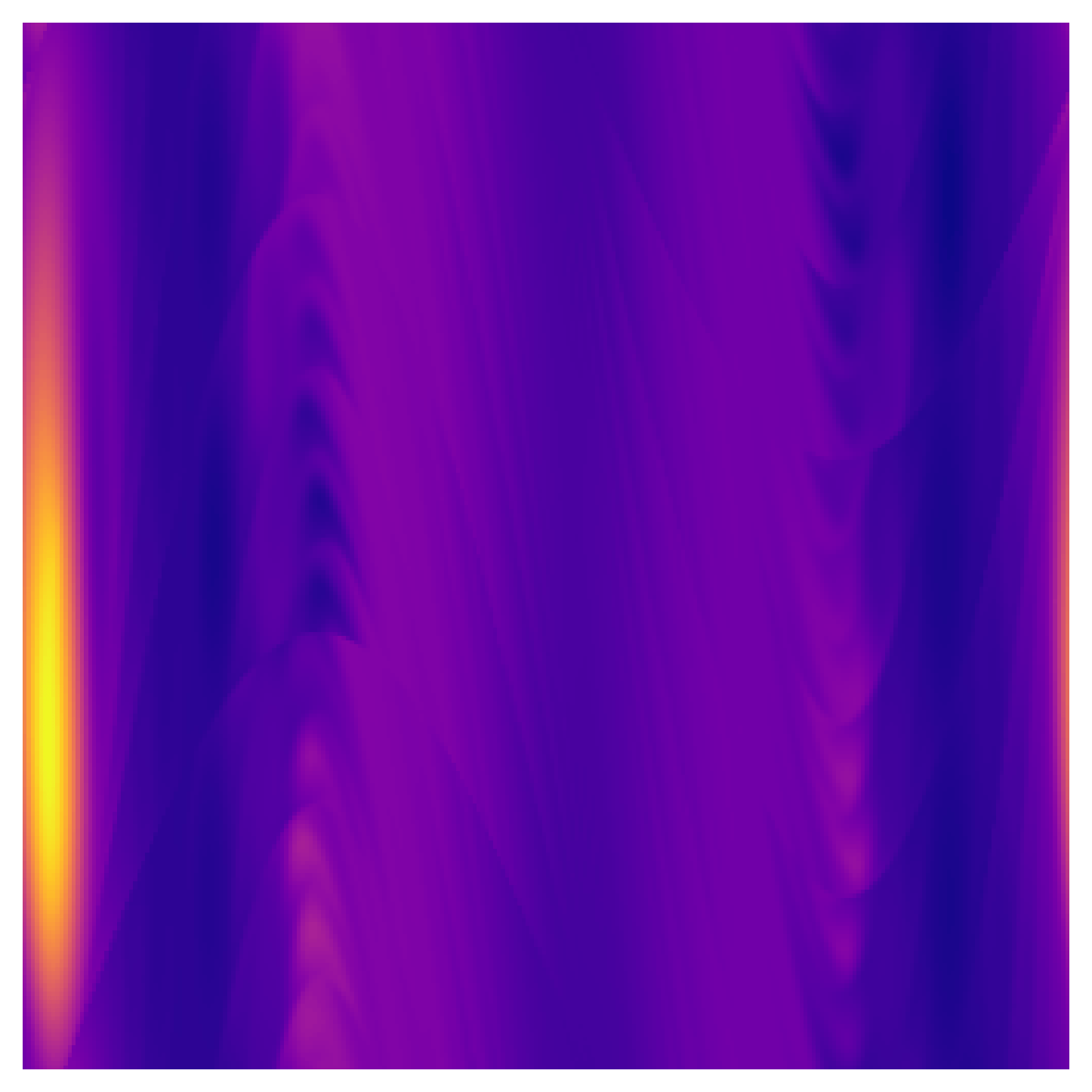} &
    \includegraphics[width=0.11\linewidth]{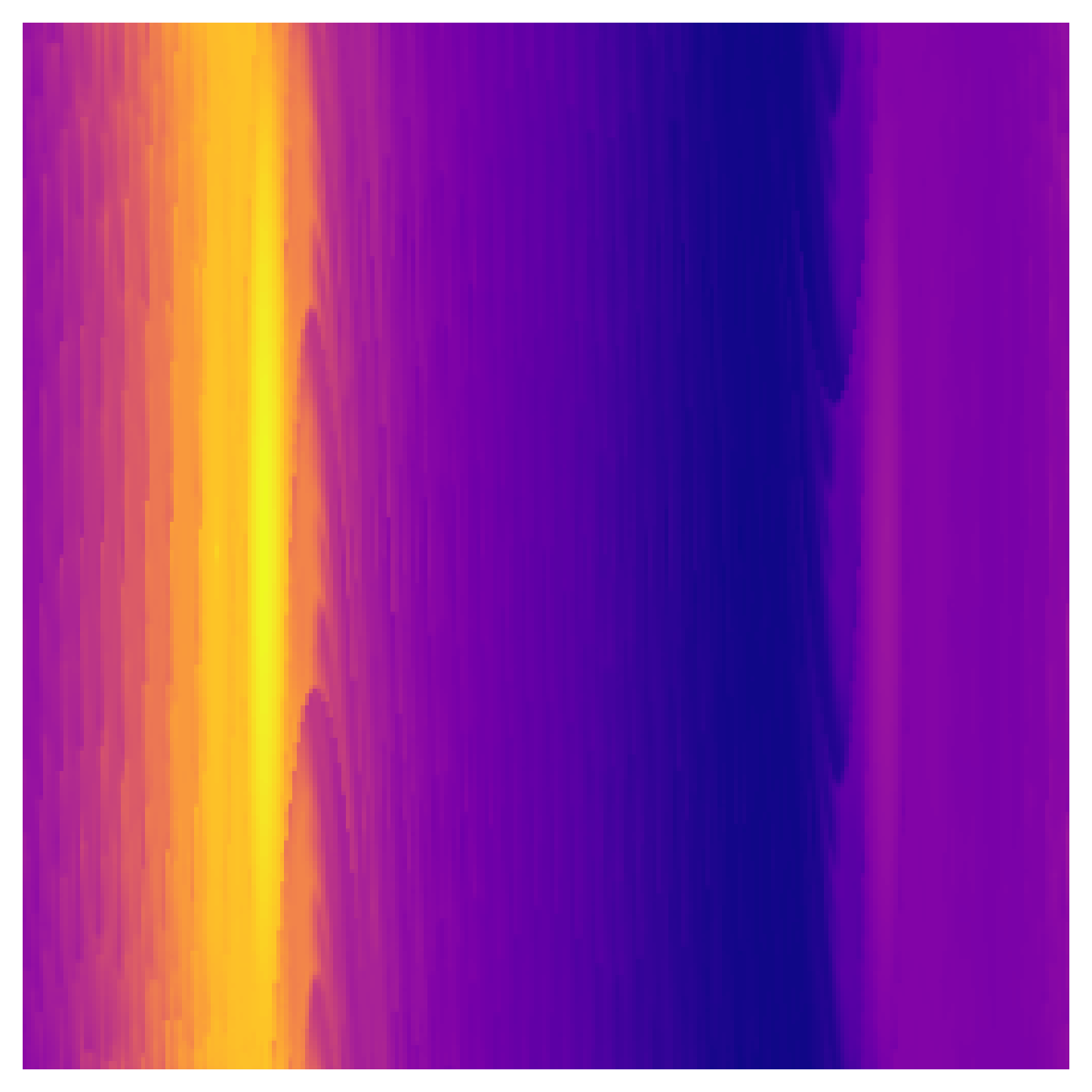} & \includegraphics[width=0.11\linewidth]{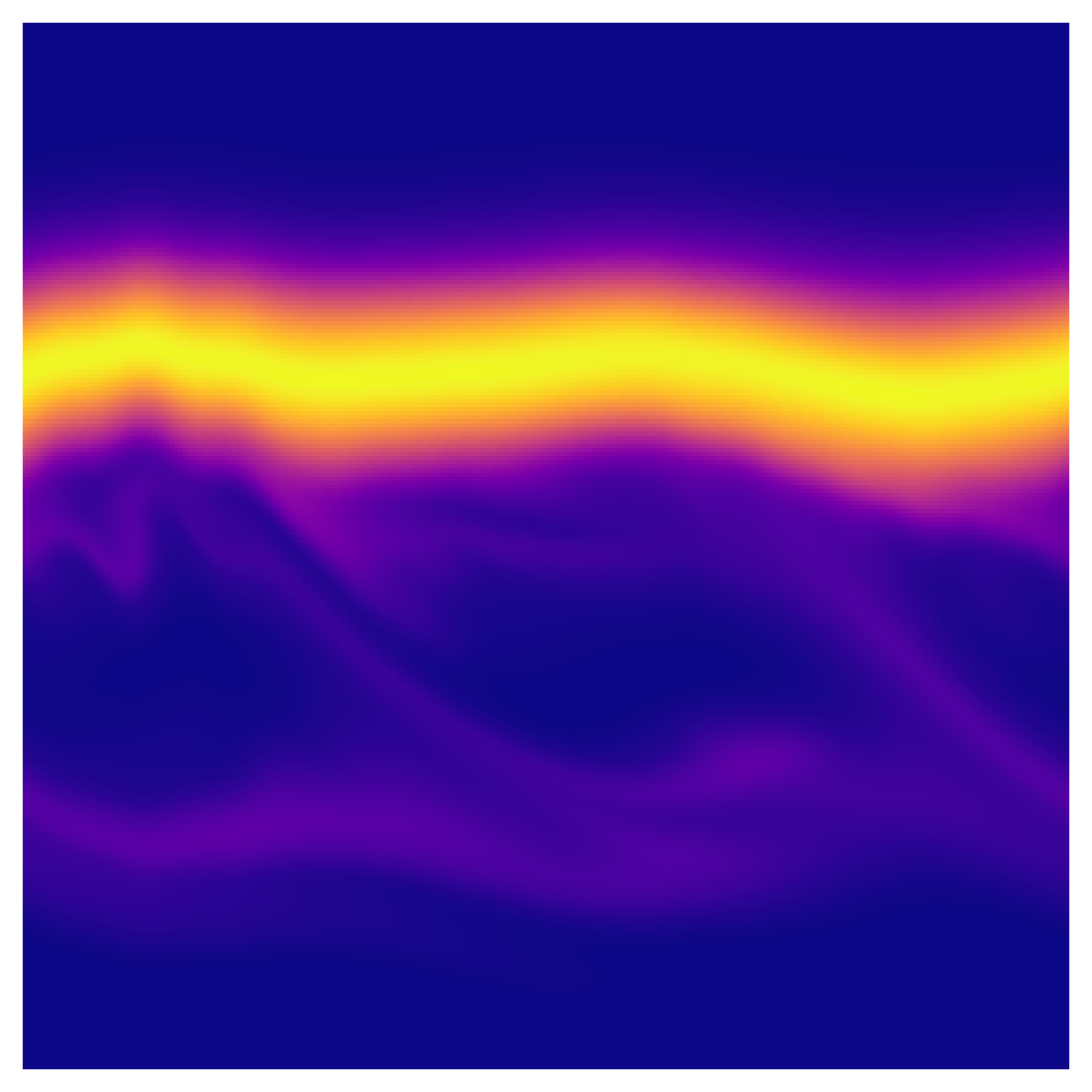} & \includegraphics[width=0.11\linewidth]{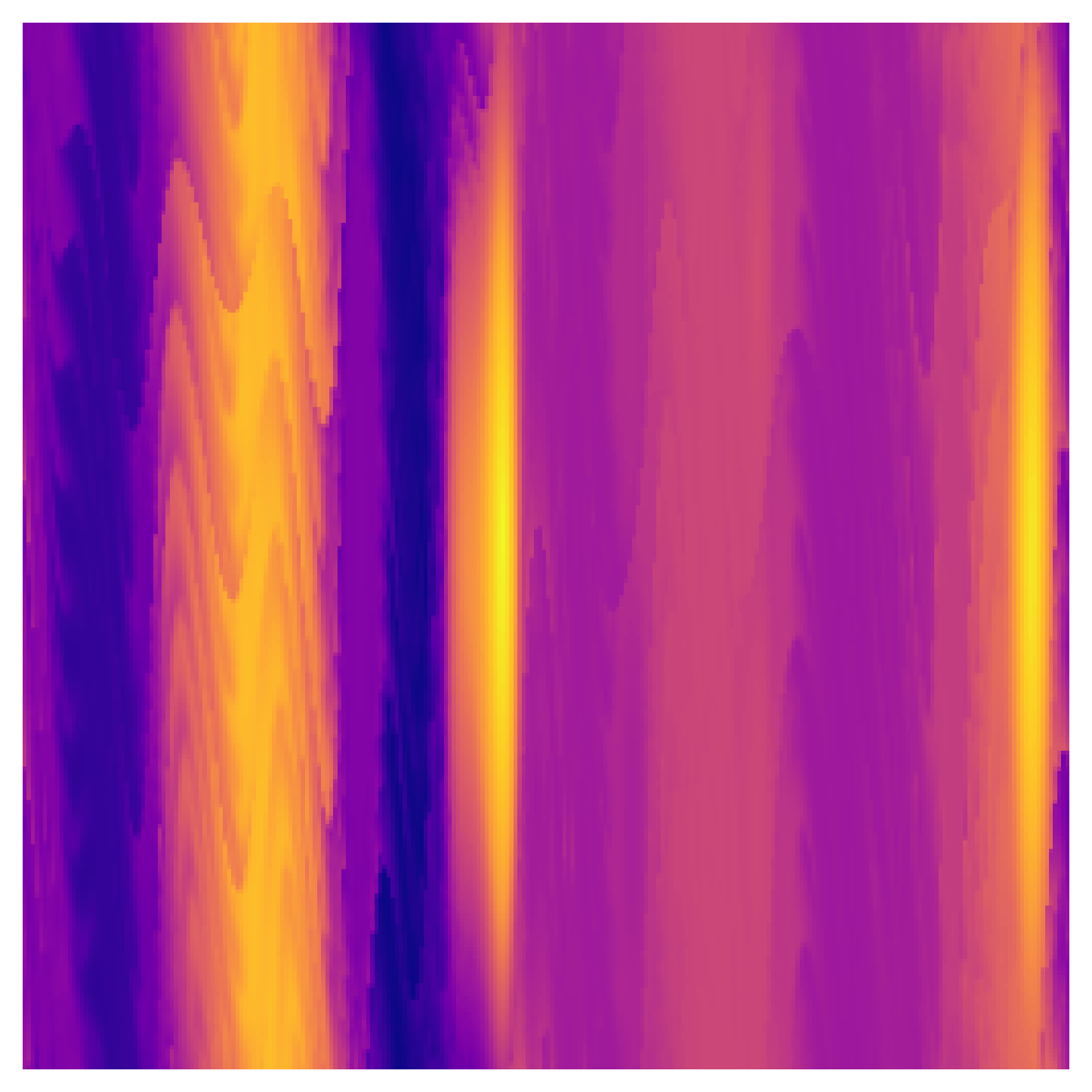} & \includegraphics[width=0.11\linewidth]{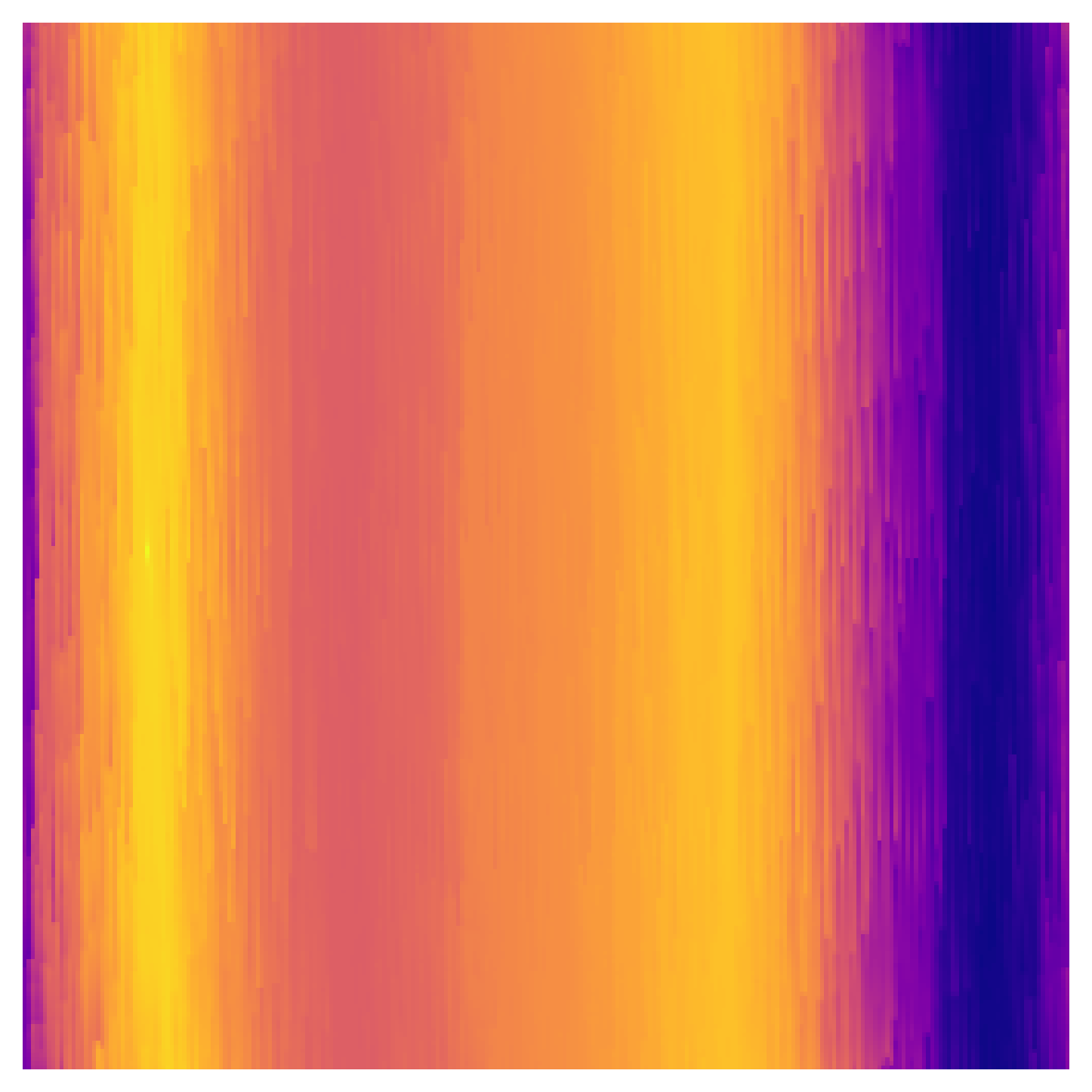} \\
    & & (Fig.~\ref{fig:BoT_KL_GDL_near_under}) & (Fig.~\ref{fig:BoT_ee_lf_GDL_near_under}) & (Fig.~\ref{fig:BoT_ee_GDL_near_under}) & (Fig.~\ref{fig:BoT_KL_GDL_near_over}) & (Fig.~\ref{fig:BoT_ee_lf_GDL_near_over}) & (Fig.~\ref{fig:BoT_ee_GDL_near_over}) \\ \cline{2-8}
    & \multirow{2}{*}[22pt]{Local} & \includegraphics[width=0.11\linewidth]{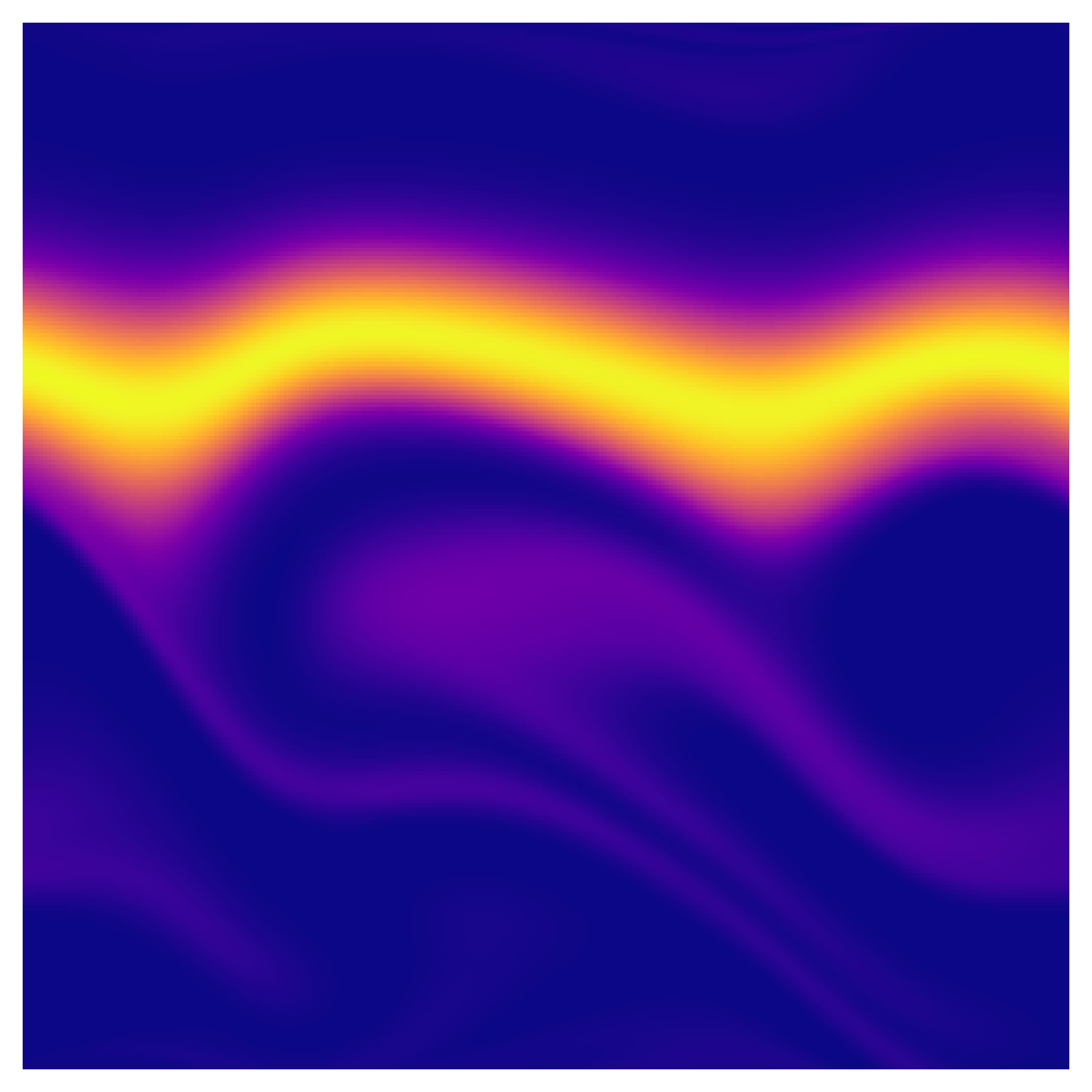} & \includegraphics[width=0.11\linewidth]{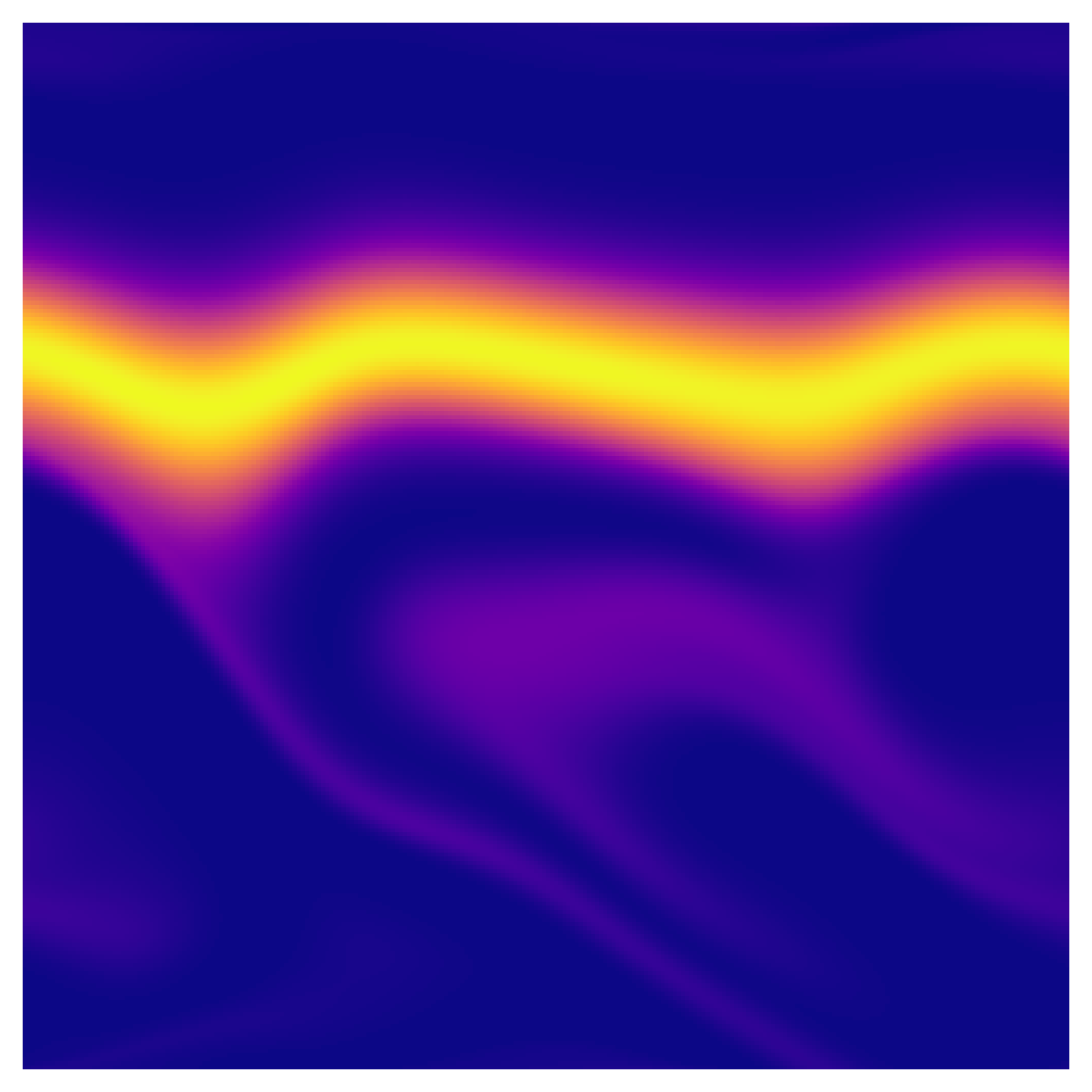} & \includegraphics[width=0.11\linewidth]{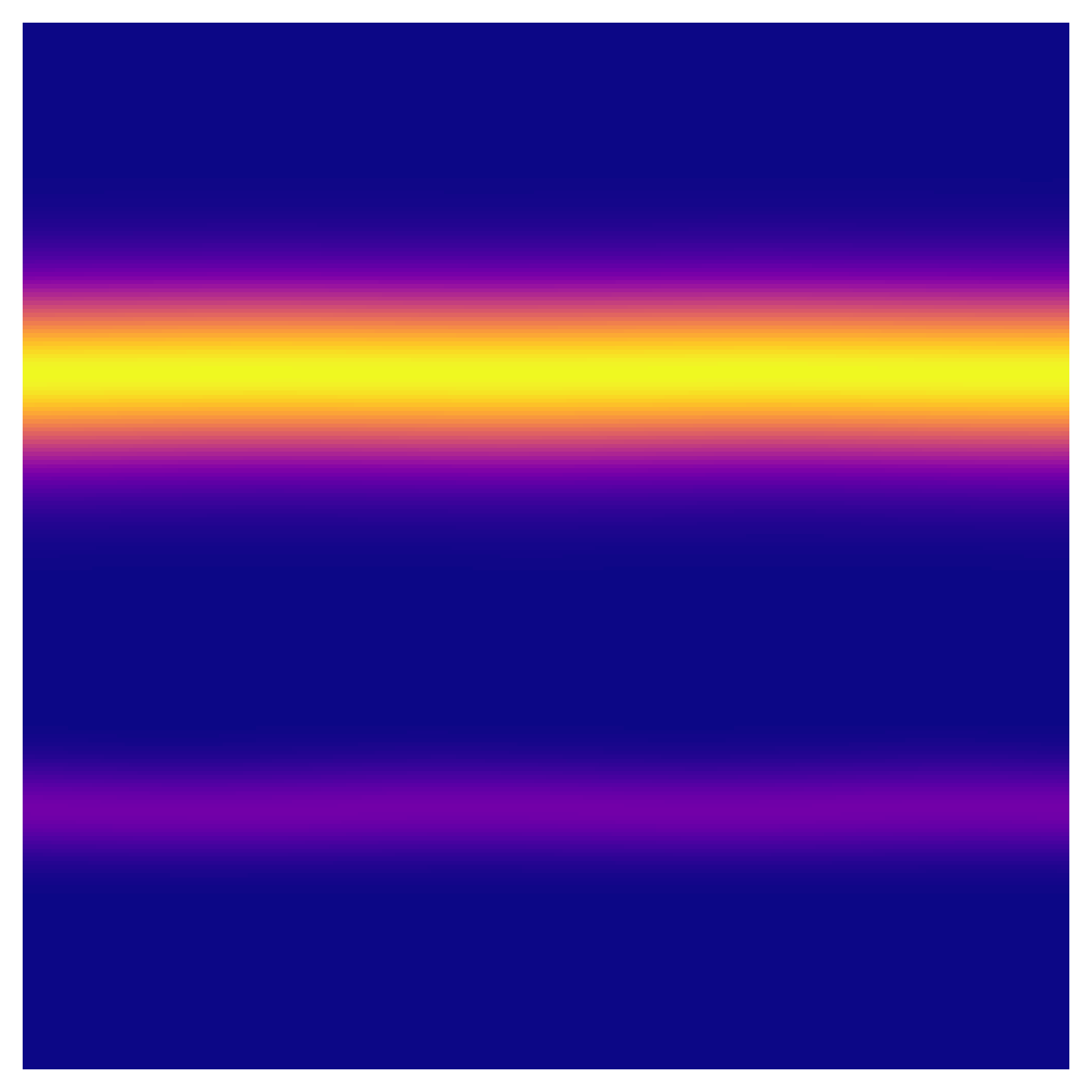} & \includegraphics[width=0.11\linewidth]{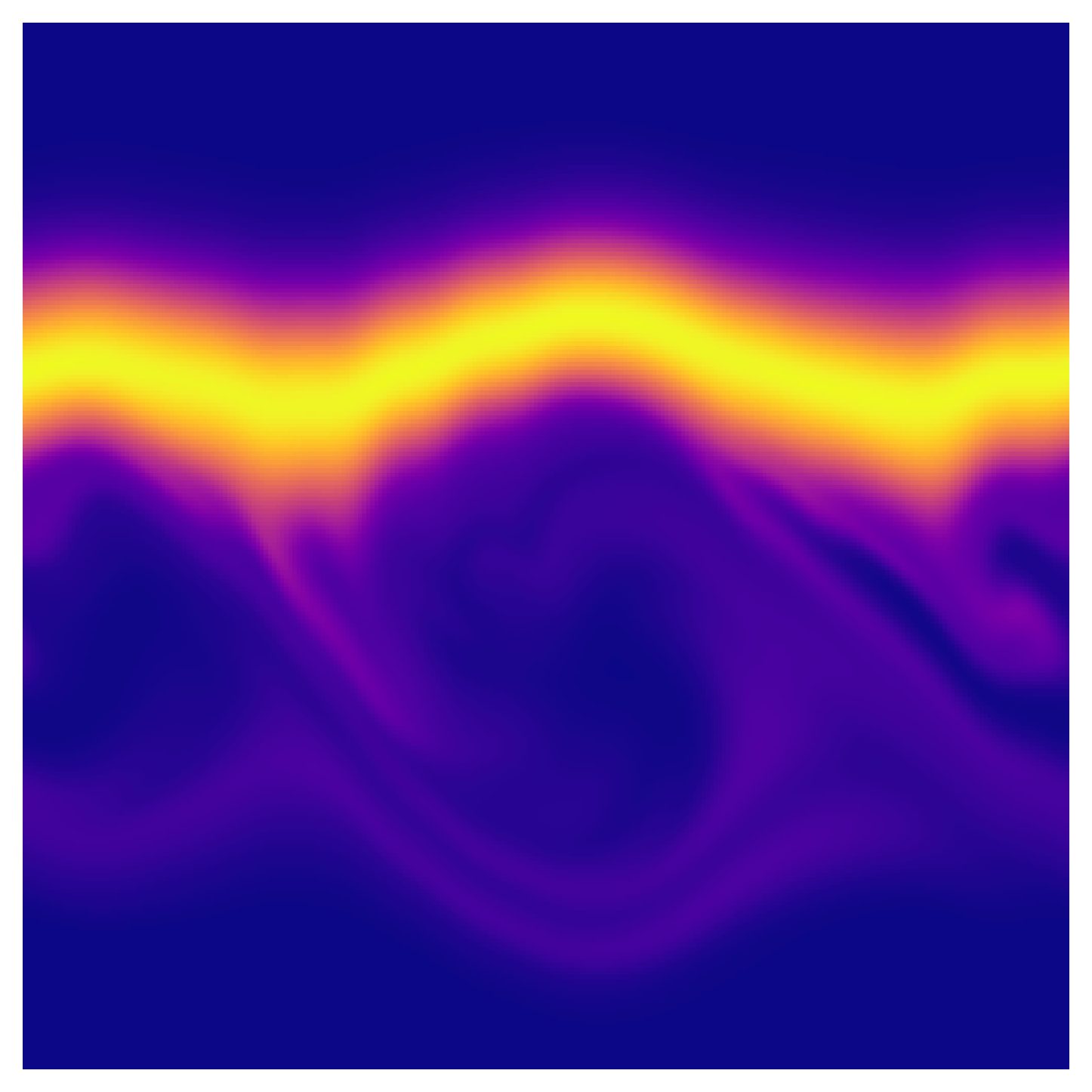} & \includegraphics[width=0.11\linewidth]{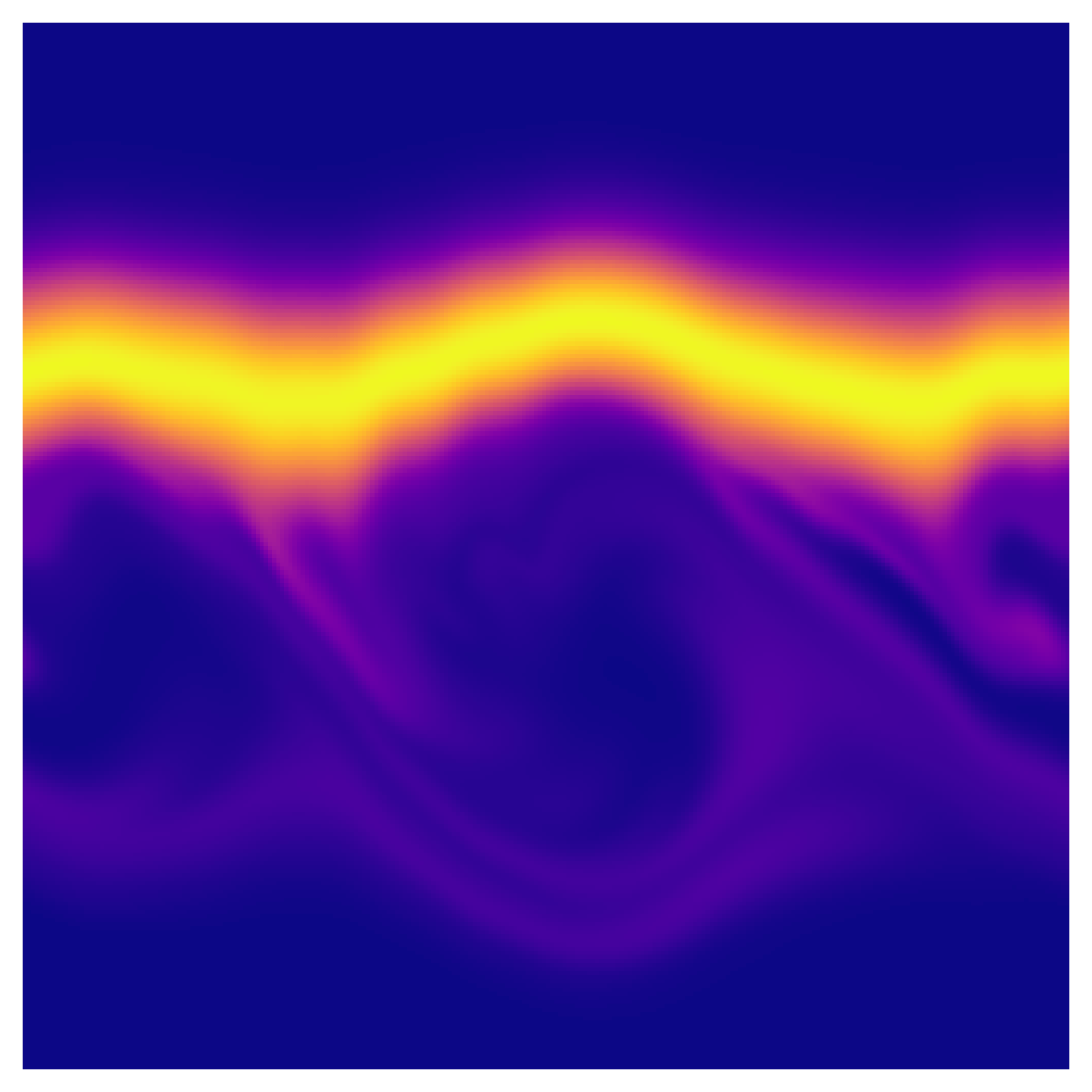} & \includegraphics[width=0.11\linewidth]{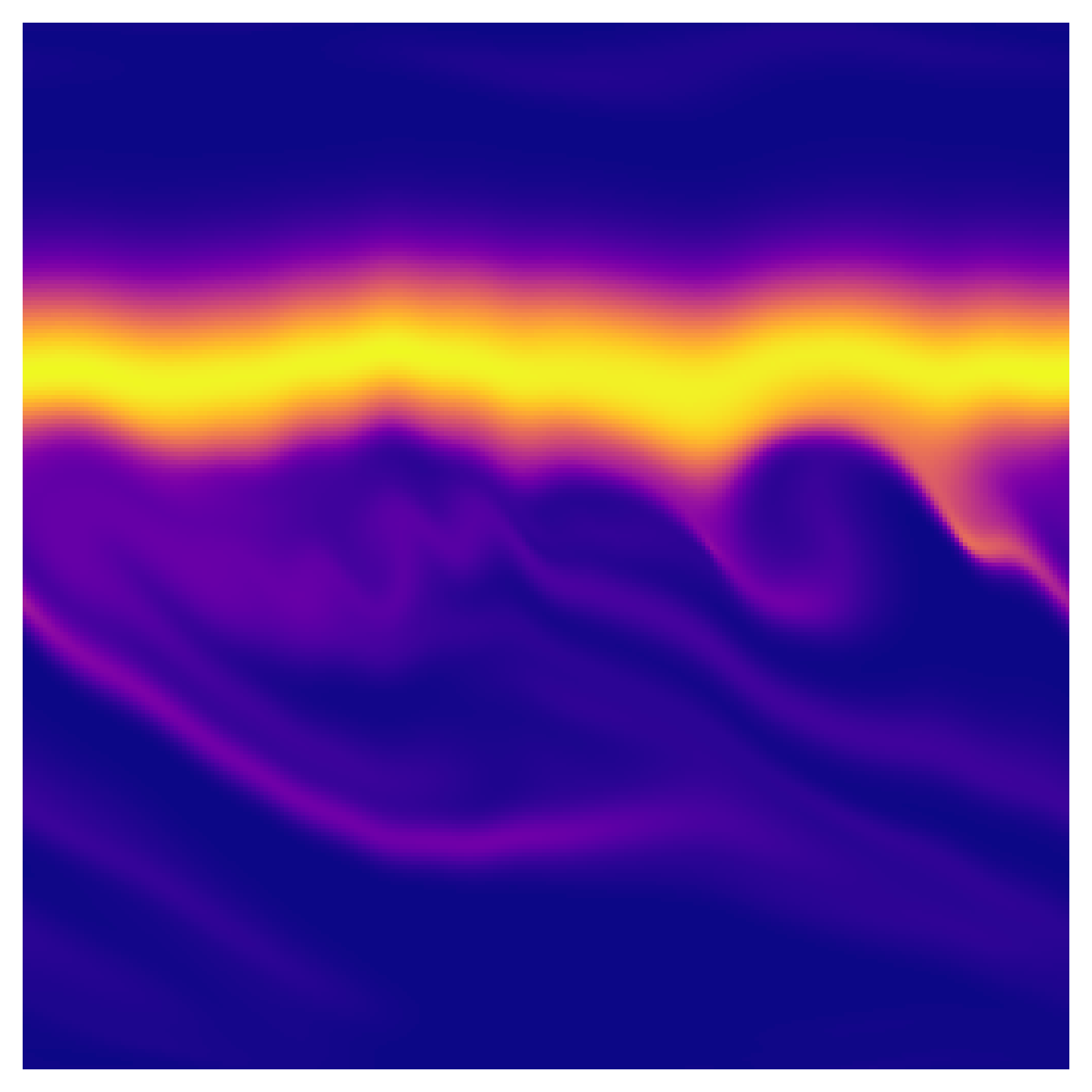} \\
    & & (Fig.~\ref{fig:BoT_KL_GD_near_under}) & (Fig.~\ref{fig:BoT_ee_lf_GD_near_under}) & 
    (Fig.~\ref{fig:BoT_ee_GD_near_under}) & (Fig.~\ref{fig:BoT_KL_GD_near_over}) & (Fig.~\ref{fig:BoT_ee_lf_GD_near_over}) & (Fig.~\ref{fig:BoT_ee_GD_near_over}) \\
    \midrule
    \multirow{4}{*}[0pt]{Near} & \multirow{2}{*}[22pt]{Adaptive} & \includegraphics[width=0.11\linewidth]{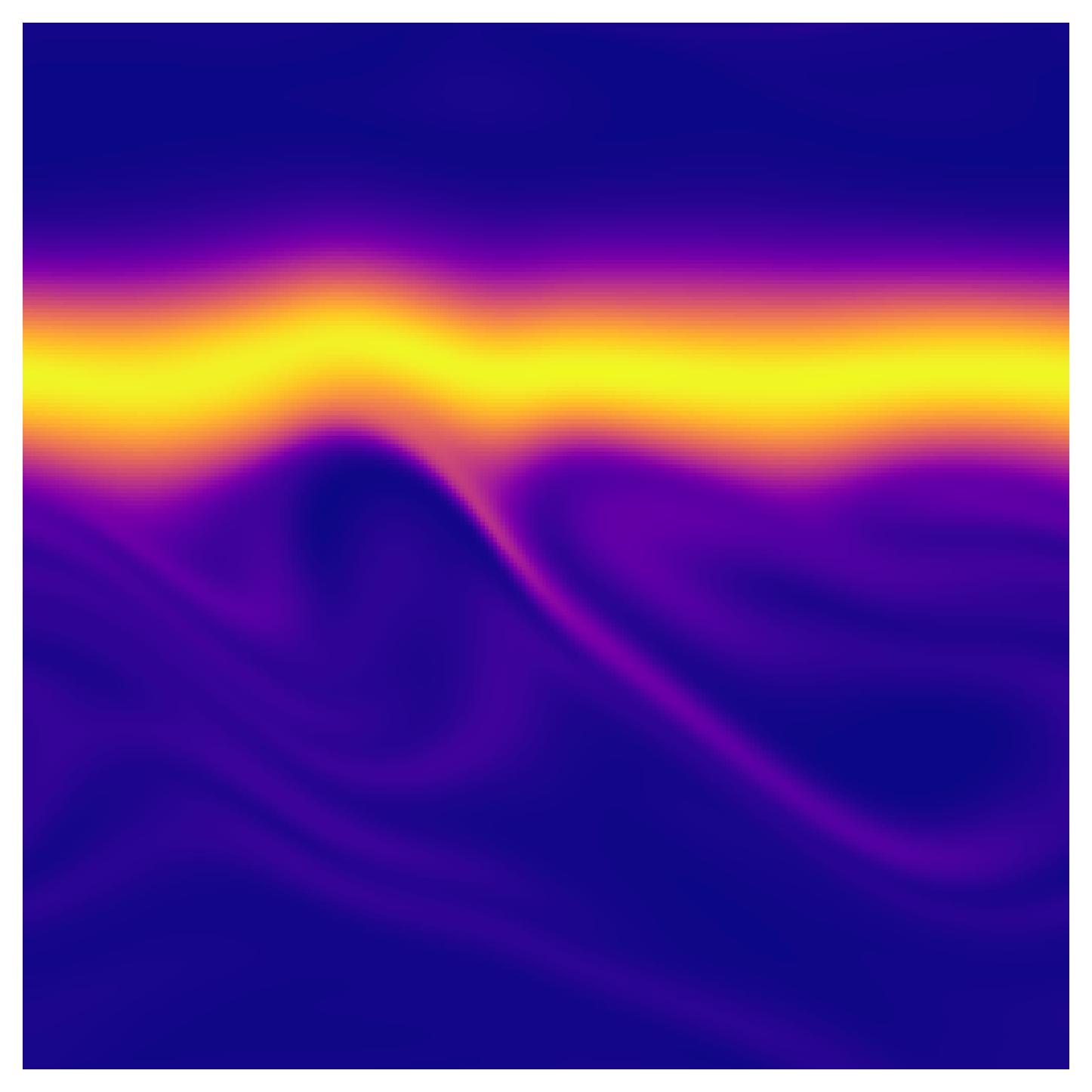}  &  \includegraphics[width=0.11\linewidth]{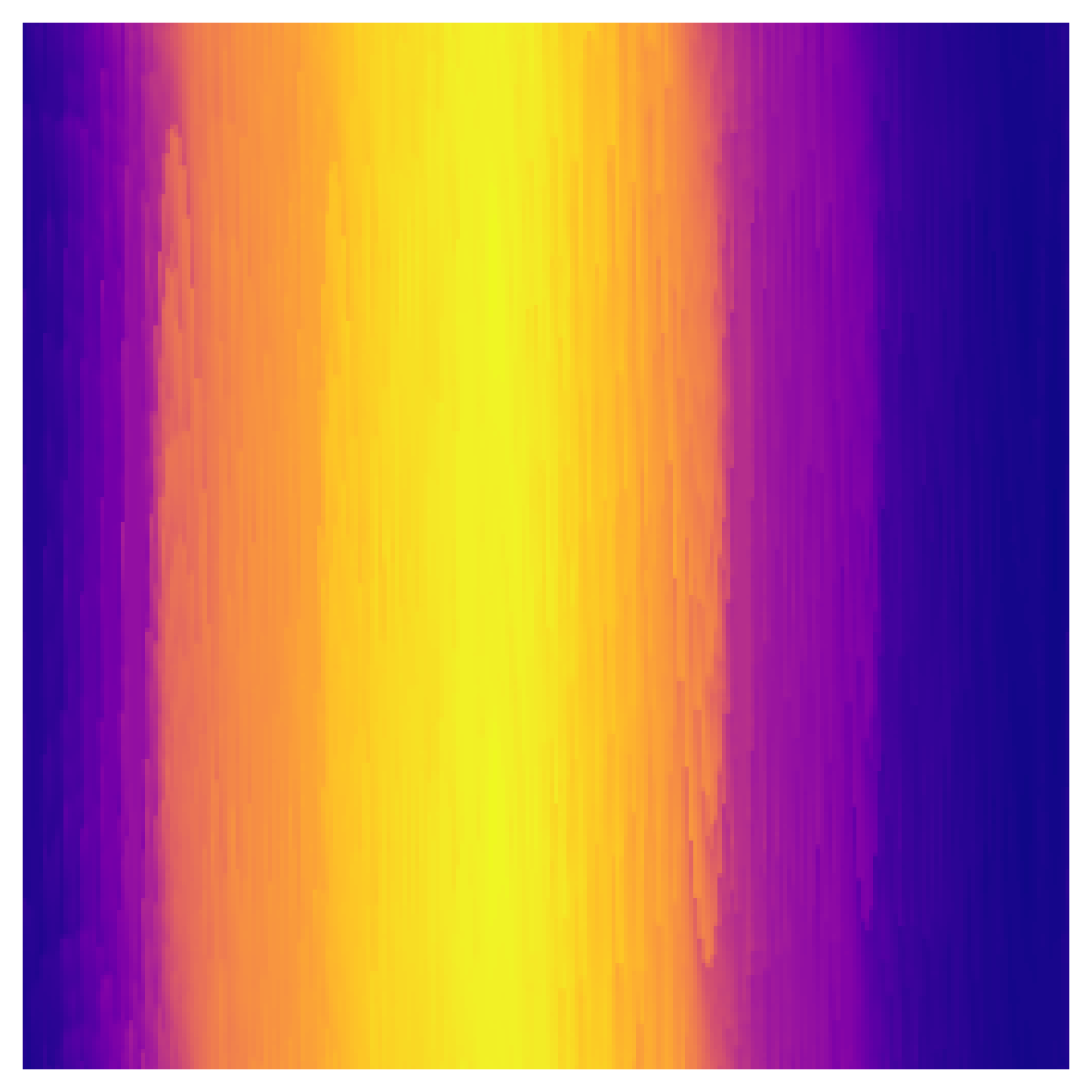}  &
    \includegraphics[width=0.11\linewidth]{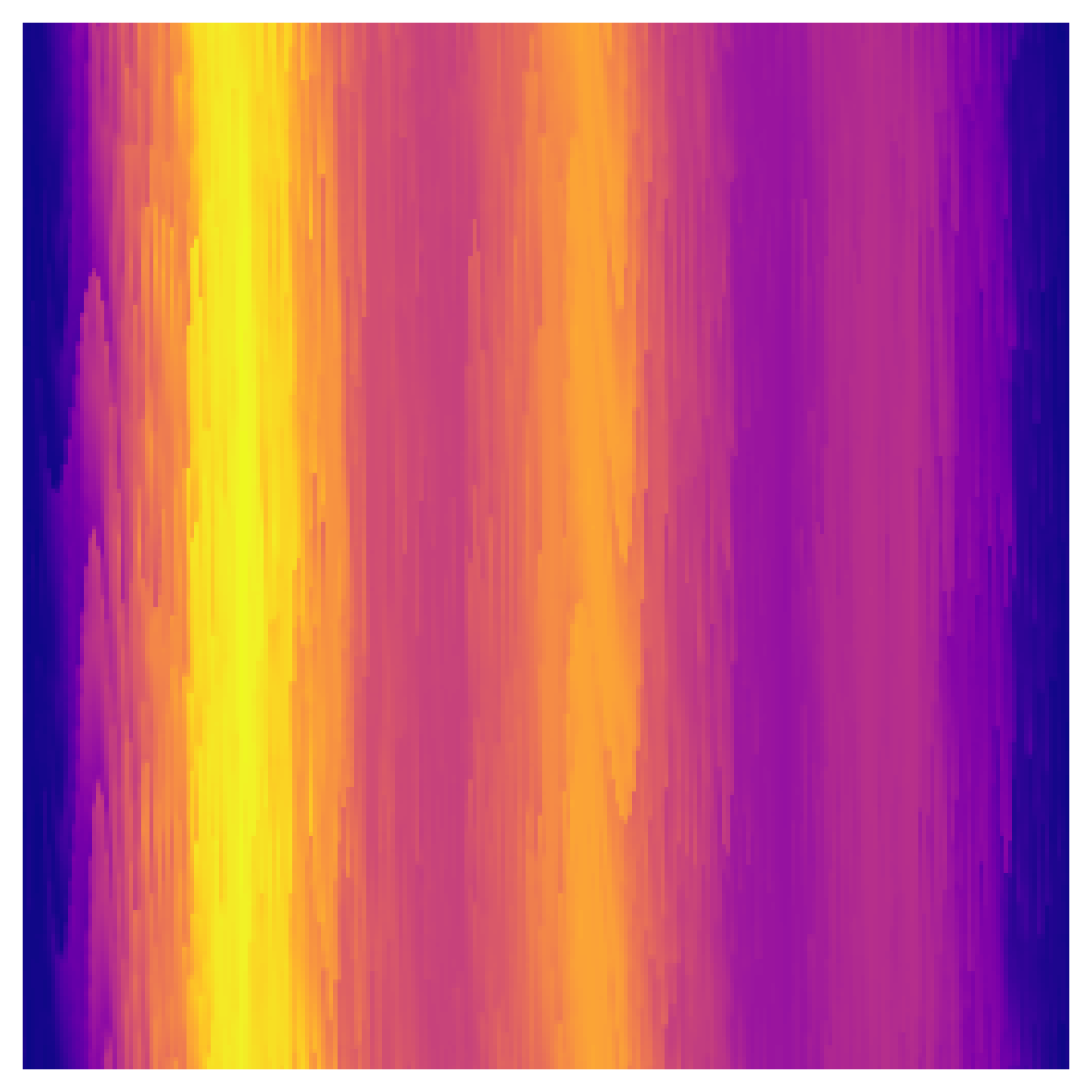}  &
    \includegraphics[width=0.11\linewidth]{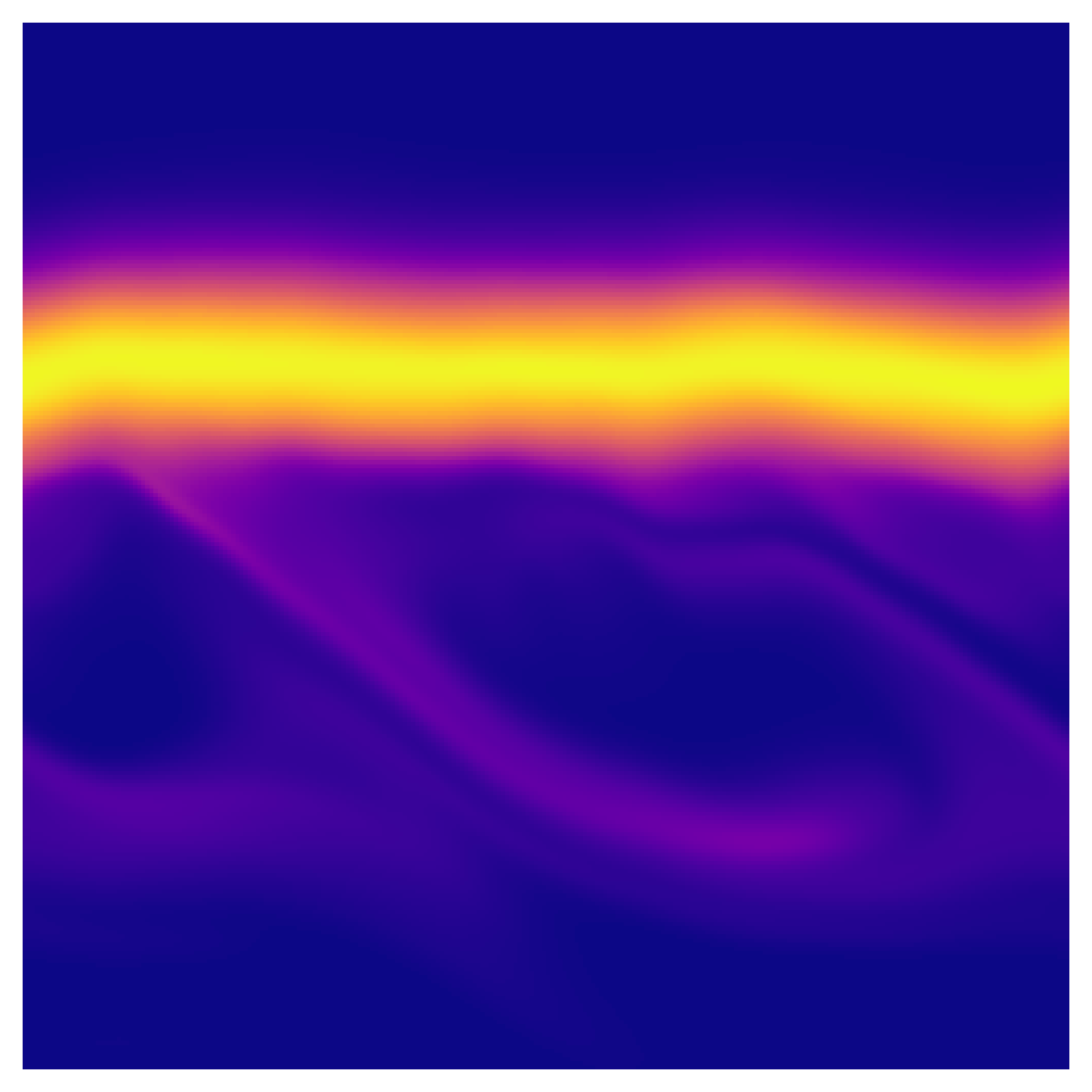} & \includegraphics[width=0.11\linewidth]{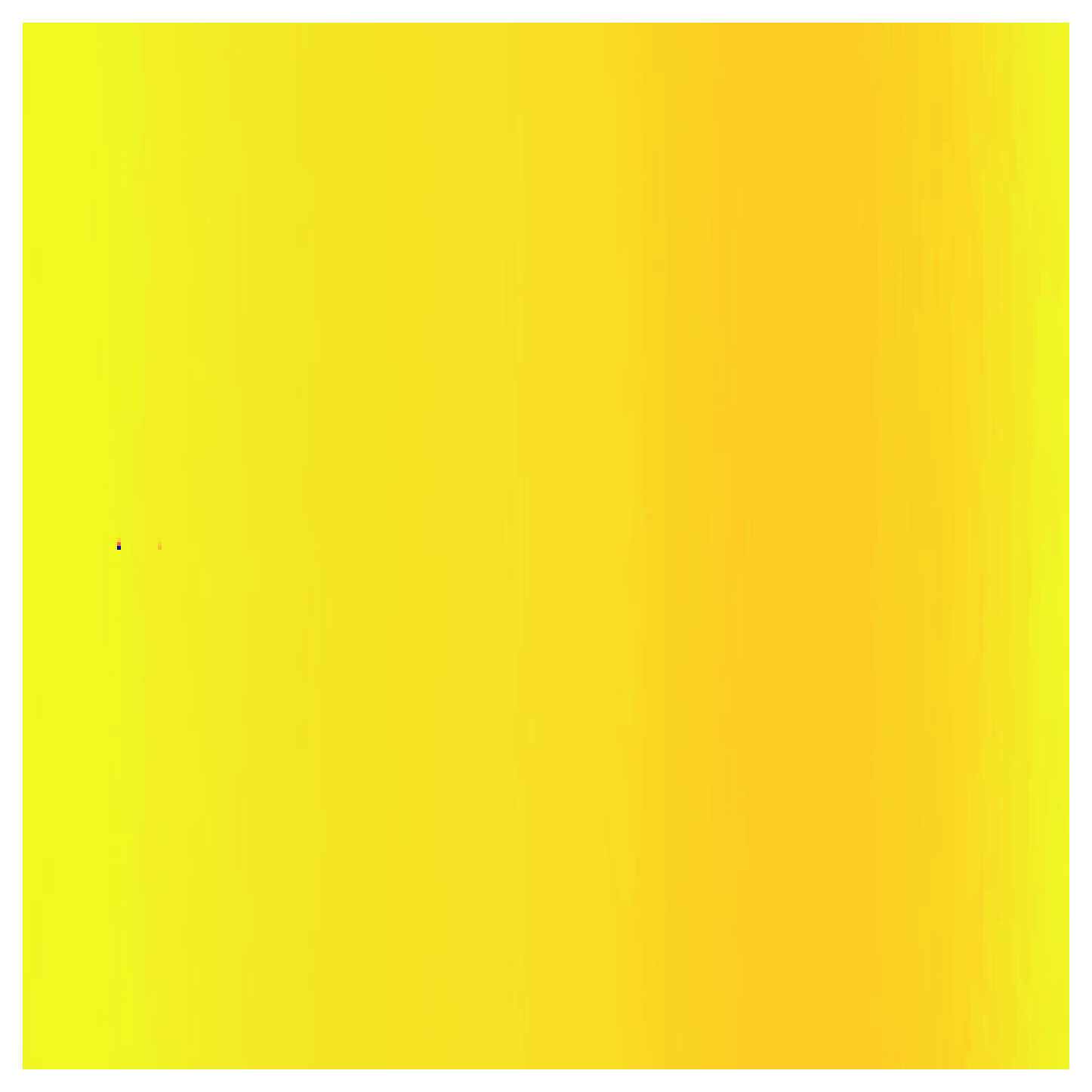} & \includegraphics[width=0.11\linewidth]{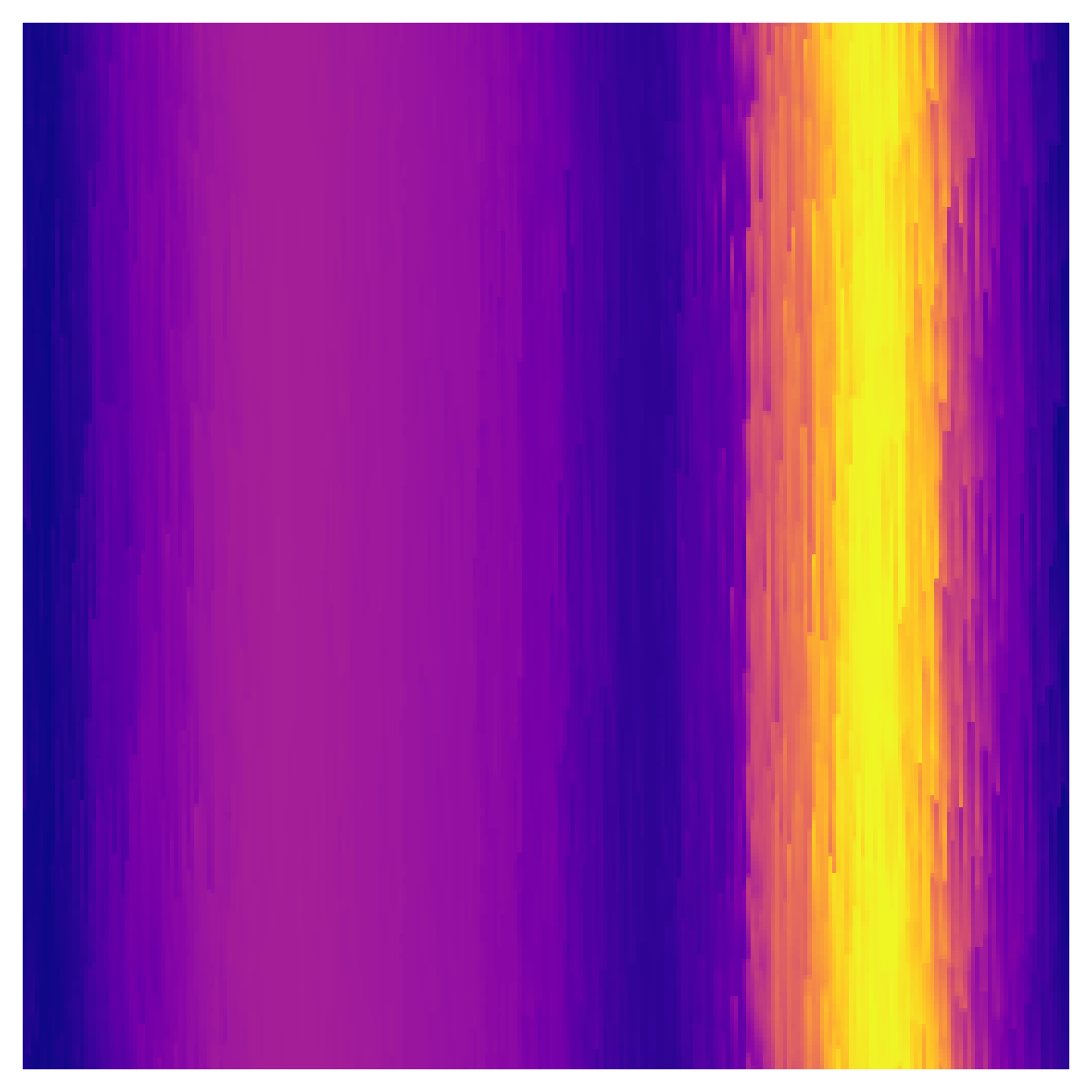} \\
    & & (Fig.~\ref{fig:BoT_KL_GDL_local_under}) & (Fig.~\ref{fig:BoT_ee_lf_GDL_local_under}) & 
    (Fig.~\ref{fig:BoT_ee_GDL_local_under}) & (Fig.~\ref{fig:BoT_KL_GDL_local_over}) & (Fig.~\ref{fig:BoT_ee_lf_GDL_local_over}) & (Fig.~\ref{fig:BoT_ee_GDL_local_over}) \\ \cline{2-8}
    & \multirow{2}{*}[22pt]{Local} & \includegraphics[width=0.11\linewidth]{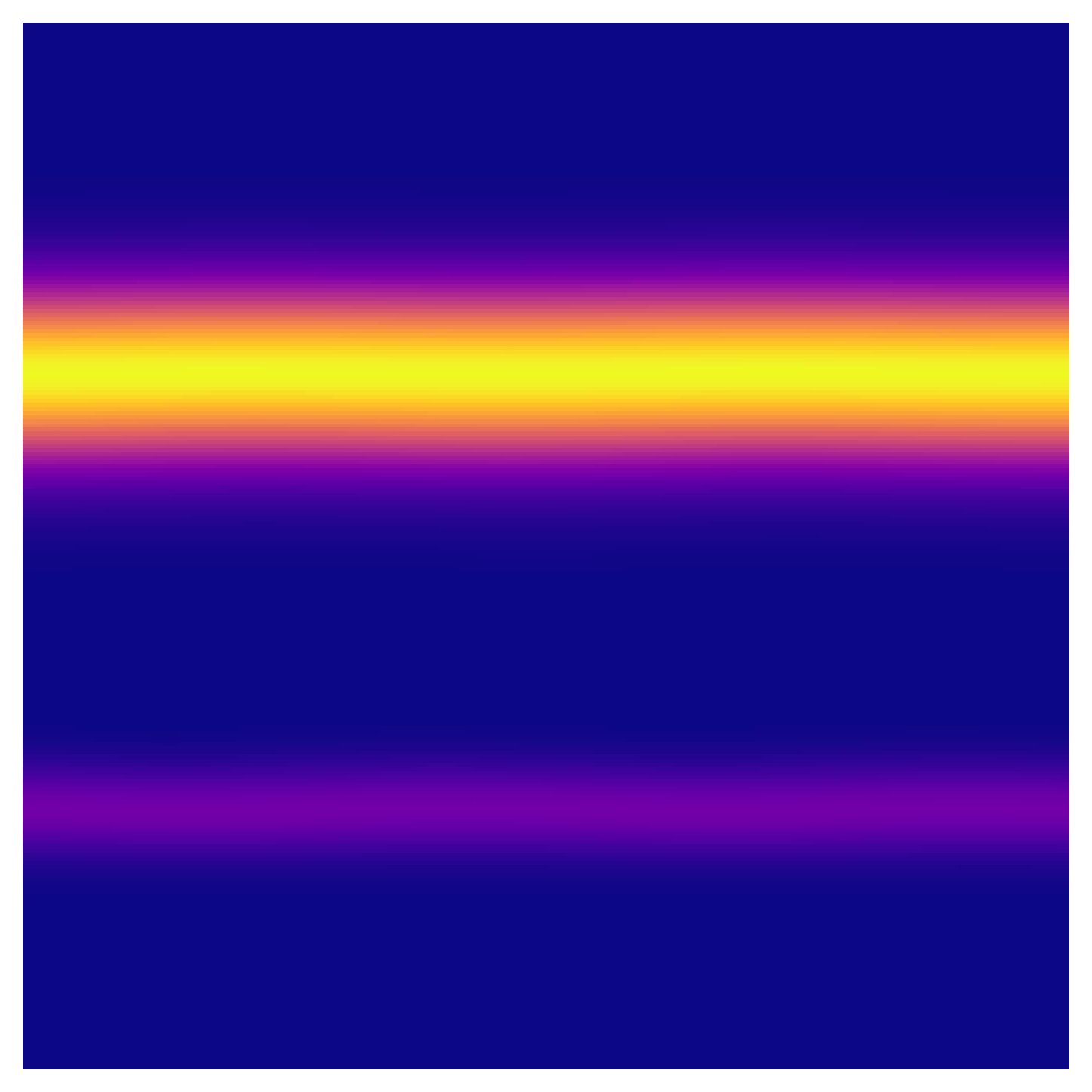} & \includegraphics[width=0.11\linewidth]{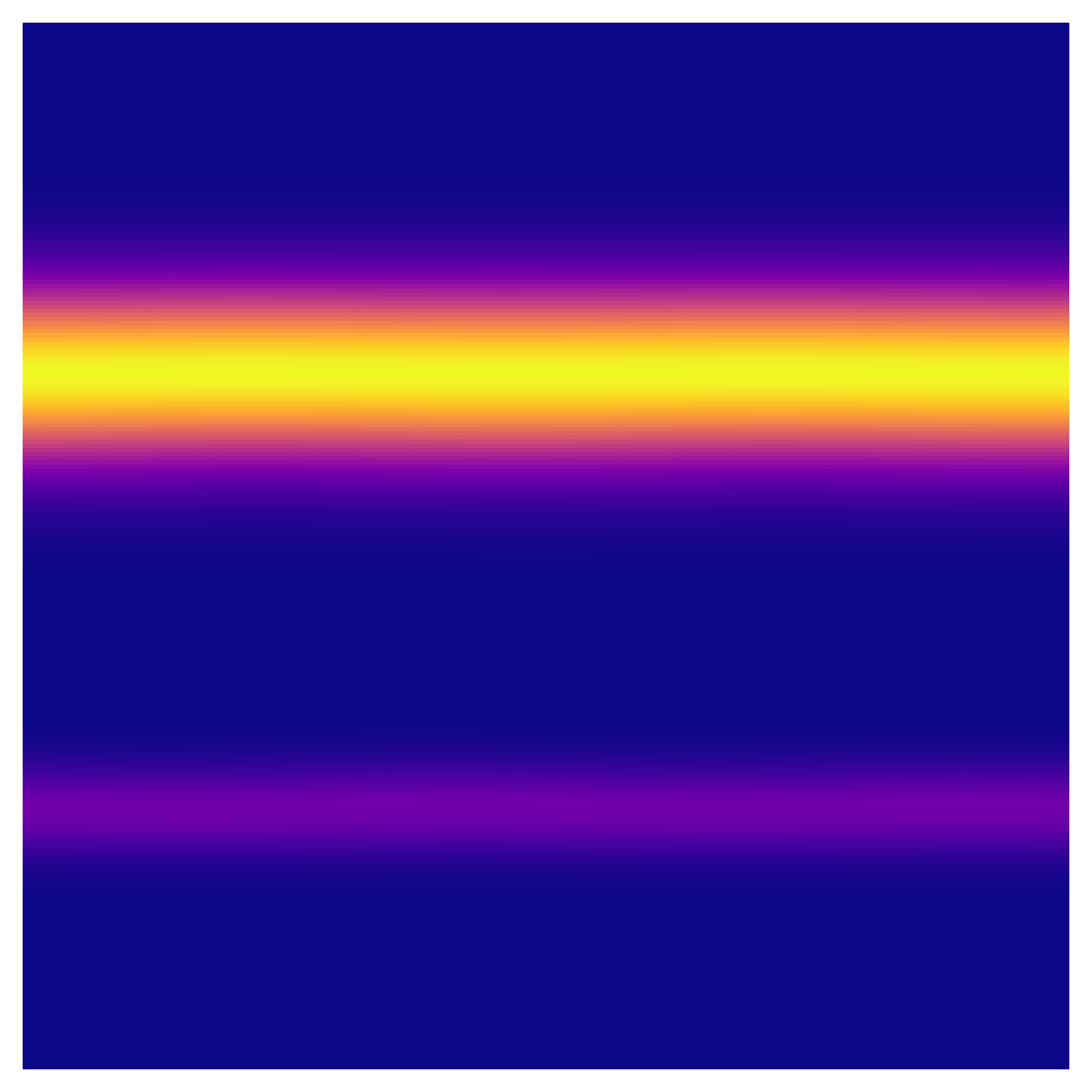} &
    \includegraphics[width=0.11\linewidth]{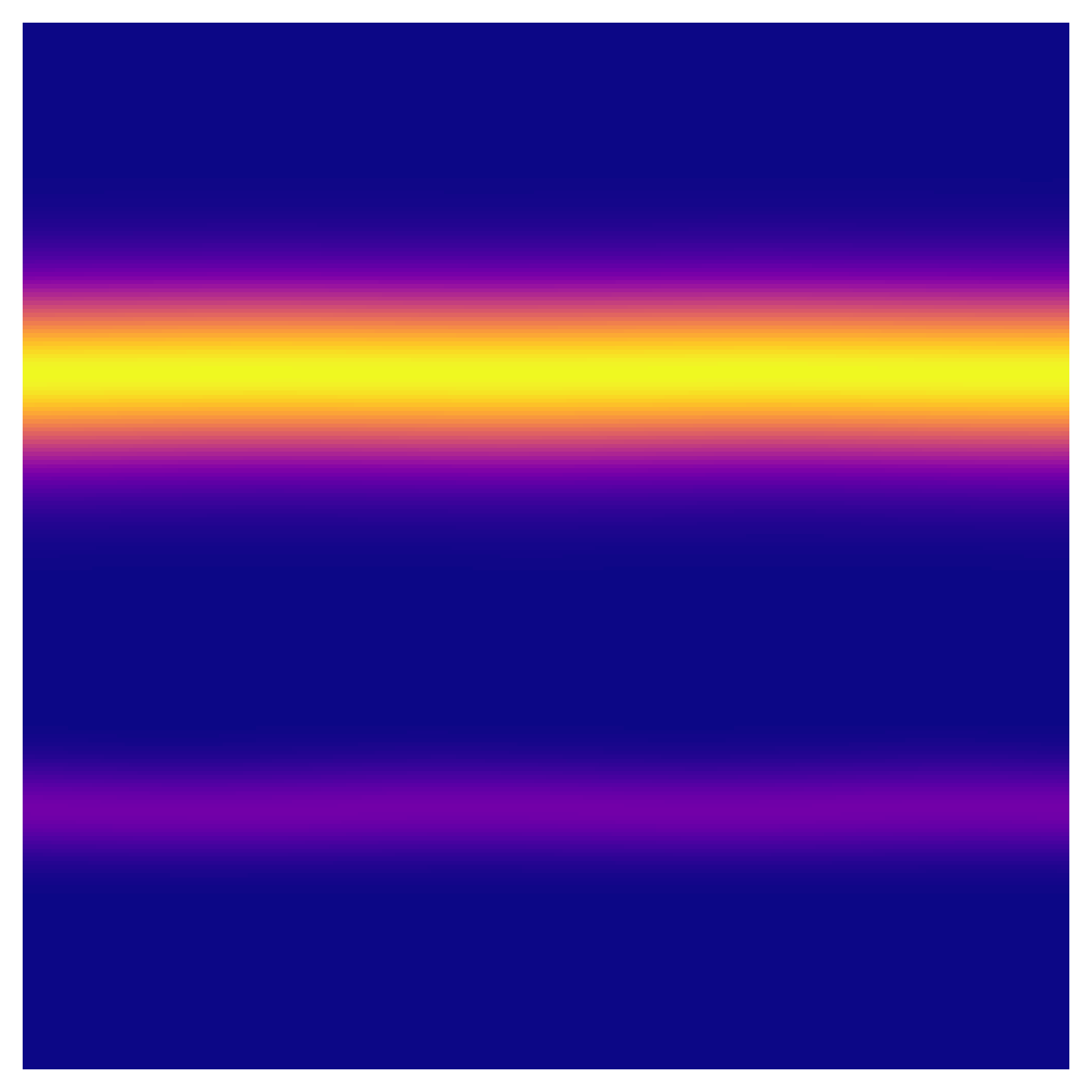} & \includegraphics[width=0.11\linewidth]{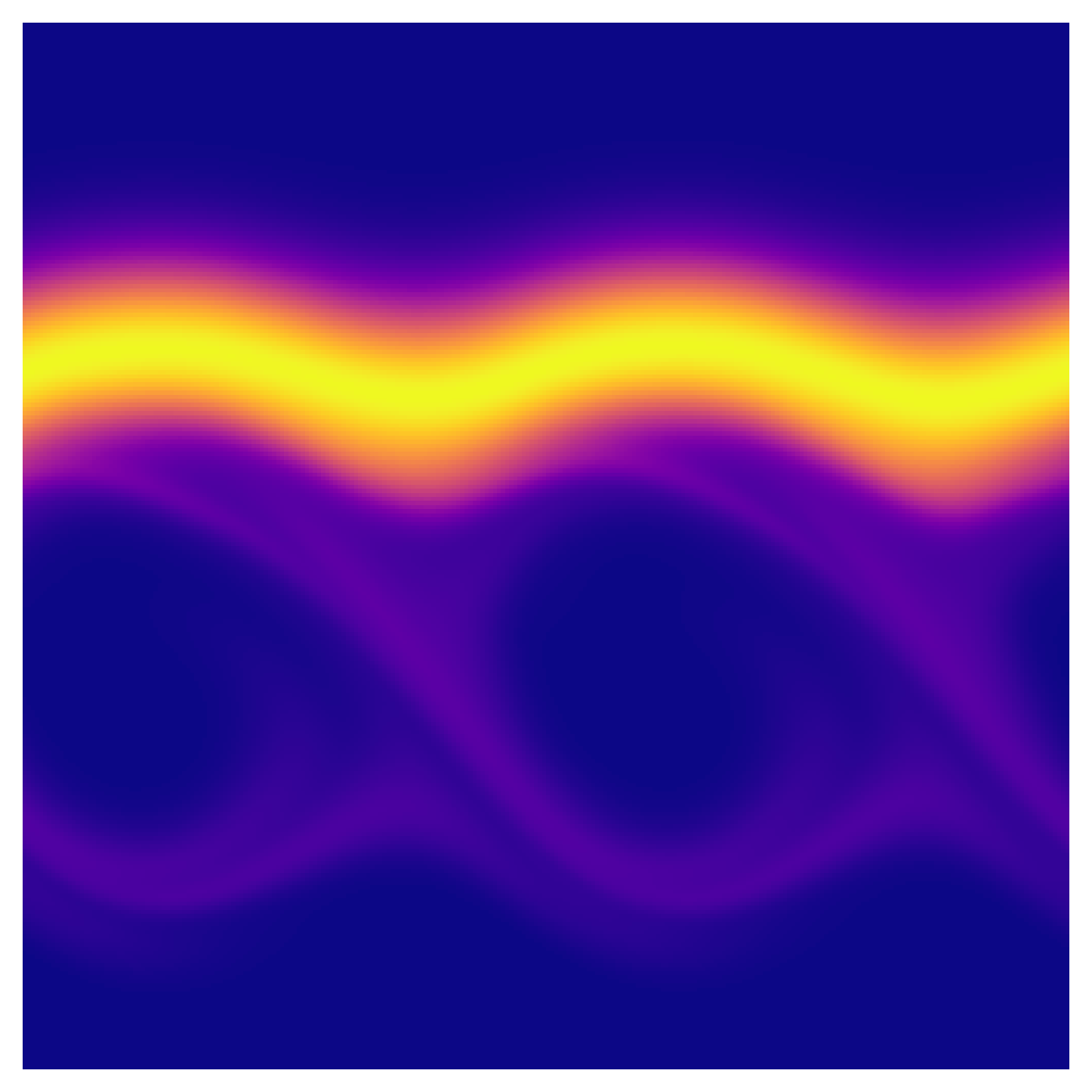} & \includegraphics[width=0.11\linewidth]{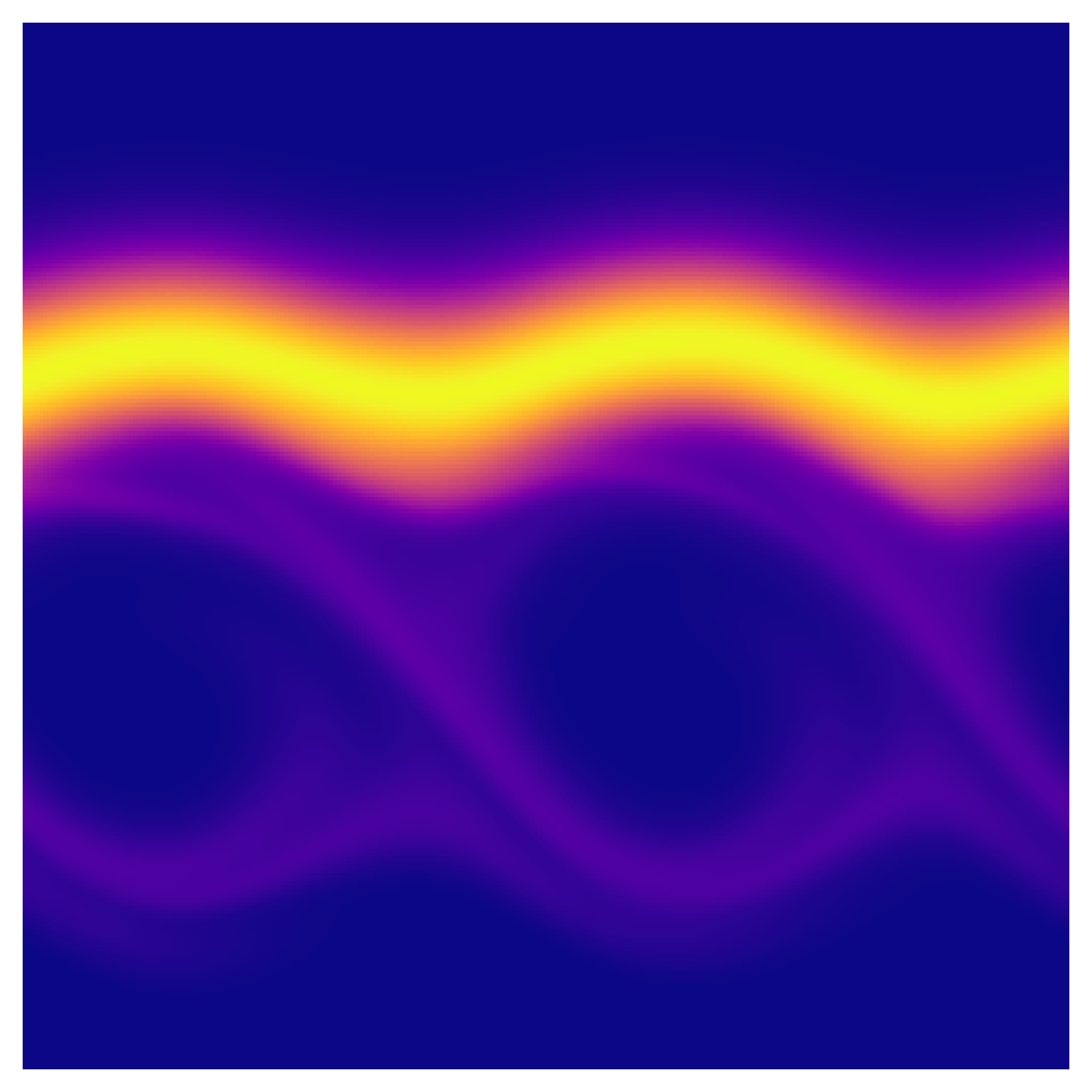} & \includegraphics[width=0.11\linewidth]{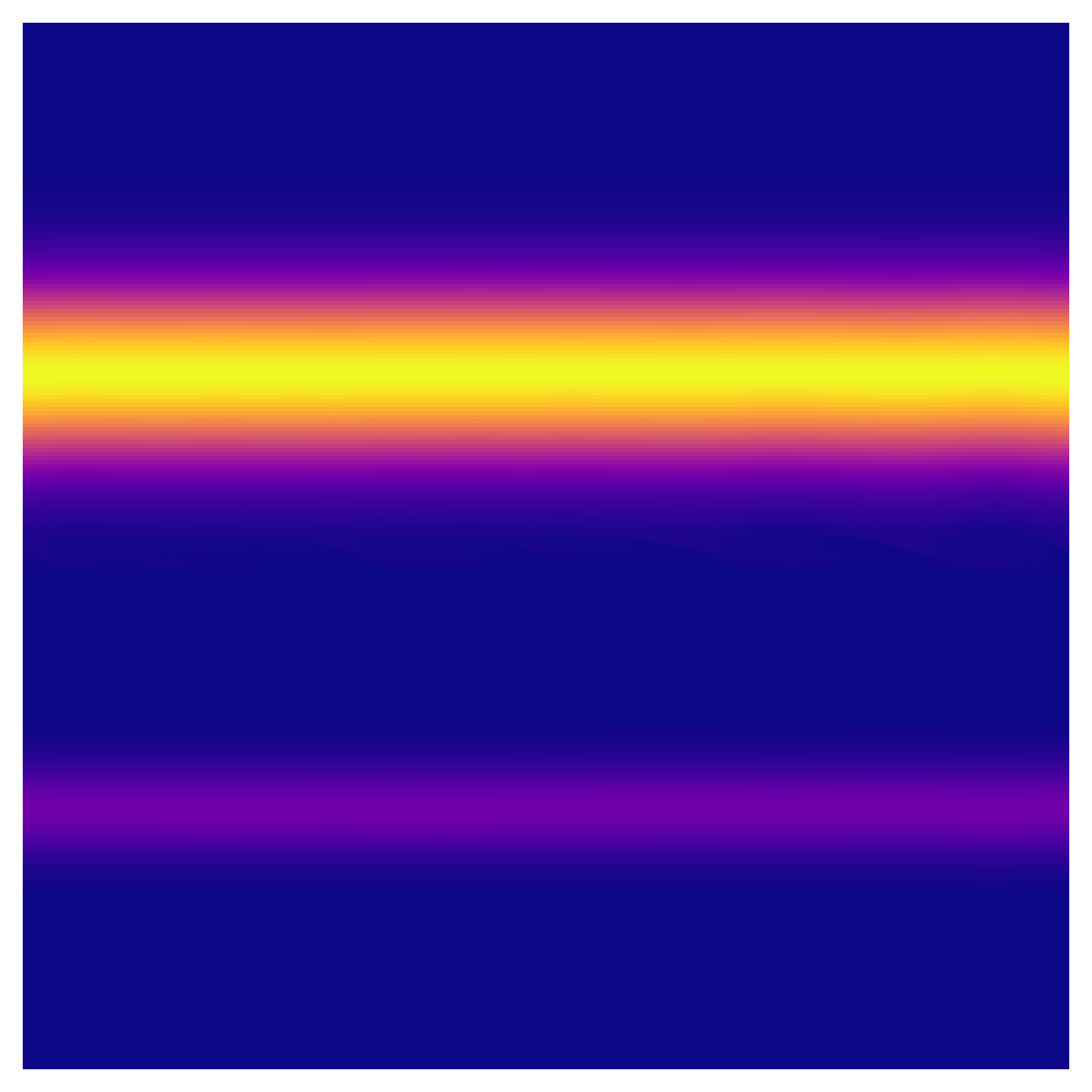} \\
    & & (Fig.~\ref{fig:BoT_KL_GD_local_under}) & (Fig.~\ref{fig:BoT_ee_lf_GD_local_under}) & (Fig.~\ref{fig:BoT_ee_GD_local_under}) & (Fig.~\ref{fig:BoT_KL_GD_local_over}) & (Fig.~\ref{fig:BoT_ee_lf_GD_local_over}) & (Fig.~\ref{fig:BoT_ee_GD_local_over}) \\
    \bottomrule
    \end{tabular}
    \caption{Summary of results for the Bump-on-Tail equilibrium.}
    \label{tab:BoT_results}
\end{table}

\section{Conclusion}

Stabilizing plasma dynamics modeled by the Vlasov–Poisson system can be formulated as an optimization problem: subject to the governing PDE, we seek external control fields that suppress dynamical instability. The choice of objective functional plays a crucial role in shaping the resulting optimization landscape, and therefore strongly influences the performance of numerical optimization methods.

In this work we investigated several objective functions and examined how their structure affects the optimization landscape. Our numerical experiments indicate that objective functions incorporating time-integrated information tend to produce smoother and more convex-like landscapes, which are more favorable for gradient-based optimization methods.

We also showed that analytical insight from the underlying PDE can guide the optimization process. In particular, dispersion-relation analysis of the linearized system allows one to identify unstable modes and construct control parameters that suppress them. Although this analysis is carried out in the linear regime, the resulting control field provides an effective initial guess for the nonlinear optimization problem.

Together, these observations suggest that combining PDE analysis with carefully designed objective functions can significantly improve the efficiency and robustness of optimization-based plasma stabilization strategies. This perspective highlights the importance of objective-function design when applying optimization methods to PDE-constrained plasma control problems. Future work will extend these ideas to more realistic plasma models and control settings.

\section*{Acknowledgements}
Q.L. and Y.Y. acknowledge support from NSF-DMS-2308440. M.G. acknowledges support from NSF-DMS-2012292.

\printbibliography

@article{Gilson2004,
  title = {Paul Trap Simulator Experiment to Model Intense-Beam Propagation in Alternating-Gradient Transport Systems},
  author = {Gilson, Erik P. and Davidson, Ronald C. and Efthimion, Philip C. and Majeski, Richard},
  journal = {Phys. Rev. Lett.},
  volume = {92},
  issue = {15},
  pages = {155002},
  numpages = {4},
  year = {2004},
  month = {4},
  publisher = {American Physical Society},
  doi = {10.1103/PhysRevLett.92.155002},
  url = {https://link.aps.org/doi/10.1103/PhysRevLett.92.155002}
}

@article{Mitchell2019,
  title = {Spectral {Galerkin} solver for intense beam {{Vlasov}} equilibria in nonlinear constant focusing channels},
  author = {Mitchell, Chad E. and Ryne, Robert D. and Hwang, Kilean},
  journal = {Phys. Rev. E},
  volume = {100},
  issue = {5},
  pages = {053308},
  numpages = {9},
  year = {2019},
  month = {11},
  publisher = {American Physical Society},
  doi = {10.1103/PhysRevE.100.053308},
  url = {https://link.aps.org/doi/10.1103/PhysRevE.100.053308}
}

@book{chen1984introduction,
  title={Introduction to plasma physics and controlled fusion},
  author={Chen, Francis F},
  volume={1},
  year={1984},
  publisher={Springer},
  doi = {10.1007/978-1-4757-5595-4}
}

@article{borzi2025optimal,
	author = {Borz\`{i}, A. and Infante, G. and Mascali, G.},
	doi = {10.1137/24M1643177},
	journal = {SIAM Journal on Applied Mathematics},
	number = {1},
	pages = {143--163},
	title = {Optimal Design of Equilibrium Solutions of the {Vlasov}--{{Poisson}} System by an External Electric Field},
	volume = {85},
	year = {2025}
}

@inproceedings{filbet2003numerical,
  title = {Numerical methods for the {Vlasov} equation},
  author = {Filbet, Francis and Sonnendr{\"u}cker, Eric},
  booktitle = {Numerical Mathematics and Advanced Applications: Proceedings of ENUMATH 2001 the 4th European Conference on Numerical Mathematics and Advanced Applications Ischia, July 2001},
  pages = {459--468},
  year = {2003},
  publisher = {Springer Milan},
  isbn = {978-88-470-2089-4},
  doi = {10.1007/978-88-470-2089-4_43}
}

@article{besse2003semi,
	author = {N. Besse and E. Sonnendr{\"u}cker},
	doi = {10.1016/S0021-9991(03)00318-8},
	issn = {0021-9991},
	journal = {Journal of Computational Physics},
	number = {2},
	pages = {341-376},
	title = {Semi-Lagrangian schemes for the {Vlasov} equation on an unstructured mesh of phase space},
	url = {https://www.sciencedirect.com/science/article/pii/S0021999103003188},
	volume = {191},
	year = {2003},
    publisher={Elsevier}
}

@article{o1971nonlinear,
  title={Nonlinear interaction of a small cold beam and a plasma},
  author={O'Neil, T.M. and Winfrey, J.H. and Malmberg, J.H.},
  journal={The Physics of Fluids},
  volume={14},
  number={6},
  pages={1204--1212},
  year={1971},
  publisher={AIP Publishing},
  doi = {10.1063/1.1693587}
}

@article{cai2021high,
	author = {Xiaofeng Cai and Sebastiano Boscarino and Jing-Mei Qiu},
	doi = {10.1016/j.jcp.2020.110036},
	issn = {0021-9991},
	journal = {Journal of Computational Physics},
	pages = {110036},
	title = {High order semi-Lagrangian discontinuous Galerkin method coupled with Runge-Kutta exponential integrators for nonlinear {Vlasov} dynamics},
	url = {https://www.sciencedirect.com/science/article/pii/S002199912030810X},
	volume = {427},
	year = {2021},
    publisher={Elsevier},
}

@article{knopf2020optimal,
	author = {Knopf, Patrik and Weber, J{\"o}rg},
	doi = {10.1007/s00245-018-9526-5},
	journal = {Applied Mathematics \& Optimization},
	number = {3},
	pages = {961--988},
	title = {Optimal Control of a {Vlasov}--{Poisson} Plasma by Fixed Magnetic Field Coils},
	volume = {81},
	year = {2020},
    publisher = {Springer}
}

@book{fitzpatrick2022plasma,
  title={Plasma physics: an introduction},
  author={Fitzpatrick, Richard},
  year={2022},
  publisher={CRC Press, Taylor \& Francis Group},
  doi = {10.1201/9781003268253},
  isbn = {978-1-4665-9426-5},
  note = {2nd edition}
}

@article{burm2012plasma,
	author = {Burm, KTAL},
	doi = {10.1007/s11090-012-9356-1},
	journal = {Plasma Chemistry and Plasma Processing},
	number = {2},
	pages = {401--407},
	title = {Plasma: The Fourth State of Matter},
	volume = {32},
	year = {2012},
    publisher={Springer}
}

@book{frank2012plasma,
  title={Plasma: the fourth state of matter},
  author={Frank-Kamenetskii, D},
  year={2012},
  publisher={Springer Science \& Business Media},
  doi = {10.1007/978-1-4684-1896-5}
}

@article{holm1985nonlinear,
	author = {Darryl D. Holm and Jerrold E. Marsden and Tudor Ratiu and Alan Weinstein},
	doi = {10.1016/0370-1573(85)90028-6},
	issn = {0370-1573},
	journal = {Physics Reports},
	number = {1},
	pages = {1-116},
	title = {Nonlinear stability of fluid and plasma equilibria},
	url = {https://www.sciencedirect.com/science/article/pii/0370157385900286},
	volume = {123},
	year = {1985},
    publisher = {Elsevier}
}

@book{nicholson1983introduction,
  title={Introduction to plasma theory},
  author={Nicholson, Dwight Roy},
  volume={1},
  year={1983},
  publisher={Wiley New York}
}

@article{ichimaru1993nuclear,
  title = {Nuclear fusion in dense plasmas},
  author = {Ichimaru, Setsuo},
  journal = {Reviews of Modern Physics},
  volume = {65},
  number = {2},
  pages = {255--299},
  year = {1993},
  publisher = {American Physical Society},
  doi = {10.1103/RevModPhys.65.255},
  url = {https://link.aps.org/doi/10.1103/RevModPhys.65.255}
}

@book{miyamoto1980plasma,
  title={Plasma Physics for Nuclear Fusion},
  author={Miyamoto, Kenro},
  year={1989},
  isbn = {978-0-26-263117-4},
  publisher={The MIT Press},
  note = {Revised Ed.}
}

@article{anderson2001,
  title={A tutorial presentation of the two stream instability and {Landau} damping},
  author={D. Anderson and R. Fedele and M. Lisak},
  journal={American Journal of Physics},
  year={2001},
  doi={10.1119/1.1407252},
url={https://pubs.aip.org/aapt/ajp/article-abstract/69/12/1262/529339/A-tutorial-presentation-of-the-two-stream}
}

@article{Einkemmer2024,
title = {Suppressing instability in a {{Vlasov}}–{{Poisson}} system by an external electric field through constrained optimization},
journal = {Journal of Computational Physics},
volume = {498},
pages = {112662},
year = {2024},
issn = {0021-9991},
doi = {10.1016/j.jcp.2023.112662},
url = {https://www.sciencedirect.com/science/article/pii/S002199912300757X},
author = {Lukas Einkemmer and Qin Li and Li Wang and Yang Yunan}
}

@article{snipes2021iter,
  title={{ITER} plasma control system final design and preparation for first plasma},
  author={Snipes, Joseph A and De Vries, PC and Gribov, Y and Henderson, MA and Hunt, R and Loarte, A and Nunes, I and Pitts, RA and Sinha, J and Zabeo, L and others},
  journal={Nuclear Fusion},
  volume={61},
  number={10},
  pages={106036},
  year={2021},
  publisher={IOP Publishing}
}

@article{hommen2014real,
  title={Real-time optical plasma boundary reconstruction for plasma position control at the {TCV Tokamak}},
  author={Hommen, G d and de Baar, Marco and Duval, BP and Andrebe, Y and Le, HB and Klop, MA and Doelman, NJ and Witvoet, G and Steinbuch, M and TCV Team and others},
  journal={Nuclear Fusion},
  volume={54},
  number={7},
  pages={073018},
  year={2014},
  publisher={IOP Publishing}
}

@article{mouhot2011Landau,
	author = {Mouhot, Cl{\'e}ment and Villani, C{\'e}dric},
	doi = {10.1007/s11511-011-0068-9},
	journal = {Acta Mathematica},
	number = {1},
	pages = {29--201},
	title = {On {Landau} damping},
	volume = {207},
	year = {2011},
}

@article{einkemmer2014convergence,
  title={Convergence analysis of a discontinuous {Galerkin/Strang splitting approximation for the Vlasov--Poisson equations}},
  author={Einkemmer, Lukas and Ostermann, Alexander},
  journal={SIAM Journal on Numerical Analysis},
  volume={52},
  number={2},
  pages={757--778},
  year={2014},
  publisher={SIAM}
}

@article{penrose1960electrostatic,
  title={Electrostatic instabilities of a uniform non-Maxwellian plasma},
  author={Penrose, Oliver},
  journal={The Physics of Fluids},
  volume={3},
  number={2},
  pages={258--265},
  year={1960},
  publisher={American Institute of Physics},
  doi = {10.1063/1.1706024}
}

@article{cheng1976integration,
	author={Cheng, Chio-Zong and Knorr, Georg},
	doi = {10.1016/0021-9991(76)90053-X},
	issn = {0021-9991},
	journal = {Journal of Computational Physics},
	number = {3},
	pages = {330-351},
	title = {The integration of the {Vlasov} equation in configuration space},
	url = {https://www.sciencedirect.com/science/article/pii/002199917690053X},
	volume = {22},
	year = {1976},
    publisher={Elsevier}
}

@article{sonnendrucker1999semi,
  title={The semi-{Lagrangian} method for the numerical resolution of the {{Vlasov}} equation},
  author={Sonnendr{\"u}cker, Eric and Roche, Jean and Bertrand, Pierre and Ghizzo, Alain},
  journal={Journal of Computational Physics},
  volume={149},
  number={2},
  pages={201--220},
  year={1999},
  publisher={Elsevier},
  doi = {10.1006/jcph.1998.6148},
  url = {https://www.sciencedirect.com/science/article/pii/S0021999198961484}
}

@article{albi2025instantaneous,
title = {Instantaneous control strategies for magnetically confined fusion plasma},
journal = {Journal of Computational Physics},
volume = {527},
pages = {113804},
year = {2025},
issn = {0021-9991},
doi = {10.1016/j.jcp.2025.113804},
url = {https://www.sciencedirect.com/science/article/pii/S0021999125000877},
author = {Giacomo Albi and Giacomo Dimarco and Federica Ferrarese and Lorenzo Pareschi}
}

@article{kawata2019dynamic,
  title={{Dynamic stabilization of plasma instability}},
  author={Kawata, S. and Karino, T. and Gu, Y.J.},
  journal={High Power Laser Science and Engineering},
  volume={7},
  pages={e3},
  year={2019},
  publisher={Cambridge University Press},
  doi = {10.1017/hpl.2018.61}
}

@article{filbet2003comparison,
	author = {Filbet, Francis and Sonnendr{\"u}cker, Eric},
	doi = {10.1016/S0010-4655(02)00694-X},
	issn = {0010-4655},
	journal = {Computer Physics Communications},
	number = {3},
	pages = {247-266},
	title = {Comparison of Eulerian {Vlasov} solvers},
	url = {https://www.sciencedirect.com/science/article/pii/S001046550200694X},
	volume = {150},
	year = {2003},
    publisher={Elsevier}
}

@article{Einkemmer2025,
title = {Control of instability in a {Vlasov}-{Poisson} system through an external electric field},
journal = {Journal of Computational Physics},
volume = {530},
pages = {113904},
year = {2025},
issn = {0021-9991},
doi = {10.1016/j.jcp.2025.113904},
url = {https://www.sciencedirect.com/science/article/pii/S0021999125001871},
author = {Lukas Einkemmer and Qin Li and Clément Mouhot and Yukun Yue}
}

@article{Engquist2022,
author = {Engquist, Bj\"orn and Yang, Yunan},
title = {Optimal Transport Based Seismic Inversion:Beyond Cycle Skipping},
journal = {Communications on Pure and Applied Mathematics},
volume = {75},
number = {10},
pages = {2201-2244},
doi = {10.1002/cpa.21990},
url = {https://onlinelibrary.wiley.com/doi/abs/10.1002/cpa.21990},
year = {2022}
}

@article{Strait2014,
    author = {Strait, E. J.},
    title = {Magnetic control of magnetohydrodynamic instabilities in tokamaks},
    journal = {Physics of Plasmas},
    volume = {22},
    number = {2},
    pages = {021803},
    year = {2014},
    month = {11},
    issn = {1070-664X},
    doi = {10.1063/1.4902126},
    url = {https://pubs.aip.org/aip/pop/article-pdf/doi/10.1063/1.4902126/16137245/021803\_1\_online.pdf}
}

@article{Bialek2001,
    author = {Bialek, James and Boozer, Allen H. and Mauel, M. E. and Navratil, G. A.},
    title = {Modeling of active control of external magnetohydrodynamic instabilities},
    journal = {Physics of Plasmas},
    volume = {8},
    number = {5},
    pages = {2170-2180},
    year = {2001},
    month = {05},
    issn = {1070-664X},
    doi = {10.1063/1.1362532},
    url = {https://pubs.aip.org/aip/pop/article-pdf/8/5/2170/19187304/2170\_1\_online.pdf}
}

@article{Klimas1994,
title = {A Splitting Algorithm for {Vlasov} Simulation with Filamentation Filtration},
journal = {Journal of Computational Physics},
volume = {110},
number = {1},
pages = {150-163},
year = {1994},
issn = {0021-9991},
doi = {10.1006/jcph.1994.1011},
url = {https://www.sciencedirect.com/science/article/pii/S0021999184710114},
author = {A.J. Klimas and W.M. Farrell}
}

@article{Klimas2018,
title={Absence of recurrence in {Fourier}–{Fourier} transformed {Vlasov}–{Poisson} simulations},
volume={84},
DOI={10.1017/S0022377818000776},
number={4}, journal={Journal of Plasma Physics},
author={Klimas, Alexander J. and Viñas, Adolfo. F.},
year={2018},
pages={905840405}
}

@article{Rossmanith2011,
title = {A positivity-preserving high-order semi-Lagrangian discontinuous Galerkin scheme for the {Vlasov}–{Poisson} equations},
journal = {Journal of Computational Physics},
volume = {230},
number = {16},
pages = {6203-6232},
year = {2011},
issn = {0021-9991},
doi = {10.1016/j.jcp.2011.04.018},
url = {https://www.sciencedirect.com/science/article/pii/S0021999111002579},
author = {James A. Rossmanith and David C. Seal}
}

@article{Verboncoeur2005,
doi = {10.1088/0741-3335/47/5A/017},
year = {2005},
volume = {47},
number = {5A},
pages = {A231},
author = {Verboncoeur, J P},
title = {Particle simulation of plasmas: review and advances},
journal = {Plasma Physics and Controlled Fusion}
}

@article{Bartsch2024,
author = {Bartsch, Jan and Knopf, Patrik and Scheurer, Stefania and Weber, J\"{o}rg},
title = {Controlling a {Vlasov}–{Poisson} Plasma by a Particle-in-Cell Method Based on a {Monte Carlo} Framework},
journal = {SIAM Journal on Control and Optimization},
volume = {62},
number = {4},
pages = {1977-2011},
year = {2024},
doi = {10.1137/23M1563852}
}

@article{Glass2003,
title = {On the controllability of the {Vlasov}–{Poisson} system},
journal = {Journal of Differential Equations},
volume = {195},
number = {2},
pages = {332-379},
year = {2003},
issn = {0022-0396},
doi = {10.1016/S0022-0396(03)00066-4},
url = {https://www.sciencedirect.com/science/article/pii/S0022039603000664},
author = {Olivier Glass}
}

@article{Glass2012,
title = {On the controllability of the {Vlasov}–{Poisson} system in the presence of external force fields},
journal = {Journal of Differential Equations},
volume = {252},
number = {10},
pages = {5453-5491},
year = {2012},
issn = {0022-0396},
doi = {10.1016/j.jde.2012.02.007},
url = {https://www.sciencedirect.com/science/article/pii/S0022039612000848},
author = {Olivier Glass and Daniel Han-Kwan}
}

@article{Knopf2018,
title = {Optimal control of a {Vlasov}–{Poisson} plasma by an external magnetic field},
author = {Knopf, Patrik},
journal = {Calculus of Variations and Partial Differential Equations},
volume = {57},
number = {5},
year = {2018},
doi = {10.1007/s00526-018-1407-x}
}

@article{Knopf2019,
author = {Knopf, Patrik},
title = {Confined steady states of a {Vlasov}-{Poisson} plasma in an infinitely long cylinder},
journal = {Mathematical Methods in the Applied Sciences},
volume = {42},
number = {18},
pages = {6369-6384},
doi = {10.1002/mma.5728},
url = {https://onlinelibrary.wiley.com/doi/abs/10.1002/mma.5728},
year = {2019}
}

@misc{jax2018github,
  author = {James Bradbury and Roy Frostig and Peter Hawkins and Matthew James Johnson and Chris Leary and Dougal Maclaurin and George Necula and Adam Paszke and Jake Vander{P}las and Skye Wanderman-{M}ilne and Qiao Zhang},
  title = {{JAX}: composable transformations of {P}ython+{N}um{P}y programs},
  url = {http://github.com/jax-ml/jax},
  note = {0.6.0},
  year = {2018},
}

@article{Symes1991,
author = {W. W. Symes and J. J. Carazzone},
title = {Velocity inversion by differential semblance optimization},
journal = {GEOPHYSICS},
volume = {56},
number = {5},
pages = {654-663},
year = {1991},
doi = {10.1190/1.1443082},
URL = {https://library.seg.org/doi/10.1190/1.1443082}
}

@article{ChenDing2023,
	author = {Chen, Shi and Ding, Zhiyan and Li, Qin and Zepeda-N\'{u}\~{n}ez, Leonardo},
	doi = {10.1137/22M147075X},
	journal = {SIAM Journal on Imaging Sciences},
	number = {1},
	pages = {111-143},
	title = {High-Frequency Limit of the Inverse Scattering Problem: Asymptotic Convergence from Inverse Helmholtz to Inverse Liouville},
	volume = {16},
	year = {2023}
}

@article{nature2022,
	author = {Degrave, Jonas and Felici, Federico and Buchli, Jonas and Neunert, Michael and Tracey, Brendan and Carpanese, Francesco and Ewalds, Timo and Hafner, Roland and Abdolmaleki, Abbas and de las Casas, Diego and Donner, Craig and Fritz, Leslie and Galperti, Cristian and Huber, Andrea and Keeling, James and Tsimpoukelli, Maria and Kay, Jackie and Merle, Antoine and Moret, Jean-Marc and Noury, Seb and Pesamosca, Federico and Pfau, David and Sauter, Olivier and Sommariva, Cristian and Coda, Stefano and Duval, Basil and Fasoli, Ambrogio and Kohli, Pushmeet and Kavukcuoglu, Koray and Hassabis, Demis and Riedmiller, Martin},
	doi = {10.1038/s41586-021-04301-9},
	journal = {Nature},
	number = {7897},
	pages = {414--419},
	title = {Magnetic control of tokamak plasmas through deep reinforcement learning},
	volume = {602},
	year = {2022}
}

@article{nature2024,
	author = {Seo, Jaemin and Kim, SangKyeun and Jalalvand, Azarakhsh and Conlin, Rory and Rothstein, Andrew and Abbate, Joseph and Erickson, Keith and Wai, Josiah and Shousha, Ricardo and Kolemen, Egemen},
	doi = {10.1038/s41586-024-07024-9},
	journal = {Nature},
	number = {8000},
	pages = {746--751},
	title = {Avoiding fusion plasma tearing instability with deep reinforcement learning},
	volume = {626},
	year = {2024}
}

@inproceedings{Seo2023,
  author={Seo, Jaemin and Conlin, Rory and Rothstein, Andrew and Kim, SangKyeun and Abbate, Joseph and Jalalvand, Azarakhsh and Kolemen, Egemen},
  booktitle={2023 International Joint Conference on Neural Networks (IJCNN)}, 
  title={Multimodal Prediction of Tearing Instabilities in a Tokamak}, 
  year={2023},
  volume={},
  number={},
  pages={1-8},
  isbn = {978-1-6654-8867-9},
  doi={10.1109/IJCNN54540.2023.10191359}
}

@inproceedings{Blum2016,
	author = {Blum, Jacques and Boulbe, C{\'e}dric and Faugeras, Blaise and Heumann, Holger},
	booktitle = {IFIP TC7 2015 - 27th Conference on System Modeling and
Optimization},
	editor = {Bociu, Lorena and D{\'e}sid{\'e}ri, Jean-Antoine and Habbal, Abderrahmane},
	isbn = {978-3-319-55795-3},
    doi = {10.1007/978-3-319-55795-3_1},
	pages = {1--20},
	publisher = {Springer International Publishing},
	title = {Control Methods for the Optimization of Plasma Scenarios in a Tokamak},
	year = {2016}
}

@misc{albi2025robustfeedbackcontrolcollisional,
      title={Robust feedback control of collisional plasma dynamics in presence of uncertainties}, 
      author={Giacomo Albi and Giacomo Dimarco and Federica Ferrarese and Lorenzo Pareschi},
      year={2025},
      eprint={2505.19992},
      archivePrefix={arXiv},
      primaryClass={math.NA},
      url={https://arxiv.org/abs/2505.19992}, 
}

@article{HanKwan2021,
	author = {Han-Kwan, Daniel and Nguyen, Toan T. and Rousset, Fr{\'e}d{\'e}ric},
	doi = {10.1007/s40818-021-00110-5},
	journal = {Annals of PDE},
	number = {2},
	pages = {Paper No. 18, 37 pp},
	title = {Asymptotic Stability of Equilibria for Screened {Vlasov}--{Poisson} Systems via Pointwise Dispersive Estimates},
	volume = {7},
	year = {2021},
}

@article{Moret1988,
    author = {Moret, J.-M. and Buhlmann, F. and Fasel, D. and Hofmann, F. and Tonetti, G.},
    title = {Magnetic measurements on the TCV Tokamak},
    journal = {Review of Scientific Instruments},
    volume = {69},
    number = {6},
    pages = {2333-2348},
    year = {1998},
    month = {06},
    issn = {0034-6748},
    doi = {10.1063/1.1148940},
    url = {https://pubs.aip.org/aip/rsi/article-pdf/69/6/2333/19105874/2333\_1\_online.pdf},
}

@book{Freidberg2007,
    title={Plasma Physics and Fusion Energy},
    publisher={Cambridge University Press},
    author={Freidberg, Jeffrey P.},
    year={2007},
    isbn = {9780511755705},
    doi = {10.1017/CBO9780511755705}
}

@book{NocedalNumOpt,
  address = {New York, NY},
  author = {Nocedal, Jorge and Wright, Stephen J.},
  note = {2nd Ed.},
  isbn = {978-0-387-30303-1},
  publisher = {Springer},
  series = {Springer Series in Operations Research and Financial Engineering},
  title = {Numerical Optimization},
  year = {2006},
  doi = {10.1007/978-0-387-40065-5}
}

\appendix

\section{Comparison between \texorpdfstring{$L^{2}$}{L2} and KL}\label{sec:L2_KL}

Another possible function is to use the $L^{2}$ disparity:
\begin{equation}\tag{$L^{2}$}\label{eq:L2_obj}
    \mathcal{J}(f[H]) = \frac{1}{2}\|f_{T}[H] - f_{\text{eq}}\|_{L^{2}(x,v)}^{2} = \frac{1}{2}\int_{-L_{v}}^{L_{v}}(f[H](T,x,v) - f_{\text{eq}}(v))^{2}\mathrm{d}x\mathrm{d}v
\end{equation}
\begin{equation}\tag{$L^{2}$T}\label{eq:L2T_obj}
    \mathcal{J}(f[H]) = \frac{1}{2}\int_{0}^{T}\int_{-L_{v}}^{L_{v}}(f[H](t,x,v) - f_{\text{eq}}(v))^{2}\mathrm{d}x\mathrm{d}v\mathrm{d}t\,.
\end{equation}

In our simulation, we find the numerical behavior of these objective functions are rather similar to those of KL (\eqref{eq:KL_obj} and \eqref{eq:KLT_obj}). We provide numerical evidence in this section.

\begin{figure}[ht]
    \centering
    \includegraphics[width=1.0\linewidth]{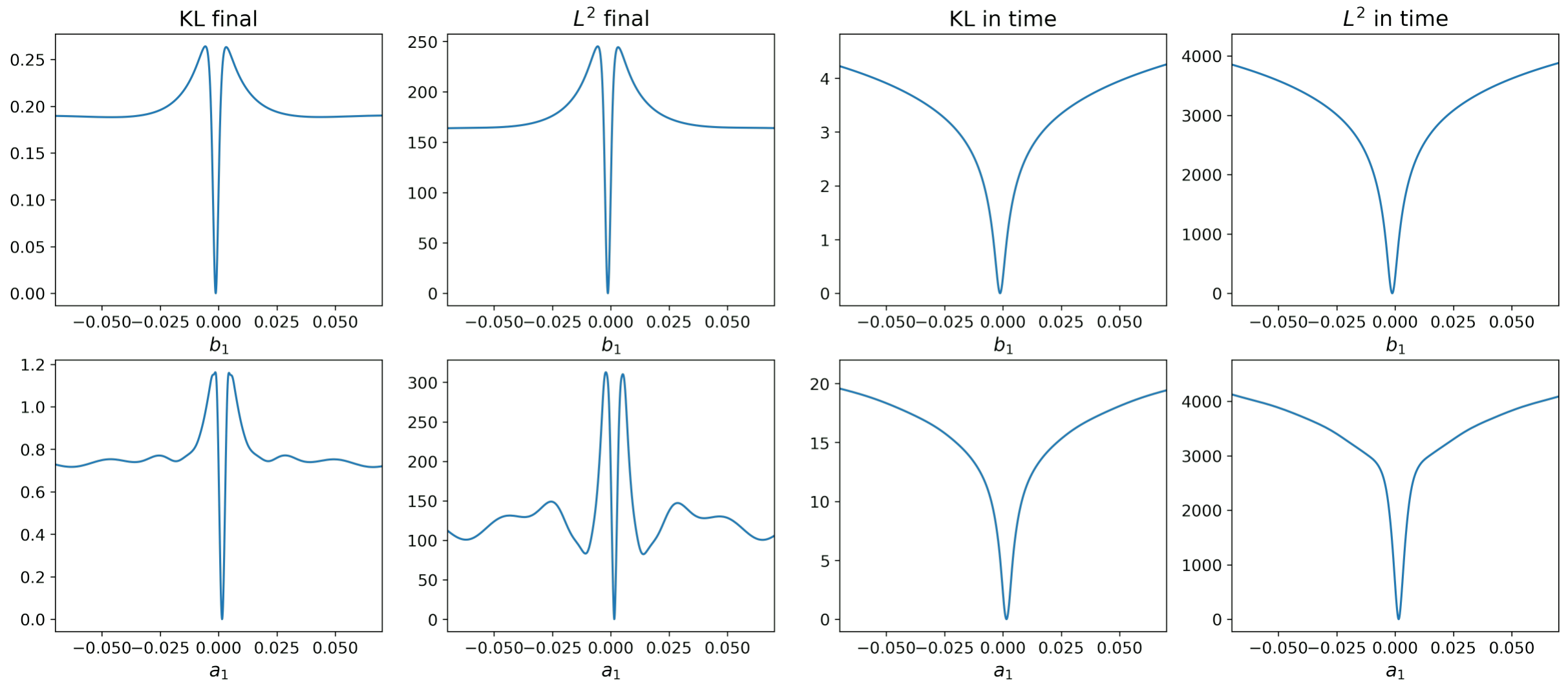}
    \caption{Landscape of Two Stream instability for $b_{1}$ for~\eqref{eq:H_sin}(top) and Bump-on-Tail instability for $a_{1}$ for~\eqref{eq:H_cos}(bottom) on the domain $[0.07,0.07]$ of the objectives~\eqref{eq:KL_obj}(left), \eqref{eq:L2_obj}(center-left), \eqref{eq:KLT_obj}(center-right) and \eqref{eq:L2T_obj}(right).}
    \label{fig:1D_L2_KL}
\end{figure}

In Figure~\ref{fig:1D_L2_KL}, we show the one-dimensional case for both examples. For the Two Stream case, both landscapes exhibit no apparent differences, for both integrated-in-time case, and the disparity computed from the last frame. The situation for the Bump-on-Tail case is similar, with landscape for \eqref{eq:L2_obj} agreeing with that for \eqref{eq:KL_obj} and \eqref{eq:L2T_obj} agreeing with \eqref{eq:KLT_obj}. We should note, however,  \eqref{eq:L2_obj} seems to show deeper local minima, making the computation more difficult to find the global optimal solution.

\section{Hessian computations of objective functions}\label{sec:Hess_computations}

\begin{figure}[H]
    \centering
    \includegraphics[width=1.0\linewidth]{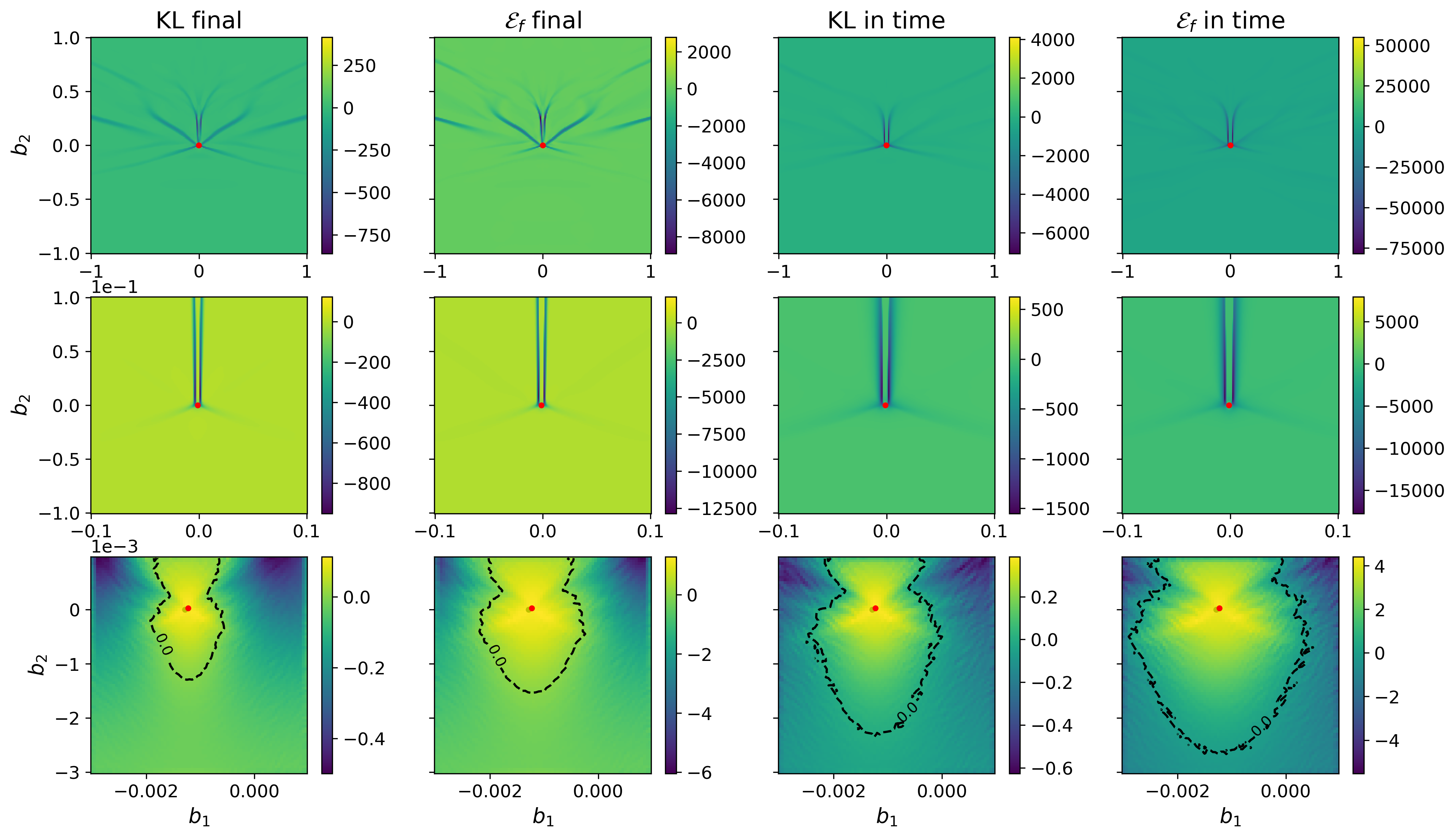}
    \caption{Landscape of the minimum eigenvalue of the Hessian of Two Stream instability with applied control over the domain \([-1,1]^{2}\)(top), \([-0.1,0.1]^{2}\)(center) and, \([-0.003,0.001]^{2}\)(bottom) for \((b_{1}, b_{2})\) for~\eqref{eq:H_sins}. The objective functions are~\eqref{eq:KL_obj}(left),~\eqref{eq:EE_obj}(center-left),~\eqref{eq:KLT_obj}(center-right) and~\eqref{eq:EET_obj}(right). Yellow dot from last row represents good initial guess from Example~\ref{ex:initial_guess}, red dot represents approximate global minimum and black dashed line represents the $0$ level set around the global minimum. Close to the global minima, Hessians are positive definite, meaning the objective functions are convex.}
    \label{fig:landscape_two_stream_2D_hess}
\end{figure}

\begin{figure}[H]
    \centering
    \includegraphics[width=1.0\linewidth]{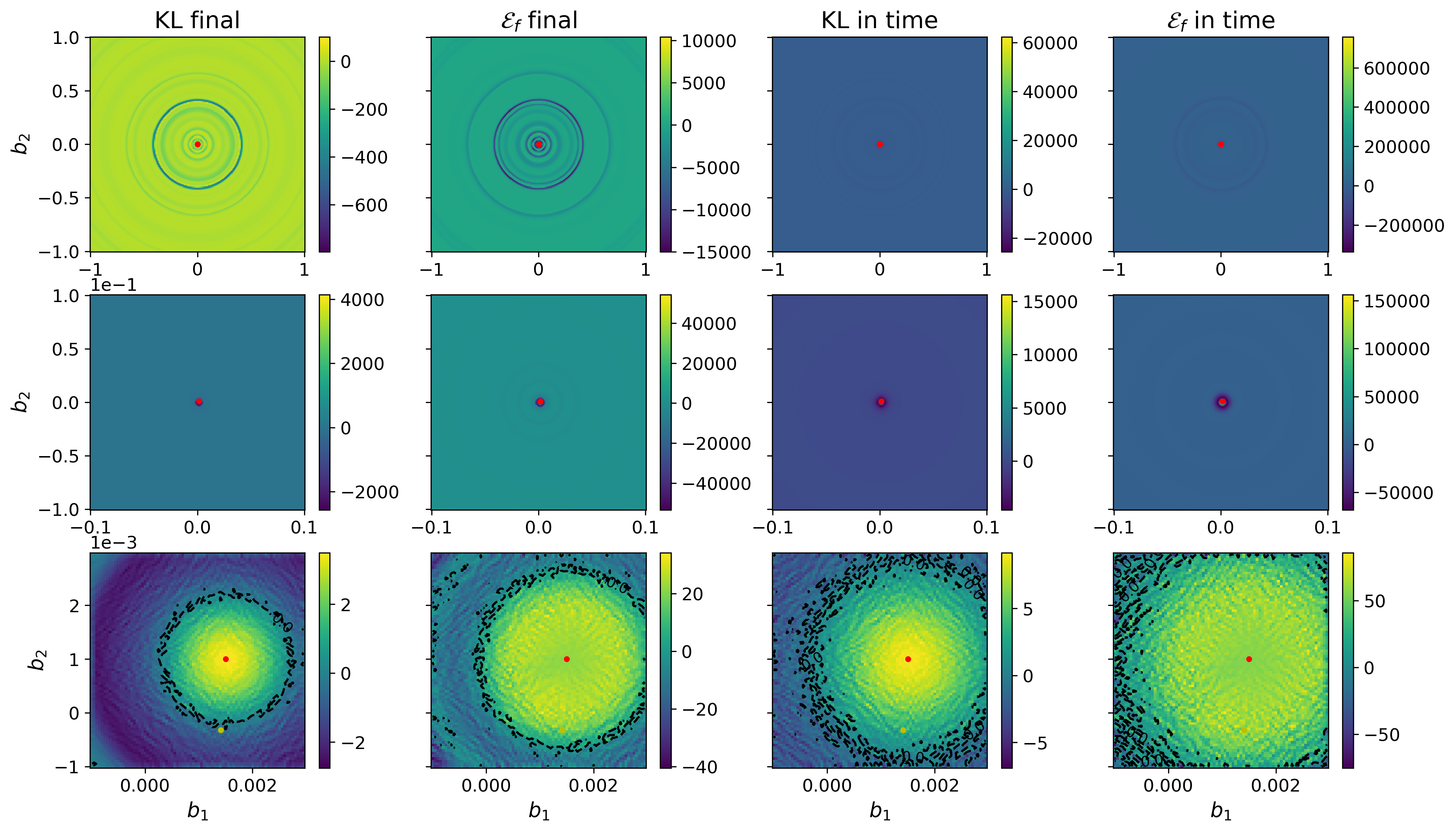}
    \caption{Landscape of the minimum eigenvalue of the Hessian of Bump-on-Tail instability with applied control over the domain \([-1,1]^{2}\)(top), \([-0.1,0.1]^{2}\)(center) and, \([-0.001,0.003]^{2}\)(bottom) for \((a_{1}, b_{1})\) for~\eqref{eq:H_cos_sin}. The objective functions are~\eqref{eq:KL_obj}(left),~\eqref{eq:EE_obj}(center-left),~\eqref{eq:KLT_obj}(center-right) and~\eqref{eq:EET_obj}(right). Yellow dot from last row represents good initial guess from Example~\ref{ex:initial_guess}, red dot represents approximate global minimum and black dashed line represents the $0$ level set around the global minimum. Close to the global minima, Hessians are positive definite, meaning the objective functions are convex.}
    \label{fig:landscape_bump-on-tail_2D_hess}
\end{figure}

In this section we present the Hessian analysis for the objective functions~\eqref{eq:KL_obj},~\eqref{eq:EE_obj},~\eqref{eq:KLT_obj}, and~\eqref{eq:EET_obj}. In particular, we compute the Hessian of each objective with respect to the control parameters and plot the minimum eigenvalue of the resulting Hessian matrices. Figures~\ref{fig:landscape_two_stream_2D_hess} and~\ref{fig:landscape_bump-on-tail_2D_hess} show the results for the Two-Stream and Bump-on-Tail instabilities, respectively.

\section{Summary of results}
In this section we present in detail the results seen in Tables~\ref{tab:TS_results}  and \ref{tab:BoT_results}.

\subsection{Two Stream example}\label{sec:TS_summary}

\begin{figure}[H]
    \centering
    \includegraphics[width=0.85\linewidth]{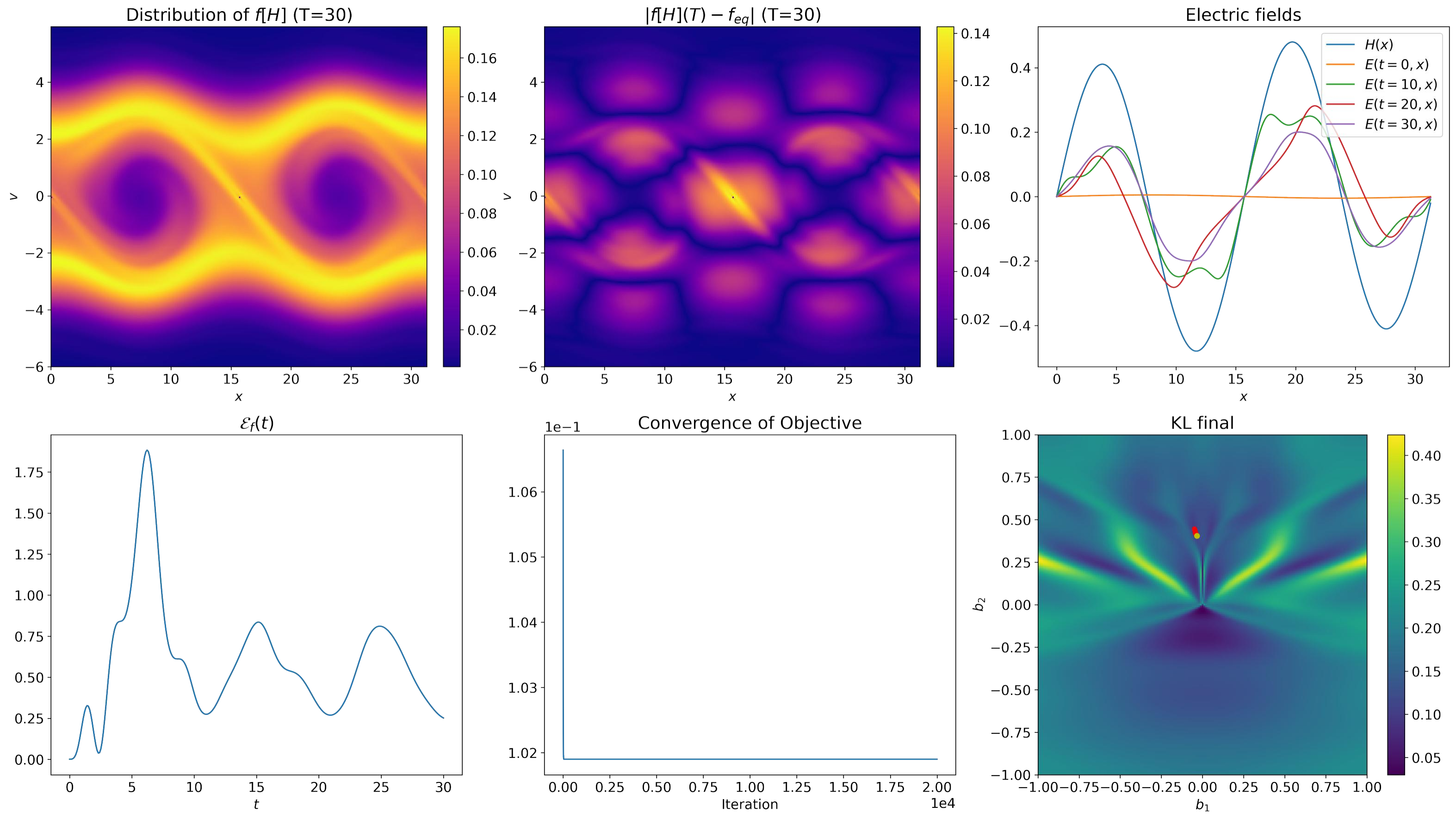}
    \caption{Simulation of~\eqref{eq:vlasov-poisoon_system_ext_1d} with under-parametrized $H$ obtained from~\eqref{eq:optimization_pb_simple} using~\eqref{eq:KL_obj} with far initialization using GD with line-search. From left to right and top to bottom: $f[H](T=30,x,v)$, $|f[H](T,x,v)-f_{\text{eq}}(v)|$, $H$ and $E_{f[H]}(t,x)$, $\mathcal{E}_{f[H]}(t)$, convergence of objective and, trajectory over the landscape of the objective (yellow dot is initial guess).}
    \label{fig:TS_KL_GDL_far_under}
\end{figure}

\begin{figure}[H]
    \centering
    \includegraphics[width=0.85\linewidth]{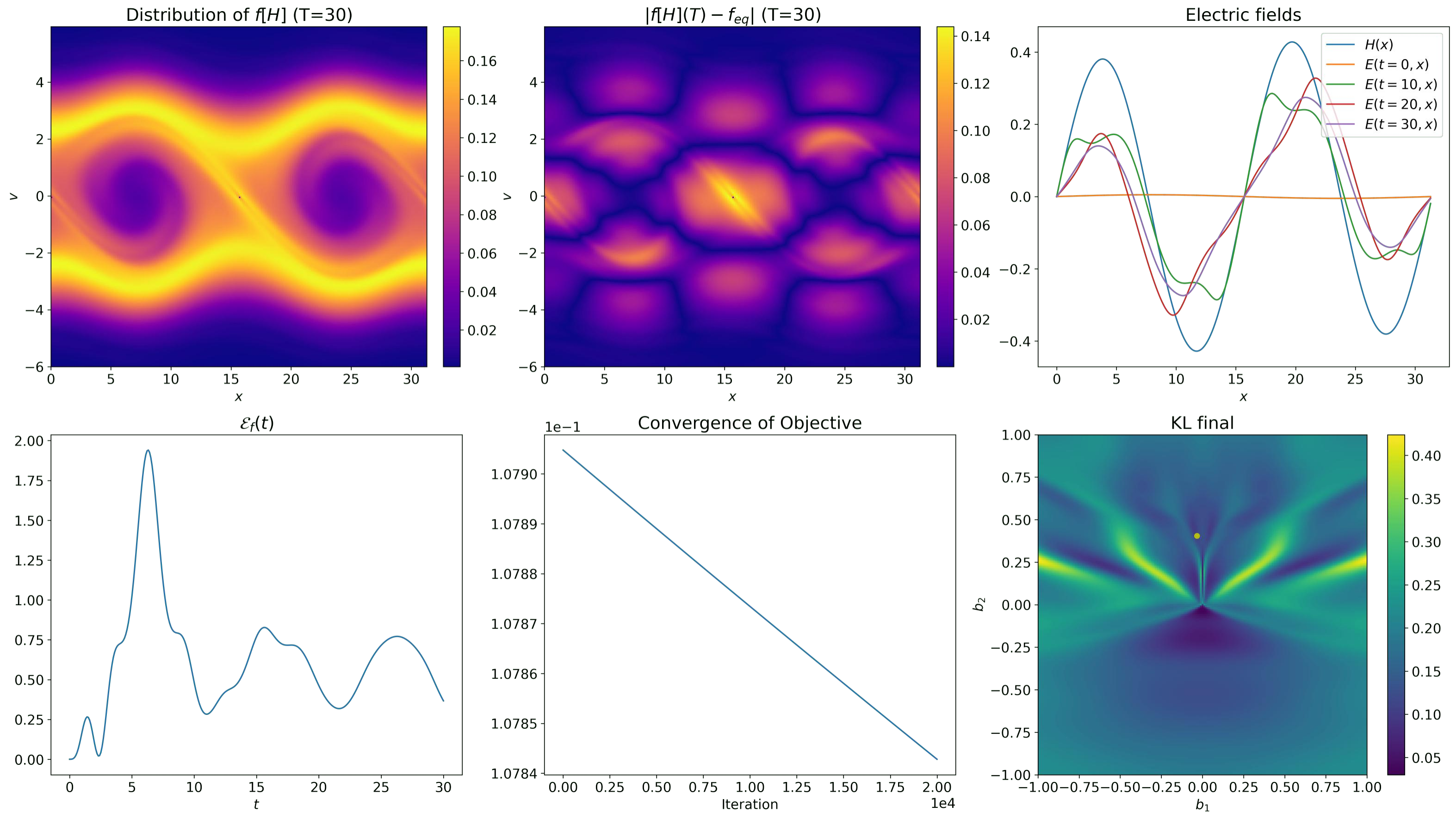}
    \caption{Simulation of~\eqref{eq:vlasov-poisoon_system_ext_1d} with under-parametrized $H$ obtained from~\eqref{eq:optimization_pb_simple} using~\eqref{eq:KL_obj} with far initialization using GD with constant stepsize. From left to right and top to bottom: $f[H](T=30,x,v)$, $|f[H](T,x,v)-f_{\text{eq}}(v)|$, $H$ and $E_{f[H]}(t,x)$, $\mathcal{E}_{f[H]}(t)$, convergence of objective and, trajectory over the landscape of the objective (yellow dot is initial guess).}
    \label{fig:TS_KL_GD_far_under}
\end{figure}

\begin{figure}[H]
    \centering
    \includegraphics[width=0.85\linewidth]{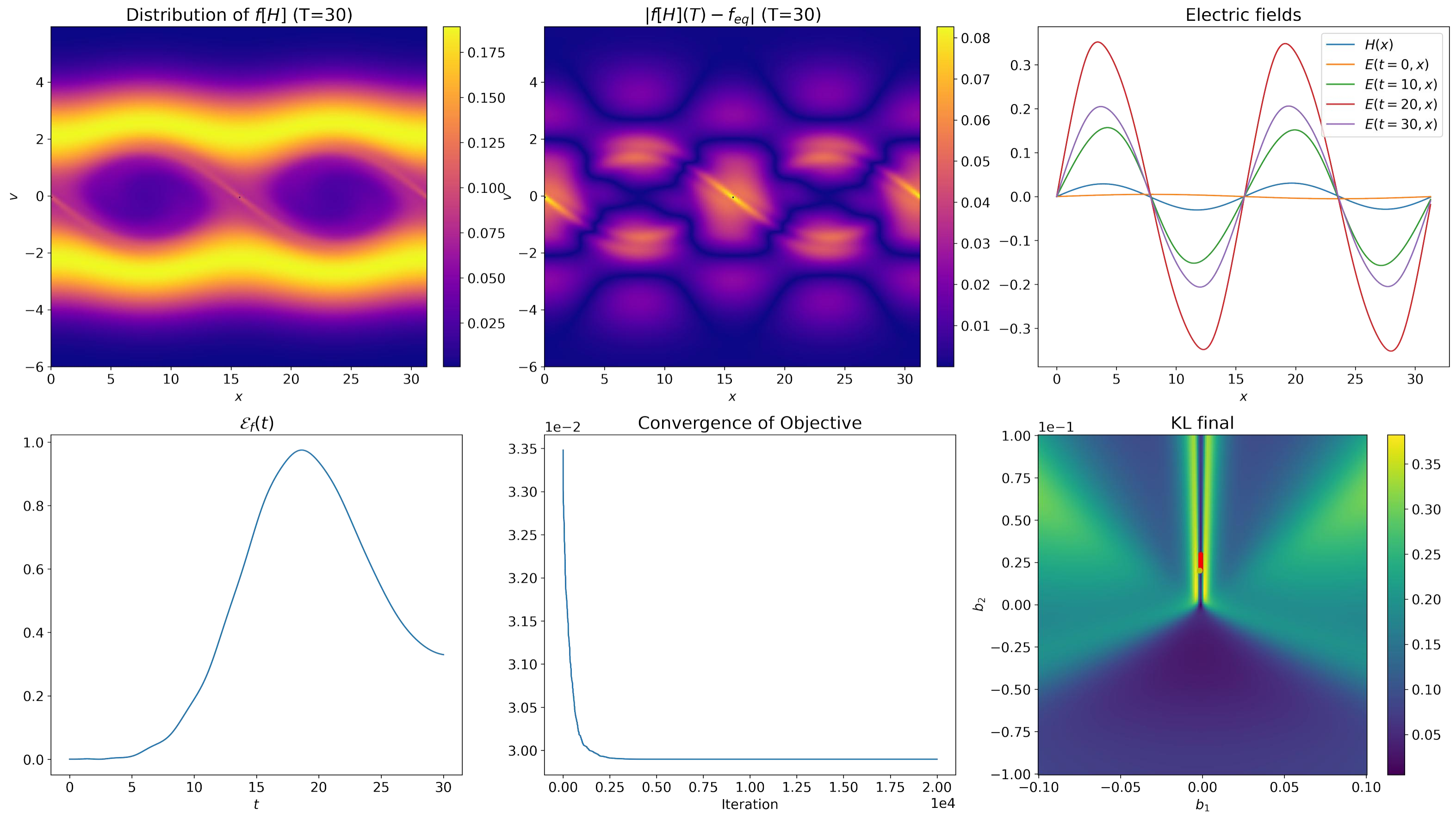}
    \caption{Simulation of~\eqref{eq:vlasov-poisoon_system_ext_1d} with under-parametrized $H$ obtained from~\eqref{eq:optimization_pb_simple} using~\eqref{eq:KL_obj} with near initialization using GD with line-search. From left to right and top to bottom: $f[H](T=30,x,v)$, $|f[H](T,x,v)-f_{\text{eq}}(v)|$, $H$ and $E_{f[H]}(t,x)$, $\mathcal{E}_{f[H]}(t)$, convergence of objective and, trajectory over the landscape of the objective (yellow dot is initial guess).}
    \label{fig:TS_KL_GDL_near_under}
\end{figure}

\begin{figure}[H]
    \centering
    \includegraphics[width=0.85\linewidth]{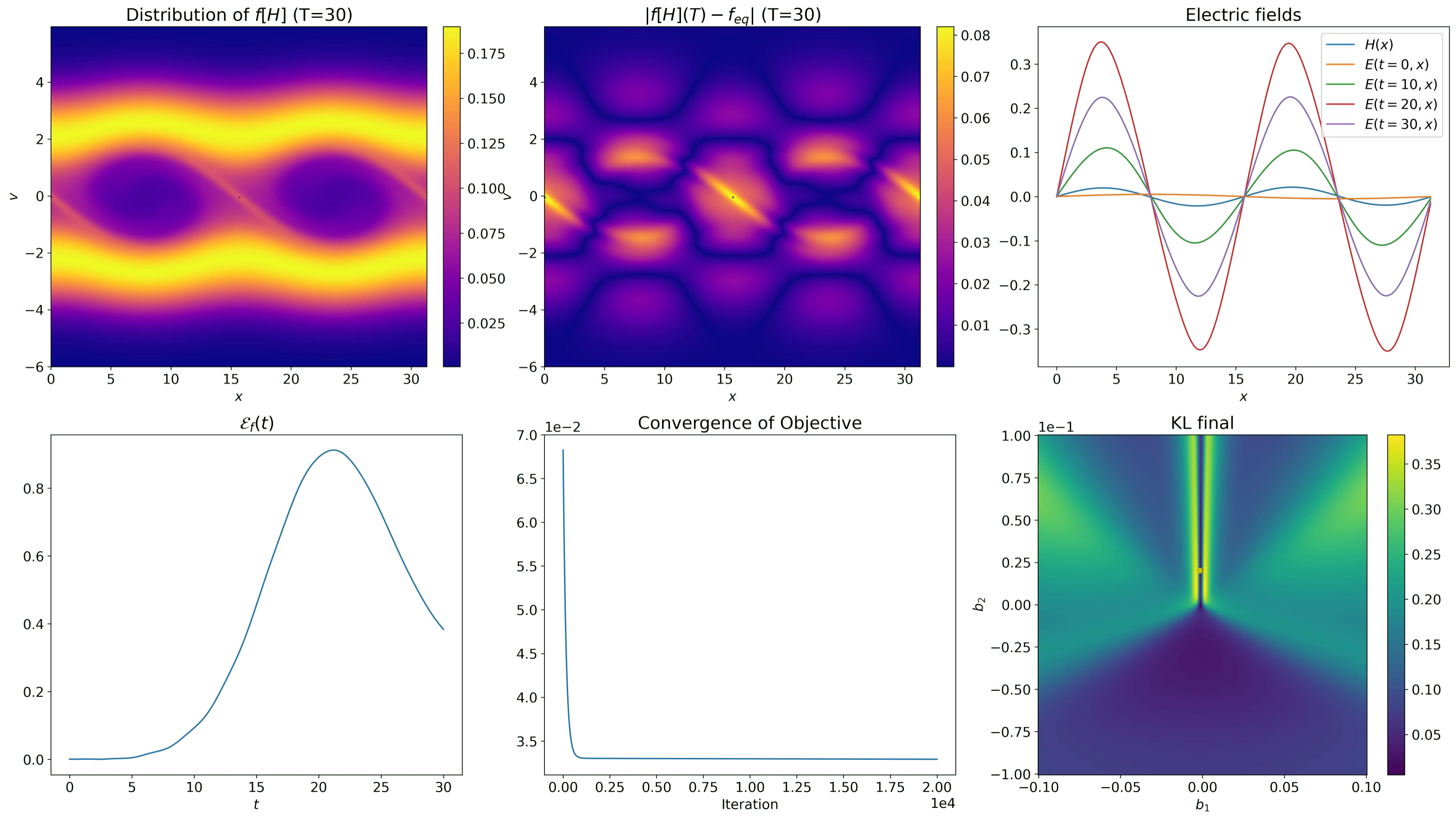}
    \caption{Simulation of~\eqref{eq:vlasov-poisoon_system_ext_1d} with under-parametrized $H$ obtained from~\eqref{eq:optimization_pb_simple} using~\eqref{eq:KL_obj} with near initialization using GD with constant stepsize. From left to right and top to bottom: $f[H](T=30,x,v)$, $|f[H](T,x,v)-f_{\text{eq}}(v)|$, $H$ and $E_{f[H]}(t,x)$, $\mathcal{E}_{f[H]}(t)$, convergence of objective and, trajectory over the landscape of the objective (yellow dot is initial guess).}
    \label{fig:TS_KL_GD_near_under}
\end{figure}

\begin{figure}[H]
    \centering
    \includegraphics[width=0.85\linewidth]{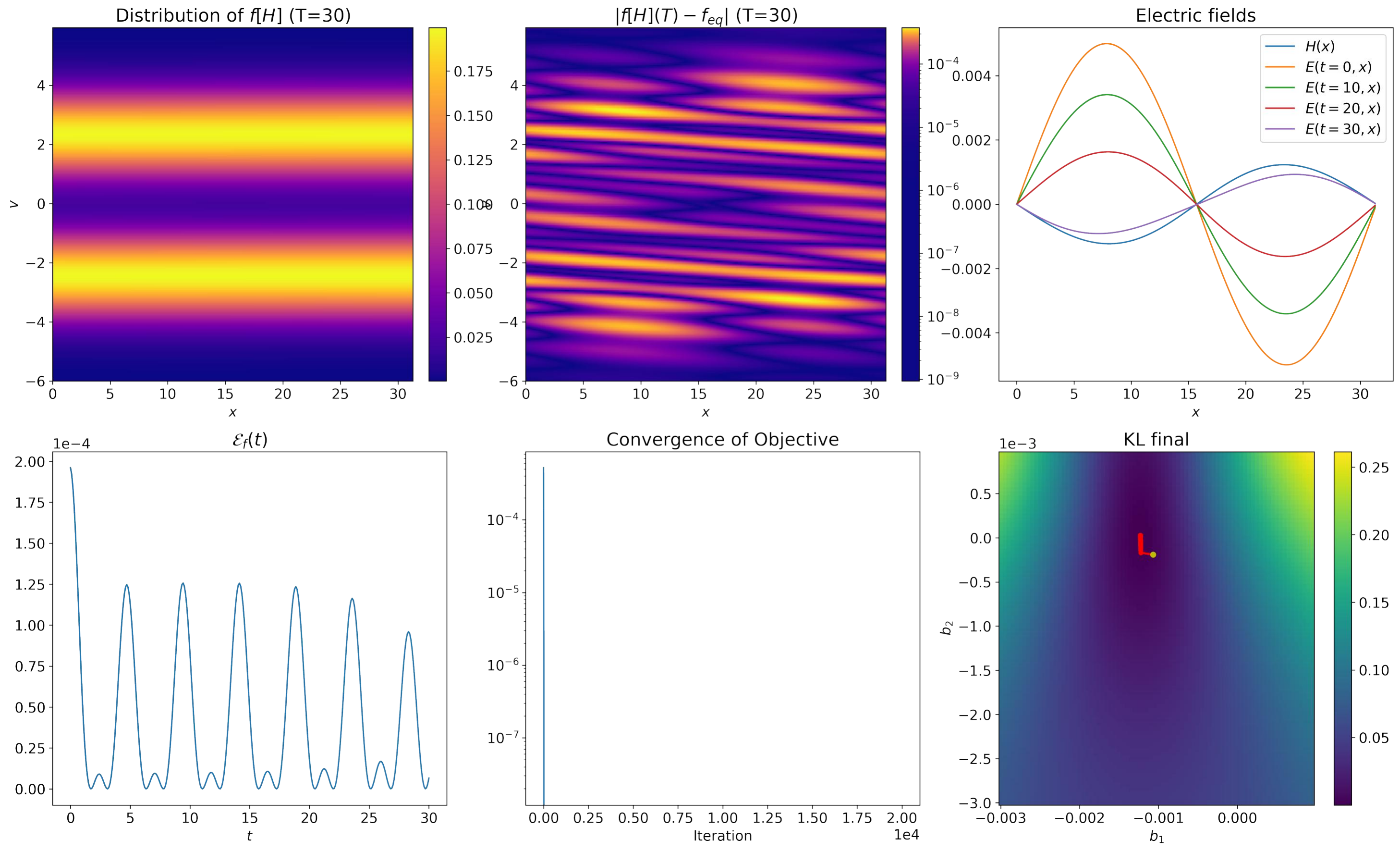}
    \caption{Simulation of~\eqref{eq:vlasov-poisoon_system_ext_1d} with under-parametrized $H$ obtained from~\eqref{eq:optimization_pb_simple} using~\eqref{eq:KL_obj} with local initialization using GD with line-search. From left to right and top to bottom: $f[H](T=30,x,v)$, $|f[H](T,x,v)-f_{\text{eq}}(v)|$,  $H$ and $E_{f[H]}(t,x)$, $\mathcal{E}_{f[H]}(t)$, convergence of objective and, trajectory over the landscape of the objective (yellow dot is initial guess).}
    \label{fig:TS_KL_GDL_local_under}
\end{figure}

\begin{figure}[H]
    \centering
    \includegraphics[width=0.85\linewidth]{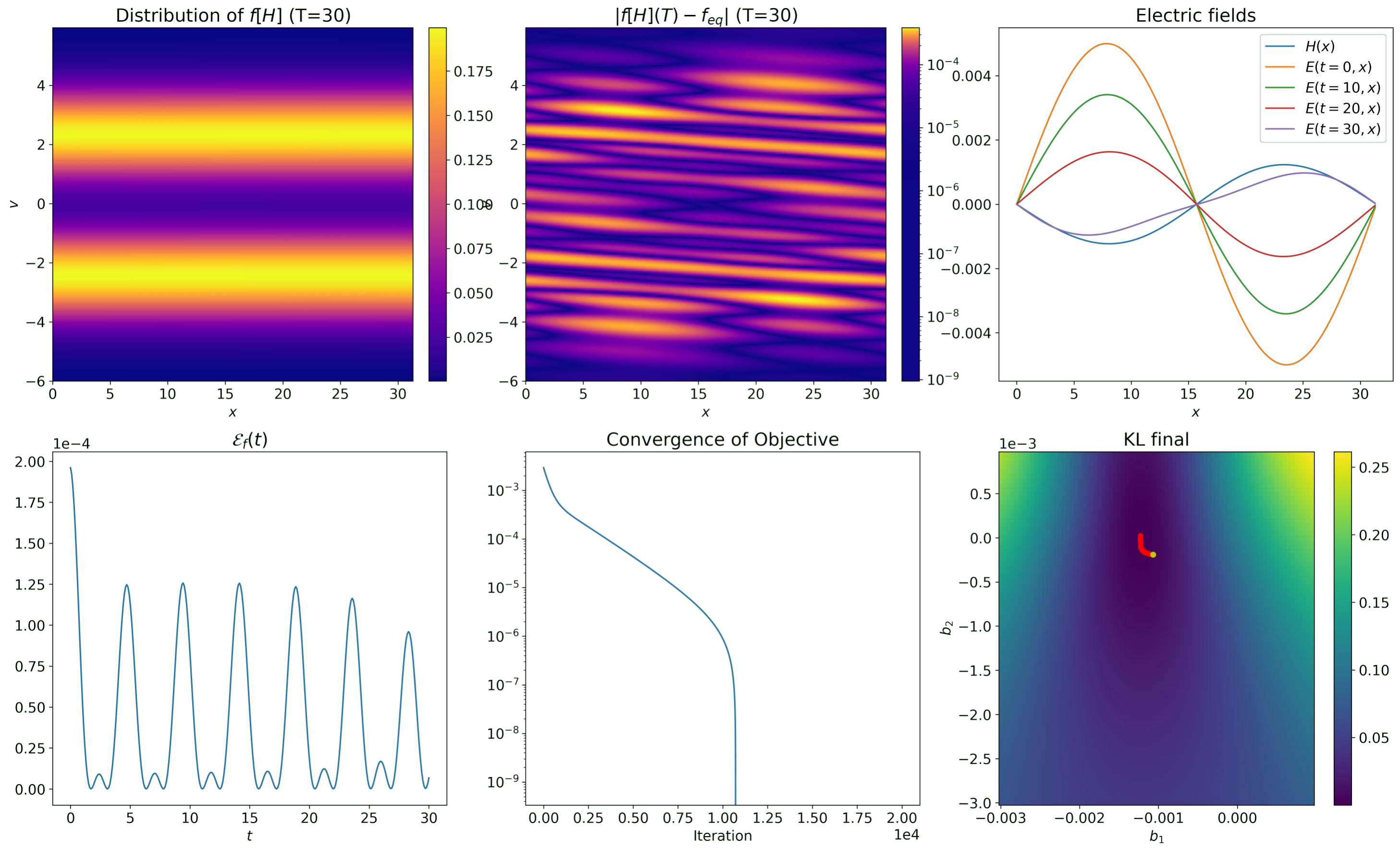}
    \caption{Simulation of~\eqref{eq:vlasov-poisoon_system_ext_1d} with under-parametrized $H$ obtained from~\eqref{eq:optimization_pb_simple} using~\eqref{eq:KL_obj} with local initialization using GD with constant stepsize. From left to right and top to bottom: $f[H](T=30,x,v)$, $|f[H](T,x,v)-f_{\text{eq}}(v)|$, $H$ and $E_{f[H]}(t,x)$, $\mathcal{E}_{f[H]}(t)$, convergence of objective and, trajectory over the landscape of the objective (yellow dot is initial guess).}
    \label{fig:TS_KL_GD_local_under}
\end{figure}

\begin{figure}[H]
    \centering
    \includegraphics[width=0.85\linewidth]{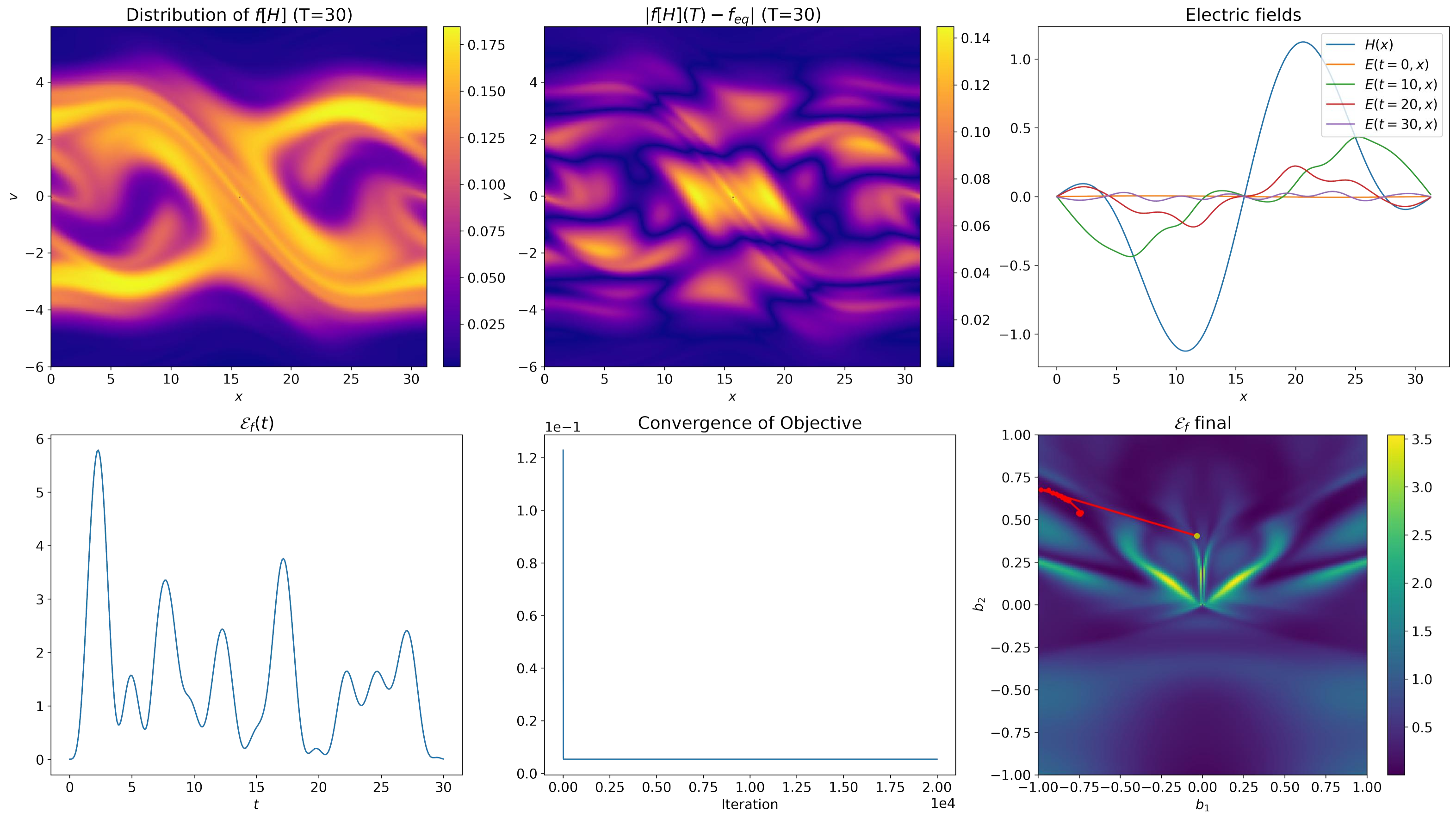}
    \caption{Simulation of~\eqref{eq:vlasov-poisoon_system_ext_1d} with under-parametrized $H$ obtained from~\eqref{eq:optimization_pb_simple} using~\eqref{eq:EE_obj} with far initialization using GD with line-search. From left to right and top to bottom: $f[H](T=30,x,v)$, $|f[H](T,x,v)-f_{\text{eq}}(v)|$, $H$ and $E_{f[H]}(t,x)$, $\mathcal{E}_{f[H]}(t)$, convergence of objective and, trajectory over the landscape of the objective (yellow dot is initial guess).}
    \label{fig:TS_ee_lf_GDL_far_under}
\end{figure}

\begin{figure}[H]
    \centering
    \includegraphics[width=0.85\linewidth]{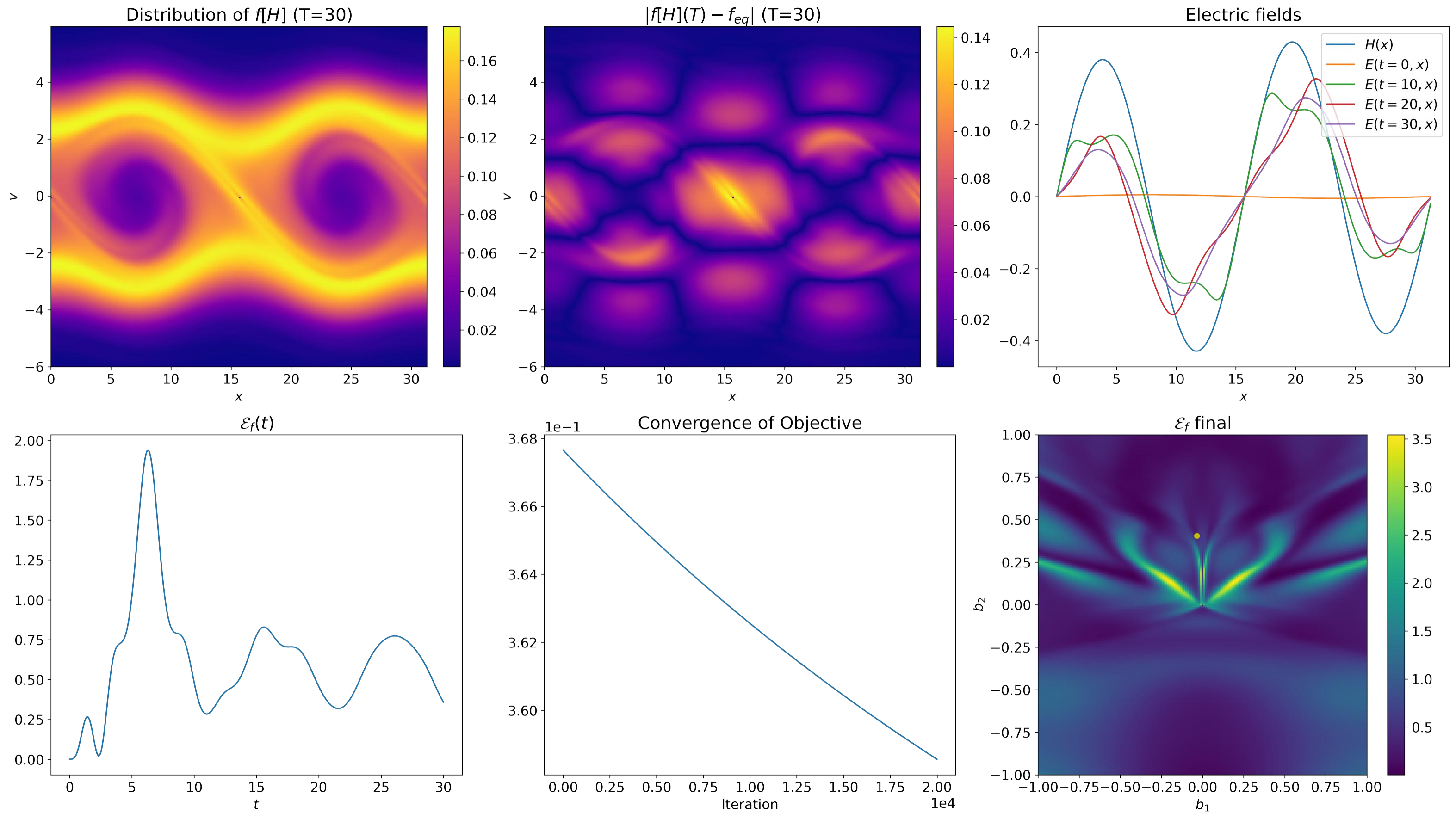}
    \caption{Simulation of~\eqref{eq:vlasov-poisoon_system_ext_1d} with under-parametrized $H$ obtained from~\eqref{eq:optimization_pb_simple} using~\eqref{eq:EE_obj} with far initialization using GD with constant stepsize. From left to right and top to bottom: $f[H](T=30,x,v)$, $|f[H](T,x,v)-f_{\text{eq}}(v)|$, $H$ and $E_{f[H]}(T,x)$, $E_{f[H]}(t,x)$, $\mathcal{E}_{f[H]}(t)$ and, convergence of objective.}
    \label{fig:TS_ee_lf_GD_far_under}
\end{figure}

\begin{figure}[H]
    \centering
    \includegraphics[width=0.85\linewidth]{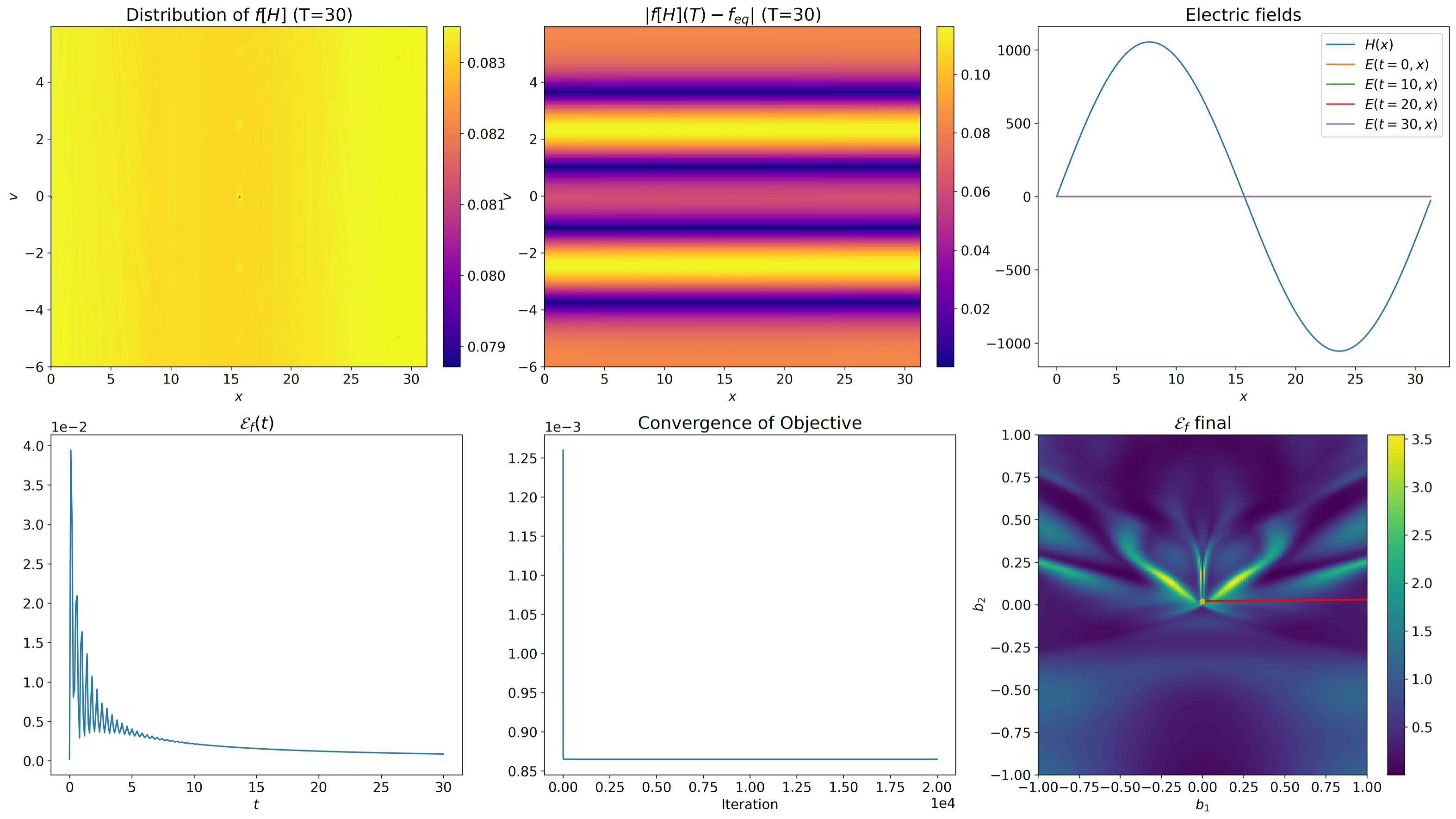}
    \caption{Simulation of~\eqref{eq:vlasov-poisoon_system_ext_1d} with under-parametrized $H$ obtained from~\eqref{eq:optimization_pb_simple} using~\eqref{eq:EE_obj} with near initialization using GD with line-search. From left to right and top to bottom: $f[H](T=30,x,v)$, $|f[H](T,x,v)-f_{\text{eq}}(v)|$, $H$ and $E_{f[H]}(t,x)$, $\mathcal{E}_{f[H]}(t)$, convergence of objective and, trajectory over the landscape of the objective (yellow dot is initial guess).}
    \label{fig:TS_ee_lf_GDL_near_under}
\end{figure}

\begin{figure}[H]
    \centering
    \includegraphics[width=0.85\linewidth]{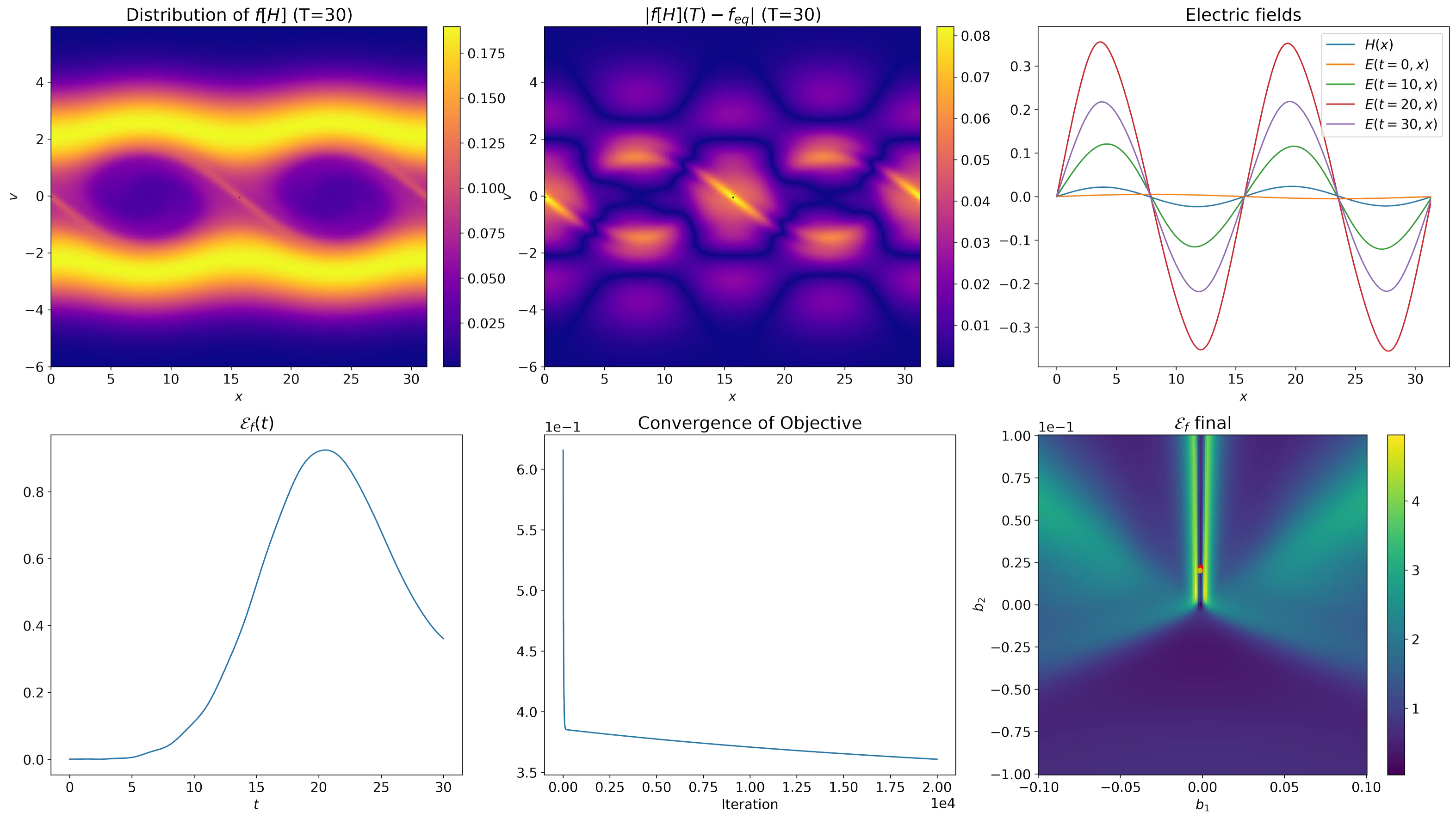}
    \caption{Simulation of~\eqref{eq:vlasov-poisoon_system_ext_1d} using~\eqref{eq:EE_obj} with under-parametrized $H$ obtained from~\eqref{eq:optimization_pb_simple} with near initialization using GD with constant stepsize. From left to right and top to bottom: $f[H](T=30,x,v)$, $|f[H](T,x,v)-f_{\text{eq}}(v)|$, $H$ and $E_{f[H]}(t,x)$, $\mathcal{E}_{f[H]}(t)$, convergence of objective and, trajectory over the landscape of the objective (yellow dot is initial guess).}
    \label{fig:TS_ee_lf_GD_near_under}
\end{figure}

\begin{figure}[H]
    \centering
    \includegraphics[width=0.85\linewidth]{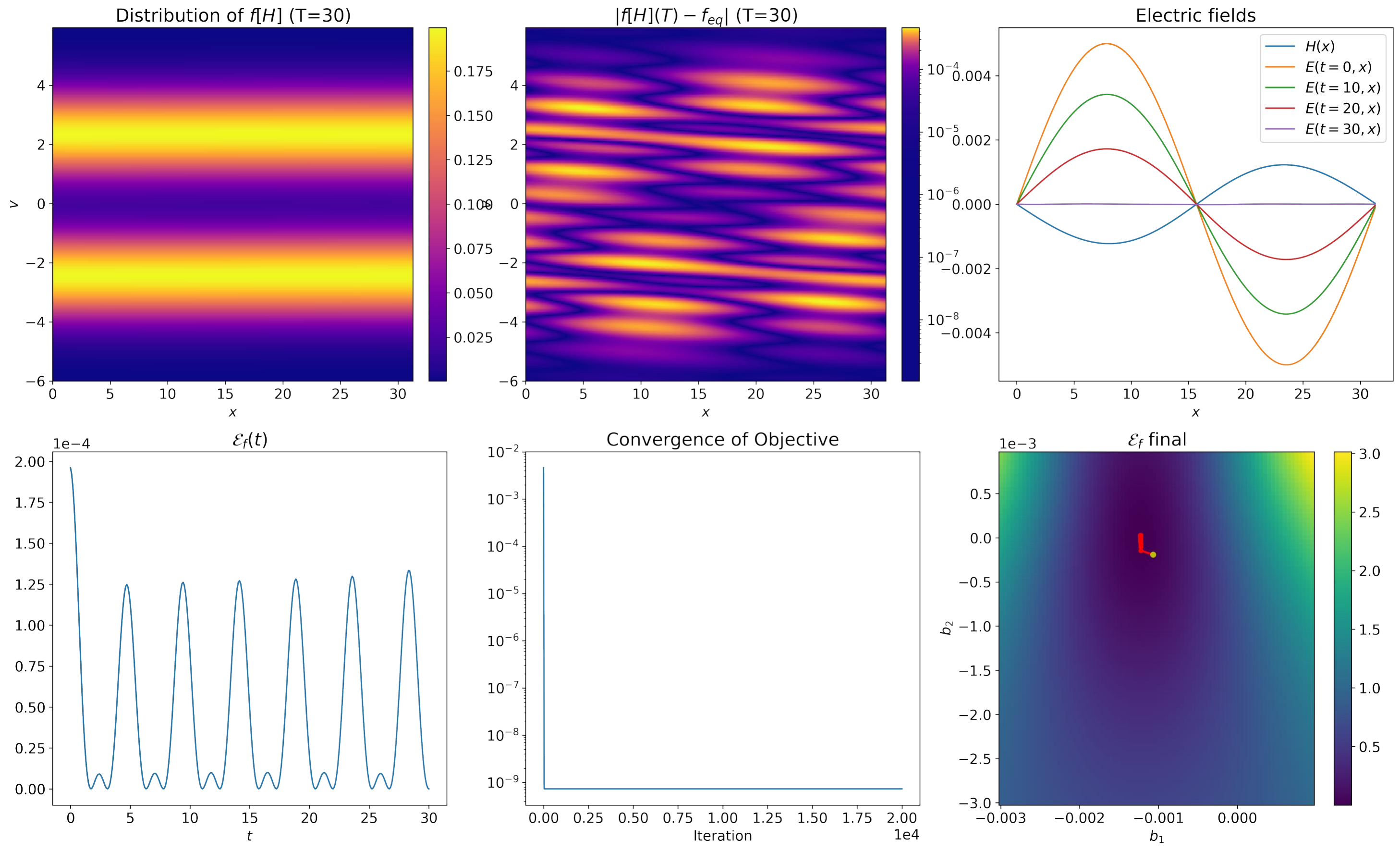}
    \caption{Simulation of~\eqref{eq:vlasov-poisoon_system_ext_1d} with under-parametrized $H$ obtained from~\eqref{eq:optimization_pb_simple} using~\eqref{eq:EE_obj} with local initialization using GD with line-search. From left to right and top to bottom: $f[H](T=30,x,v)$, $|f[H](T,x,v)-f_{\text{eq}}(v)|$, $H$ and $E_{f[H]}(t,x)$, $\mathcal{E}_{f[H]}(t)$, convergence of objective and, trajectory over the landscape of the objective (yellow dot is initial guess).}
    \label{fig:TS_ee_lf_GDL_local_under}
\end{figure}

\begin{figure}[H]
    \centering
    \includegraphics[width=0.85\linewidth]{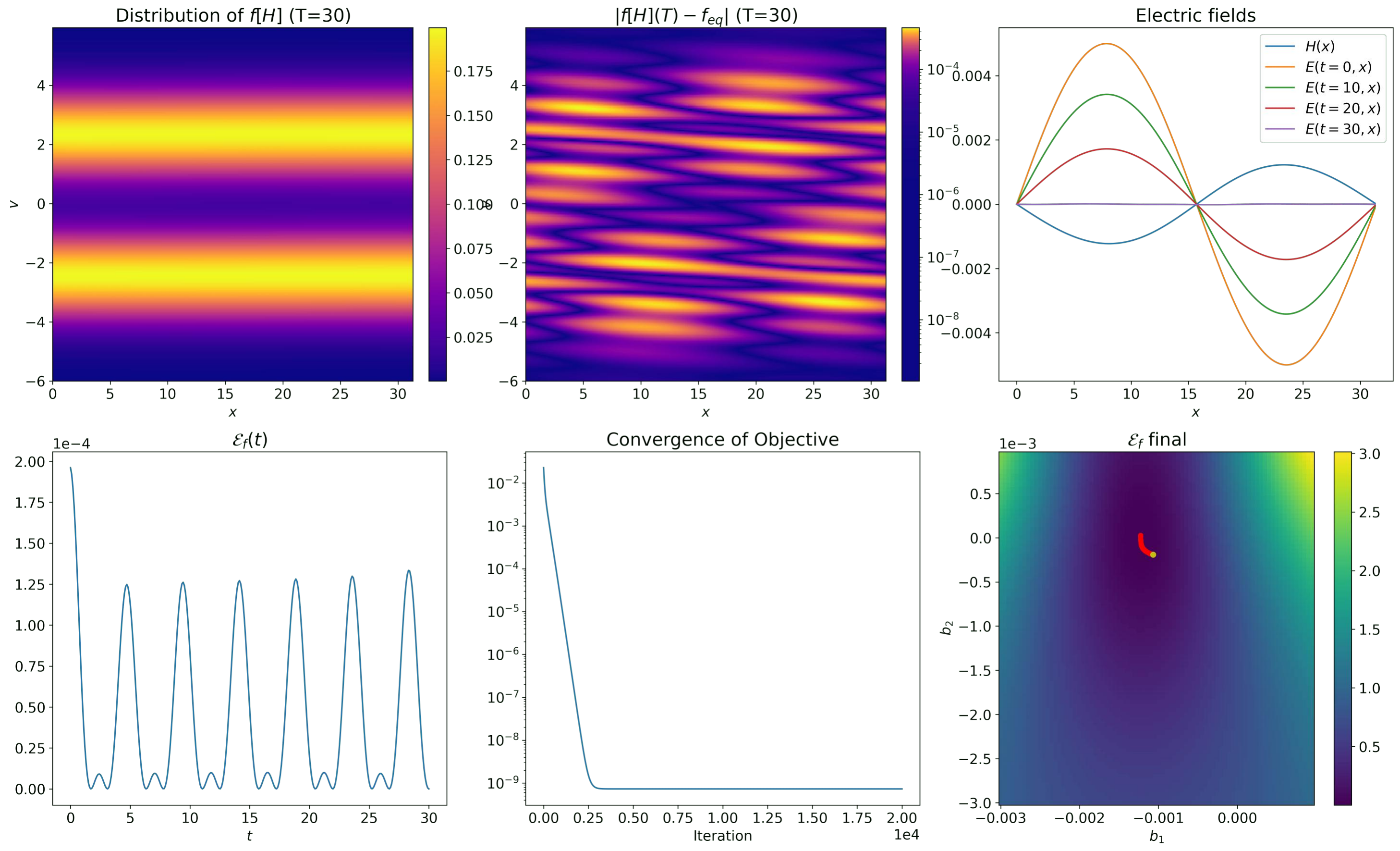}
    \caption{Simulation of~\eqref{eq:vlasov-poisoon_system_ext_1d} using~\eqref{eq:EE_obj} with under-parametrized $H$ obtained from~\eqref{eq:optimization_pb_simple} with local initialization using GD with constant stepsize. From left to right and top to bottom: $f[H](T=30,x,v)$, $|f[H](T,x,v)-f_{\text{eq}}(v)|$, $H$ and $E_{f[H]}(t,x)$, $\mathcal{E}_{f[H]}(t)$, convergence of objective and, trajectory over the landscape of the objective (yellow dot is initial guess).}
    \label{fig:TS_ee_lf_GD_local_under}
\end{figure}

\begin{figure}[H]
    \centering
    \includegraphics[width=0.85\linewidth]{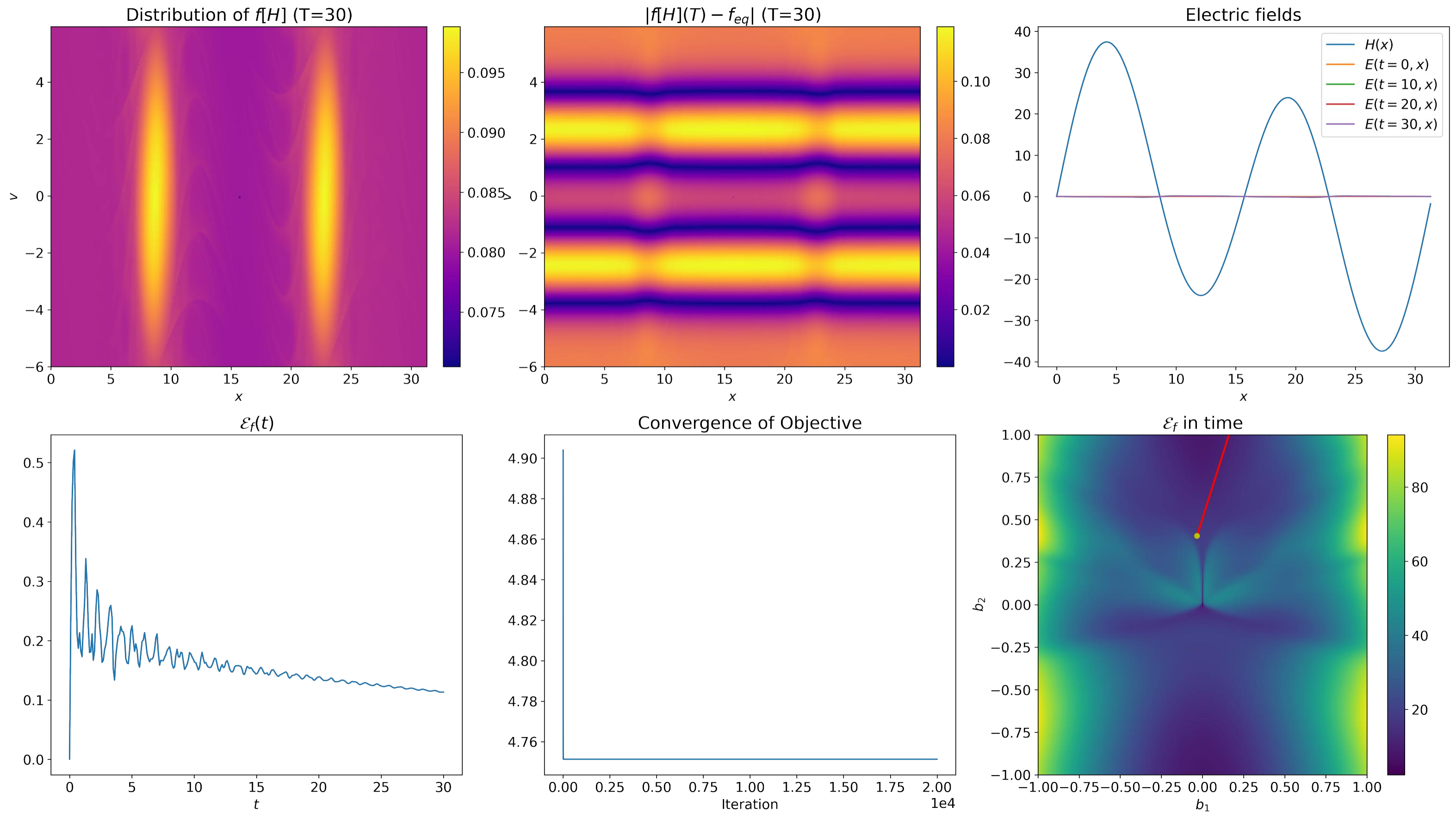}
    \caption{Simulation of~\eqref{eq:vlasov-poisoon_system_ext_1d} with under-parametrized $H$ obtained from~\eqref{eq:optimization_pb_simple} using~\eqref{eq:EET_obj} with far initialization using GD with line-search. From left to right and top to bottom: $f[H](T=30,x,v)$, $|f[H](T,x,v)-f_{\text{eq}}(v)|$, $H$ and $E_{f[H]}(t,x)$, $\mathcal{E}_{f[H]}(t)$, convergence of objective and, trajectory over the landscape of the objective (yellow dot is initial guess).}
    \label{fig:TS_ee_GDL_far_under}
\end{figure}

\begin{figure}[H]
    \centering
    \includegraphics[width=0.85\linewidth]{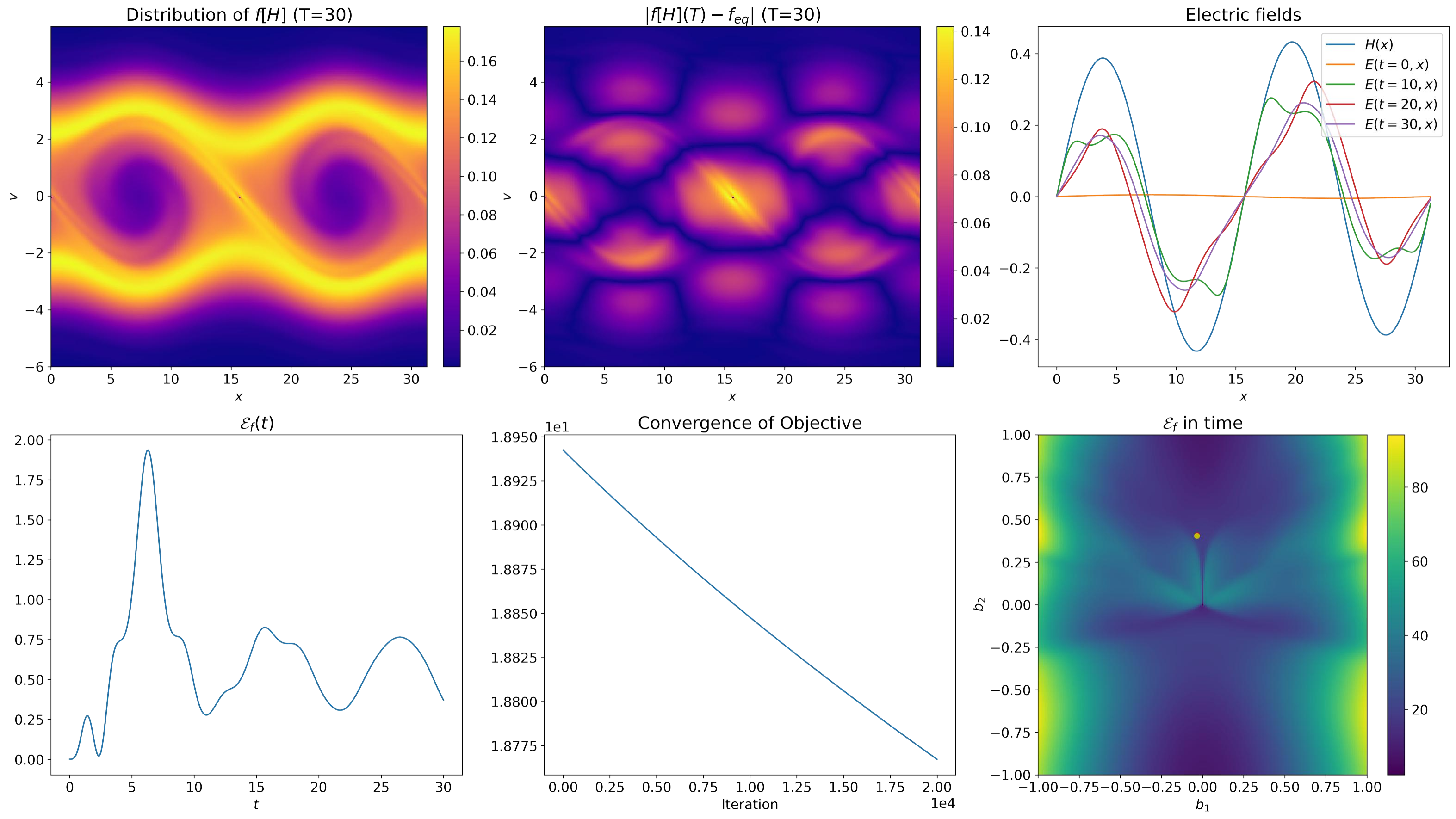}
    \caption{Simulation of~\eqref{eq:vlasov-poisoon_system_ext_1d} with under-parametrized $H$ obtained from~\eqref{eq:optimization_pb_simple} using~\eqref{eq:EET_obj} with far initialization using GD with constant stepsize. From left to right and top to bottom: $f[H](T=30,x,v)$, $|f[H](T,x,v)-f_{\text{eq}}(v)|$, $H$ and $E_{f[H]}(T,x)$, $E_{f[H]}(t,x)$, $\mathcal{E}_{f[H]}(t)$ and, convergence of objective.}
    \label{fig:TS_ee_GD_far_under}
\end{figure}

\begin{figure}[H]
    \centering
    \includegraphics[width=0.85\linewidth]{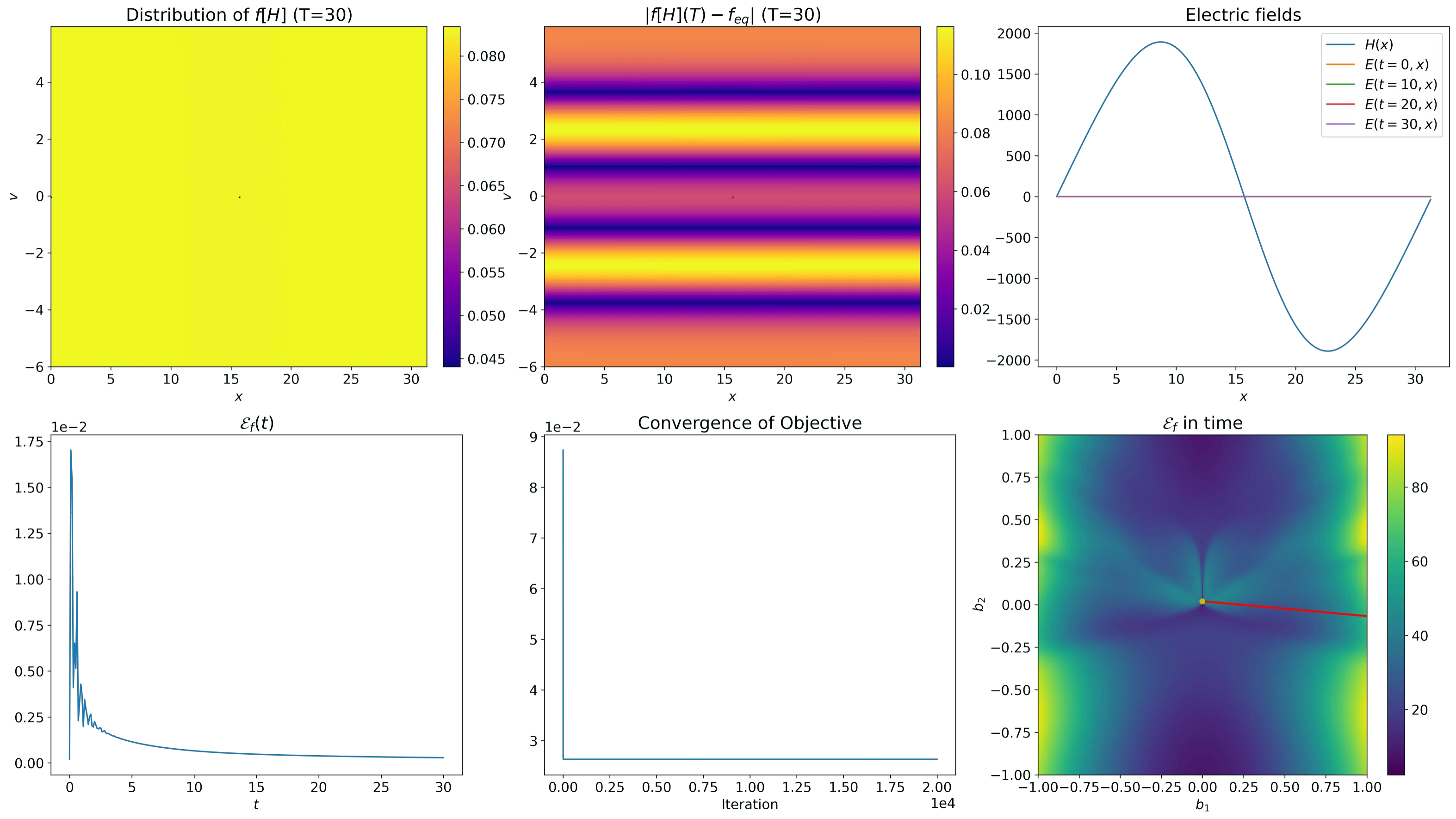}
    \caption{Simulation of~\eqref{eq:vlasov-poisoon_system_ext_1d} with under-parametrized $H$ obtained from~\eqref{eq:optimization_pb_simple} using~\eqref{eq:EET_obj} with near initialization using GD with line-search. From left to right and top to bottom: $f[H](T=30,x,v)$, $|f[H](T,x,v)-f_{\text{eq}}(v)|$, $H$ and $E_{f[H]}(t,x)$, $\mathcal{E}_{f[H]}(t)$, convergence of objective and, trajectory over the landscape of the objective (yellow dot is initial guess).}
    \label{fig:TS_ee_GDL_near_under}
\end{figure}

\begin{figure}[H]
    \centering
    \includegraphics[width=0.85\linewidth]{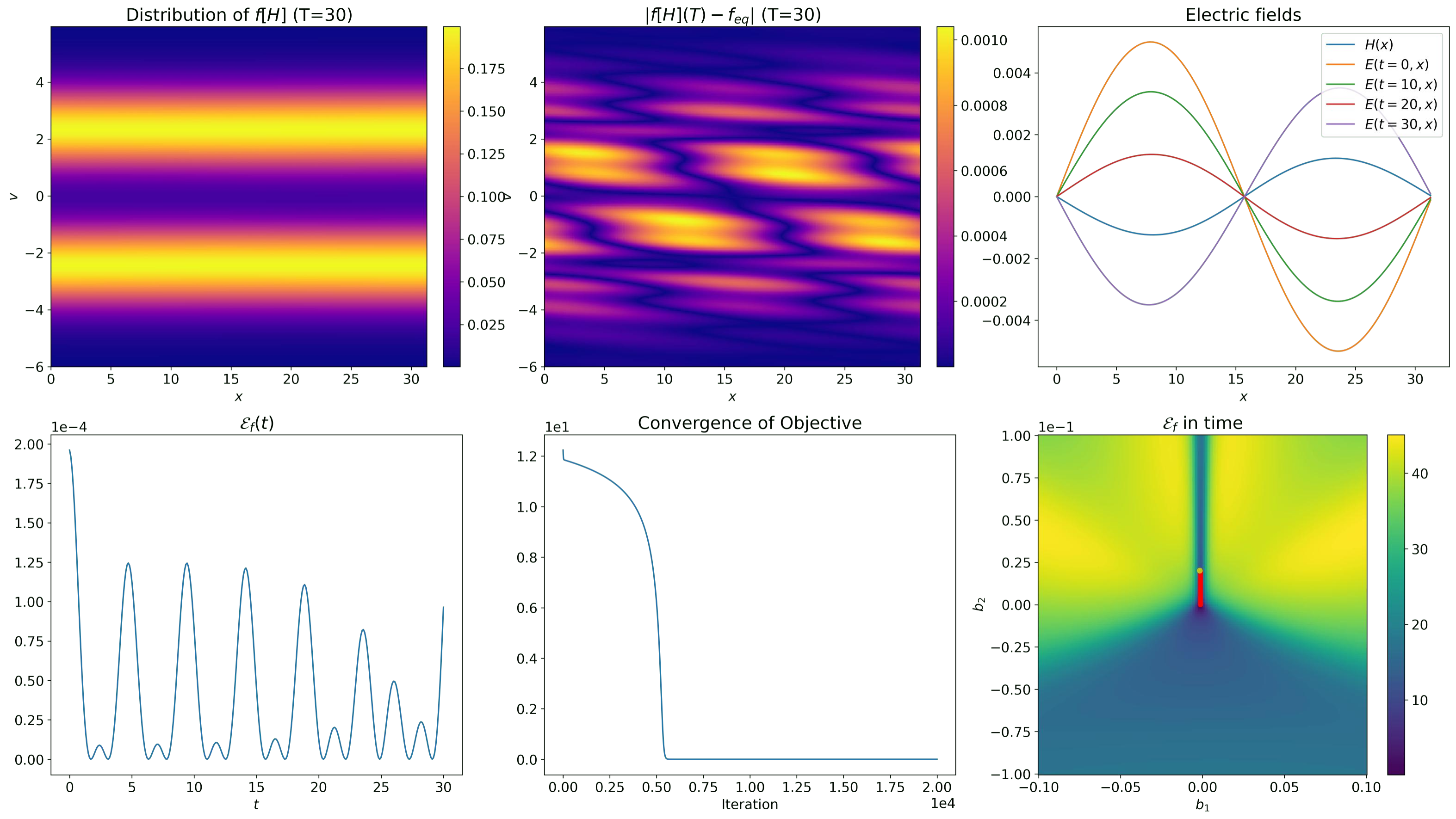}
    \caption{Simulation of~\eqref{eq:vlasov-poisoon_system_ext_1d} using~\eqref{eq:EET_obj} with under-parametrized $H$ obtained from~\eqref{eq:optimization_pb_simple} with near initialization using GD with constant stepsize. From left to right and top to bottom: $f[H](T=30,x,v)$, $|f[H](T,x,v)-f_{\text{eq}}(v)|$, $H$ and $E_{f[H]}(t,x)$, $\mathcal{E}_{f[H]}(t)$, convergence of objective and, trajectory over the landscape of the objective (yellow dot is initial guess).}
    \label{fig:TS_ee_GD_near_under}
\end{figure}

\begin{figure}[H]
    \centering
    \includegraphics[width=0.85\linewidth]{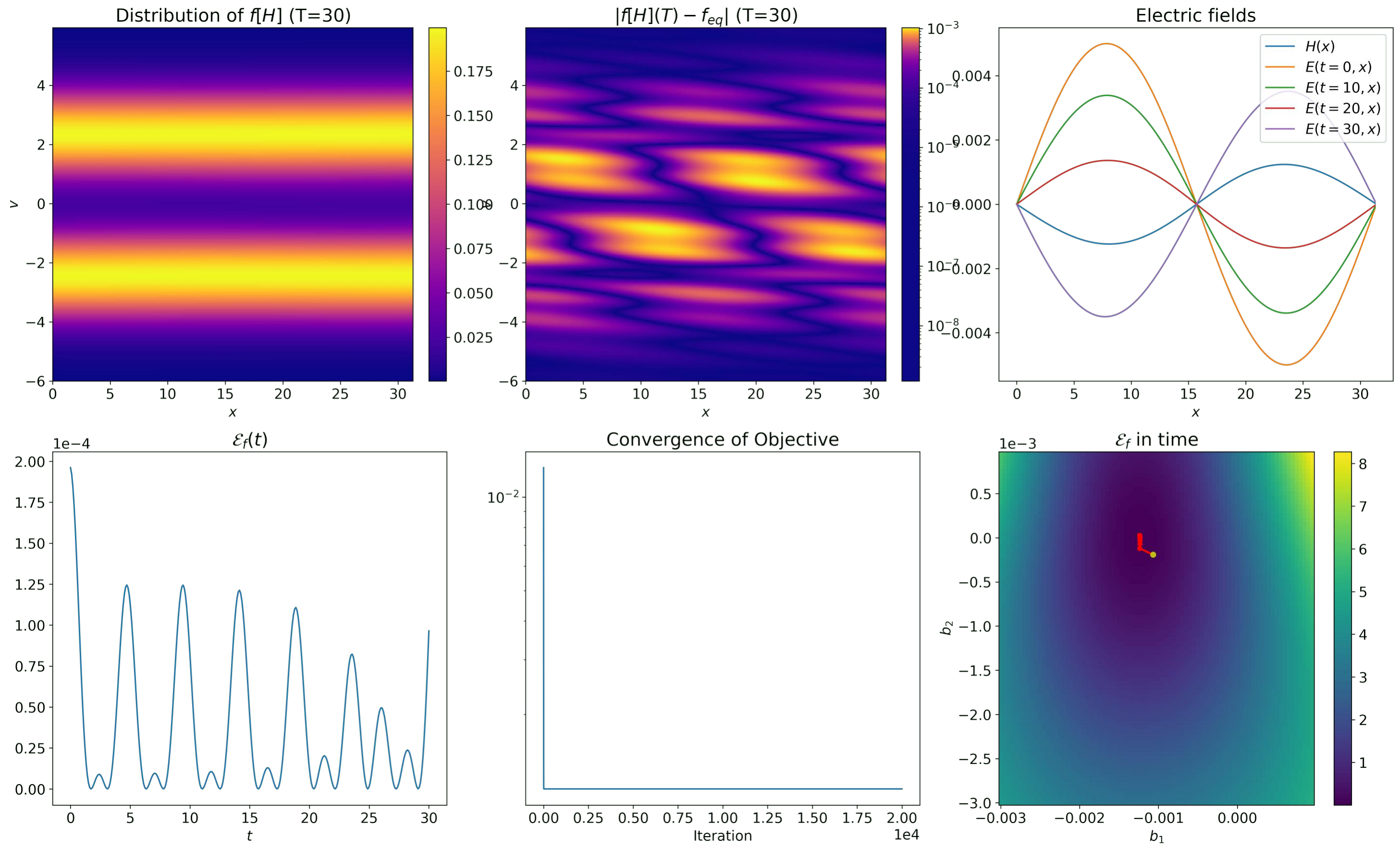}
    \caption{Simulation of~\eqref{eq:vlasov-poisoon_system_ext_1d} with under-parametrized $H$ obtained from~\eqref{eq:optimization_pb_simple} using~\eqref{eq:EET_obj} with local initialization using GD with line-search. From left to right and top to bottom: $f[H](T=30,x,v)$, $|f[H](T,x,v)-f_{\text{eq}}(v)|$, $H$ and $E_{f[H]}(t,x)$, $\mathcal{E}_{f[H]}(t)$, convergence of objective and, trajectory over the landscape of the objective (yellow dot is initial guess).}
    \label{fig:TS_ee_GDL_local_under}
\end{figure}

\begin{figure}[H]
    \centering
    \includegraphics[width=0.85\linewidth]{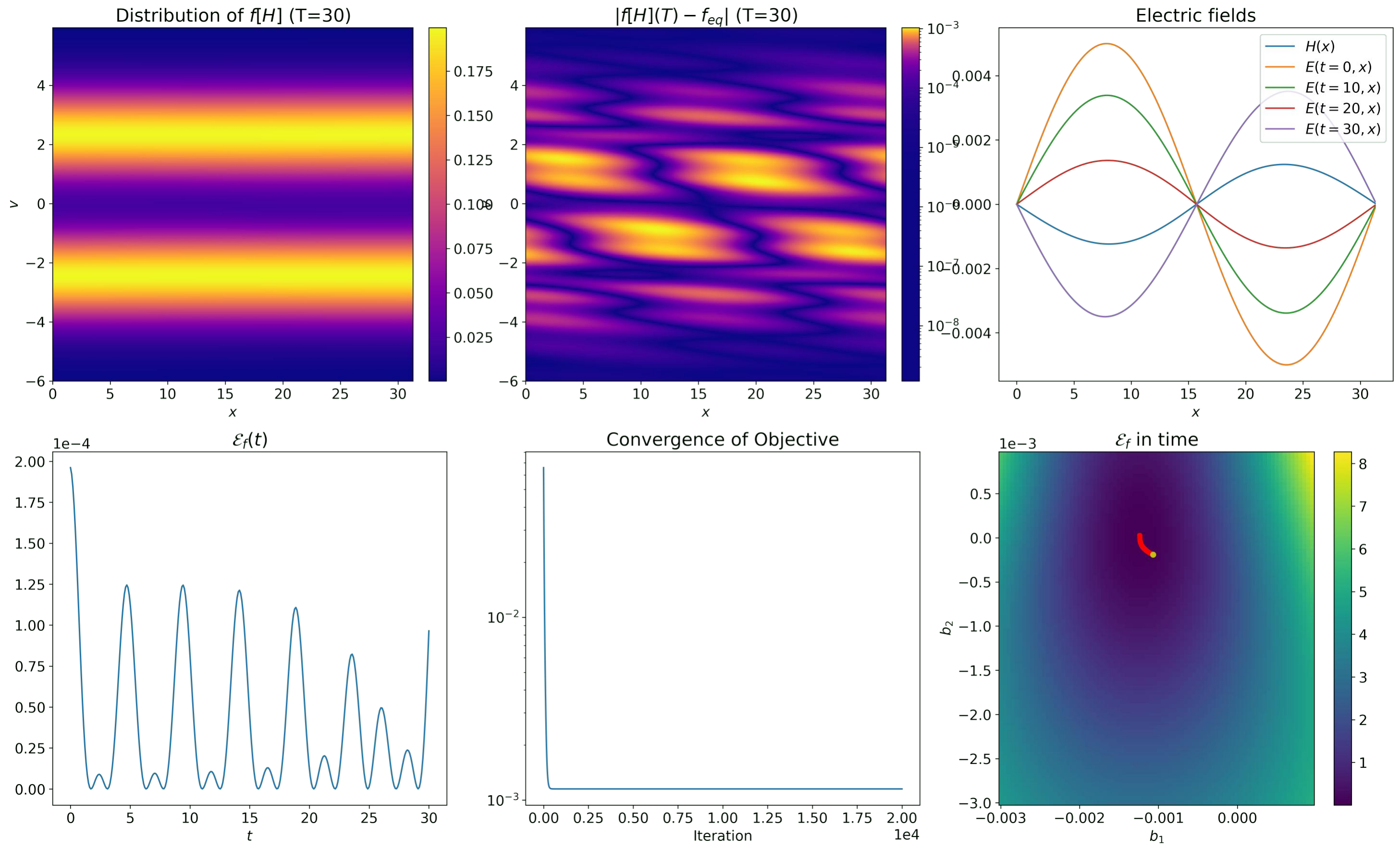}
    \caption{Simulation of~\eqref{eq:vlasov-poisoon_system_ext_1d} using~\eqref{eq:EET_obj} with under-parametrized $H$ obtained from~\eqref{eq:optimization_pb_simple} with local initialization using GD with constant stepsize. From left to right and top to bottom: $f[H](T=30,x,v)$, $|f[H](T,x,v)-f_{\text{eq}}(v)|$, $H$ and $E_{f[H]}(t,x)$, $\mathcal{E}_{f[H]}(t)$, convergence of objective and, trajectory over the landscape of the objective (yellow dot is initial guess).}
    \label{fig:TS_ee_GD_local_under}
\end{figure}

\begin{figure}[H]
    \centering
    \includegraphics[width=0.85\linewidth]{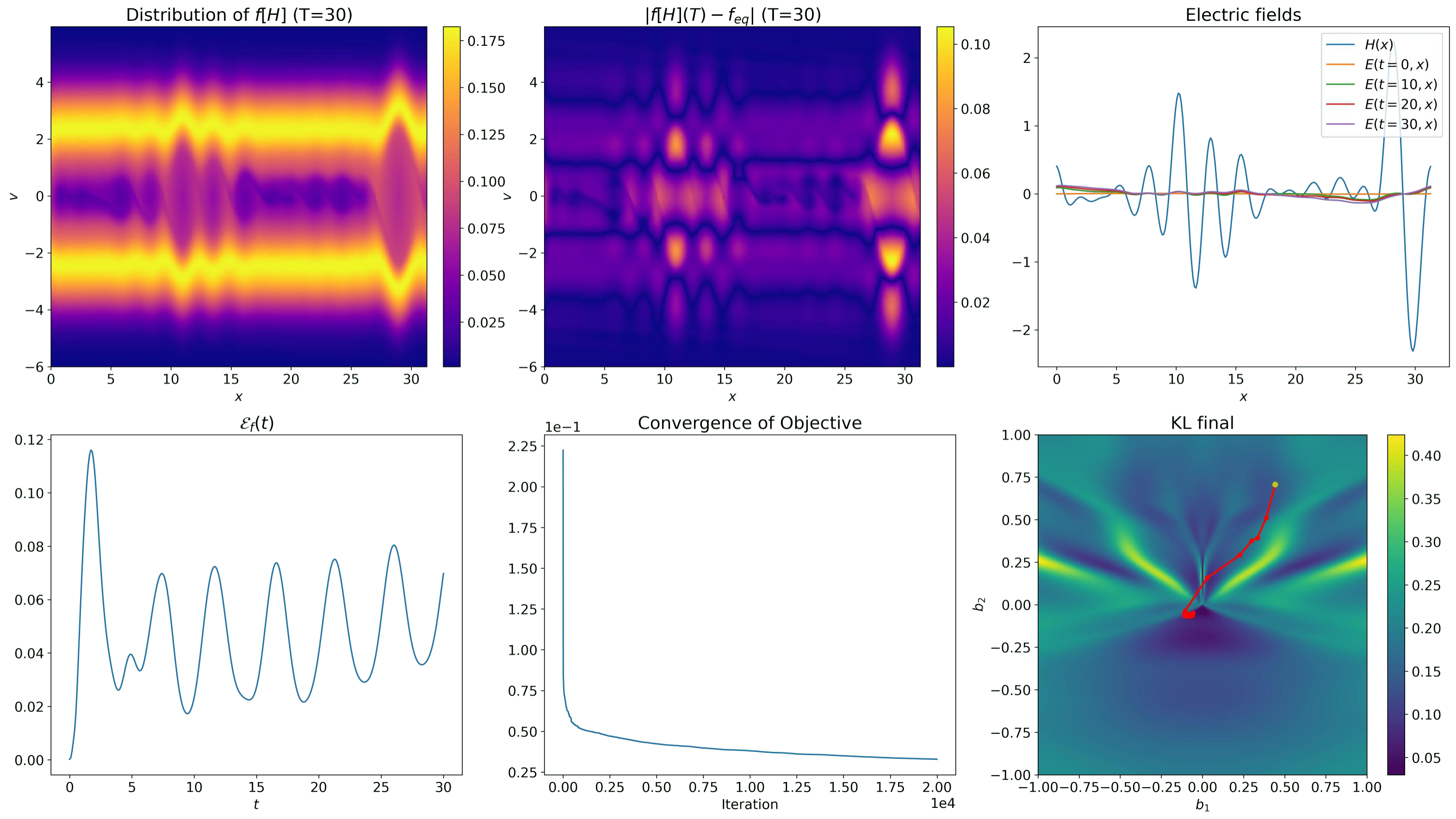}
    \caption{Simulation of~\eqref{eq:vlasov-poisoon_system_ext_1d} with over-parametrized $H$ obtained from~\eqref{eq:optimization_pb_simple} using~\eqref{eq:KL_obj} with far initialization using GD with line-search. From left to right and top to bottom: $f[H](T=30,x,v)$, $|f[H](T,x,v)-f_{\text{eq}}(v)|$, $H$ and $E_{f[H]}(t,x)$, $\mathcal{E}_{f[H]}(t)$, convergence of objective and, trajectory over the landscape of the objective (yellow dot is initial guess).}
    \label{fig:TS_KL_GDL_far_over}
\end{figure}

\begin{figure}[H]
    \centering
    \includegraphics[width=0.85\linewidth]{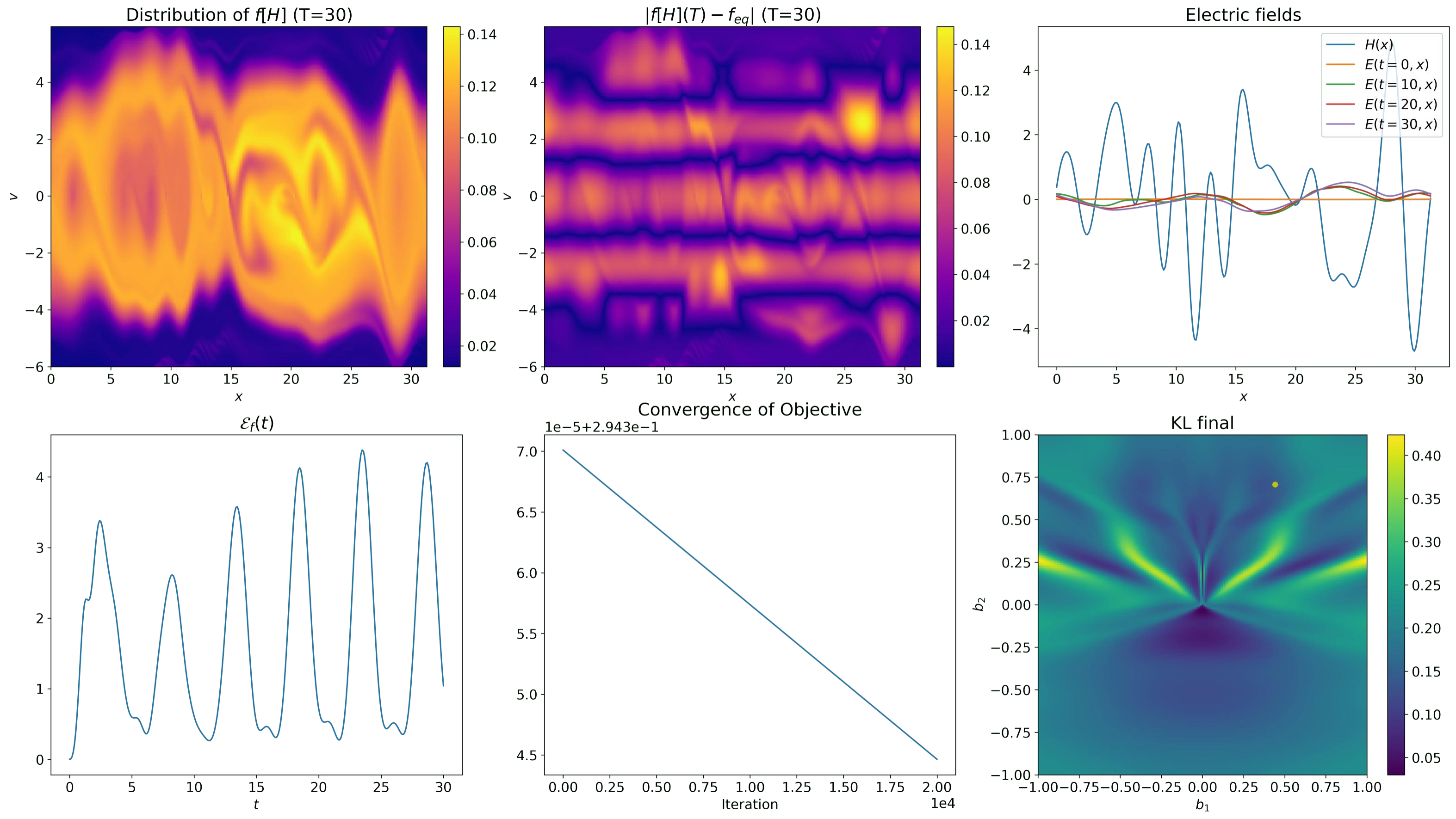}
    \caption{Simulation of~\eqref{eq:vlasov-poisoon_system_ext_1d} with over-parametrized $H$ obtained from~\eqref{eq:optimization_pb_simple} using~\eqref{eq:KL_obj} with far initialization using GD with constant stepsize. From left to right and top to bottom: $f[H](T=30,x,v)$, $|f[H](T,x,v)-f_{\text{eq}}(v)|$, $H$ and $E_{f[H]}(t,x)$, $\mathcal{E}_{f[H]}(t)$, convergence of objective and, trajectory over the landscape of the objective (yellow dot is initial guess).}
    \label{fig:TS_KL_GD_far_over}
\end{figure}

\begin{figure}[H]
    \centering
    \includegraphics[width=0.85\linewidth]{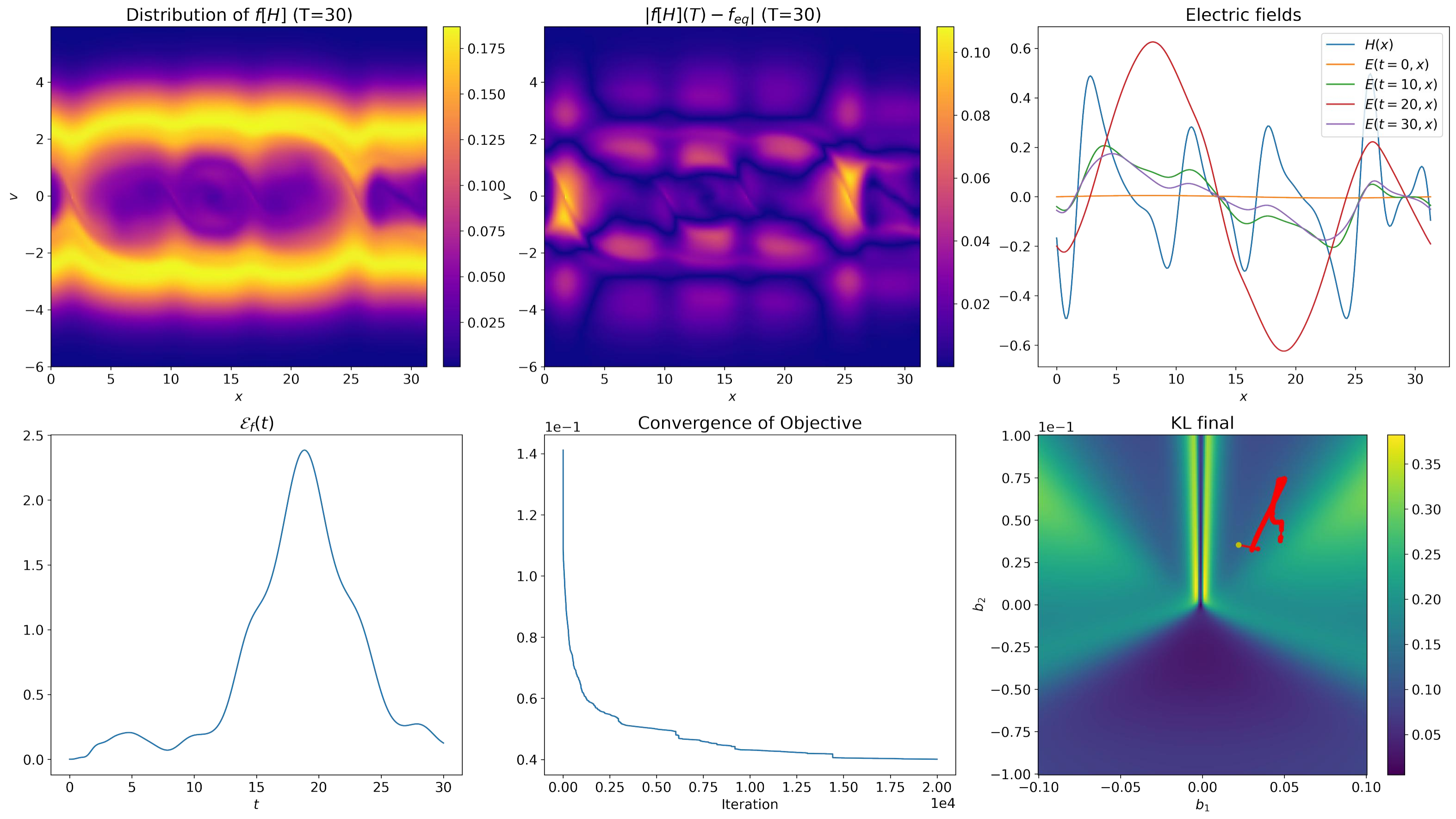}
    \caption{Simulation of~\eqref{eq:vlasov-poisoon_system_ext_1d} with over-parametrized $H$ obtained from~\eqref{eq:optimization_pb_simple} using~\eqref{eq:KL_obj} with near initialization using GD with line-search. From left to right and top to bottom: $f[H](T=30,x,v)$, $|f[H](T,x,v)-f_{\text{eq}}(v)|$, $H$ and $E_{f[H]}(t,x)$, $\mathcal{E}_{f[H]}(t)$, convergence of objective and, trajectory over the landscape of the objective (yellow dot is initial guess).}
    \label{fig:TS_KL_GDL_near_over}
\end{figure}

\begin{figure}[H]
    \centering
    \includegraphics[width=0.85\linewidth]{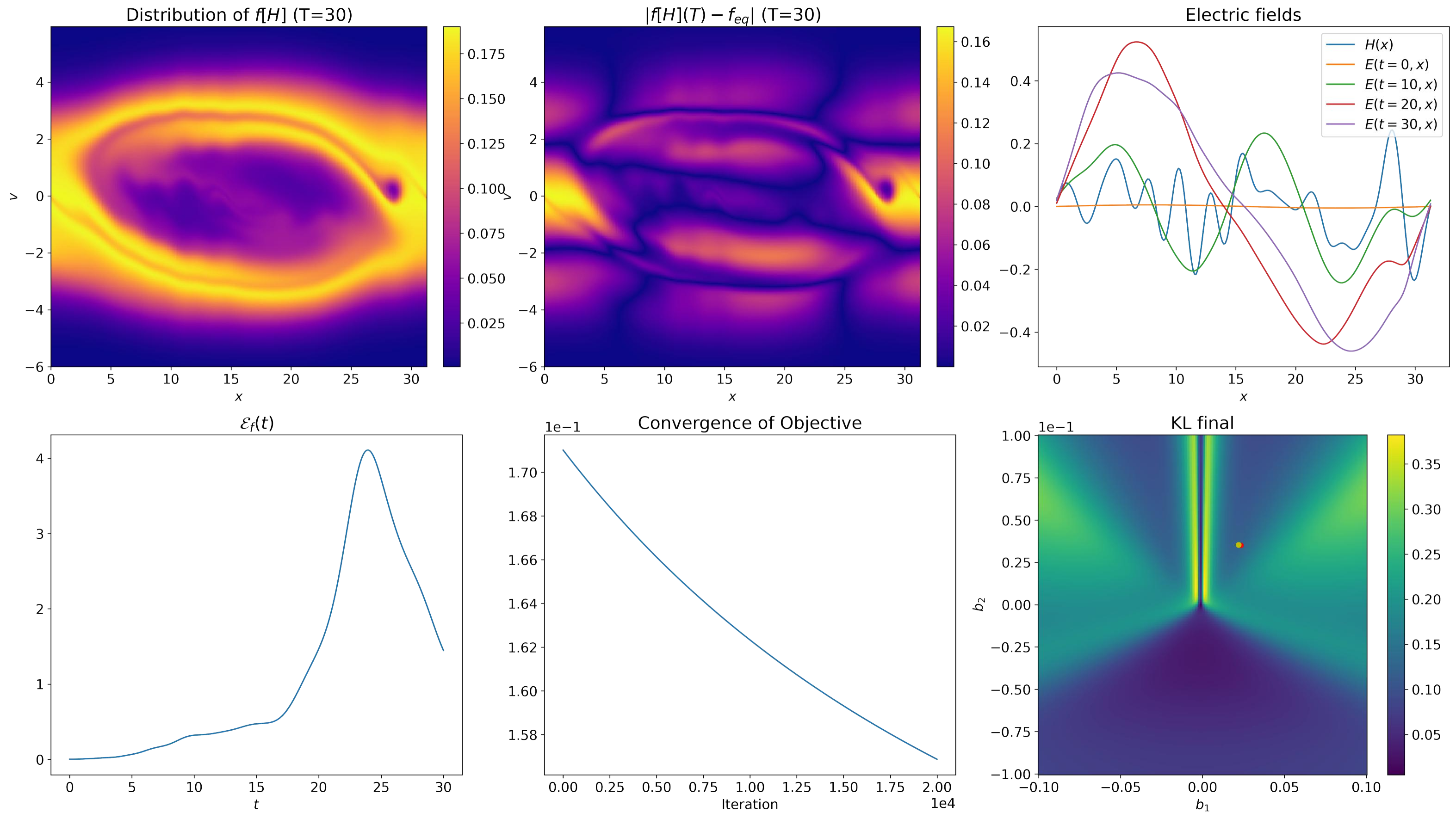}
    \caption{Simulation of~\eqref{eq:vlasov-poisoon_system_ext_1d} with over-parametrized $H$ obtained from~\eqref{eq:optimization_pb_simple} using~\eqref{eq:KL_obj} with near initialization using GD with constant stepsize. From left to right and top to bottom: $f[H](T=30,x,v)$, $|f[H](T,x,v)-f_{\text{eq}}(v)|$, $H$ and $E_{f[H]}(t,x)$, $\mathcal{E}_{f[H]}(t)$, convergence of objective and, trajectory over the landscape of the objective (yellow dot is initial guess).}
    \label{fig:TS_KL_GD_near_over}
\end{figure}

\begin{figure}[H]
    \centering
    \includegraphics[width=0.85\linewidth]{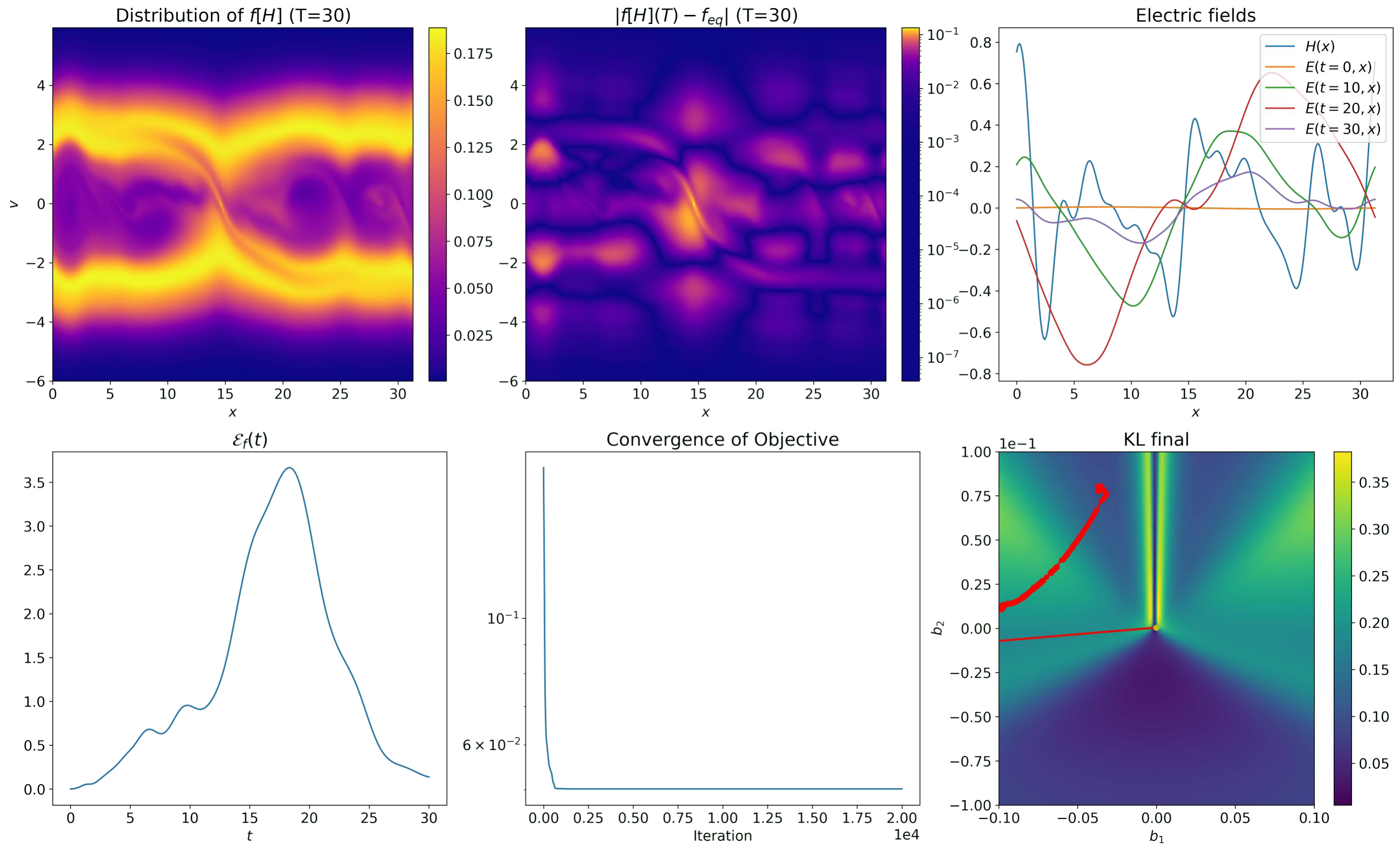}
    \caption{Simulation of~\eqref{eq:vlasov-poisoon_system_ext_1d} with over-parametrized $H$ obtained from~\eqref{eq:optimization_pb_simple} using~\eqref{eq:KL_obj} with local initialization using GD with line-search. From left to right and top to bottom: $f[H](T=30,x,v)$, $|f[H](T,x,v)-f_{\text{eq}}(v)|$, $H$ and $E_{f[H]}(t,x)$, $\mathcal{E}_{f[H]}(t)$, convergence of objective and, trajectory over the landscape of the objective (yellow dot is initial guess).}
    \label{fig:TS_KL_GDL_local_over}
\end{figure}

\begin{figure}[H]
    \centering
    \includegraphics[width=0.85\linewidth]{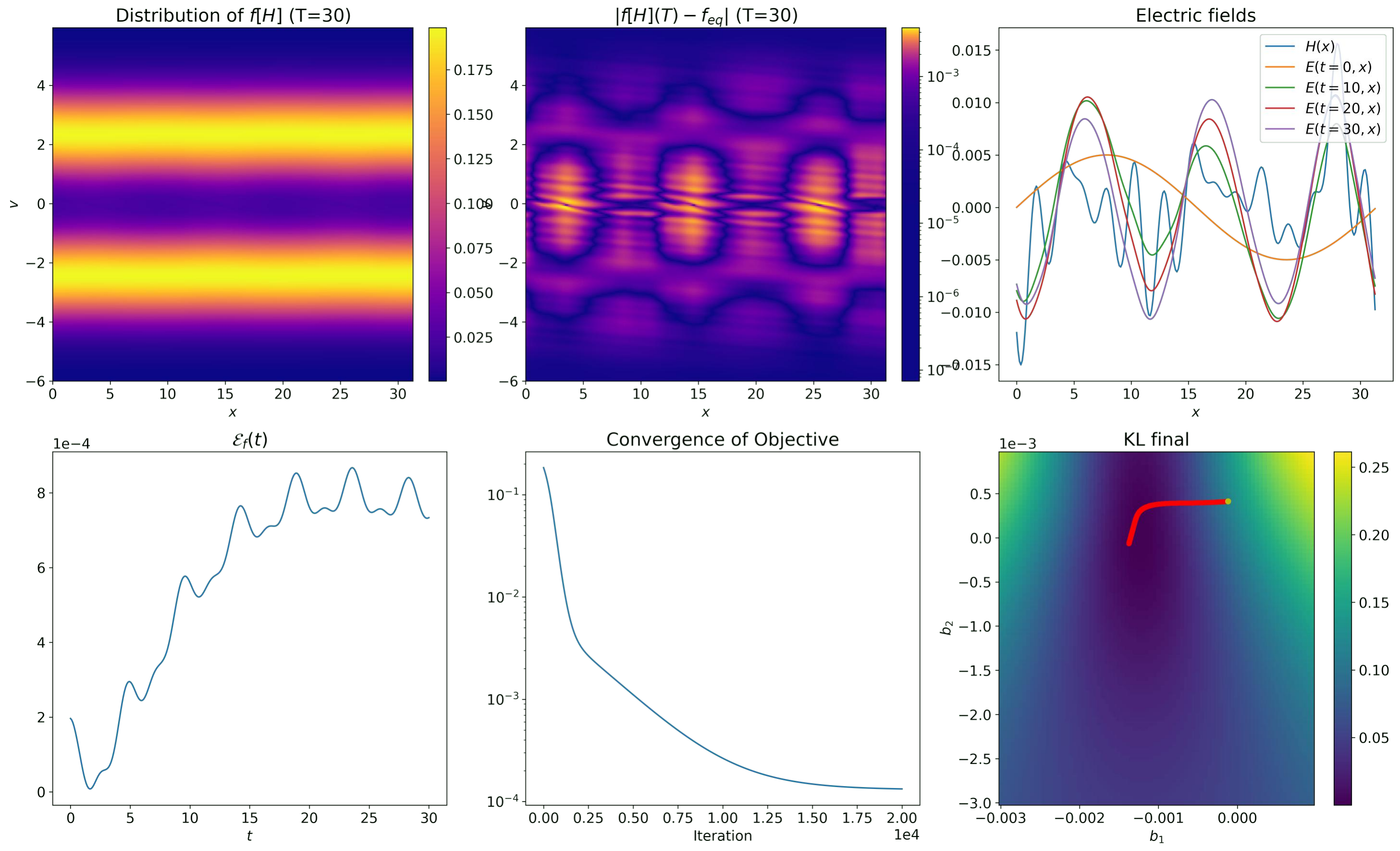}
    \caption{Simulation of~\eqref{eq:vlasov-poisoon_system_ext_1d} with over-parametrized $H$ obtained from~\eqref{eq:optimization_pb_simple} using~\eqref{eq:KL_obj} with local initialization using GD with constant stepsize. From left to right and top to bottom: $f[H](T=30,x,v)$, $|f[H](T,x,v)-f_{\text{eq}}(v)|$, $H$ and $E_{f[H]}(t,x)$, $\mathcal{E}_{f[H]}(t)$, convergence of objective and, trajectory over the landscape of the objective (yellow dot is initial guess).}
    \label{fig:TS_KL_GD_local_over}
\end{figure}

\begin{figure}[H]
    \centering
    \includegraphics[width=0.85\linewidth]{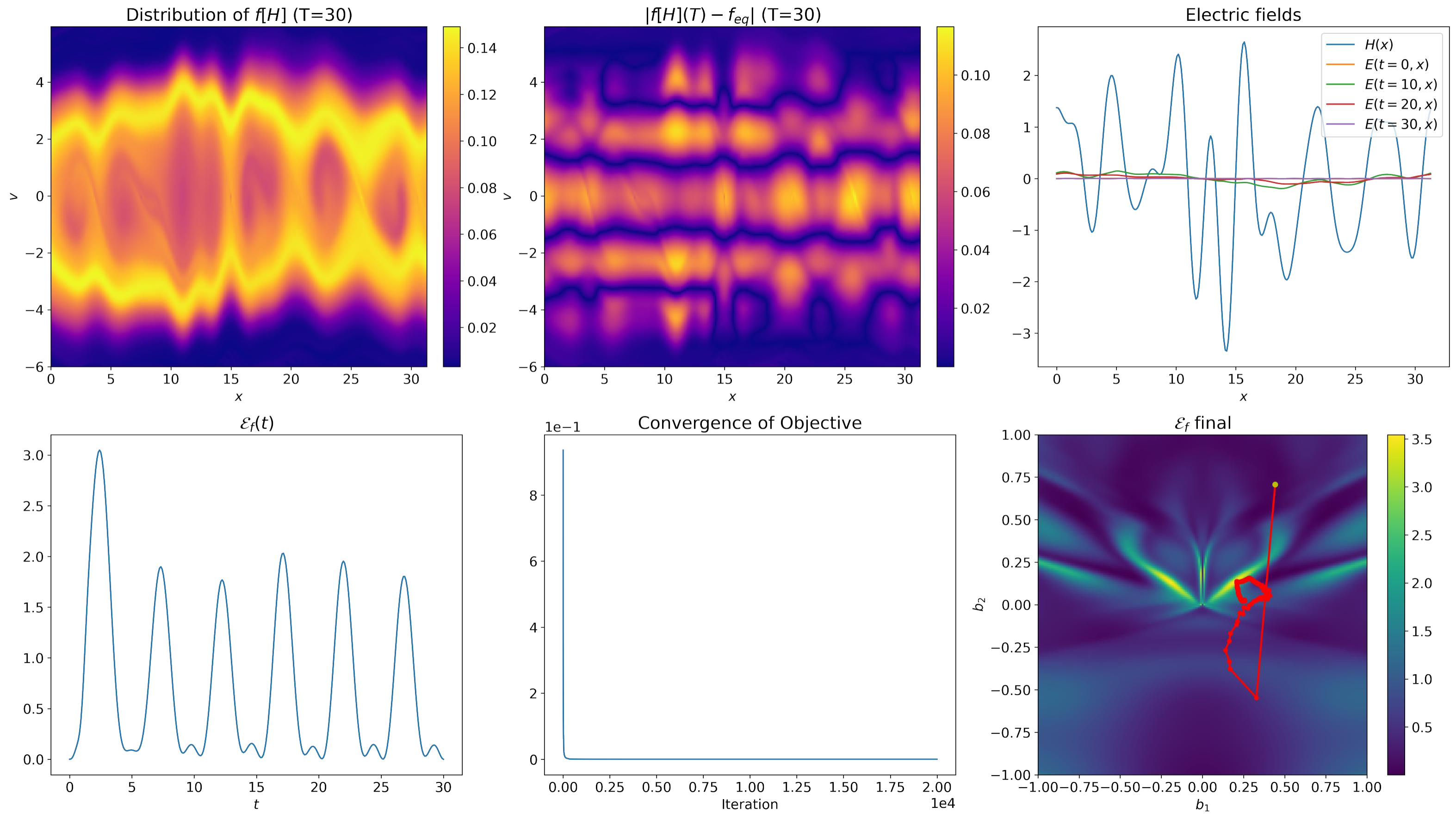}
    \caption{Simulation of~\eqref{eq:vlasov-poisoon_system_ext_1d} with under-parametrized $H$ obtained from~\eqref{eq:optimization_pb_simple} using~\eqref{eq:EE_obj} with far initialization using GD with line-search. From left to right and top to bottom: $f[H](T=30,x,v)$, $|f[H](T,x,v)-f_{\text{eq}}(v)|$, $H$ and $E_{f[H]}(t,x)$, $\mathcal{E}_{f[H]}(t)$, convergence of objective and, trajectory over the landscape of the objective (yellow dot is initial guess).}
    \label{fig:TS_ee_lf_GDL_far_over}
\end{figure}

\begin{figure}[H]
    \centering
    \includegraphics[width=0.85\linewidth]{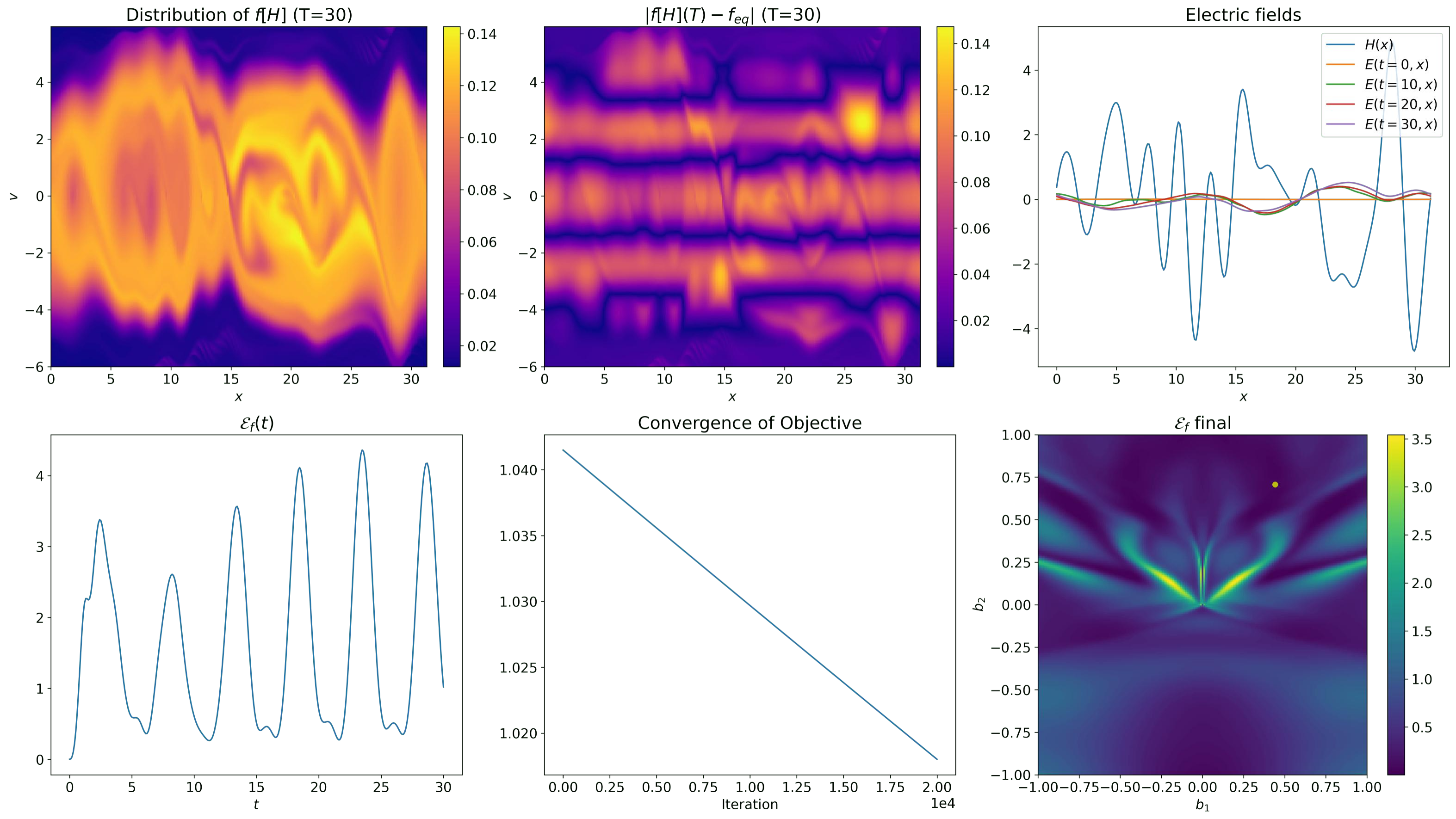}
    \caption{Simulation of~\eqref{eq:vlasov-poisoon_system_ext_1d} with over-parametrized $H$ obtained from~\eqref{eq:optimization_pb_simple} using~\eqref{eq:EE_obj} with far initialization using GD with constant stepsize. From left to right and top to bottom: $f[H](T=30,x,v)$, $|f[H](T,x,v)-f_{\text{eq}}(v)|$, $H$ and $E_{f[H]}(t,x)$, $\mathcal{E}_{f[H]}(t)$, convergence of objective and, trajectory over the landscape of the objective (yellow dot is initial guess).}
    \label{fig:TS_ee_lf_GD_far_over}
\end{figure}

\begin{figure}[H]
    \centering
    \includegraphics[width=0.85\linewidth]{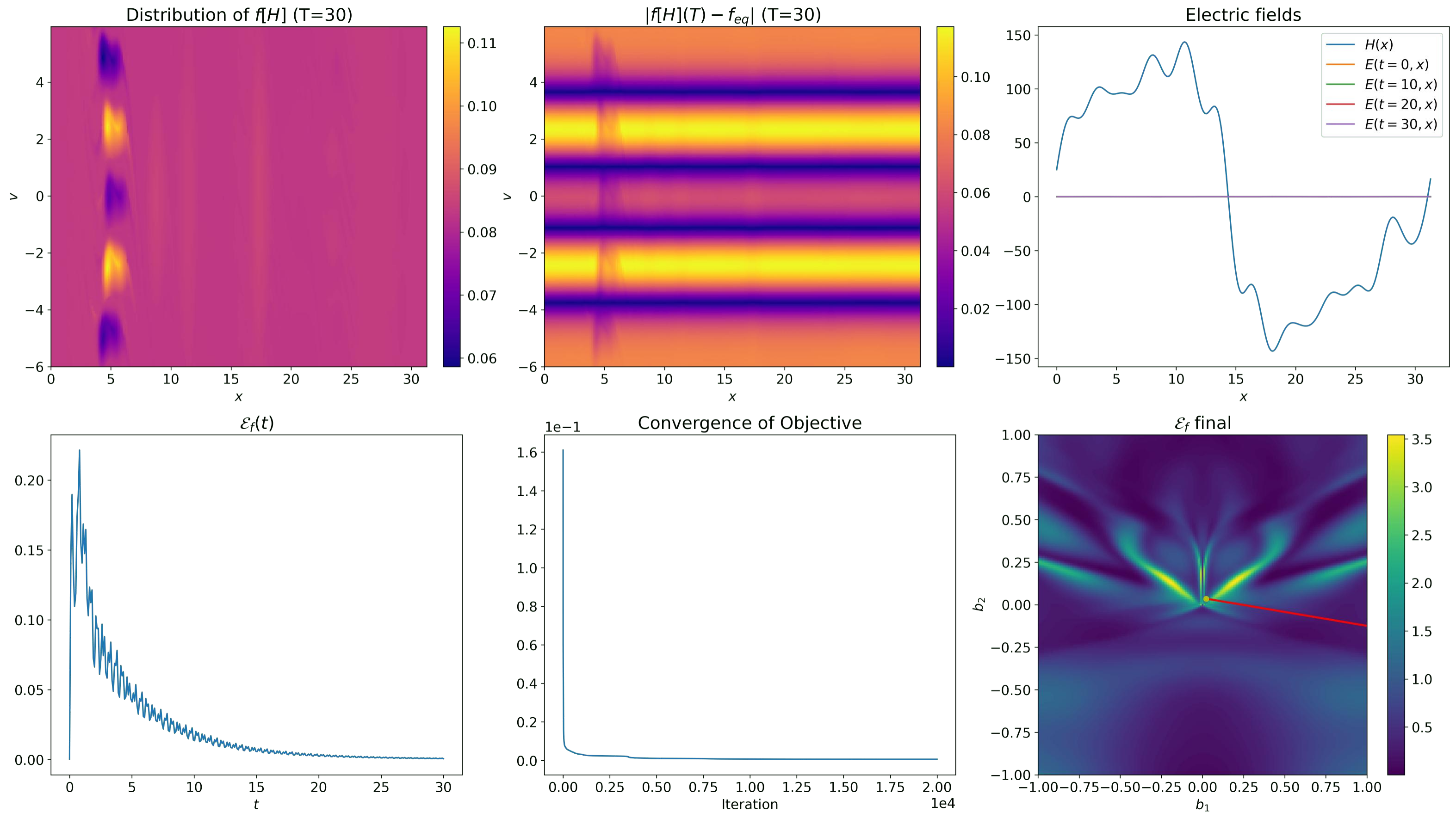}
    \caption{Simulation of~\eqref{eq:vlasov-poisoon_system_ext_1d} with over-parametrized $H$ obtained from~\eqref{eq:optimization_pb_simple} using~\eqref{eq:EE_obj} with near initialization using GD with line-search. From left to right and top to bottom: $f[H](T=30,x,v)$, $|f[H](T,x,v)-f_{\text{eq}}(v)|$, $H$ and $E_{f[H]}(t,x)$, $\mathcal{E}_{f[H]}(t)$, convergence of objective and, trajectory over the landscape of the objective (yellow dot is initial guess).}
    \label{fig:TS_ee_lf_GDL_near_over}
\end{figure}

\begin{figure}[H]
    \centering
    \includegraphics[width=0.85\linewidth]{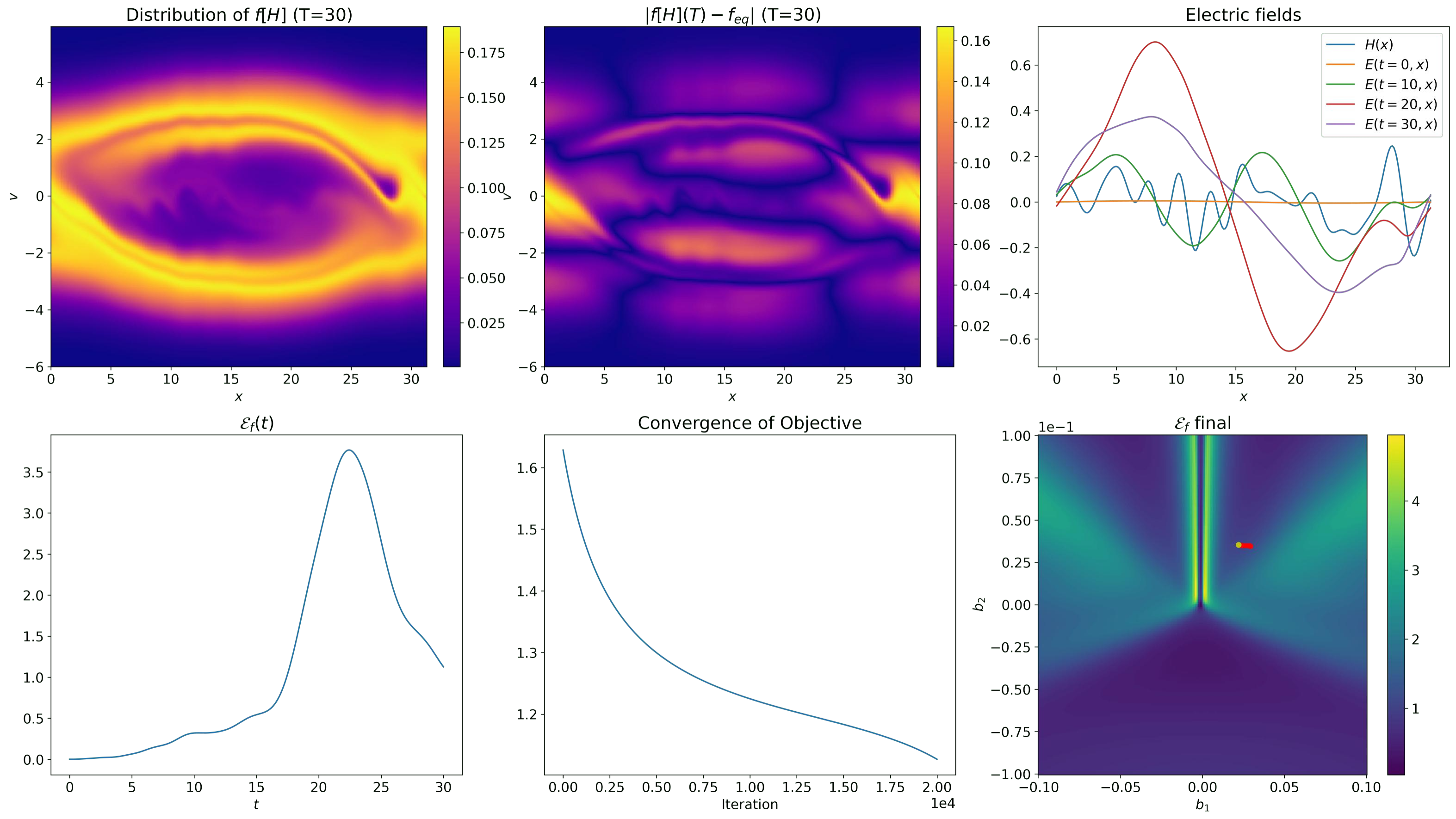}
    \caption{Simulation of~\eqref{eq:vlasov-poisoon_system_ext_1d} using~\eqref{eq:EE_obj} with over-parametrized $H$ obtained from~\eqref{eq:optimization_pb_simple} with near initialization using GD with constant stepsize. From left to right and top to bottom: $f[H](T=30,x,v)$, $|f[H](T,x,v)-f_{\text{eq}}(v)|$, $H$ and $E_{f[H]}(t,x)$, $\mathcal{E}_{f[H]}(t)$, convergence of objective and, trajectory over the landscape of the objective (yellow dot is initial guess).}
    \label{fig:TS_ee_lf_GD_near_over}
\end{figure}

\begin{figure}[H]
    \centering
    \includegraphics[width=0.85\linewidth]{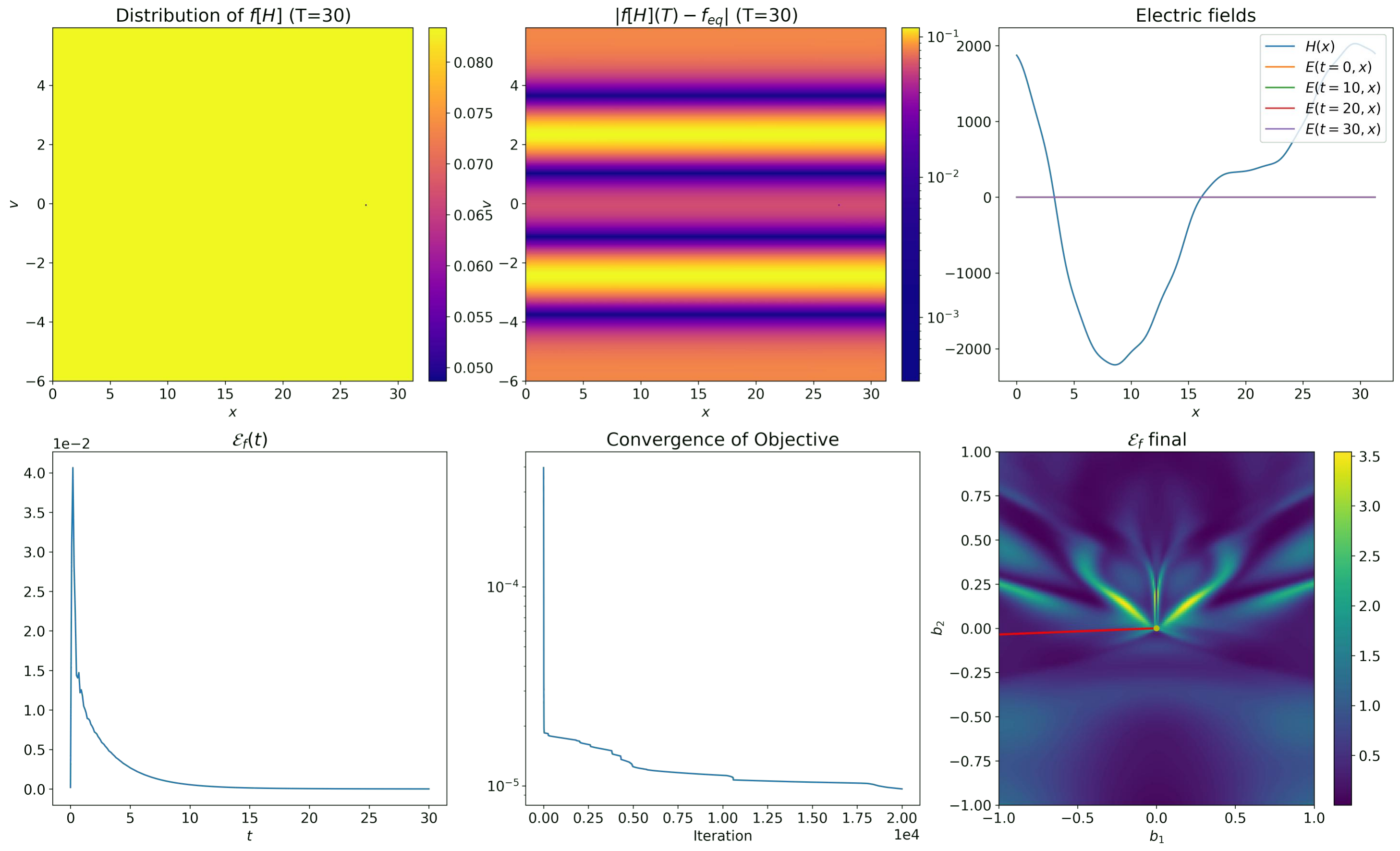}
    \caption{Simulation of~\eqref{eq:vlasov-poisoon_system_ext_1d} with over-parametrized $H$ obtained from~\eqref{eq:optimization_pb_simple} using~\eqref{eq:EE_obj} with local initialization using GD with line-search. From left to right and top to bottom: $f[H](T=30,x,v)$, $|f[H](T,x,v)-f_{\text{eq}}(v)|$, $H$ and $E_{f[H]}(t,x)$, $\mathcal{E}_{f[H]}(t)$, convergence of objective and, trajectory over the landscape of the objective (yellow dot is initial guess).}
    \label{fig:TS_ee_lf_GDL_local_over}
\end{figure}

\begin{figure}[H]
    \centering
    \includegraphics[width=0.85\linewidth]{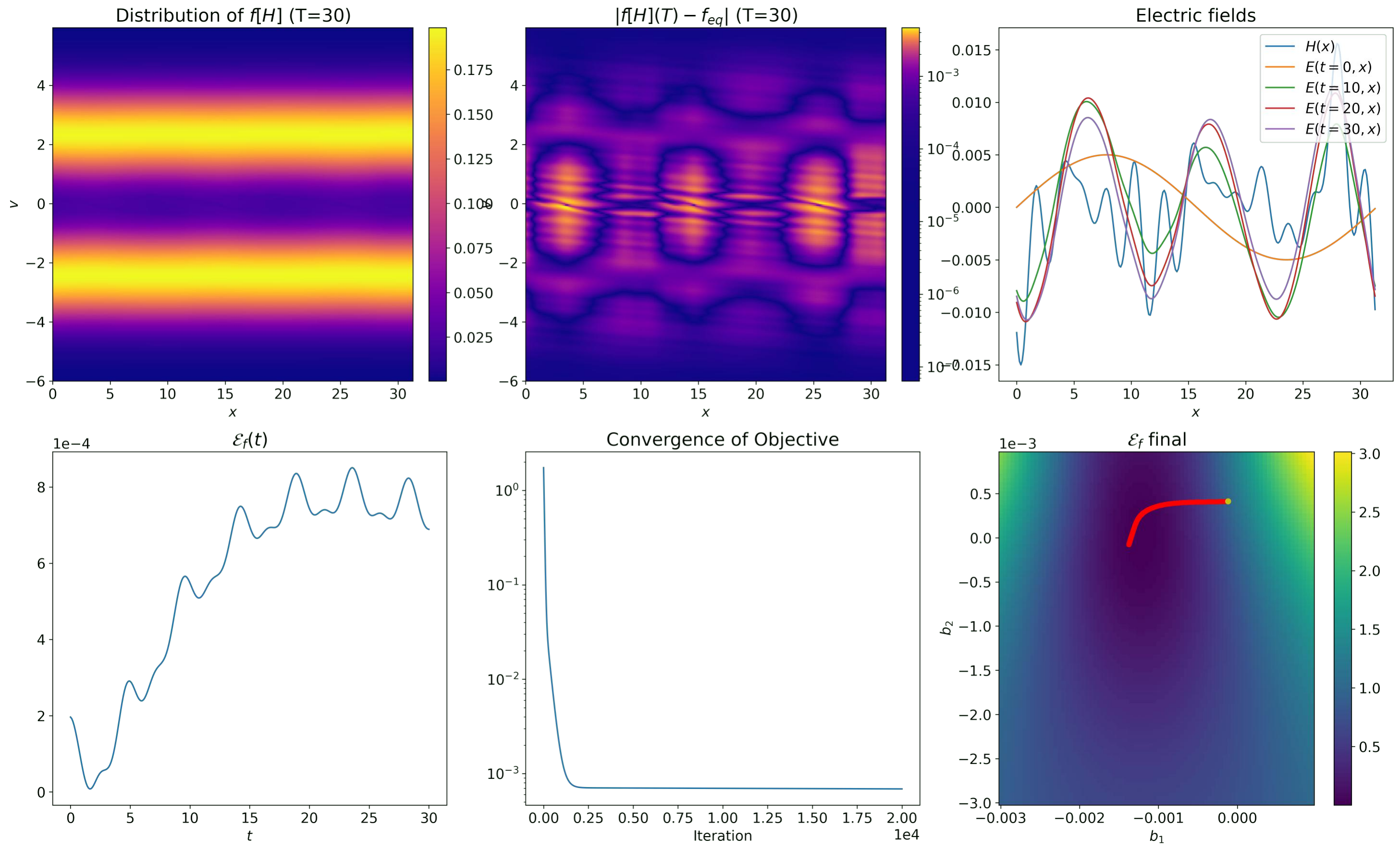}
    \caption{Simulation of~\eqref{eq:vlasov-poisoon_system_ext_1d} using~\eqref{eq:EE_obj} with over-parametrized $H$ obtained from~\eqref{eq:optimization_pb_simple} with local initialization using GD with constant stepsize. From left to right and top to bottom: $f[H](T=30,x,v)$, $|f[H](T,x,v)-f_{\text{eq}}(v)|$, $H$ and $E_{f[H]}(t,x)$, $\mathcal{E}_{f[H]}(t)$, convergence of objective and, trajectory over the landscape of the objective (yellow dot is initial guess).}
    \label{fig:TS_ee_lf_GD_local_over}
\end{figure}

\begin{figure}[H]
    \centering
    \includegraphics[width=0.85\linewidth]{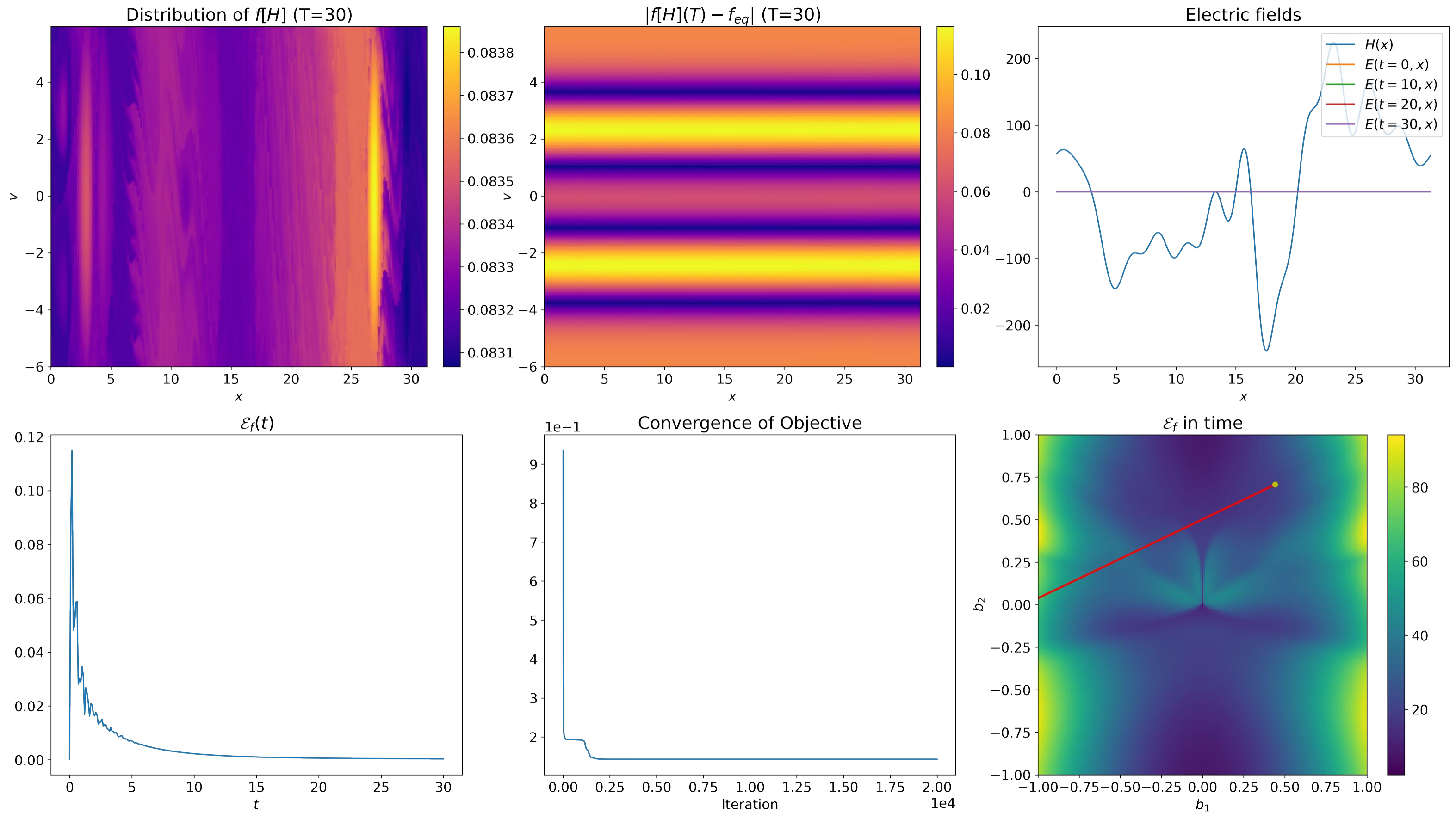}
    \caption{Simulation of~\eqref{eq:vlasov-poisoon_system_ext_1d} with under-parametrized $H$ obtained from~\eqref{eq:optimization_pb_simple} using~\eqref{eq:EET_obj} with far initialization using GD with line-search. From left to right and top to bottom: $f[H](T=30,x,v)$, $|f[H](T,x,v)-f_{\text{eq}}(v)|$, $H$ and $E_{f[H]}(t,x)$, $\mathcal{E}_{f[H]}(t)$, convergence of objective and, trajectory over the landscape of the objective (yellow dot is initial guess).}
    \label{fig:TS_ee_GDL_far_over}
\end{figure}

\begin{figure}[H]
    \centering
    \includegraphics[width=0.85\linewidth]{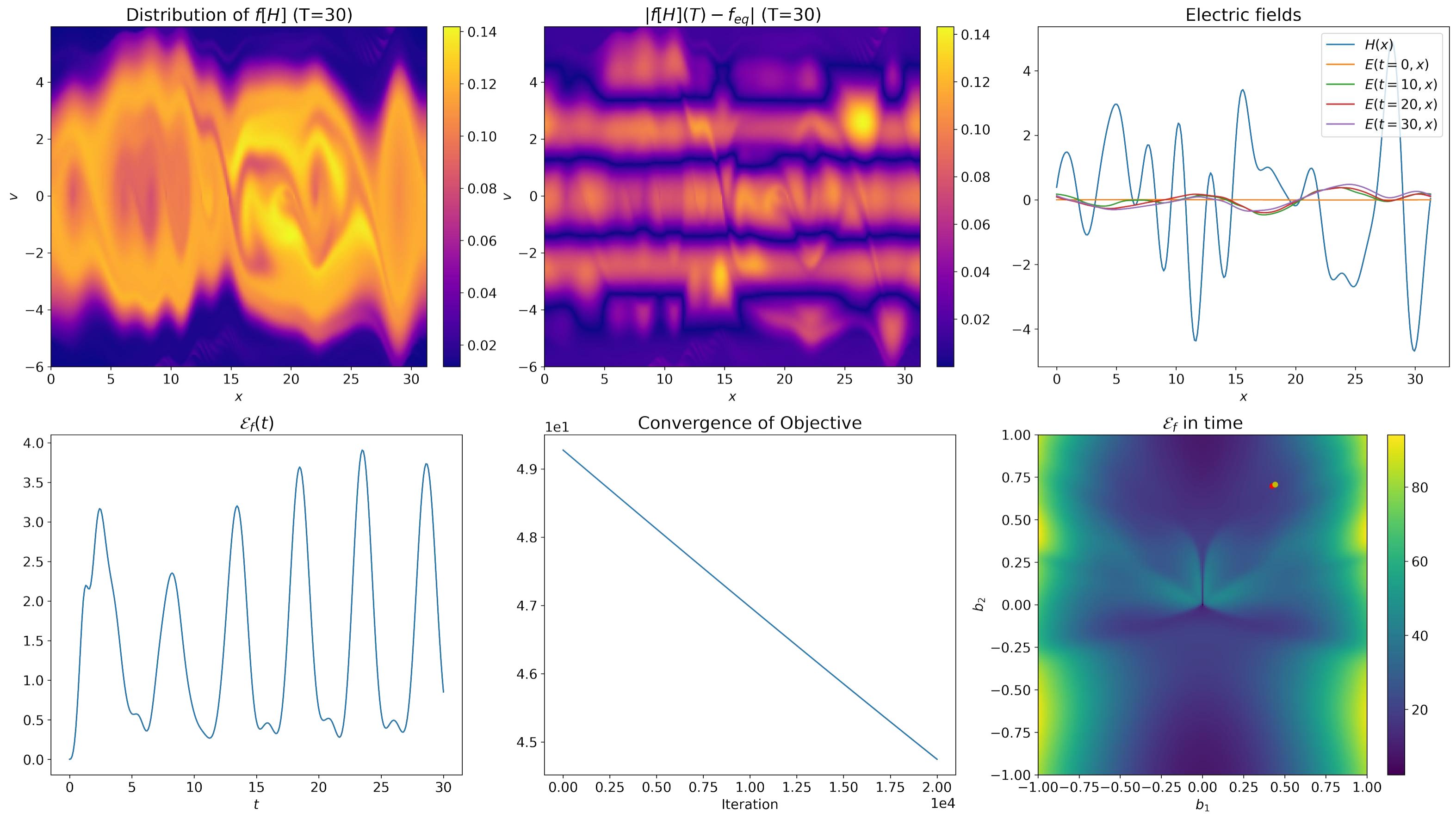}
    \caption{Simulation of~\eqref{eq:vlasov-poisoon_system_ext_1d} with over-parametrized $H$ obtained from~\eqref{eq:optimization_pb_simple} using~\eqref{eq:EET_obj} with far initialization using GD with constant stepsize. From left to right and top to bottom: $f[H](T=30,x,v)$, $|f[H](T,x,v)-f_{\text{eq}}(v)|$, $H$ and $E_{f[H]}(t,x)$, $\mathcal{E}_{f[H]}(t)$, convergence of objective and, trajectory over the landscape of the objective (yellow dot is initial guess).}
    \label{fig:TS_ee_GD_far_over}
\end{figure}

\begin{figure}[H]
    \centering
    \includegraphics[width=0.85\linewidth]{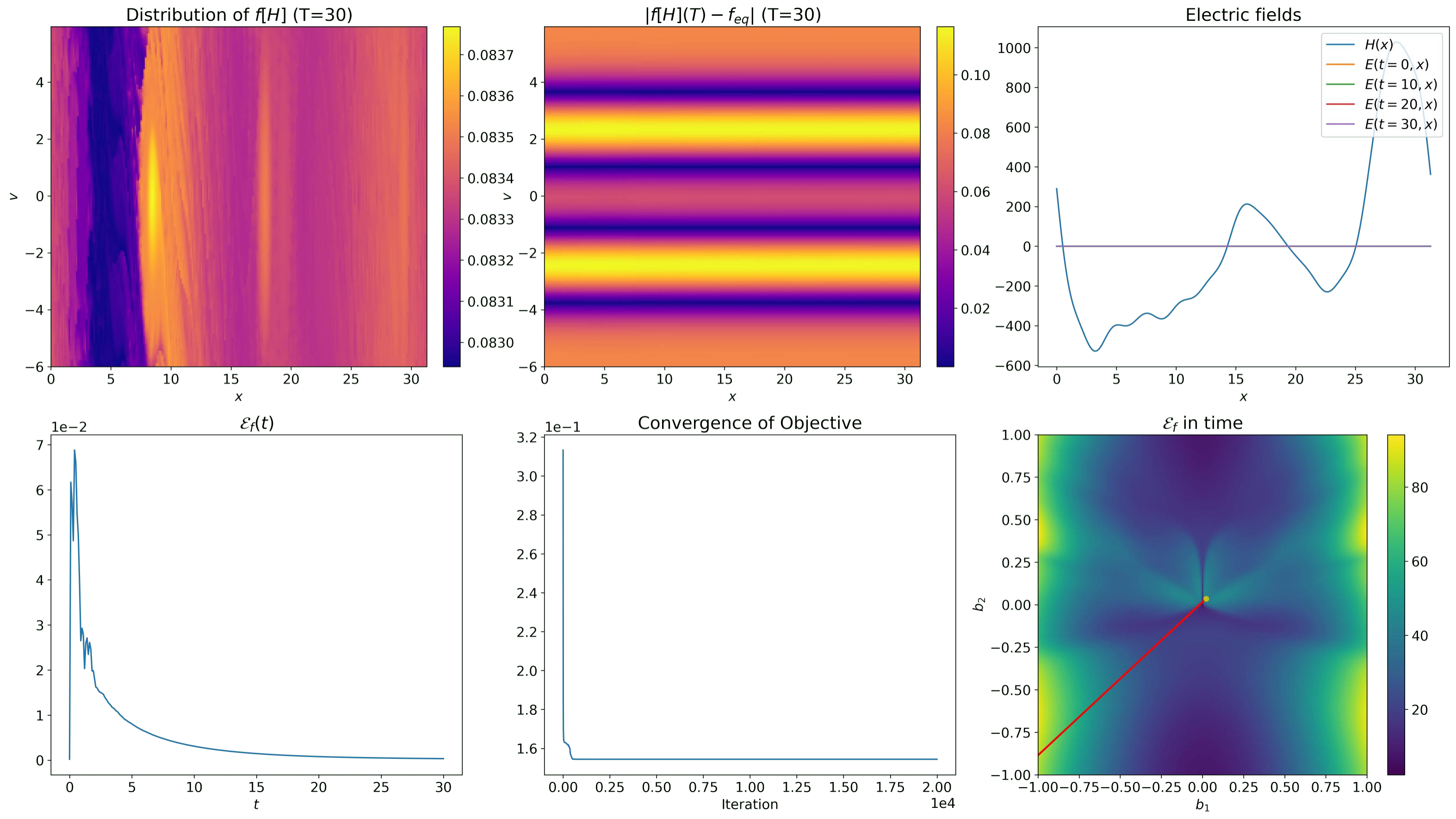}
    \caption{Simulation of~\eqref{eq:vlasov-poisoon_system_ext_1d} with over-parametrized $H$ obtained from~\eqref{eq:optimization_pb_simple} using~\eqref{eq:EET_obj} with near initialization using GD with line-search. From left to right and top to bottom: $f[H](T=30,x,v)$, $|f[H](T,x,v)-f_{\text{eq}}(v)|$, $H$ and $E_{f[H]}(t,x)$, $\mathcal{E}_{f[H]}(t)$, convergence of objective and, trajectory over the landscape of the objective (yellow dot is initial guess).}
    \label{fig:TS_ee_GDL_near_over}
\end{figure}

\begin{figure}[H]
    \centering
    \includegraphics[width=0.85\linewidth]{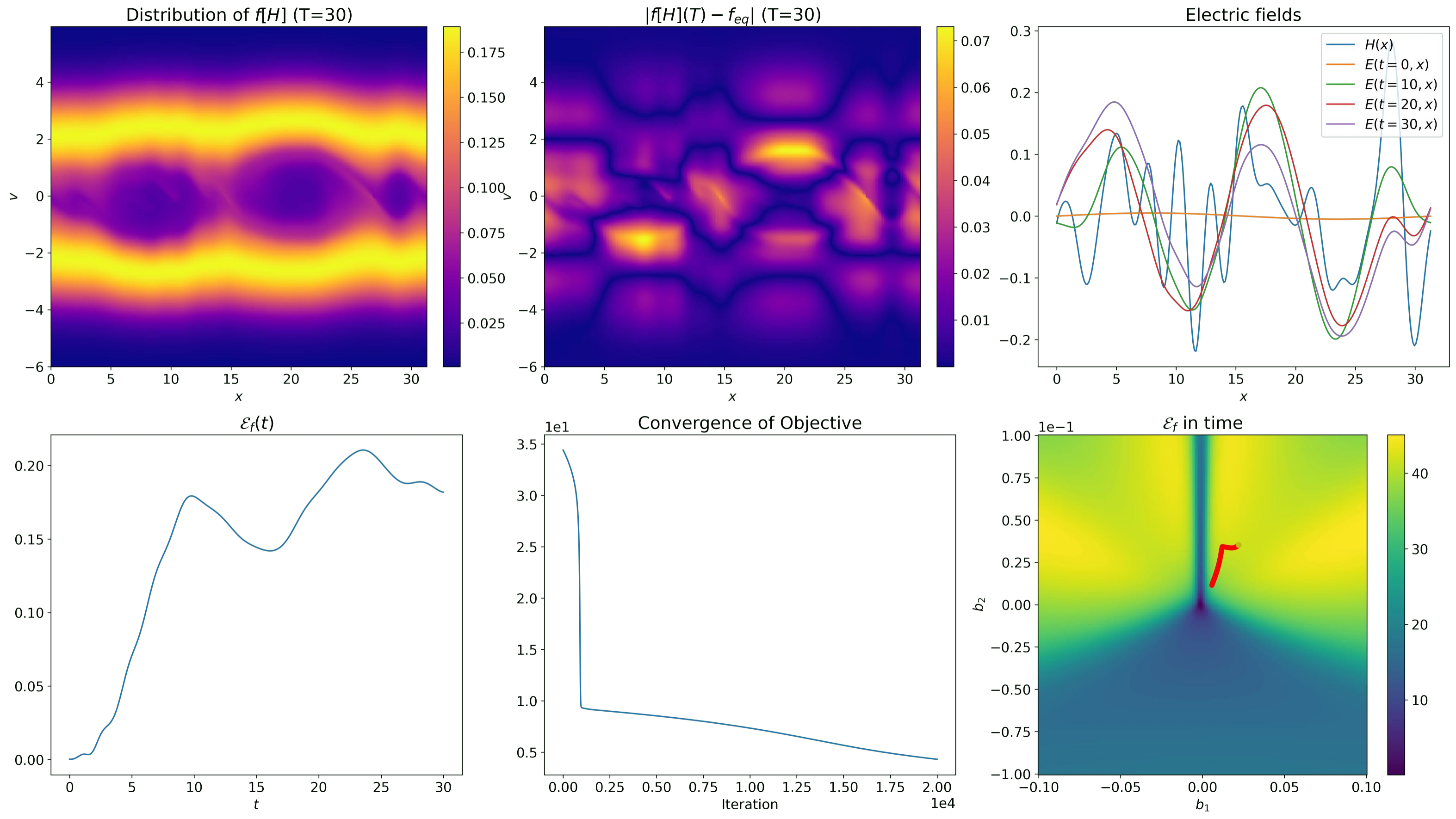}
    \caption{Simulation of~\eqref{eq:vlasov-poisoon_system_ext_1d} using~\eqref{eq:EET_obj} with over-parametrized $H$ obtained from~\eqref{eq:optimization_pb_simple} with near initialization using GD with constant stepsize. From left to right and top to bottom: $f[H](T=30,x,v)$, $|f[H](T,x,v)-f_{\text{eq}}(v)|$, $H$ and $E_{f[H]}(t,x)$, $\mathcal{E}_{f[H]}(t)$, convergence of objective and, trajectory over the landscape of the objective (yellow dot is initial guess).}
    \label{fig:TS_ee_GD_near_over}
\end{figure}

\begin{figure}[H]
    \centering
    \includegraphics[width=0.85\linewidth]{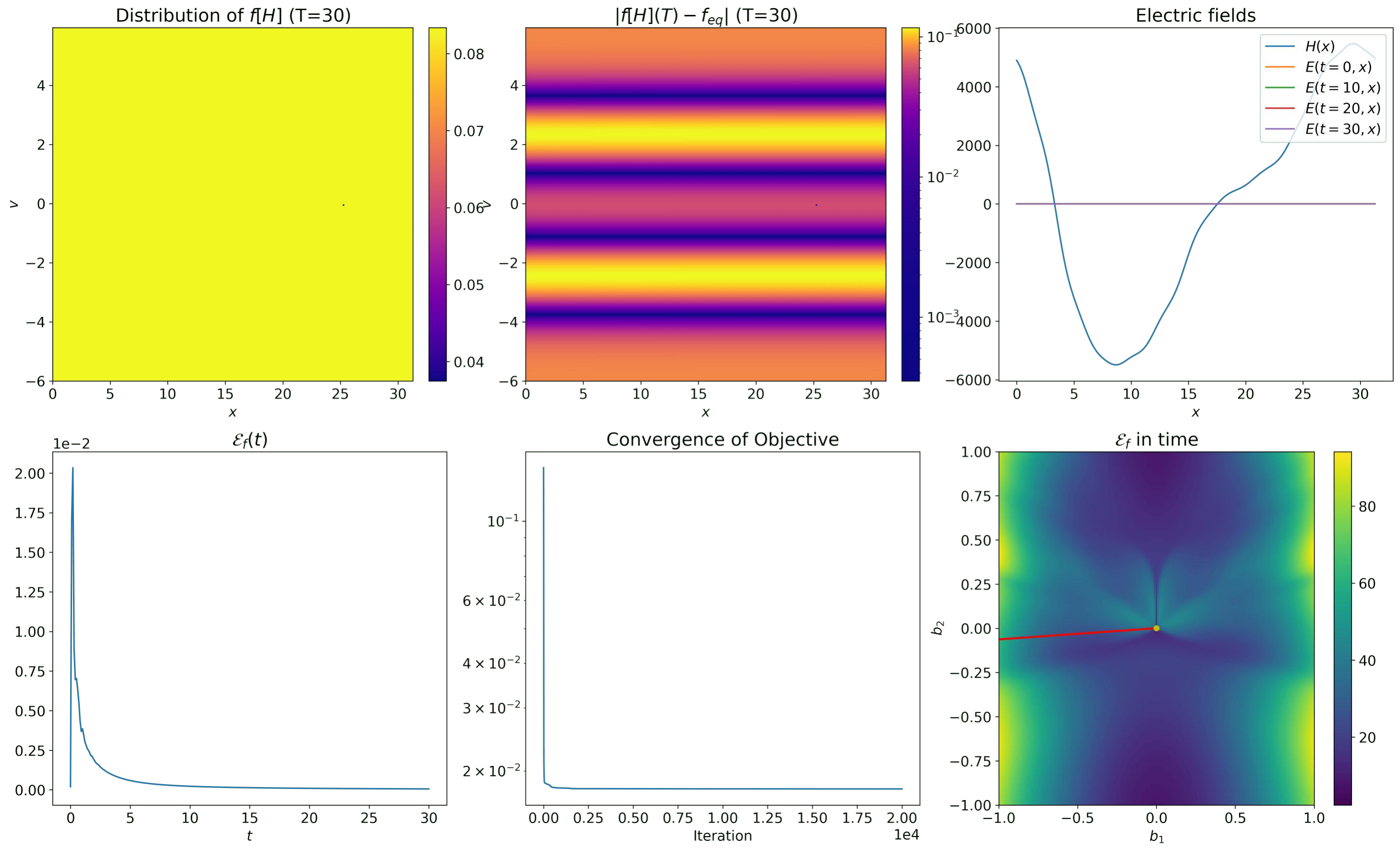}
    \caption{Simulation of~\eqref{eq:vlasov-poisoon_system_ext_1d} with over-parametrized $H$ obtained from~\eqref{eq:optimization_pb_simple} using~\eqref{eq:EET_obj} with local initialization using GD with line-search. From left to right and top to bottom: $f[H](T=30,x,v)$, $|f[H](T,x,v)-f_{\text{eq}}(v)|$, $H$ and $E_{f[H]}(t,x)$, $\mathcal{E}_{f[H]}(t)$, convergence of objective and, trajectory over the landscape of the objective (yellow dot is initial guess).}
    \label{fig:TS_ee_GDL_local_over}
\end{figure}

\begin{figure}[H]
    \centering
    \includegraphics[width=0.85\linewidth]{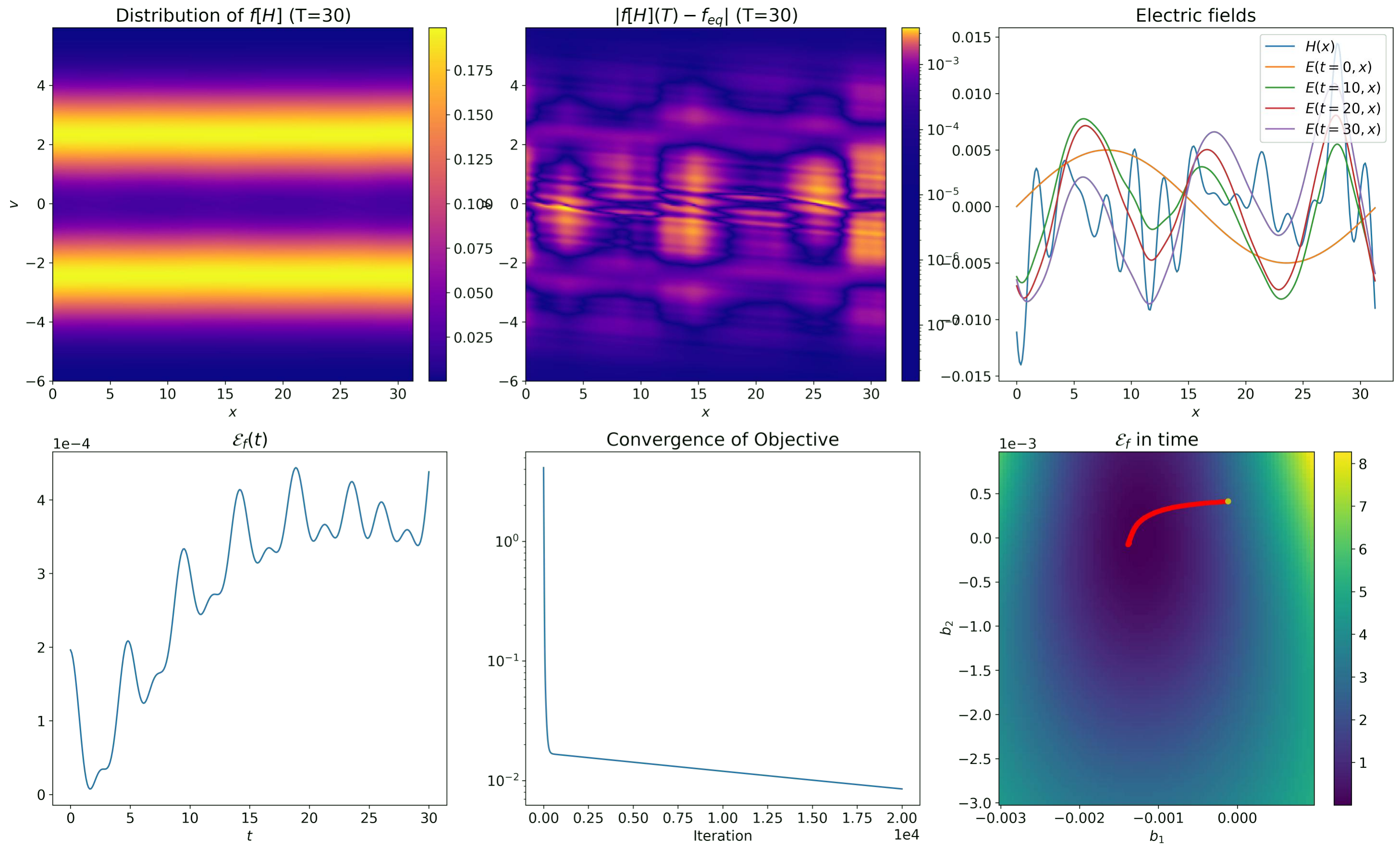}
    \caption{Simulation of~\eqref{eq:vlasov-poisoon_system_ext_1d} using~\eqref{eq:EET_obj} with over-parametrized $H$ obtained from~\eqref{eq:optimization_pb_simple} with local initialization using GD with constant stepsize. From left to right and top to bottom: $f[H](T=30,x,v)$, $|f[H](T,x,v)-f_{\text{eq}}(v)|$, $H$ and $E_{f[H]}(t,x)$, $\mathcal{E}_{f[H]}(t)$, convergence of objective and, trajectory over the landscape of the objective (yellow dot is initial guess).}
    \label{fig:TS_ee_GD_local_over}
\end{figure}

\subsection{Bump-on-Tail example}\label{sec:BoT_summary}

\begin{figure}[H]
    \centering
    \includegraphics[width=0.85\linewidth]{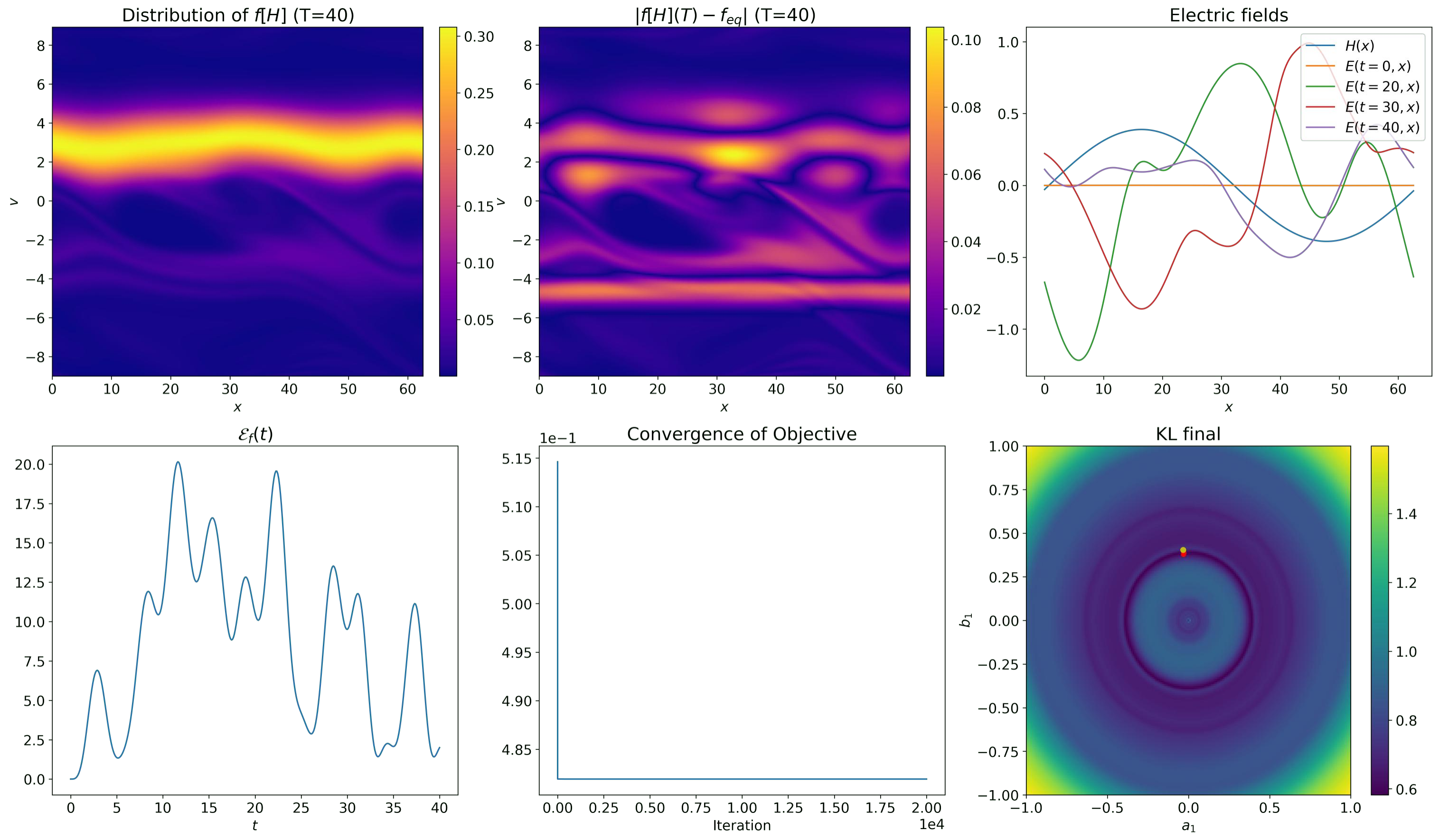}
    \caption{Simulation of~\eqref{eq:vlasov-poisoon_system_ext_1d} with under-parametrized $H$ obtained from~\eqref{eq:optimization_pb_simple} using~\eqref{eq:KL_obj} with far initialization using GD with line-search. From left to right and top to bottom: $f[H](T=30,x,v)$, $|f[H](T,x,v)-f_{\text{eq}}(v)|$, $H$ and $E_{f[H]}(t,x)$, $\mathcal{E}_{f[H]}(t)$, convergence of objective and, trajectory over the landscape of the objective (yellow dot is initial guess).}
    \label{fig:BoT_KL_GDL_far_under}
\end{figure}

\begin{figure}[H]
    \centering
    \includegraphics[width=0.85\linewidth]{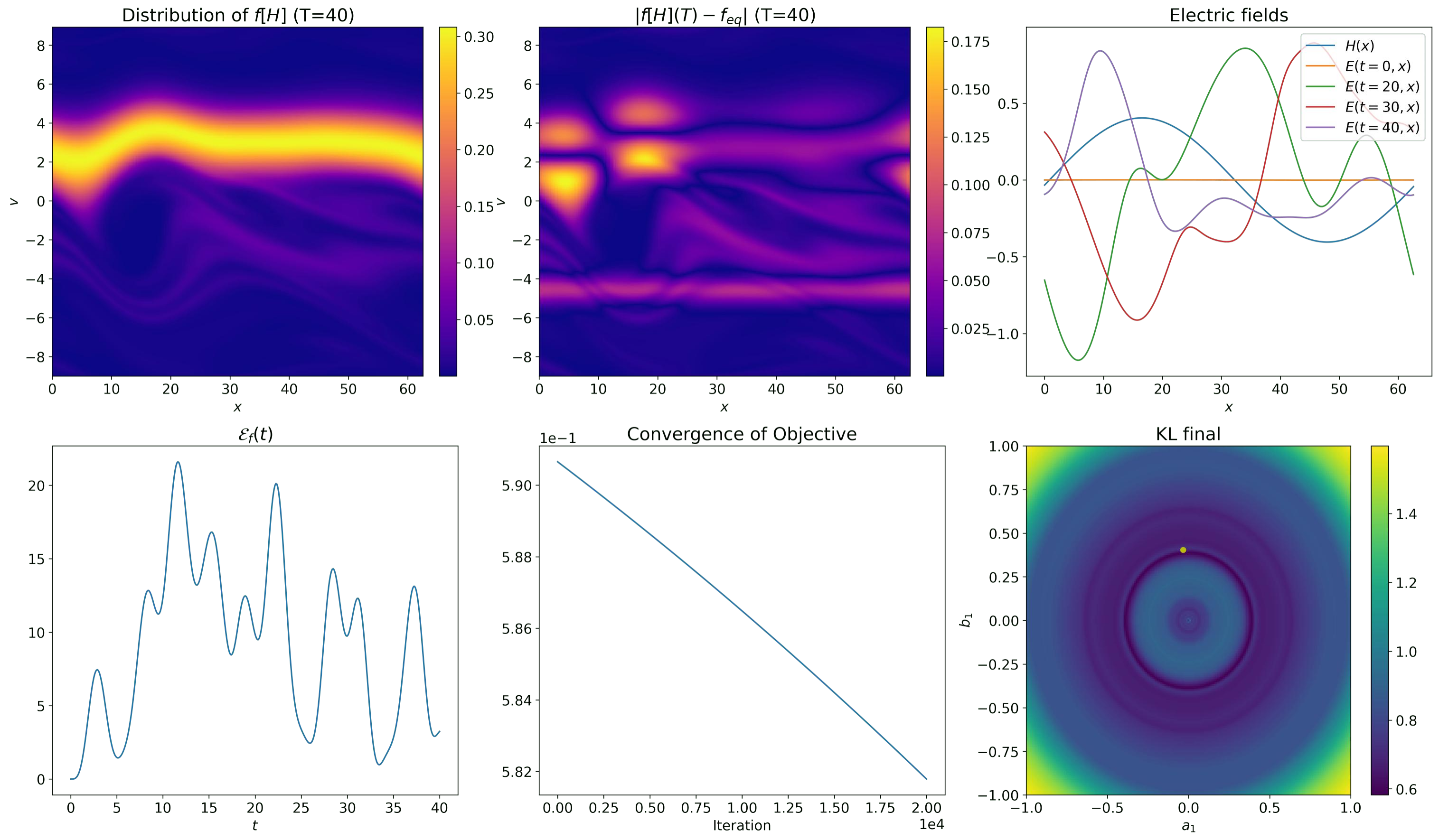}
    \caption{Simulation of~\eqref{eq:vlasov-poisoon_system_ext_1d} with under-parametrized $H$ obtained from~\eqref{eq:optimization_pb_simple} using~\eqref{eq:KL_obj} with far initialization using GD with constant stepsize. From left to right and top to bottom: $f[H](T=30,x,v)$, $|f[H](T,x,v)-f_{\text{eq}}(v)|$, $H$ and $E_{f[H]}(t,x)$, $\mathcal{E}_{f[H]}(t)$, convergence of objective and, trajectory over the landscape of the objective (yellow dot is initial guess).}
    \label{fig:BoT_KL_GD_far_under}
\end{figure}

\begin{figure}[H]
    \centering
    \includegraphics[width=0.85\linewidth]{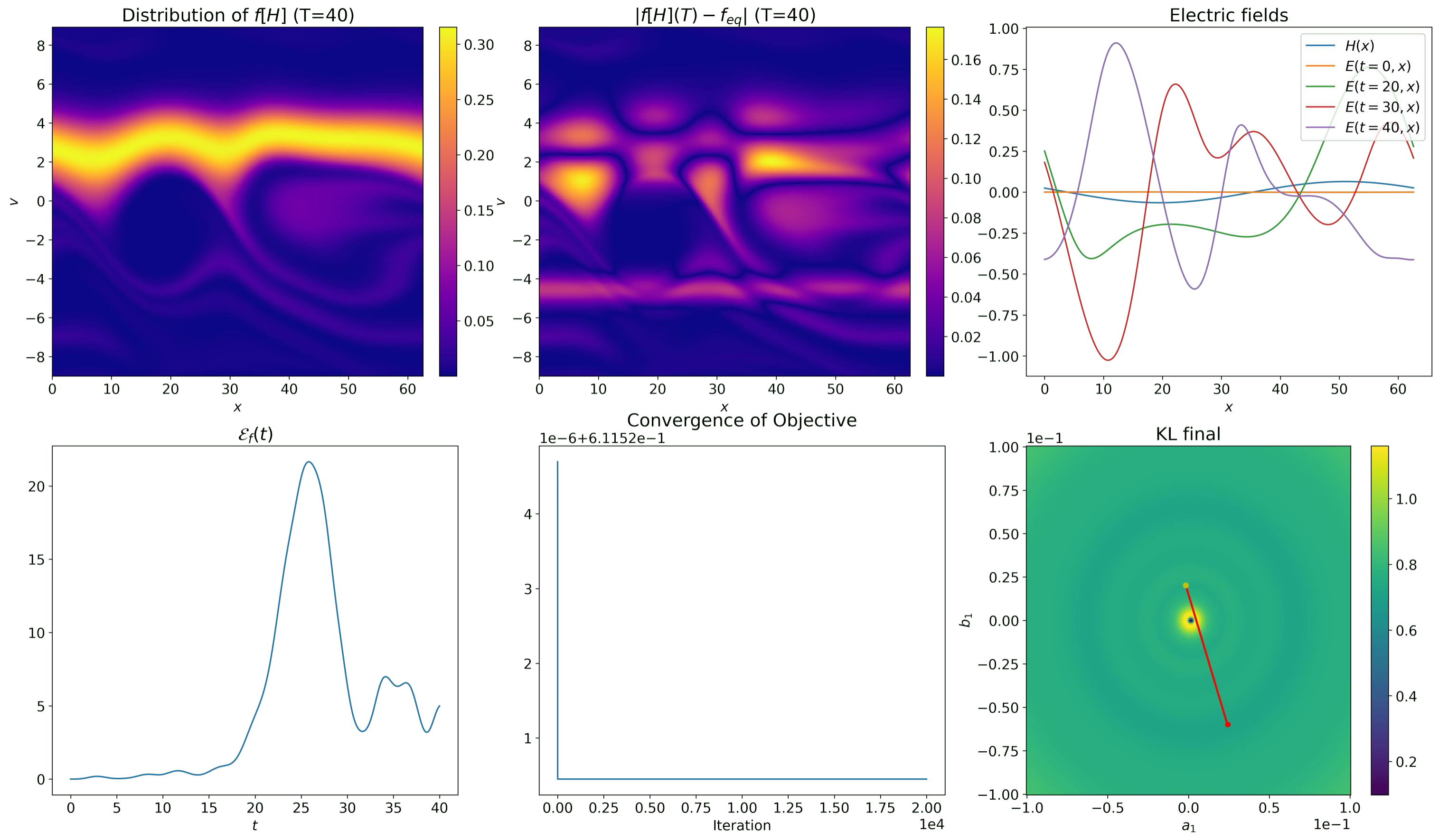}
    \caption{Simulation of~\eqref{eq:vlasov-poisoon_system_ext_1d} with under-parametrized $H$ obtained from~\eqref{eq:optimization_pb_simple} using~\eqref{eq:KL_obj} with near initialization using GD with line-search. From left to right and top to bottom: $f[H](T=30,x,v)$, $|f[H](T,x,v)-f_{\text{eq}}(v)|$, $H$ and $E_{f[H]}(t,x)$, $\mathcal{E}_{f[H]}(t)$, convergence of objective and, trajectory over the landscape of the objective (yellow dot is initial guess).}
    \label{fig:BoT_KL_GDL_near_under}
\end{figure}

\begin{figure}[H]
    \centering
    \includegraphics[width=0.85\linewidth]{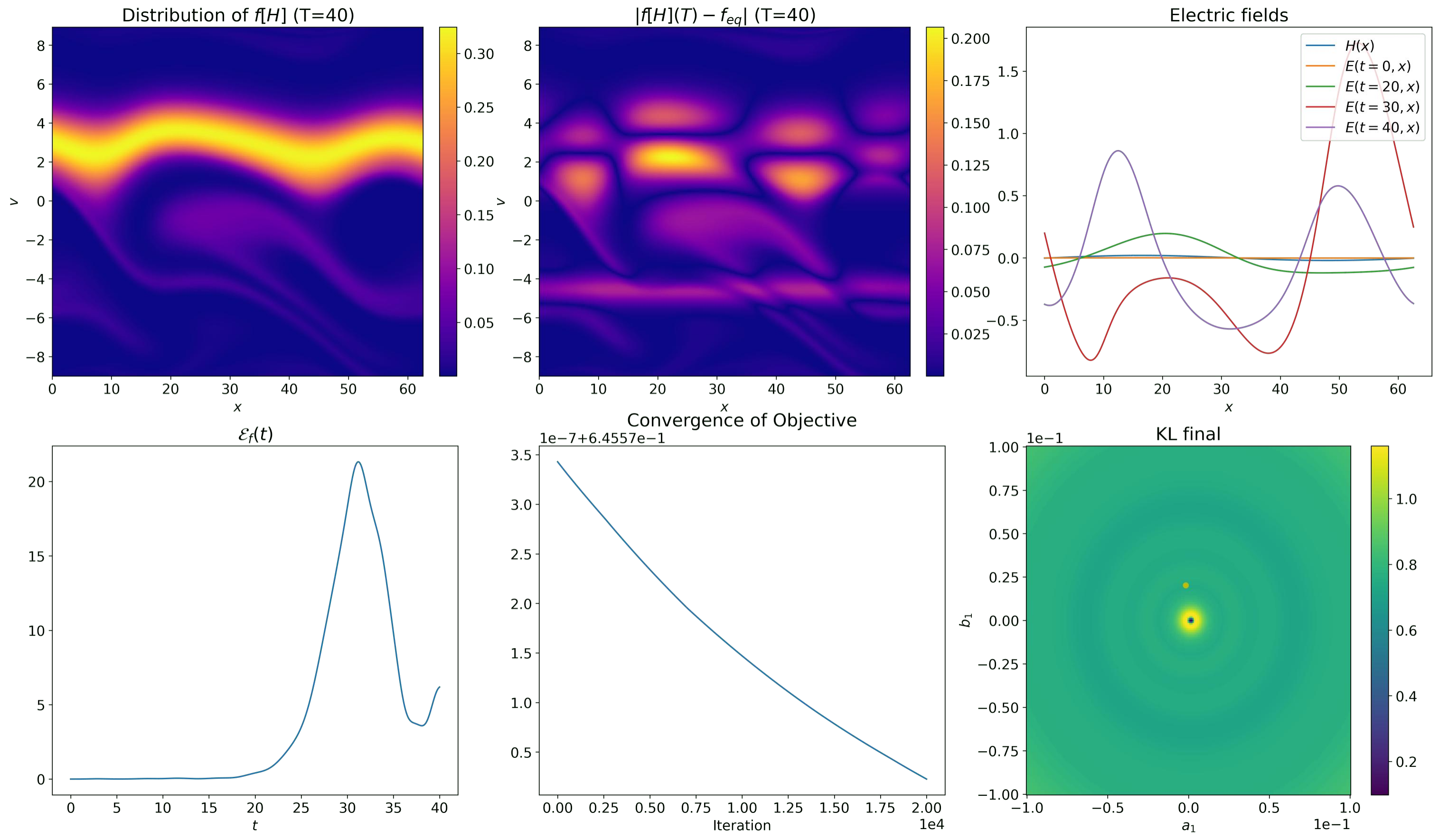}
    \caption{Simulation of~\eqref{eq:vlasov-poisoon_system_ext_1d} with under-parametrized $H$ obtained from~\eqref{eq:optimization_pb_simple} using~\eqref{eq:KL_obj} with near initialization using GD with constant stepsize. From left to right and top to bottom: $f[H](T=30,x,v)$, $|f[H](T,x,v)-f_{\text{eq}}(v)|$, $H$ and $E_{f[H]}(t,x)$, $\mathcal{E}_{f[H]}(t)$, convergence of objective and, trajectory over the landscape of the objective (yellow dot is initial guess).}
    \label{fig:BoT_KL_GD_near_under}
\end{figure}

\begin{figure}[H]
    \centering
    \includegraphics[width=0.85\linewidth]{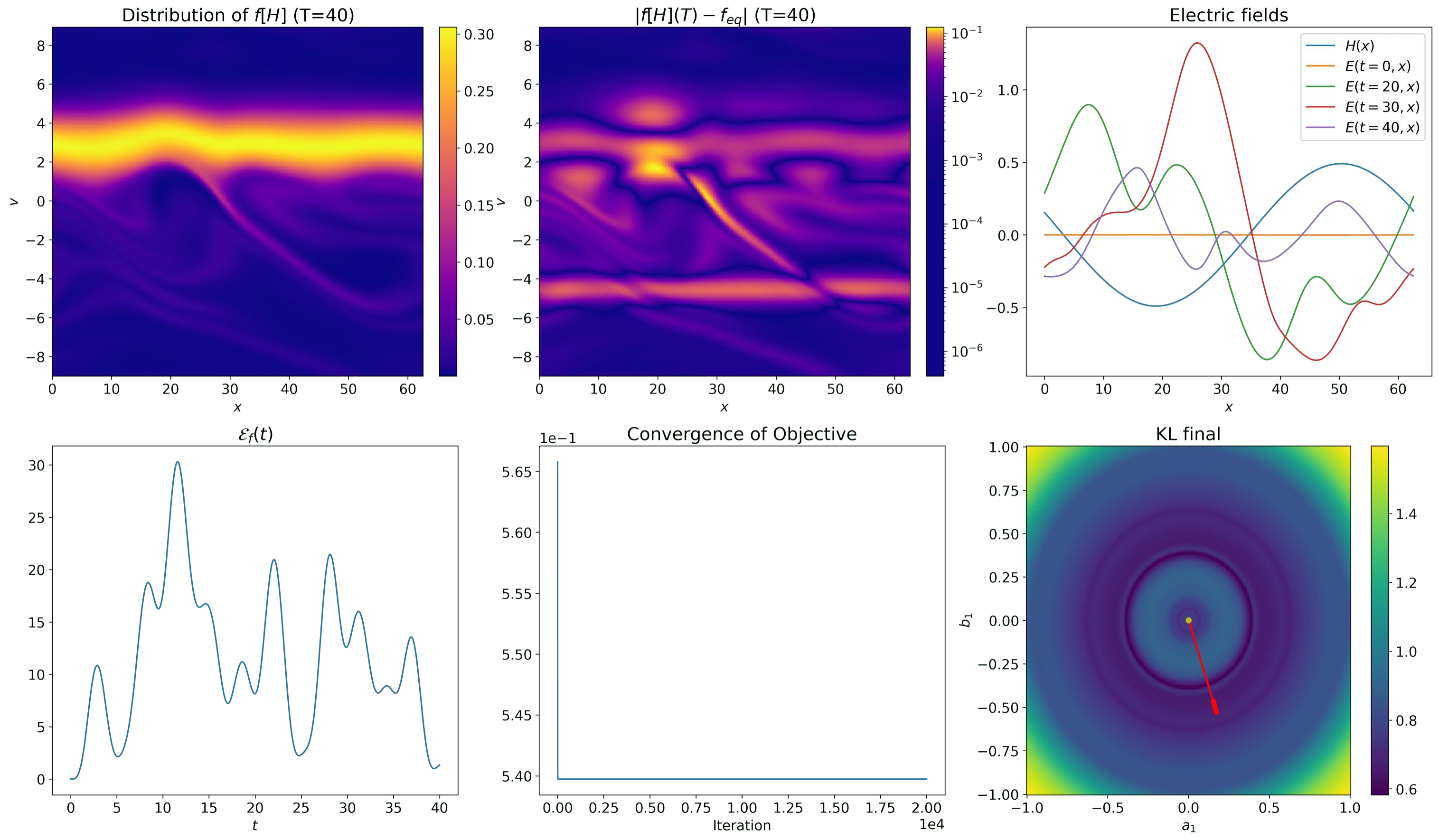}
    \caption{Simulation of~\eqref{eq:vlasov-poisoon_system_ext_1d} with under-parametrized $H$ obtained from~\eqref{eq:optimization_pb_simple} using~\eqref{eq:KL_obj} with local initialization using GD with line-search. From left to right and top to bottom: $f[H](T=30,x,v)$, $|f[H](T,x,v)-f_{\text{eq}}(v)|$, $H$ and $E_{f[H]}(t,x)$, $\mathcal{E}_{f[H]}(t)$, convergence of objective and, trajectory over the landscape of the objective (yellow dot is initial guess).}
    \label{fig:BoT_KL_GDL_local_under}
\end{figure}

\begin{figure}[H]
    \centering
    \includegraphics[width=0.85\linewidth]{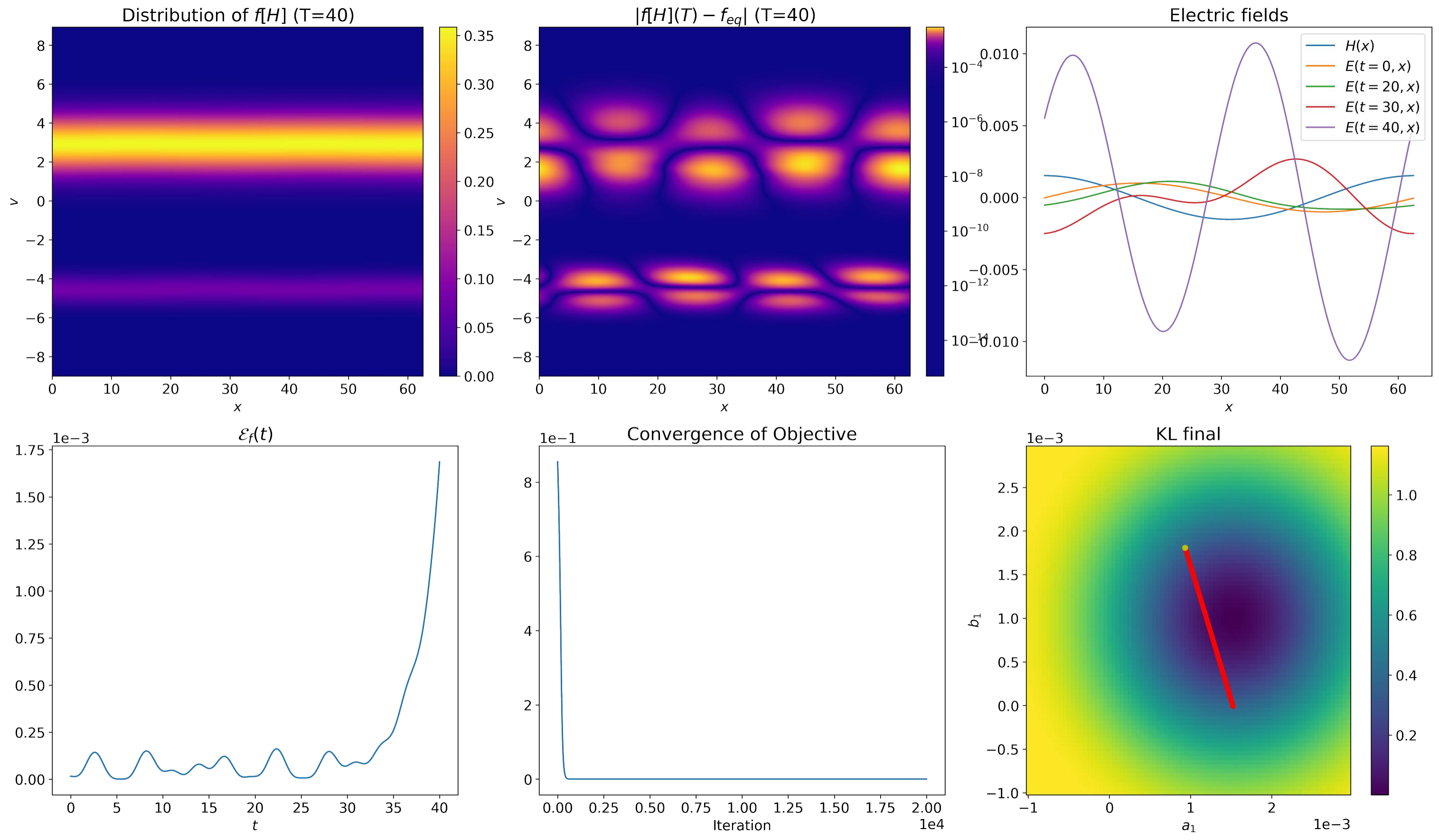}
    \caption{Simulation of~\eqref{eq:vlasov-poisoon_system_ext_1d} with under-parametrized $H$ obtained from~\eqref{eq:optimization_pb_simple} using~\eqref{eq:KL_obj} with local initialization using GD with constant stepsize. From left to right and top to bottom: $f[H](T=30,x,v)$, $|f[H](T,x,v)-f_{\text{eq}}(v)|$, $H$ and $E_{f[H]}(t,x)$, $\mathcal{E}_{f[H]}(t)$, convergence of objective and, trajectory over the landscape of the objective (yellow dot is initial guess).}
    \label{fig:BoT_KL_GD_local_under}
\end{figure}

\begin{figure}[H]
    \centering
    \includegraphics[width=0.85\linewidth]{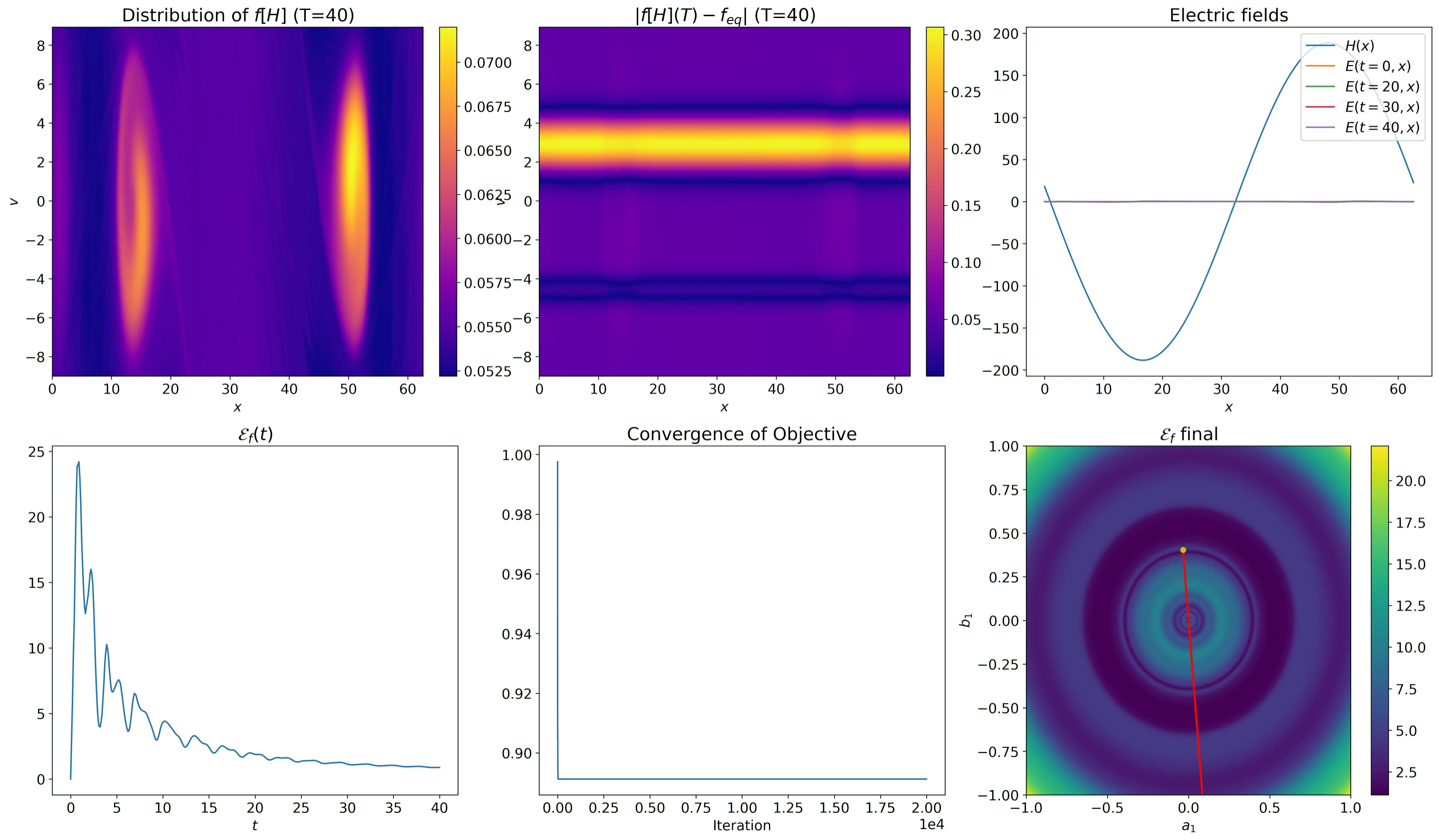}
    \caption{Simulation of~\eqref{eq:vlasov-poisoon_system_ext_1d} with under-parametrized $H$ obtained from~\eqref{eq:optimization_pb_simple} using~\eqref{eq:EE_obj} with far initialization using GD with line-search. From left to right and top to bottom: $f[H](T=30,x,v)$, $|f[H](T,x,v)-f_{\text{eq}}(v)|$, $H$ and $E_{f[H]}(t,x)$, $\mathcal{E}_{f[H]}(t)$, convergence of objective and, trajectory over the landscape of the objective (yellow dot is initial guess).}
    \label{fig:BoT_ee_lf_GDL_far_under}
\end{figure}

\begin{figure}[H]
    \centering
    \includegraphics[width=0.85\linewidth]{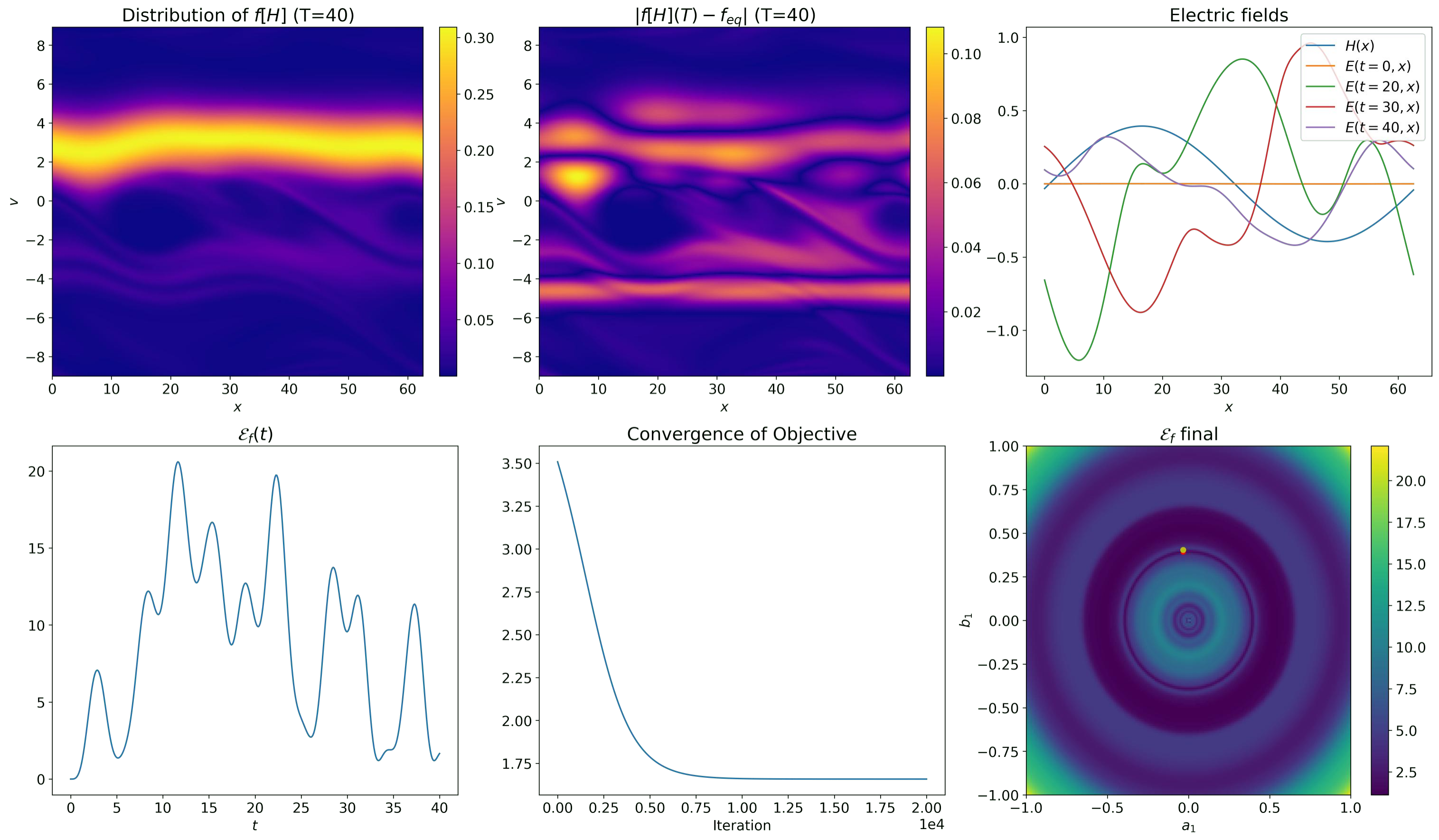}
    \caption{Simulation of~\eqref{eq:vlasov-poisoon_system_ext_1d} with under-parametrized $H$ obtained from~\eqref{eq:optimization_pb_simple} using~\eqref{eq:EE_obj} with far initialization using GD with constant stepsize. From left to right and top to bottom: $f[H](T=30,x,v)$, $|f[H](T,x,v)-f_{\text{eq}}(v)|$, $H$ and $E_{f[H]}(t,x)$, $\mathcal{E}_{f[H]}(t)$, convergence of objective and, trajectory over the landscape of the objective (yellow dot is initial guess).}
    \label{fig:BoT_ee_lf_GD_far_under}
\end{figure}

\begin{figure}[H]
    \centering
    \includegraphics[width=0.85\linewidth]{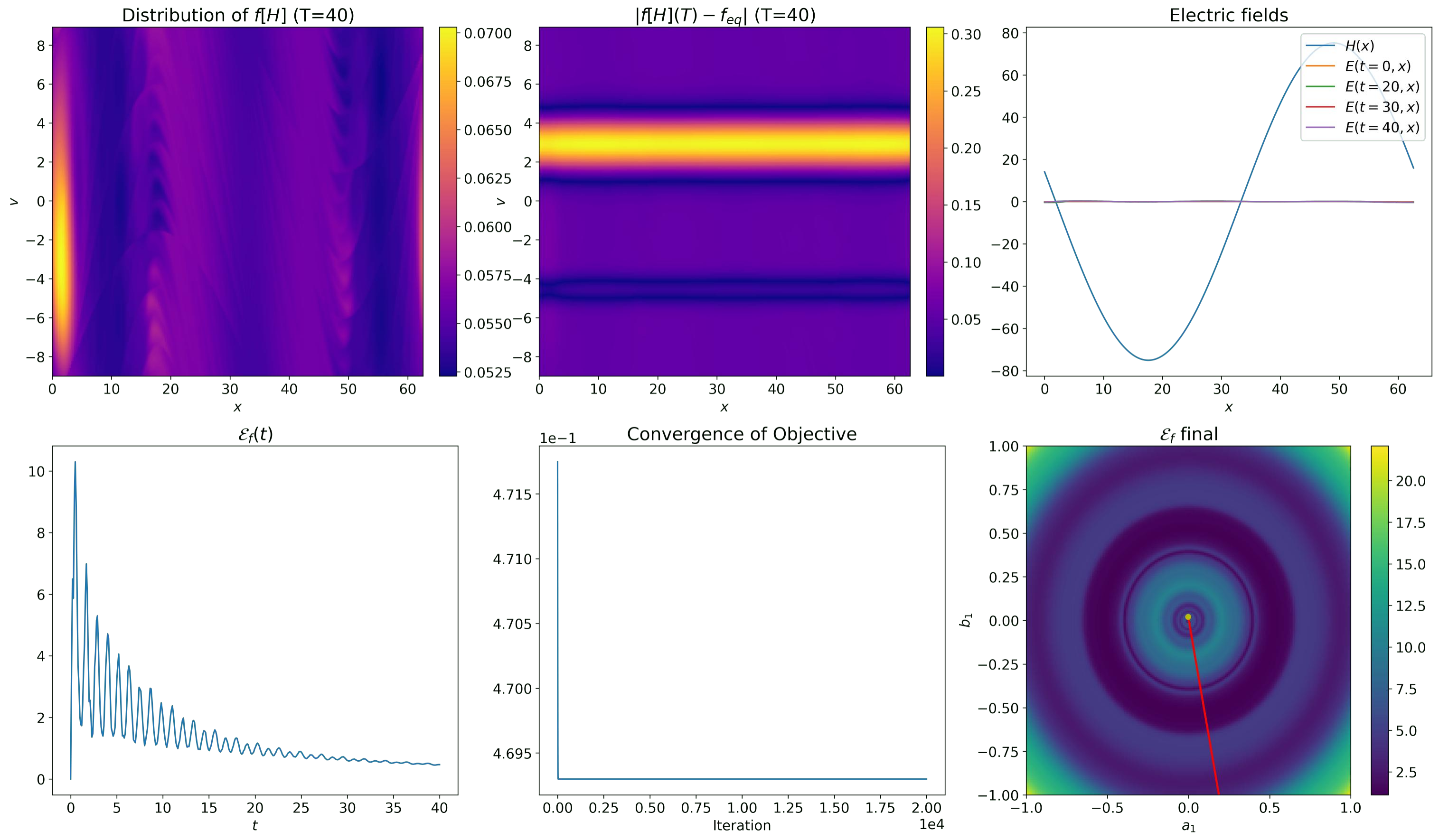}
    \caption{Simulation of~\eqref{eq:vlasov-poisoon_system_ext_1d} with under-parametrized $H$ obtained from~\eqref{eq:optimization_pb_simple} using~\eqref{eq:EE_obj} with near initialization using GD with line-search. From left to right and top to bottom: $f[H](T=30,x,v)$, $|f[H](T,x,v)-f_{\text{eq}}(v)|$, $H$ and $E_{f[H]}(t,x)$, $\mathcal{E}_{f[H]}(t)$, convergence of objective and, trajectory over the landscape of the objective (yellow dot is initial guess).}
    \label{fig:BoT_ee_lf_GDL_near_under}
\end{figure}

\begin{figure}[H]
    \centering
    \includegraphics[width=0.85\linewidth]{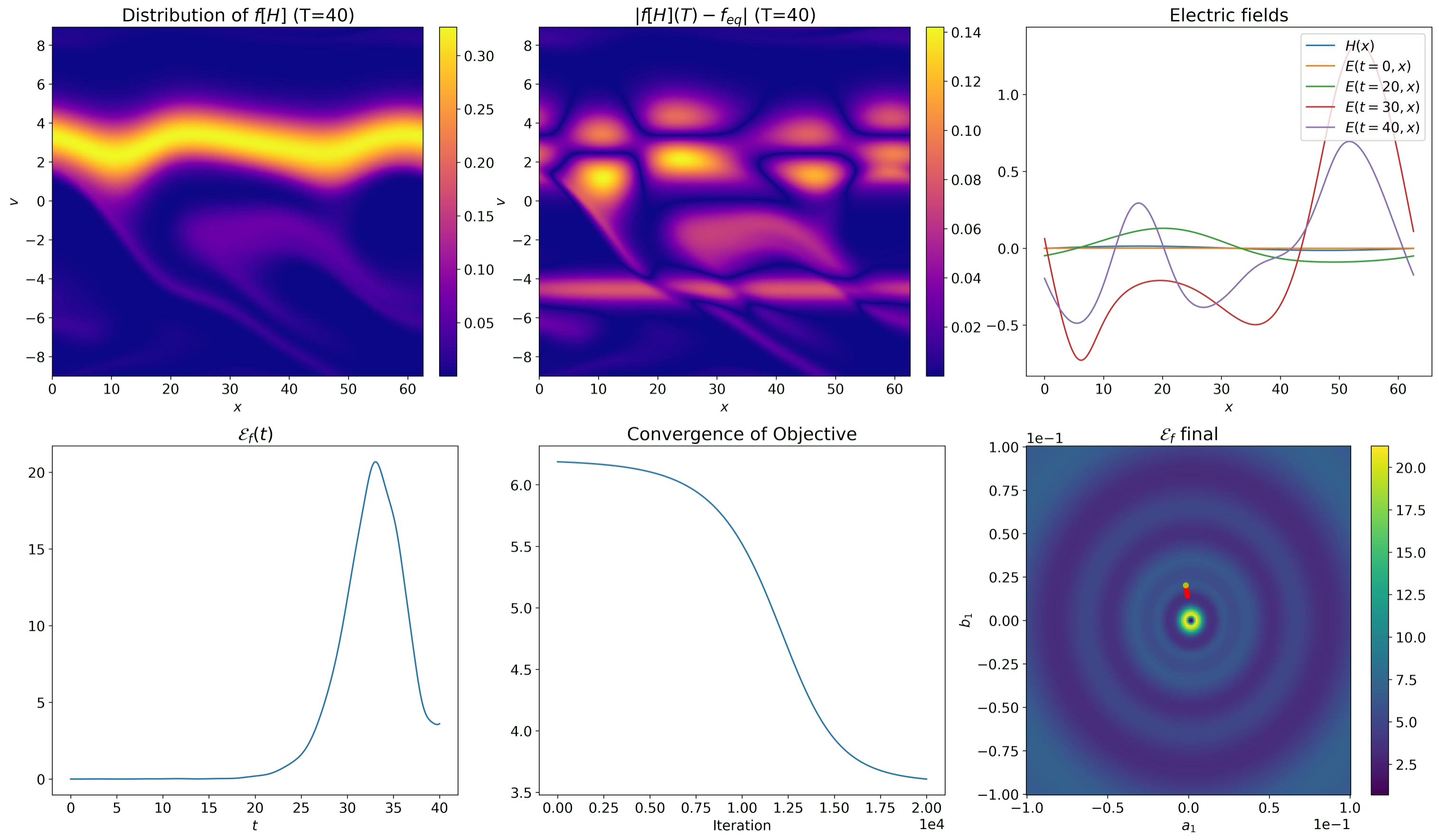}
    \caption{Simulation of~\eqref{eq:vlasov-poisoon_system_ext_1d} using~\eqref{eq:EE_obj} with under-parametrized $H$ obtained from~\eqref{eq:optimization_pb_simple} with near initialization using GD with constant stepsize. From left to right and top to bottom: $f[H](T=30,x,v)$, $|f[H](T,x,v)-f_{\text{eq}}(v)|$, $H$ and $E_{f[H]}(t,x)$, $\mathcal{E}_{f[H]}(t)$, convergence of objective and, trajectory over the landscape of the objective (yellow dot is initial guess).}
    \label{fig:BoT_ee_lf_GD_near_under}
\end{figure}

\begin{figure}[H]
    \centering
    \includegraphics[width=0.85\linewidth]{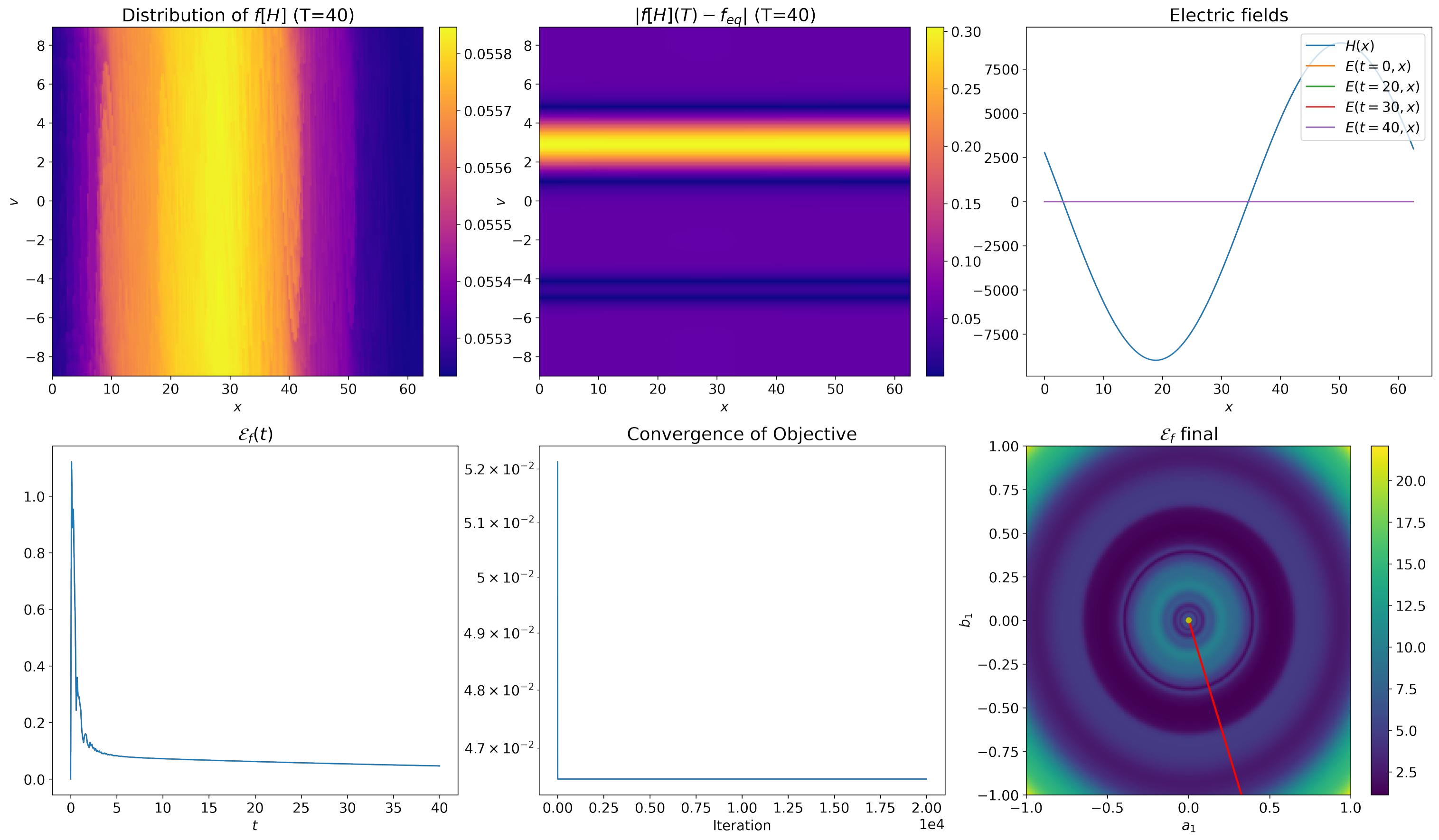}
    \caption{Simulation of~\eqref{eq:vlasov-poisoon_system_ext_1d} with under-parametrized $H$ obtained from~\eqref{eq:optimization_pb_simple} using~\eqref{eq:EE_obj} with local initialization using GD with line-search. From left to right and top to bottom: $f[H](T=30,x,v)$, $|f[H](T,x,v)-f_{\text{eq}}(v)|$,$H$ and $E_{f[H]}(t,x)$, $\mathcal{E}_{f[H]}(t)$, convergence of objective and, trajectory over the landscape of the objective (yellow dot is initial guess).}
    \label{fig:BoT_ee_lf_GDL_local_under}
\end{figure}

\begin{figure}[H]
    \centering
    \includegraphics[width=0.85\linewidth]{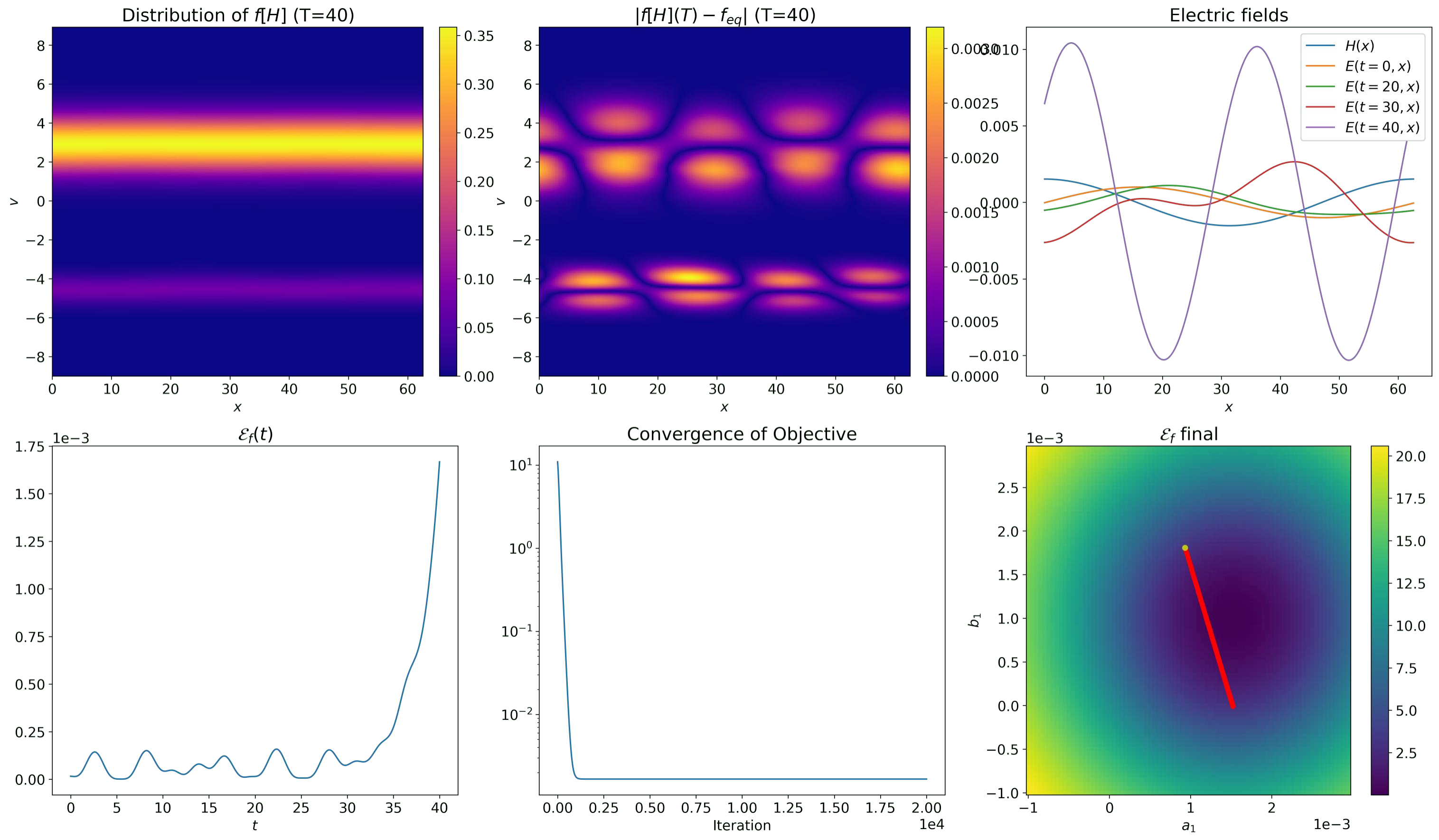}
    \caption{Simulation of~\eqref{eq:vlasov-poisoon_system_ext_1d} using~\eqref{eq:EE_obj} with under-parametrized $H$ obtained from~\eqref{eq:optimization_pb_simple} with local initialization using GD with constant stepsize. From left to right and top to bottom: $f[H](T=30,x,v)$, $|f[H](T,x,v)-f_{\text{eq}}(v)|$, $H$ and $E_{f[H]}(t,x)$, $\mathcal{E}_{f[H]}(t)$, convergence of objective and, trajectory over the landscape of the objective (yellow dot is initial guess).}
    \label{fig:BoT_ee_lf_GD_local_under}
\end{figure}

\begin{figure}[H]
    \centering
    \includegraphics[width=0.85\linewidth]{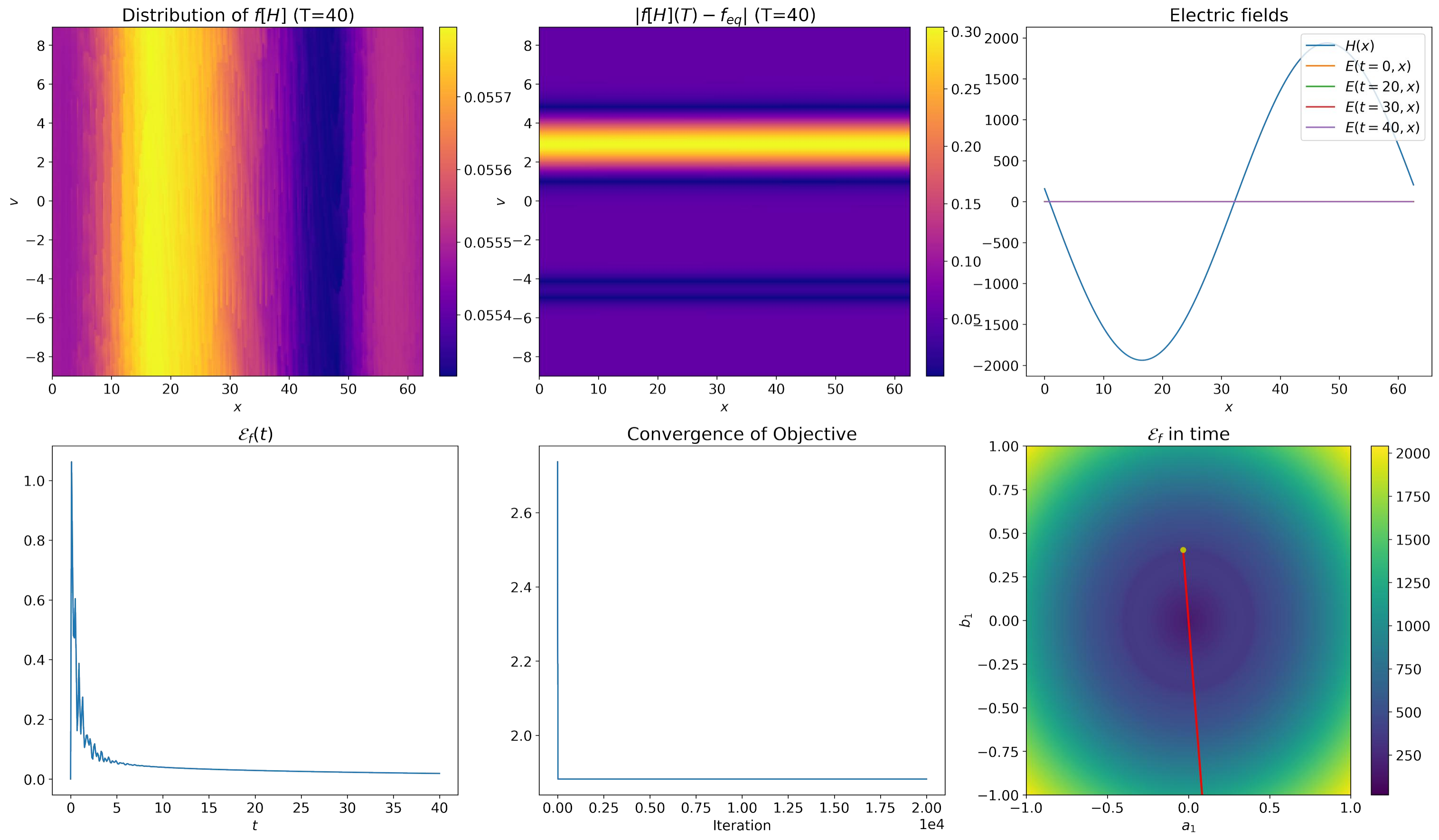}
    \caption{Simulation of~\eqref{eq:vlasov-poisoon_system_ext_1d} with under-parametrized $H$ obtained from~\eqref{eq:optimization_pb_simple} using~\eqref{eq:EET_obj} with far initialization using GD with line-search. From left to right and top to bottom: $f[H](T=30,x,v)$, $|f[H](T,x,v)-f_{\text{eq}}(v)|$, $H$ and $E_{f[H]}(t,x)$, $\mathcal{E}_{f[H]}(t)$, convergence of objective and, trajectory over the landscape of the objective (yellow dot is initial guess).}
    \label{fig:BoT_ee_GDL_far_under}
\end{figure}

\begin{figure}[H]
    \centering
    \includegraphics[width=0.85\linewidth]{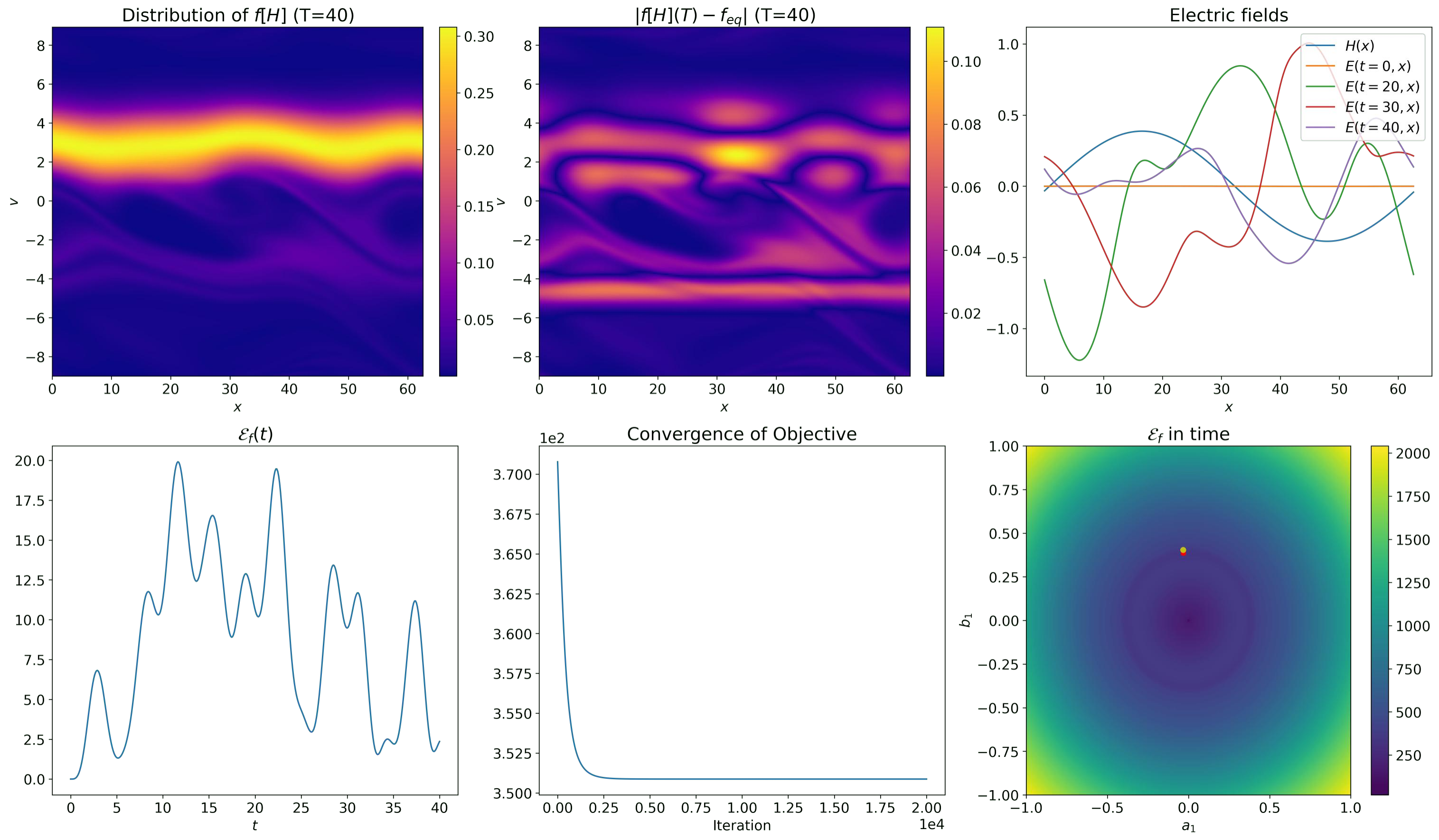}
    \caption{Simulation of~\eqref{eq:vlasov-poisoon_system_ext_1d} with under-parametrized $H$ obtained from~\eqref{eq:optimization_pb_simple} using~\eqref{eq:EET_obj} with far initialization using GD with constant stepsize. From left to right and top to bottom: $f[H](T=30,x,v)$, $|f[H](T,x,v)-f_{\text{eq}}(v)|$, $H$ and $E_{f[H]}(t,x)$, $\mathcal{E}_{f[H]}(t)$, convergence of objective and, trajectory over the landscape of the objective (yellow dot is initial guess).}
    \label{fig:BoT_ee_GD_far_under}
\end{figure}

\begin{figure}[H]
    \centering
    \includegraphics[width=0.85\linewidth]{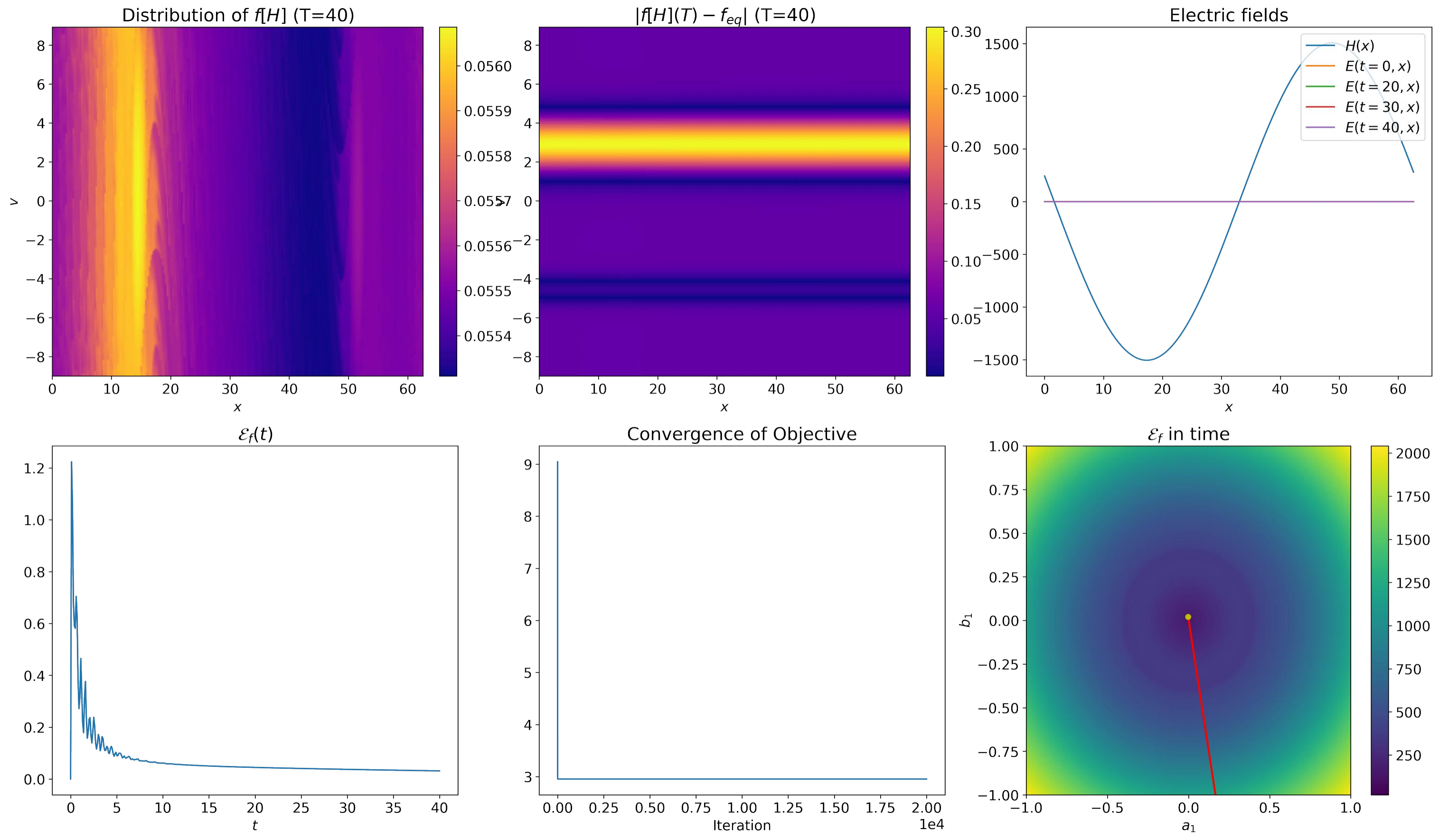}
    \caption{Simulation of~\eqref{eq:vlasov-poisoon_system_ext_1d} with under-parametrized $H$ obtained from~\eqref{eq:optimization_pb_simple} using~\eqref{eq:EET_obj} with near initialization using GD with line-search. From left to right and top to bottom: $f[H](T=30,x,v)$, $|f[H](T,x,v)-f_{\text{eq}}(v)|$, $H$ and $E_{f[H]}(t,x)$, $\mathcal{E}_{f[H]}(t)$, convergence of objective and, trajectory over the landscape of the objective (yellow dot is initial guess).}
    \label{fig:BoT_ee_GDL_near_under}
\end{figure}

\begin{figure}[H]
    \centering
    \includegraphics[width=0.85\linewidth]{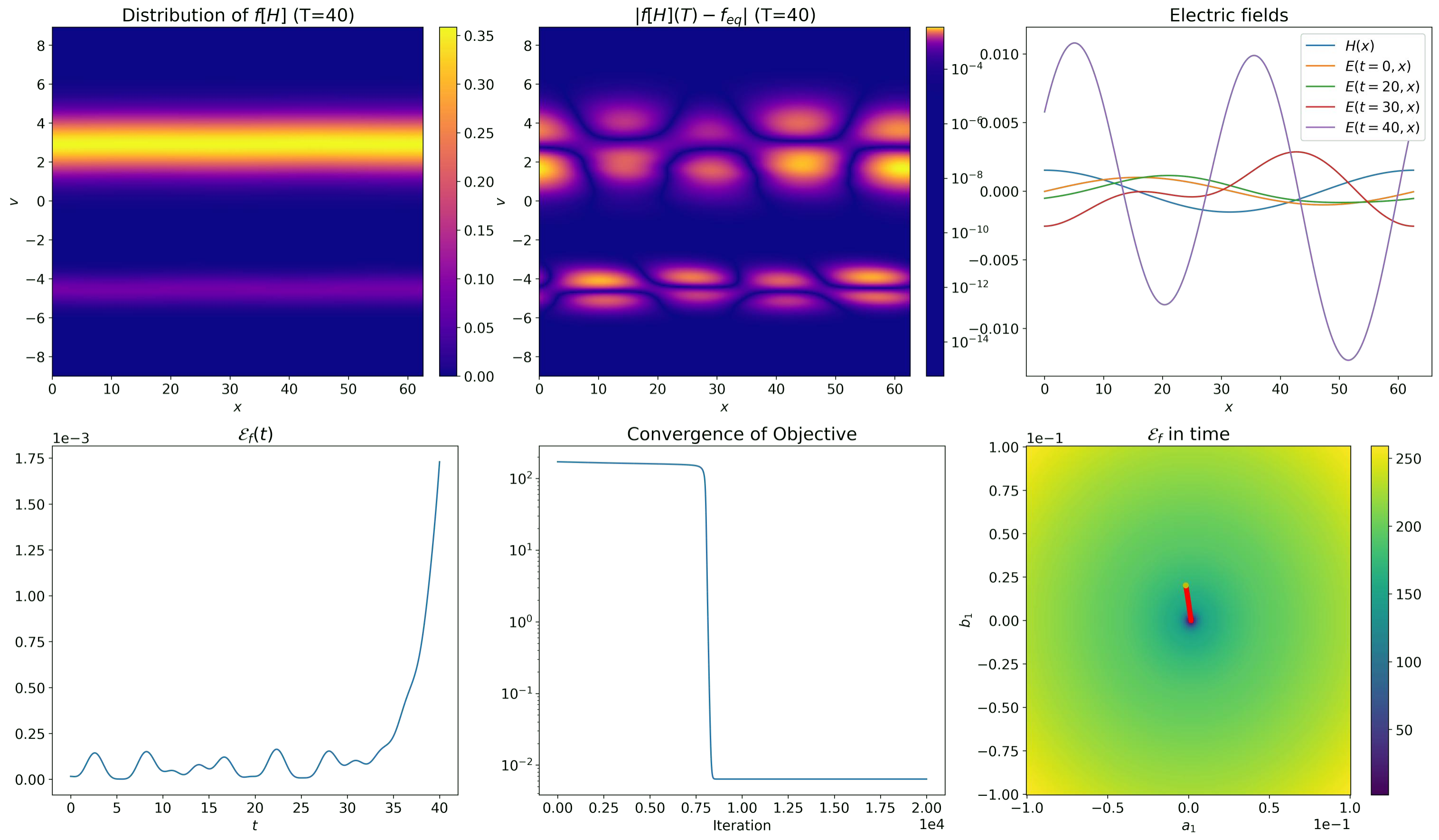}
    \caption{Simulation of~\eqref{eq:vlasov-poisoon_system_ext_1d} using~\eqref{eq:EET_obj} with under-parametrized $H$ obtained from~\eqref{eq:optimization_pb_simple} with near initialization using GD with constant stepsize. From left to right and top to bottom: $f[H](T=30,x,v)$, $|f[H](T,x,v)-f_{\text{eq}}(v)|$, $H$ and $E_{f[H]}(t,x)$, $\mathcal{E}_{f[H]}(t)$, convergence of objective and, trajectory over the landscape of the objective (yellow dot is initial guess).}
    \label{fig:BoT_ee_GD_near_under}
\end{figure}

\begin{figure}[H]
    \centering
    \includegraphics[width=0.85\linewidth]{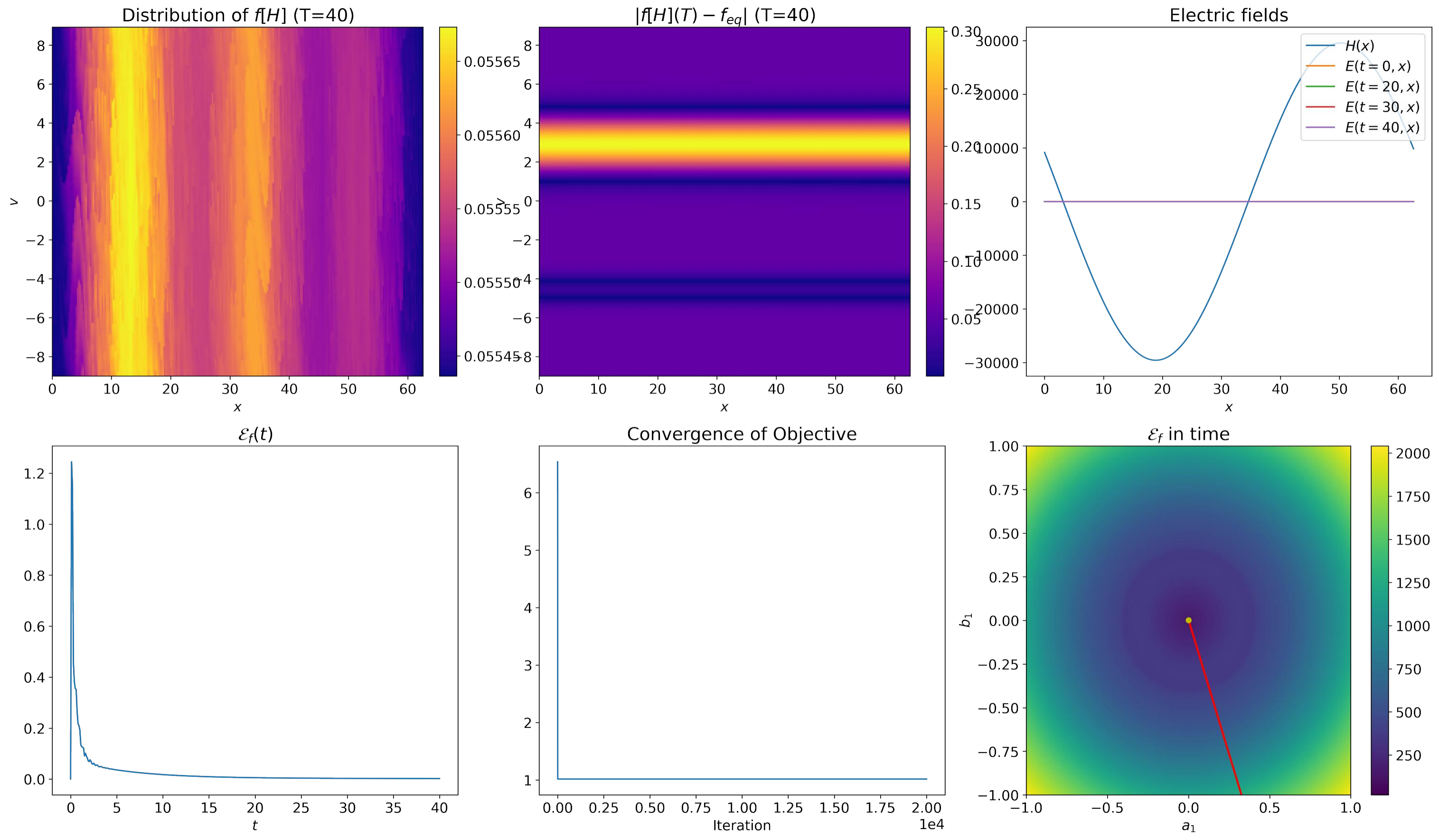}
    \caption{Simulation of~\eqref{eq:vlasov-poisoon_system_ext_1d} with under-parametrized $H$ obtained from~\eqref{eq:optimization_pb_simple} using~\eqref{eq:EET_obj} with local initialization using GD with line-search. From left to right and top to bottom: $f[H](T=30,x,v)$, $|f[H](T,x,v)-f_{\text{eq}}(v)|$,$H$ and $E_{f[H]}(t,x)$, $\mathcal{E}_{f[H]}(t)$, convergence of objective and, trajectory over the landscape of the objective (yellow dot is initial guess).}
    \label{fig:BoT_ee_GDL_local_under}
\end{figure}

\begin{figure}[H]
    \centering
    \includegraphics[width=0.85\linewidth]{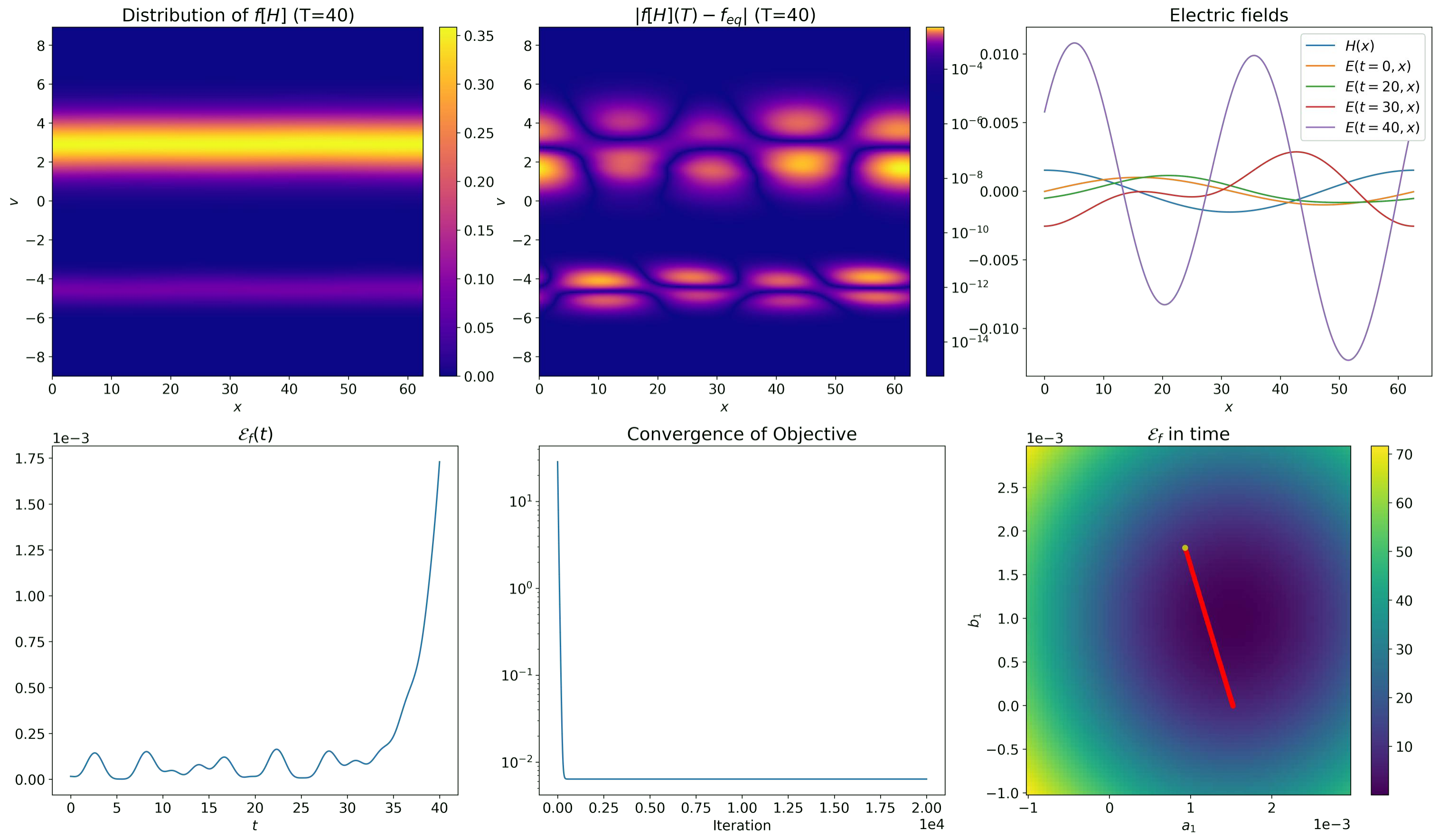}
    \caption{Simulation of~\eqref{eq:vlasov-poisoon_system_ext_1d} using~\eqref{eq:EET_obj} with under-parametrized $H$ obtained from~\eqref{eq:optimization_pb_simple} with local initialization using GD with constant stepsize. From left to right and top to bottom: $f[H](T=30,x,v)$, $|f[H](T,x,v)-f_{\text{eq}}(v)|$, $H$ and $E_{f[H]}(t,x)$, $\mathcal{E}_{f[H]}(t)$, convergence of objective and, trajectory over the landscape of the objective (yellow dot is initial guess).}
    \label{fig:BoT_ee_GD_local_under}
\end{figure}

\begin{figure}[H]
    \centering
    \includegraphics[width=0.85\linewidth]{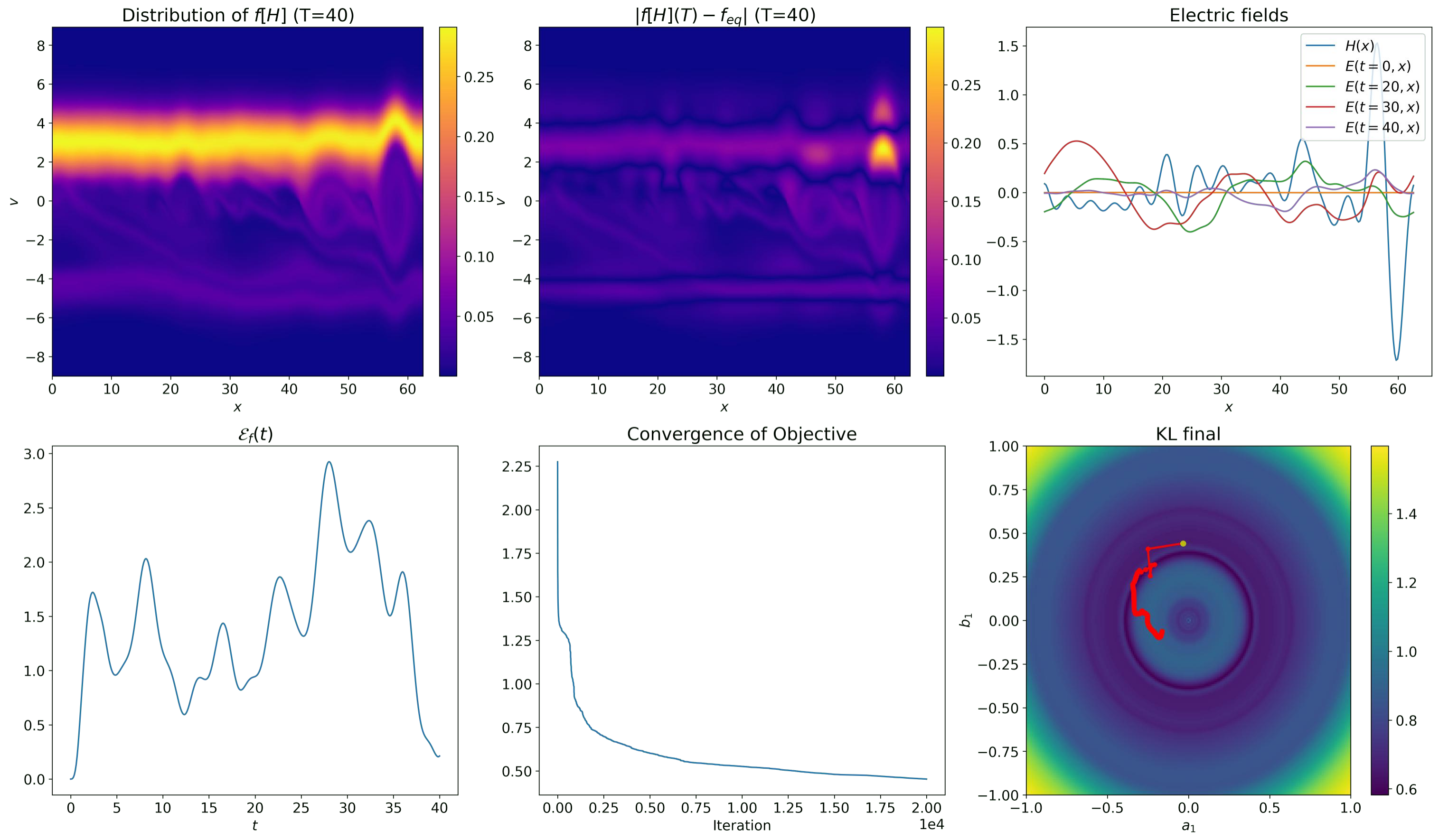}
    \caption{Simulation of~\eqref{eq:vlasov-poisoon_system_ext_1d} with over-parametrized $H$ obtained from~\eqref{eq:optimization_pb_simple} using~\eqref{eq:KL_obj} with far initialization using GD with line-search. From left to right and top to bottom: $f[H](T=30,x,v)$, $|f[H](T,x,v)-f_{\text{eq}}(v)|$, $H$ and $E_{f[H]}(t,x)$, $\mathcal{E}_{f[H]}(t)$, convergence of objective and, trajectory over the landscape of the objective (yellow dot is initial guess).}
    \label{fig:BoT_KL_GDL_far_over}
\end{figure}

\begin{figure}[H]
    \centering
    \includegraphics[width=0.85\linewidth]{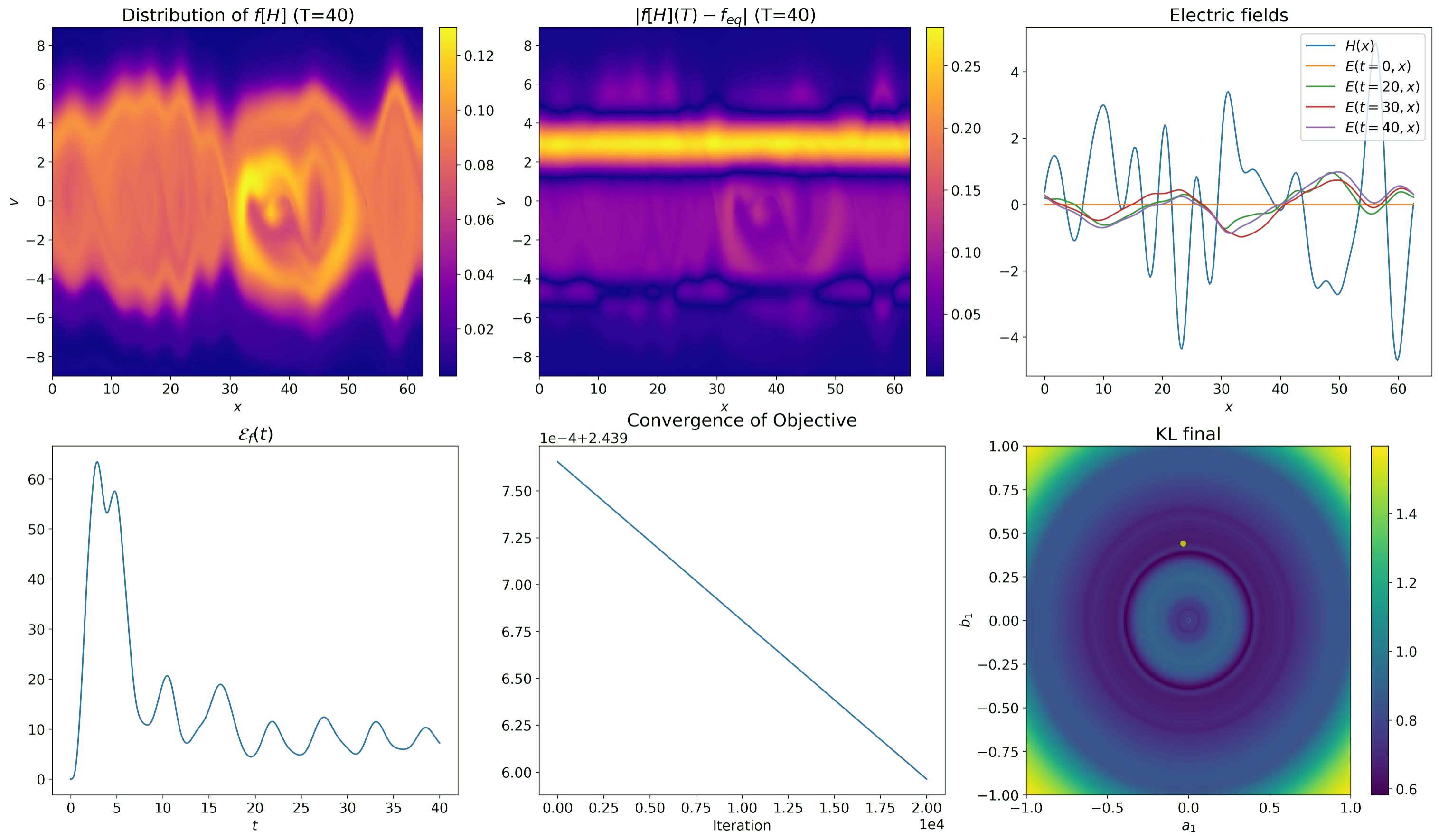}
    \caption{Simulation of~\eqref{eq:vlasov-poisoon_system_ext_1d} with over-parametrized $H$ obtained from~\eqref{eq:optimization_pb_simple} using~\eqref{eq:KL_obj} with far initialization using GD with constant stepsize. From left to right and top to bottom: $f[H](T=30,x,v)$, $|f[H](T,x,v)-f_{\text{eq}}(v)|$, $H$ and $E_{f[H]}(t,x)$, $\mathcal{E}_{f[H]}(t)$, convergence of objective and, trajectory over the landscape of the objective (yellow dot is initial guess).}
    \label{fig:BoT_KL_GD_far_over}
\end{figure}

\begin{figure}[H]
    \centering
    \includegraphics[width=0.85\linewidth]{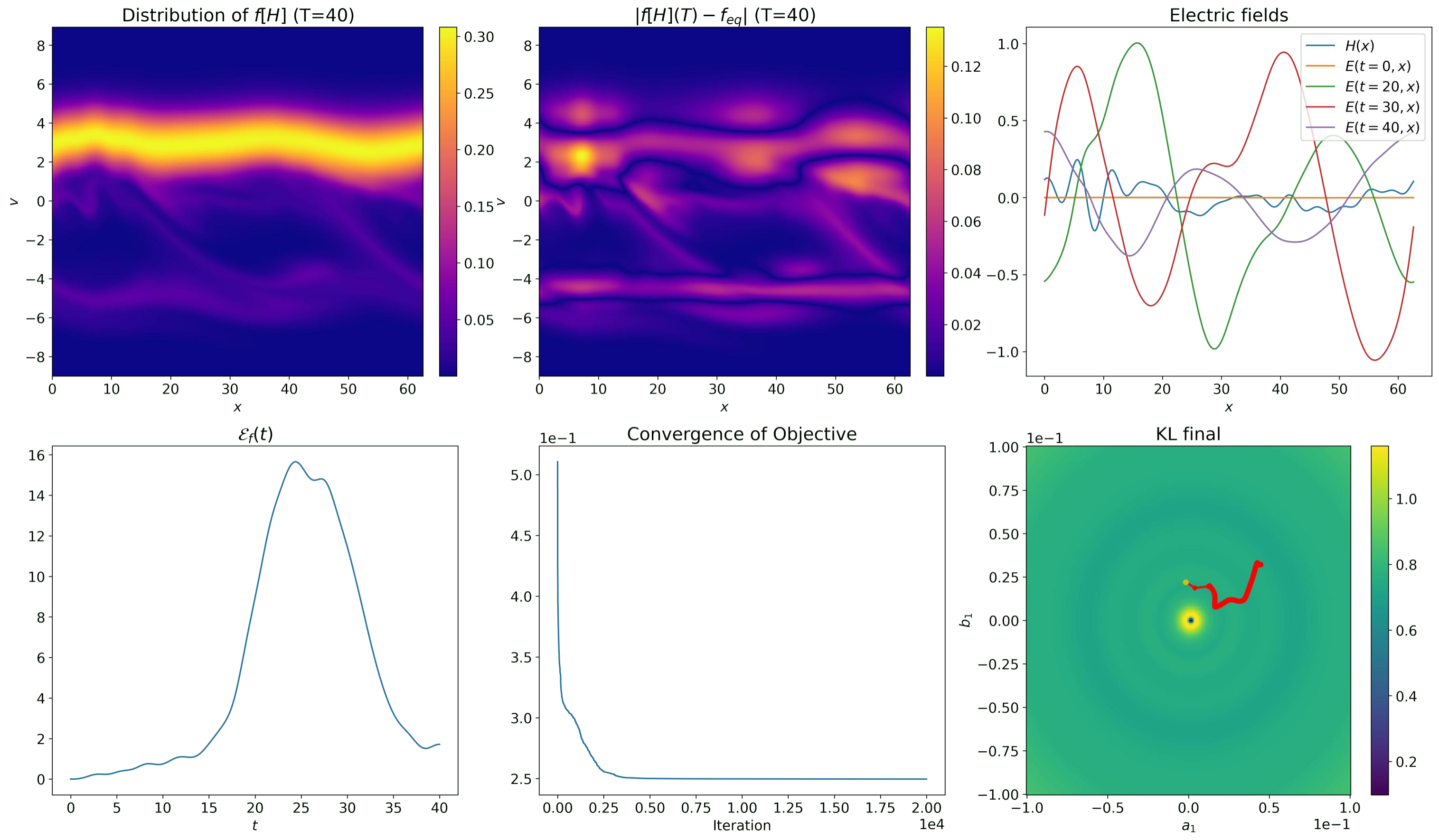}
    \caption{Simulation of~\eqref{eq:vlasov-poisoon_system_ext_1d} with over-parametrized $H$ obtained from~\eqref{eq:optimization_pb_simple} using~\eqref{eq:KL_obj} with near initialization using GD with line-search. From left to right and top to bottom: $f[H](T=30,x,v)$, $|f[H](T,x,v)-f_{\text{eq}}(v)|$, $H$ and $E_{f[H]}(t,x)$, $\mathcal{E}_{f[H]}(t)$, convergence of objective and, trajectory over the landscape of the objective (yellow dot is initial guess).}
    \label{fig:BoT_KL_GDL_near_over}
\end{figure}

\begin{figure}[H]
    \centering
    \includegraphics[width=0.85\linewidth]{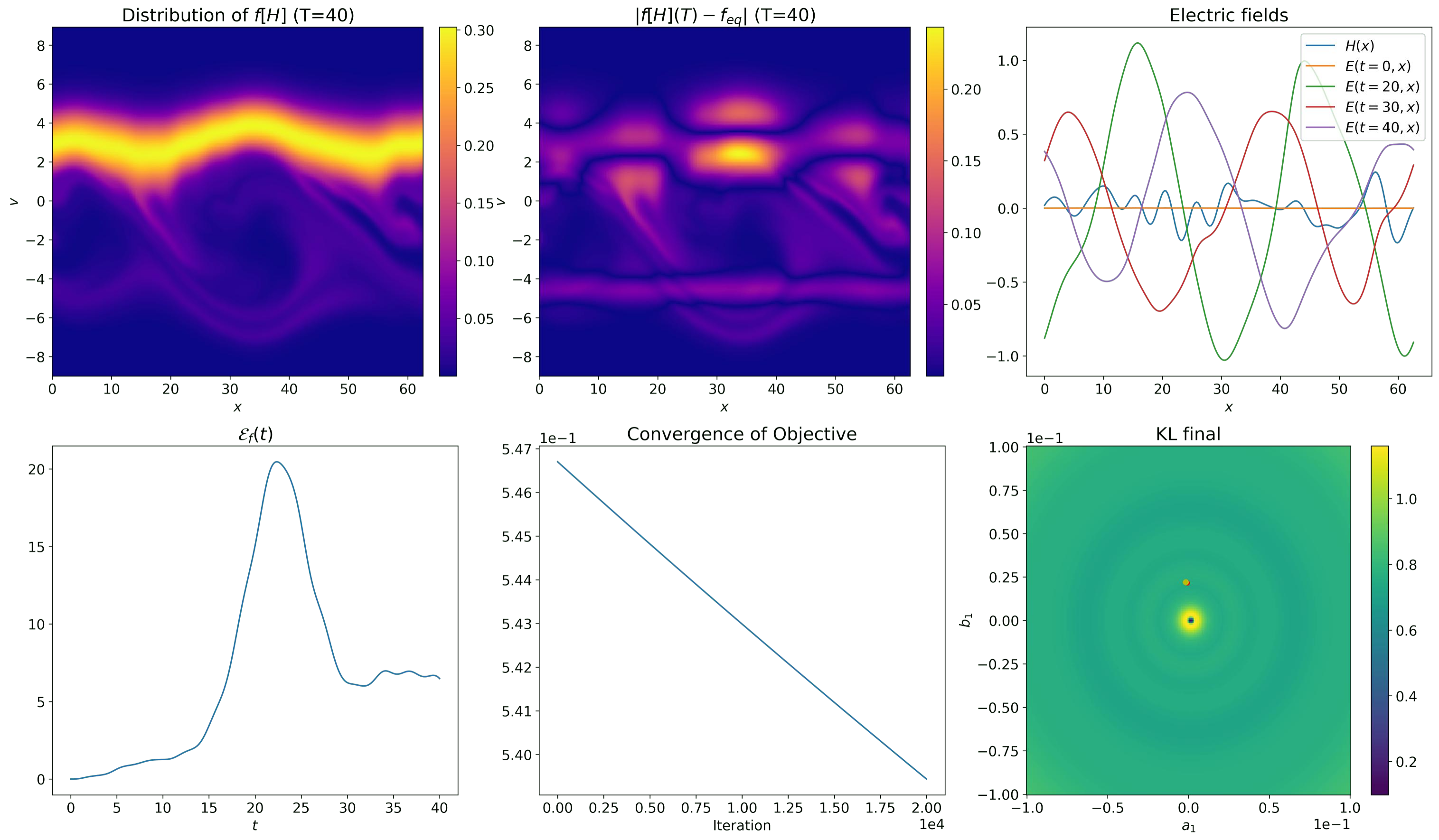}
    \caption{Simulation of~\eqref{eq:vlasov-poisoon_system_ext_1d} with over-parametrized $H$ obtained from~\eqref{eq:optimization_pb_simple} using~\eqref{eq:KL_obj} with near initialization using GD with constant stepsize. From left to right and top to bottom: $f[H](T=30,x,v)$, $|f[H](T,x,v)-f_{\text{eq}}(v)|$, $H$ and $E_{f[H]}(t,x)$, $\mathcal{E}_{f[H]}(t)$, convergence of objective and, trajectory over the landscape of the objective (yellow dot is initial guess).}
    \label{fig:BoT_KL_GD_near_over}
\end{figure}

\begin{figure}[H]
    \centering
    \includegraphics[width=0.85\linewidth]{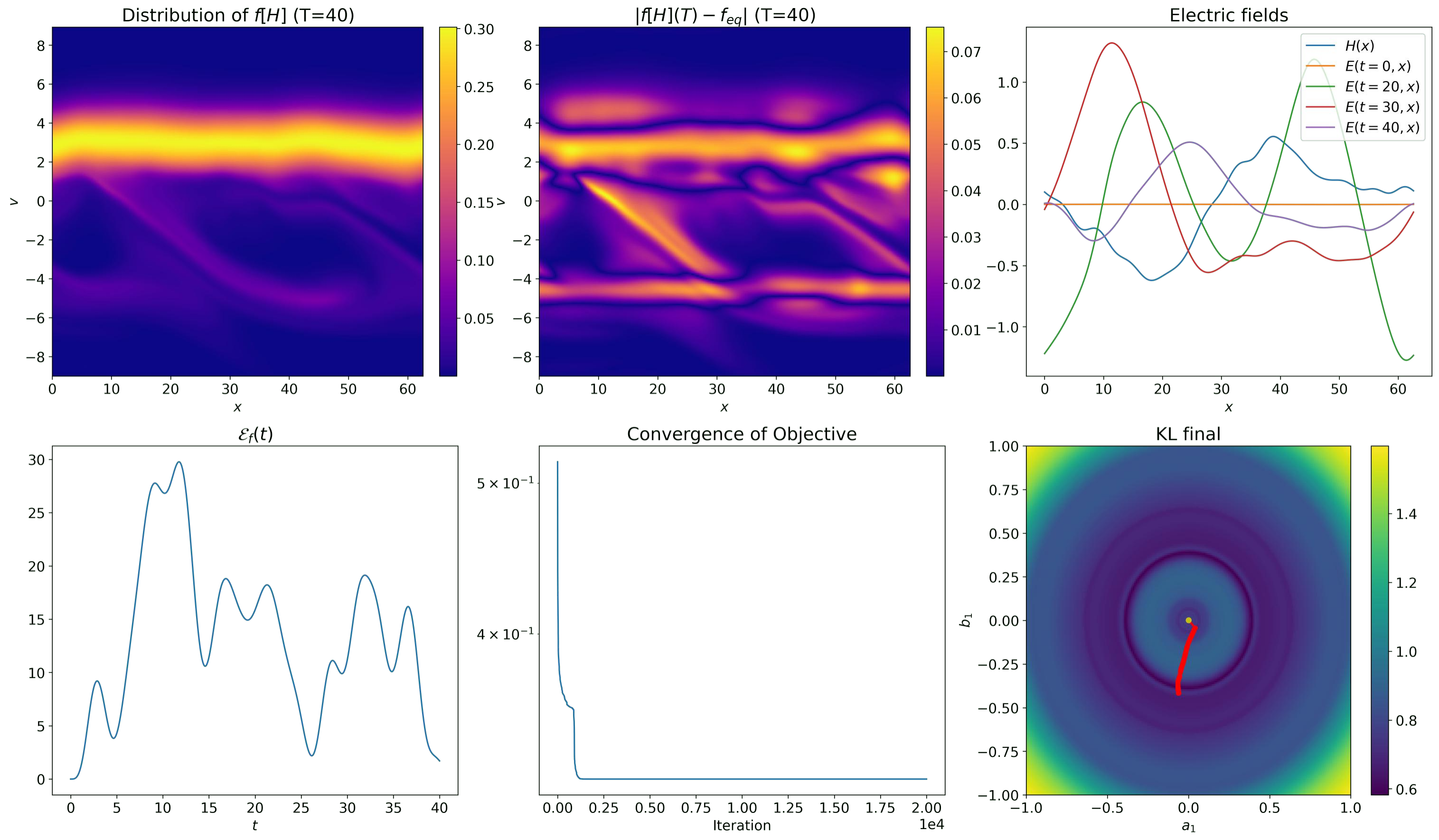}
    \caption{Simulation of~\eqref{eq:vlasov-poisoon_system_ext_1d} with over-parametrized $H$ obtained from~\eqref{eq:optimization_pb_simple} using~\eqref{eq:KL_obj} with local initialization using GD with line-search. From left to right and top to bottom: $f[H](T=30,x,v)$, $|f[H](T,x,v)-f_{\text{eq}}(v)|$, $H$ and $E_{f[H]}(t,x)$, $\mathcal{E}_{f[H]}(t)$, convergence of objective and, trajectory over the landscape of the objective (yellow dot is initial guess).}
    \label{fig:BoT_KL_GDL_local_over}
\end{figure}

\begin{figure}[H]
    \centering
    \includegraphics[width=0.85\linewidth]{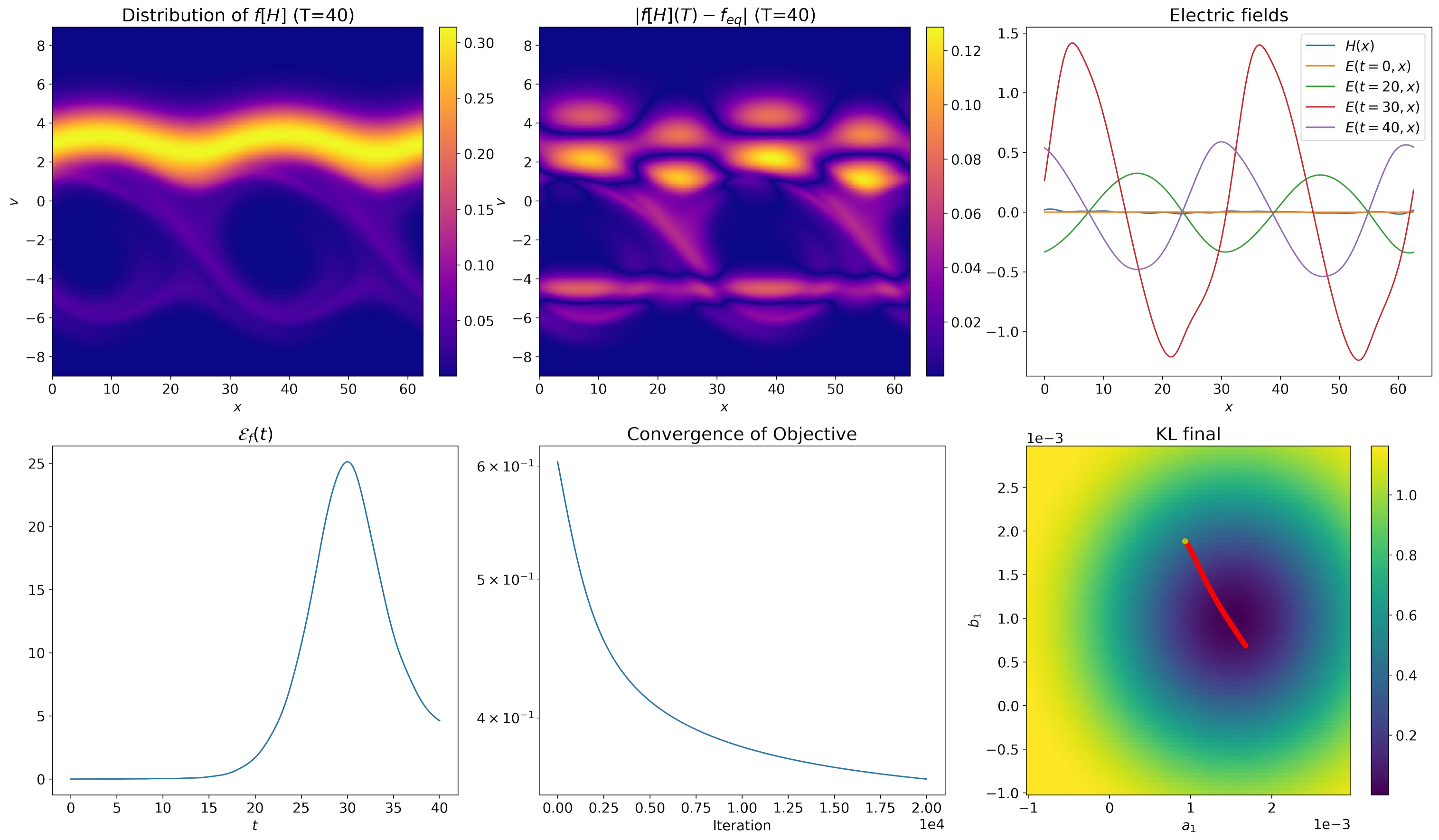}
    \caption{Simulation of~\eqref{eq:vlasov-poisoon_system_ext_1d} with over-parametrized $H$ obtained from~\eqref{eq:optimization_pb_simple} using~\eqref{eq:KL_obj} with local initialization using GD with constant stepsize. From left to right and top to bottom: $f[H](T=30,x,v)$, $|f[H](T,x,v)-f_{\text{eq}}(v)|$, $H$ and $E_{f[H]}(t,x)$, $\mathcal{E}_{f[H]}(t)$, convergence of objective and, trajectory over the landscape of the objective (yellow dot is initial guess).}
    \label{fig:BoT_KL_GD_local_over}
\end{figure}

\begin{figure}[H]
    \centering
    \includegraphics[width=0.85\linewidth]{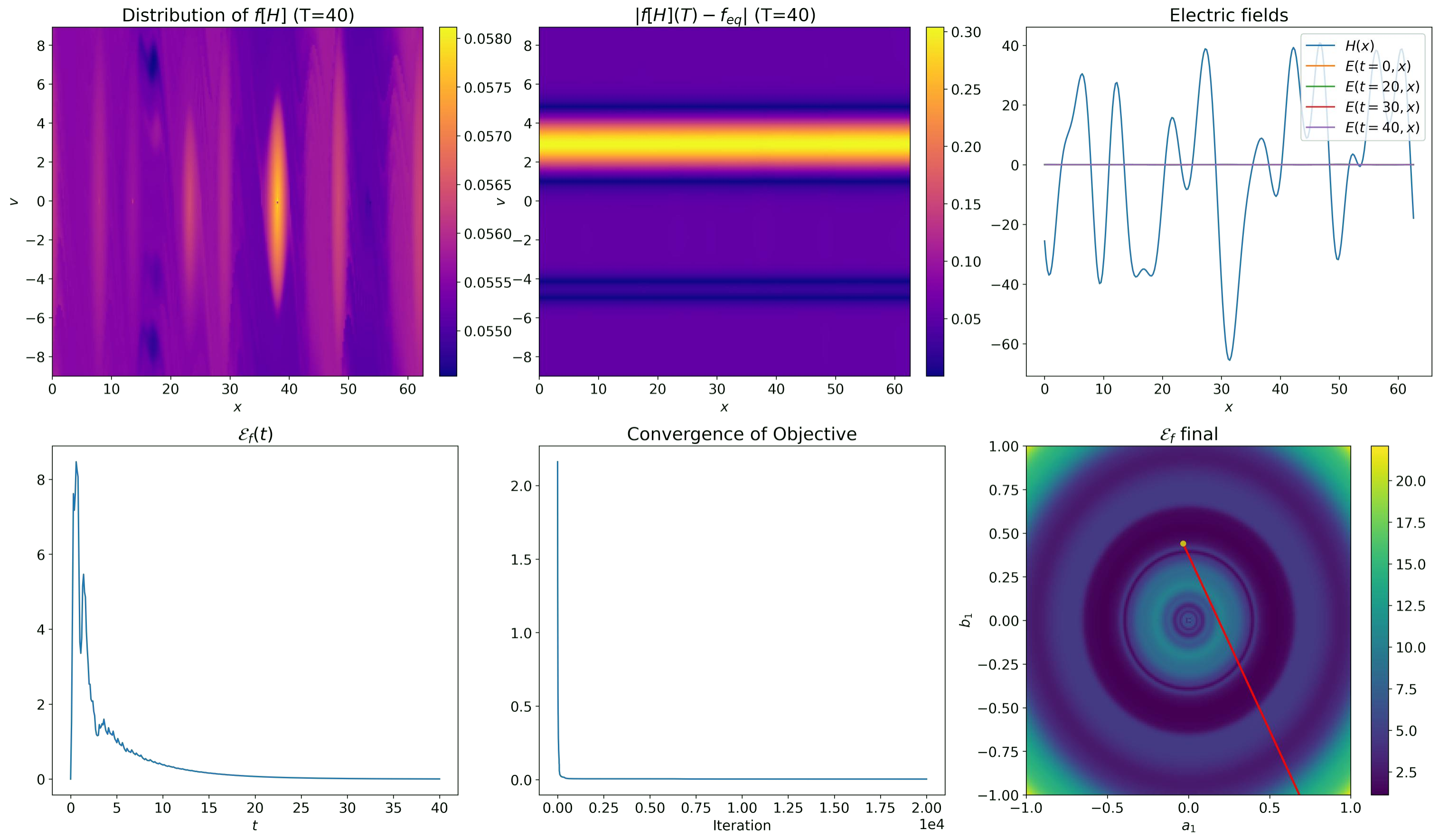}
    \caption{Simulation of~\eqref{eq:vlasov-poisoon_system_ext_1d} with under-parametrized $H$ obtained from~\eqref{eq:optimization_pb_simple} using~\eqref{eq:EE_obj} with far initialization using GD with line-search. From left to right and top to bottom: $f[H](T=30,x,v)$, $|f[H](T,x,v)-f_{\text{eq}}(v)|$, $H$ and $E_{f[H]}(T,x)$, $E_{f[H]}(t,x)$, $\mathcal{E}_{f[H]}(t)$ and, convergence of objective.}
    \label{fig:BoT_ee_lf_GDL_far_over}
\end{figure}

\begin{figure}[H]
    \centering
    \includegraphics[width=0.85\linewidth]{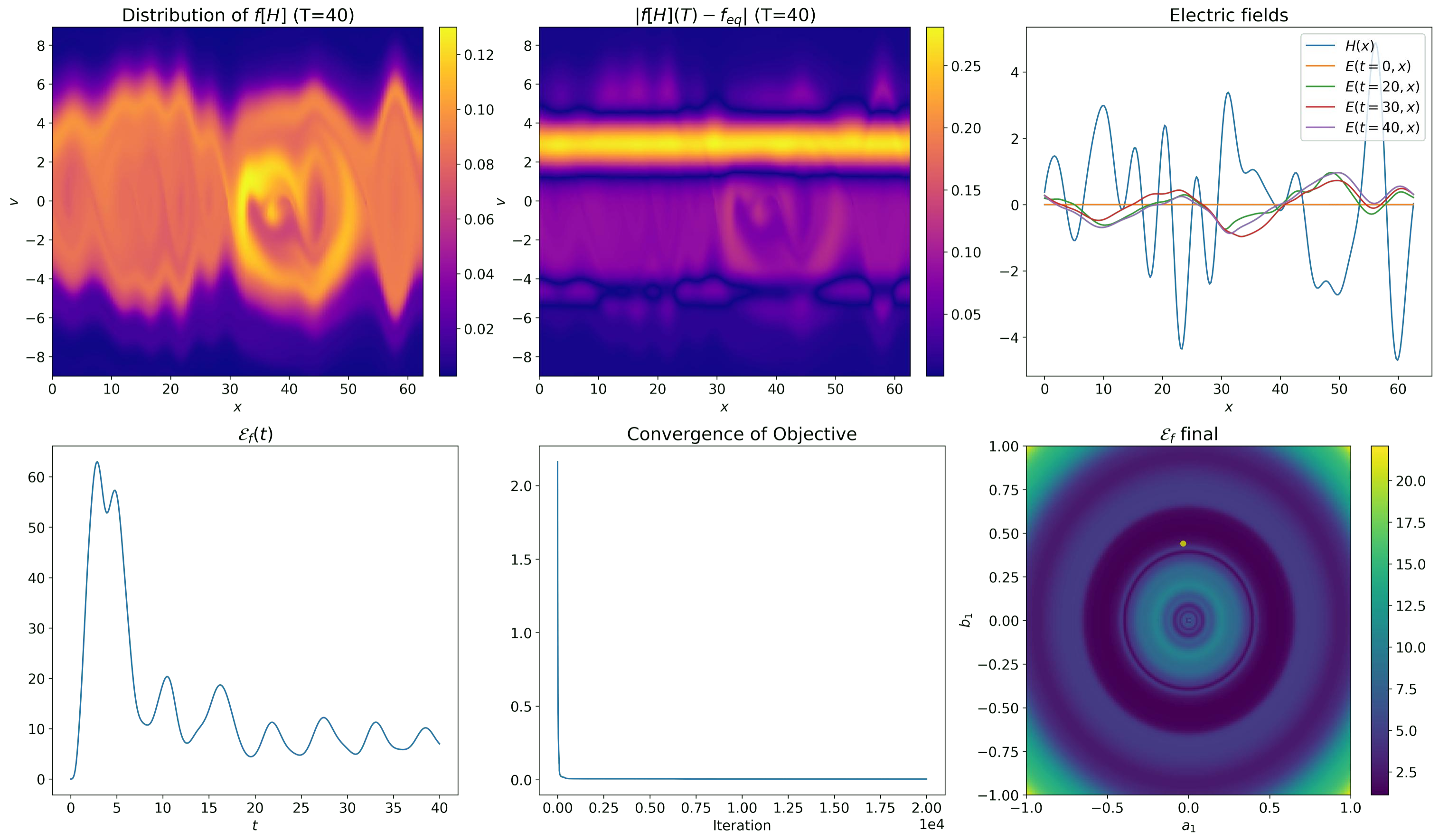}
    \caption{Simulation of~\eqref{eq:vlasov-poisoon_system_ext_1d} with over-parametrized $H$ obtained from~\eqref{eq:optimization_pb_simple} using~\eqref{eq:EE_obj} with far initialization using GD with constant stepsize. From left to right and top to bottom: $f[H](T=30,x,v)$, $|f[H](T,x,v)-f_{\text{eq}}(v)|$, $H$ and $E_{f[H]}(t,x)$, $\mathcal{E}_{f[H]}(t)$, convergence of objective and, trajectory over the landscape of the objective (yellow dot is initial guess).}
    \label{fig:BoT_ee_lf_GD_far_over}
\end{figure}

\begin{figure}[H]
    \centering
    \includegraphics[width=0.85\linewidth]{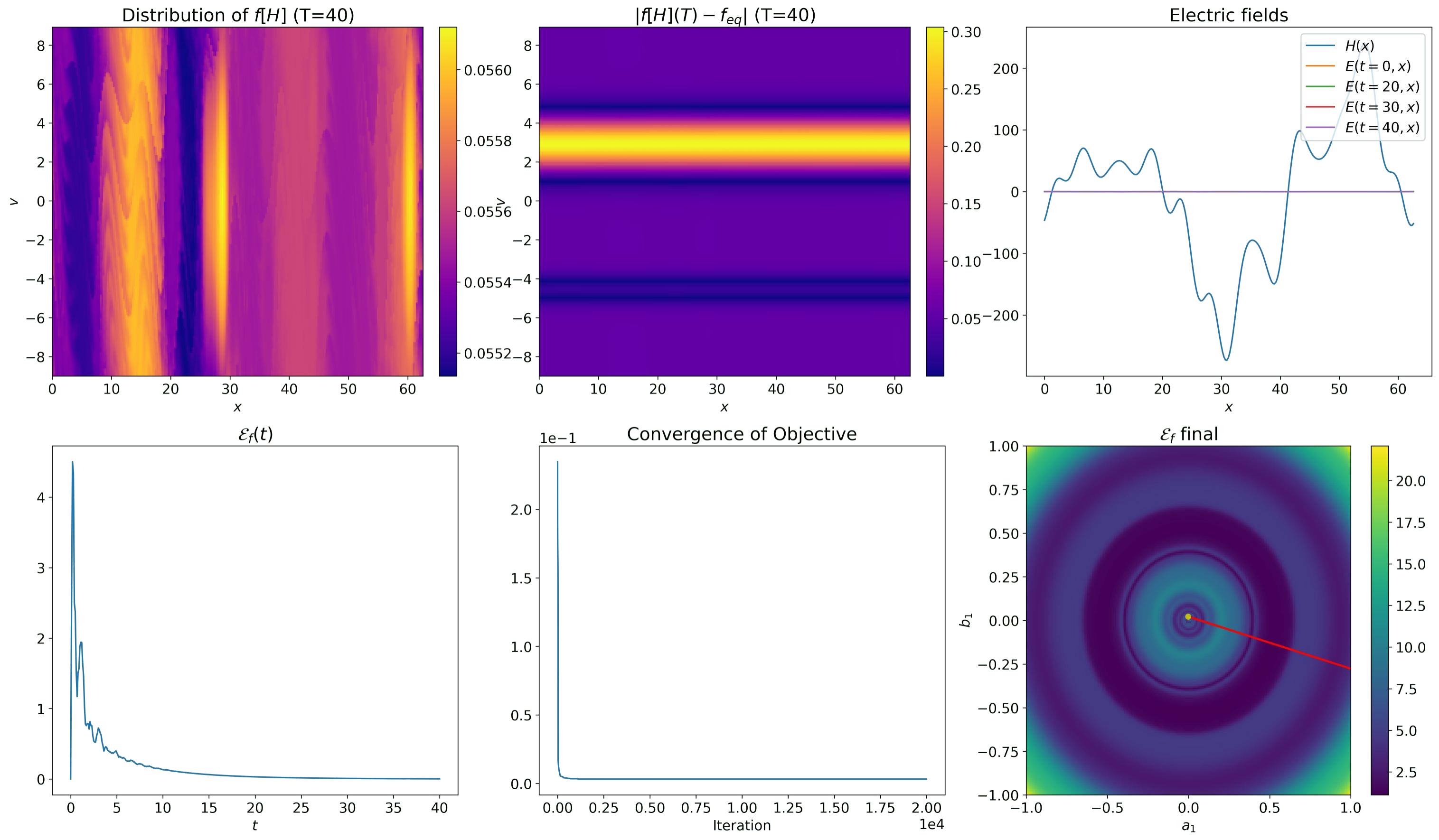}
    \caption{Simulation of~\eqref{eq:vlasov-poisoon_system_ext_1d} with over-parametrized $H$ obtained from~\eqref{eq:optimization_pb_simple} using~\eqref{eq:EE_obj} with near initialization using GD with line-search. From left to right and top to bottom: $f[H](T=30,x,v)$, $|f[H](T,x,v)-f_{\text{eq}}(v)|$, $H$ and $E_{f[H]}(t,x)$, $\mathcal{E}_{f[H]}(t)$, convergence of objective and, trajectory over the landscape of the objective (yellow dot is initial guess).}
    \label{fig:BoT_ee_lf_GDL_near_over}
\end{figure}

\begin{figure}[H]
    \centering
    \includegraphics[width=0.85\linewidth]{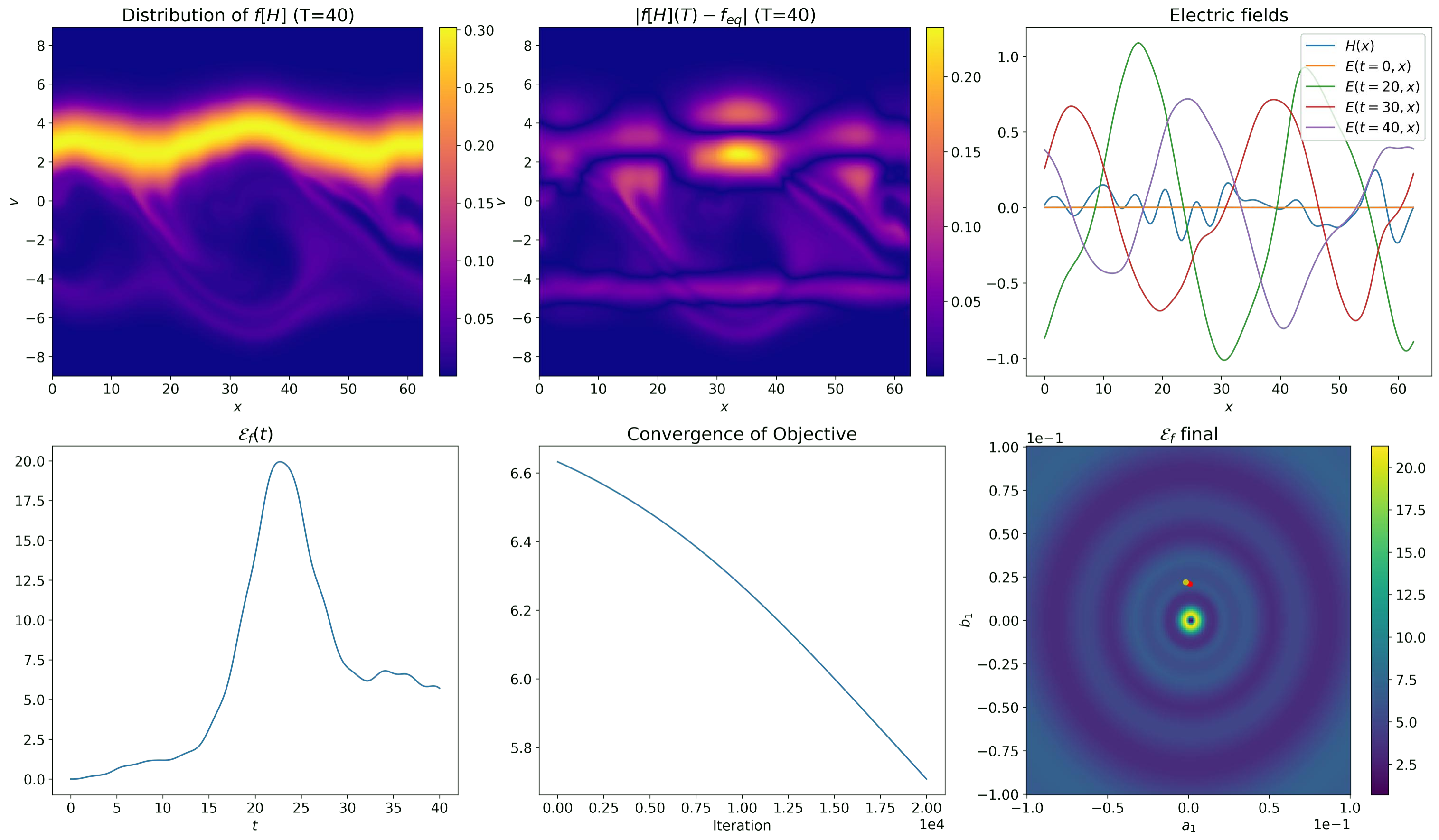}
    \caption{Simulation of~\eqref{eq:vlasov-poisoon_system_ext_1d} using~\eqref{eq:EE_obj} with over-parametrized $H$ obtained from~\eqref{eq:optimization_pb_simple} with near initialization using GD with constant stepsize. From left to right and top to bottom: $f[H](T=30,x,v)$, $|f[H](T,x,v)-f_{\text{eq}}(v)|$, $H$ and $E_{f[H]}(t,x)$, $\mathcal{E}_{f[H]}(t)$, convergence of objective and, trajectory over the landscape of the objective (yellow dot is initial guess).}
    \label{fig:BoT_ee_lf_GD_near_over}
\end{figure}

\begin{figure}[H]
    \centering
    \includegraphics[width=0.85\linewidth]{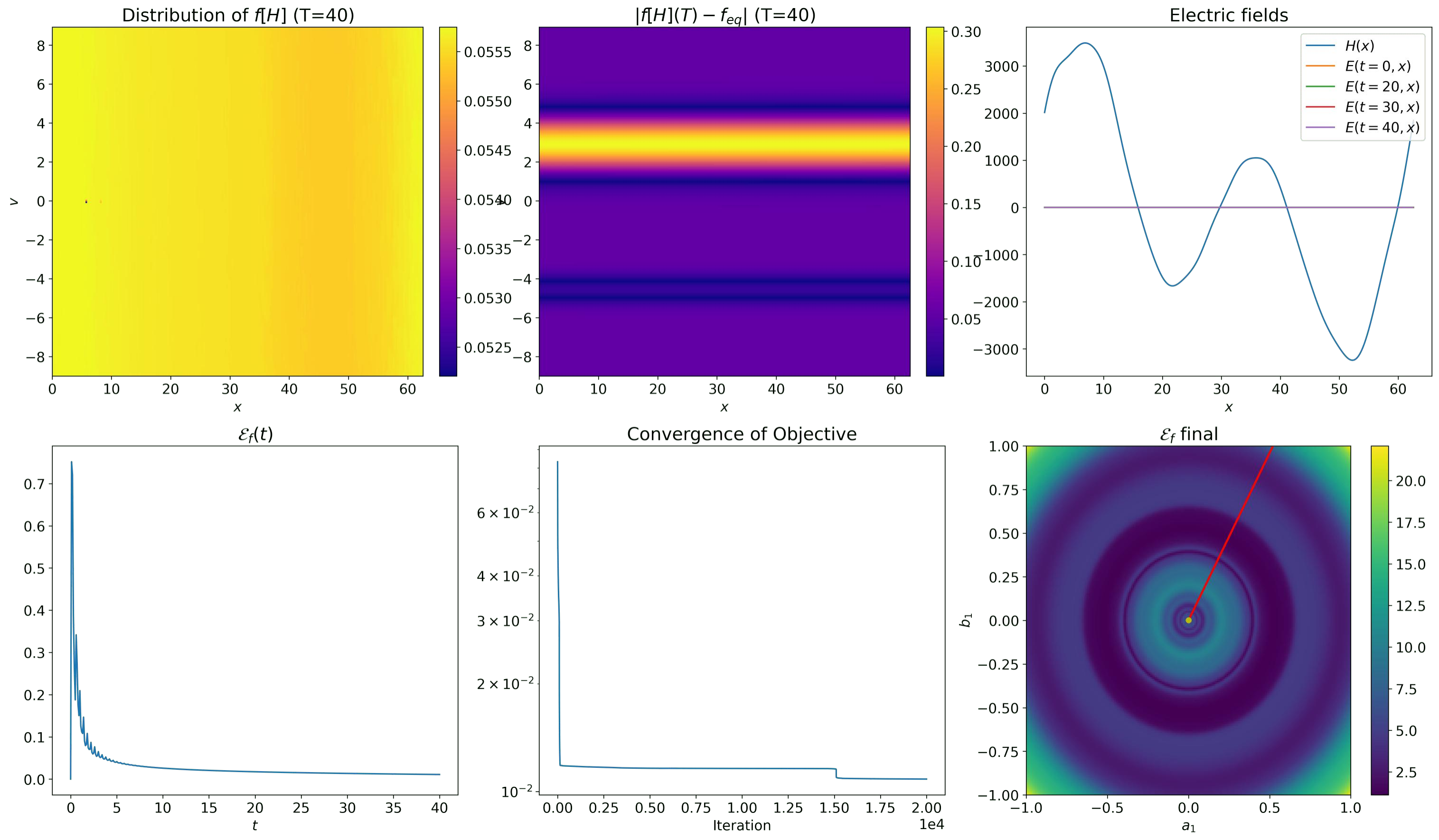}
    \caption{Simulation of~\eqref{eq:vlasov-poisoon_system_ext_1d} with over-parametrized $H$ obtained from~\eqref{eq:optimization_pb_simple} using~\eqref{eq:EE_obj} with local initialization using GD with line-search. From left to right and top to bottom: $f[H](T=30,x,v)$, $|f[H](T,x,v)-f_{\text{eq}}(v)|$, $H$ and $E_{f[H]}(t,x)$, $\mathcal{E}_{f[H]}(t)$, convergence of objective and, trajectory over the landscape of the objective (yellow dot is initial guess).}
    \label{fig:BoT_ee_lf_GDL_local_over}
\end{figure}

\begin{figure}[H]
    \centering
    \includegraphics[width=0.85\linewidth]{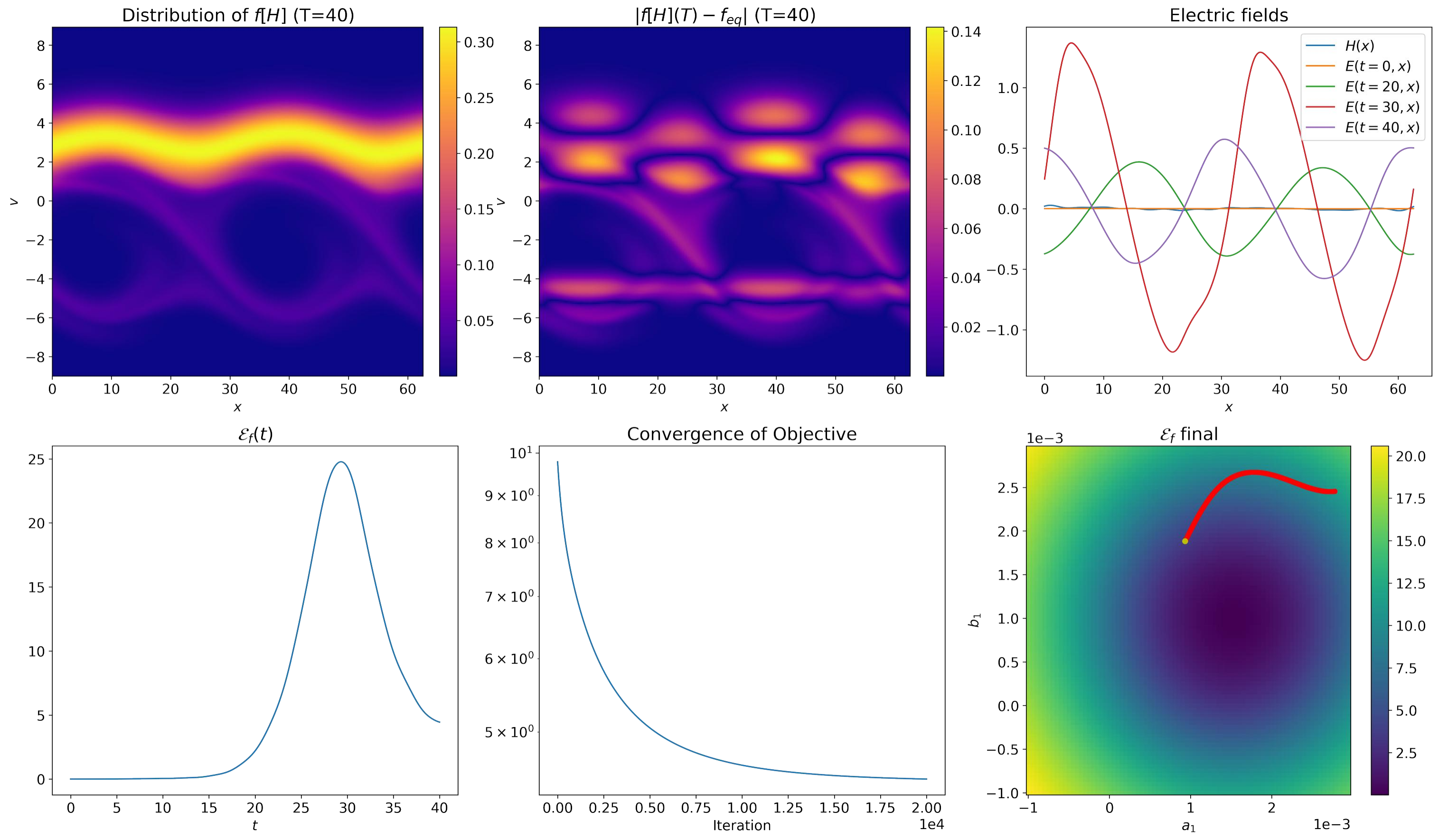}
    \caption{Simulation of~\eqref{eq:vlasov-poisoon_system_ext_1d} using~\eqref{eq:EE_obj} with over-parametrized $H$ obtained from~\eqref{eq:optimization_pb_simple} with local initialization using GD with constant stepsize. From left to right and top to bottom: $f[H](T=30,x,v)$, $|f[H](T,x,v)-f_{\text{eq}}(v)|$, $H$ and $E_{f[H]}(t,x)$, $\mathcal{E}_{f[H]}(t)$, convergence of objective and, trajectory over the landscape of the objective (yellow dot is initial guess).}
    \label{fig:BoT_ee_lf_GD_local_over}
\end{figure}

\begin{figure}[H]
    \centering
    \includegraphics[width=0.85\linewidth]{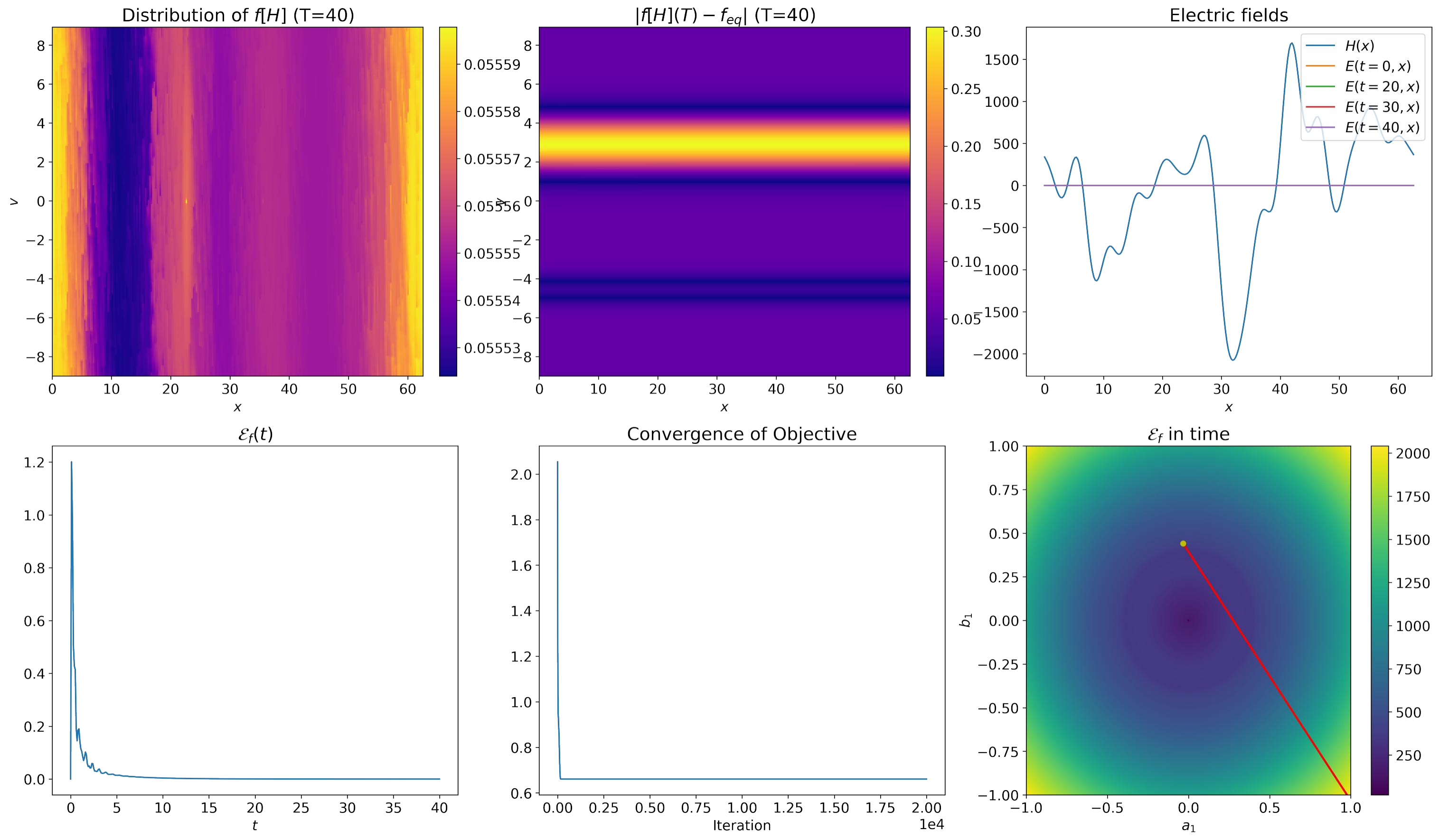}
    \caption{Simulation of~\eqref{eq:vlasov-poisoon_system_ext_1d} with under-parametrized $H$ obtained from~\eqref{eq:optimization_pb_simple} using~\eqref{eq:EET_obj} with far initialization using GD with line-search. From left to right and top to bottom: $f[H](T=30,x,v)$, $|f[H](T,x,v)-f_{\text{eq}}(v)|$, $H$ and $E_{f[H]}(T,x)$, $E_{f[H]}(t,x)$, $\mathcal{E}_{f[H]}(t)$ and, convergence of objective.}
    \label{fig:BoT_ee_GDL_far_over}
\end{figure}

\begin{figure}[H]
    \centering
    \includegraphics[width=0.85\linewidth]{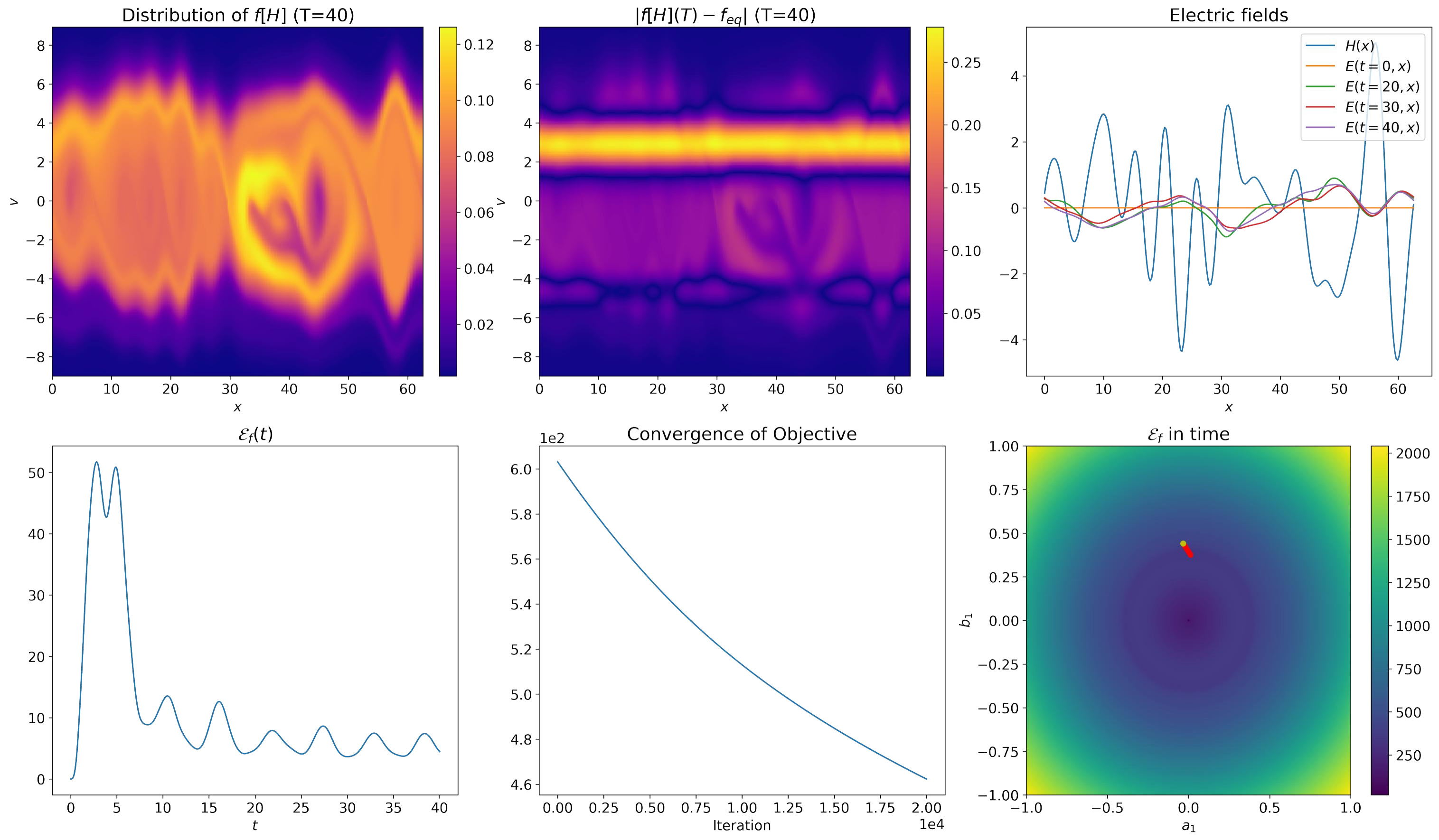}
    \caption{Simulation of~\eqref{eq:vlasov-poisoon_system_ext_1d} with over-parametrized $H$ obtained from~\eqref{eq:optimization_pb_simple} using~\eqref{eq:EET_obj} with far initialization using GD with constant stepsize. From left to right and top to bottom: $f[H](T=30,x,v)$, $|f[H](T,x,v)-f_{\text{eq}}(v)|$, $H$ and $E_{f[H]}(t,x)$, $\mathcal{E}_{f[H]}(t)$, convergence of objective and, trajectory over the landscape of the objective (yellow dot is initial guess).}
    \label{fig:BoT_ee_GD_far_over}
\end{figure}

\begin{figure}[H]
    \centering
    \includegraphics[width=0.85\linewidth]{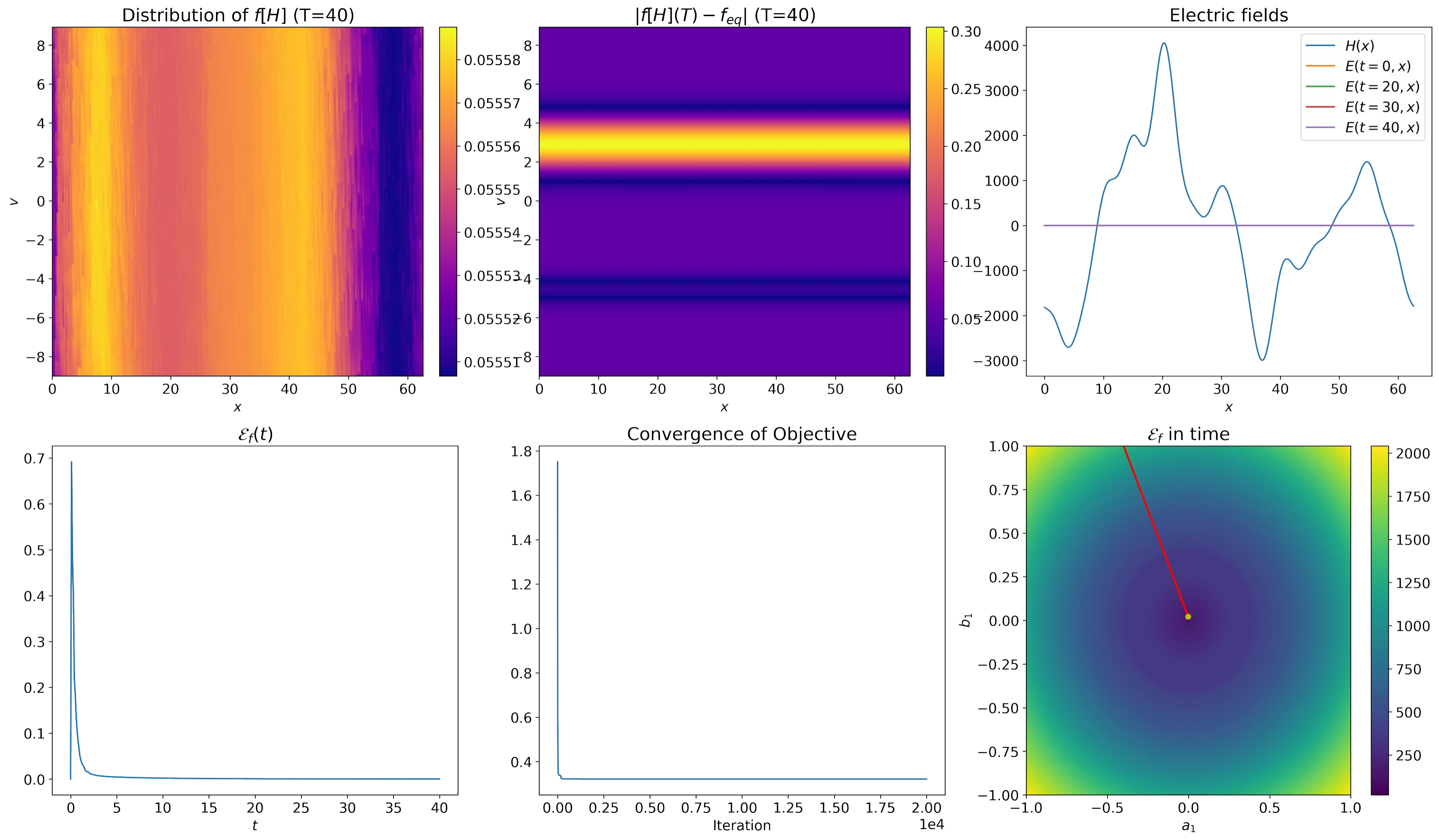}
    \caption{Simulation of~\eqref{eq:vlasov-poisoon_system_ext_1d} with over-parametrized $H$ obtained from~\eqref{eq:optimization_pb_simple} using~\eqref{eq:EET_obj} with near initialization using GD with line-search. From left to right and top to bottom: $f[H](T=30,x,v)$, $|f[H](T,x,v)-f_{\text{eq}}(v)|$, $H$ and $E_{f[H]}(t,x)$, $\mathcal{E}_{f[H]}(t)$, convergence of objective and, trajectory over the landscape of the objective (yellow dot is initial guess).}
    \label{fig:BoT_ee_GDL_near_over}
\end{figure}

\begin{figure}[H]
    \centering
    \includegraphics[width=0.85\linewidth]{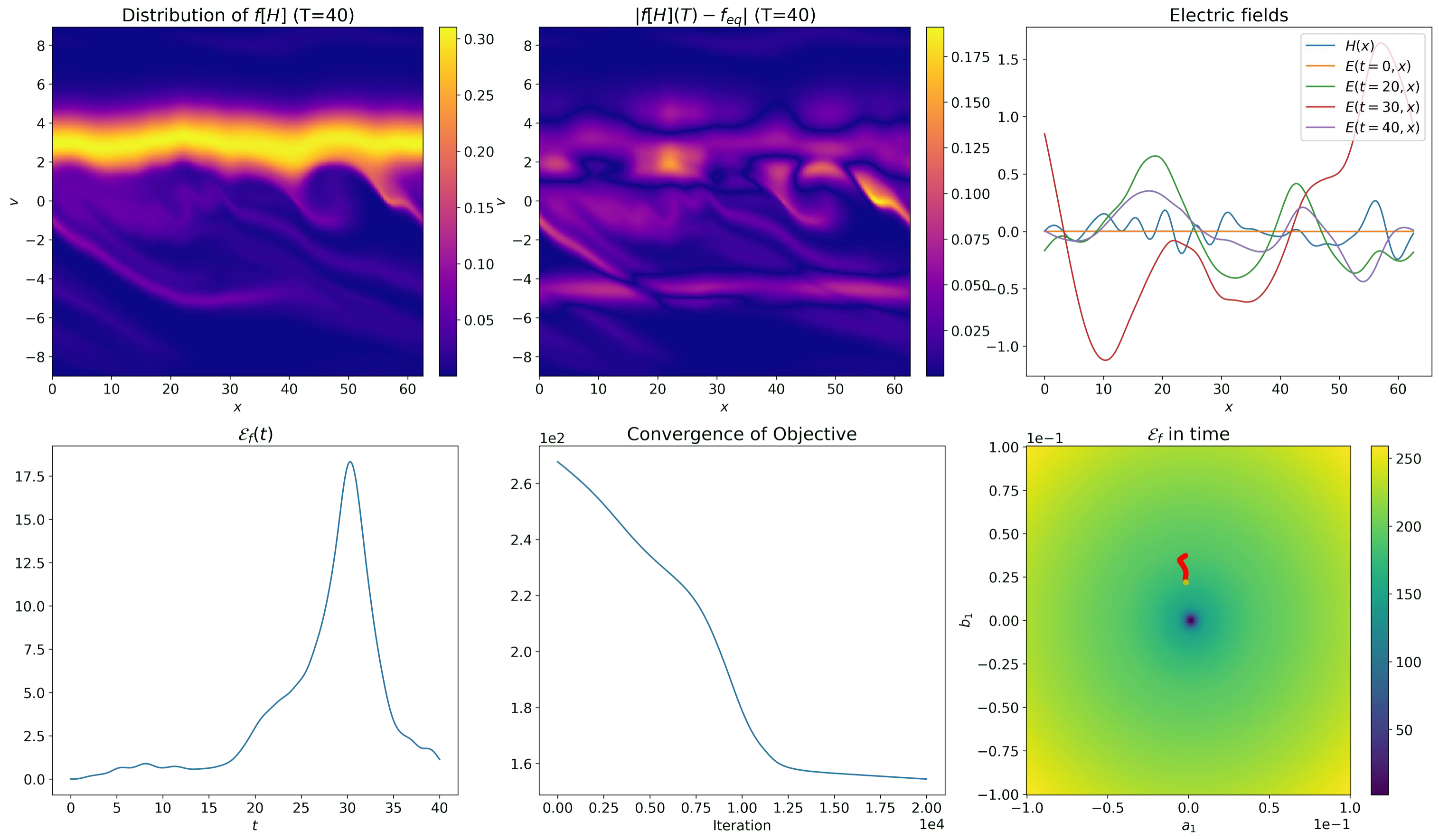}
    \caption{Simulation of~\eqref{eq:vlasov-poisoon_system_ext_1d} using~\eqref{eq:EET_obj} with over-parametrized $H$ obtained from~\eqref{eq:optimization_pb_simple} with near initialization using GD with constant stepsize. From left to right and top to bottom: $f[H](T=30,x,v)$, $|f[H](T,x,v)-f_{\text{eq}}(v)|$, $H$ and $E_{f[H]}(t,x)$, $\mathcal{E}_{f[H]}(t)$, convergence of objective and, trajectory over the landscape of the objective (yellow dot is initial guess).}
    \label{fig:BoT_ee_GD_near_over}
\end{figure}

\begin{figure}[H]
    \centering
    \includegraphics[width=0.85\linewidth]{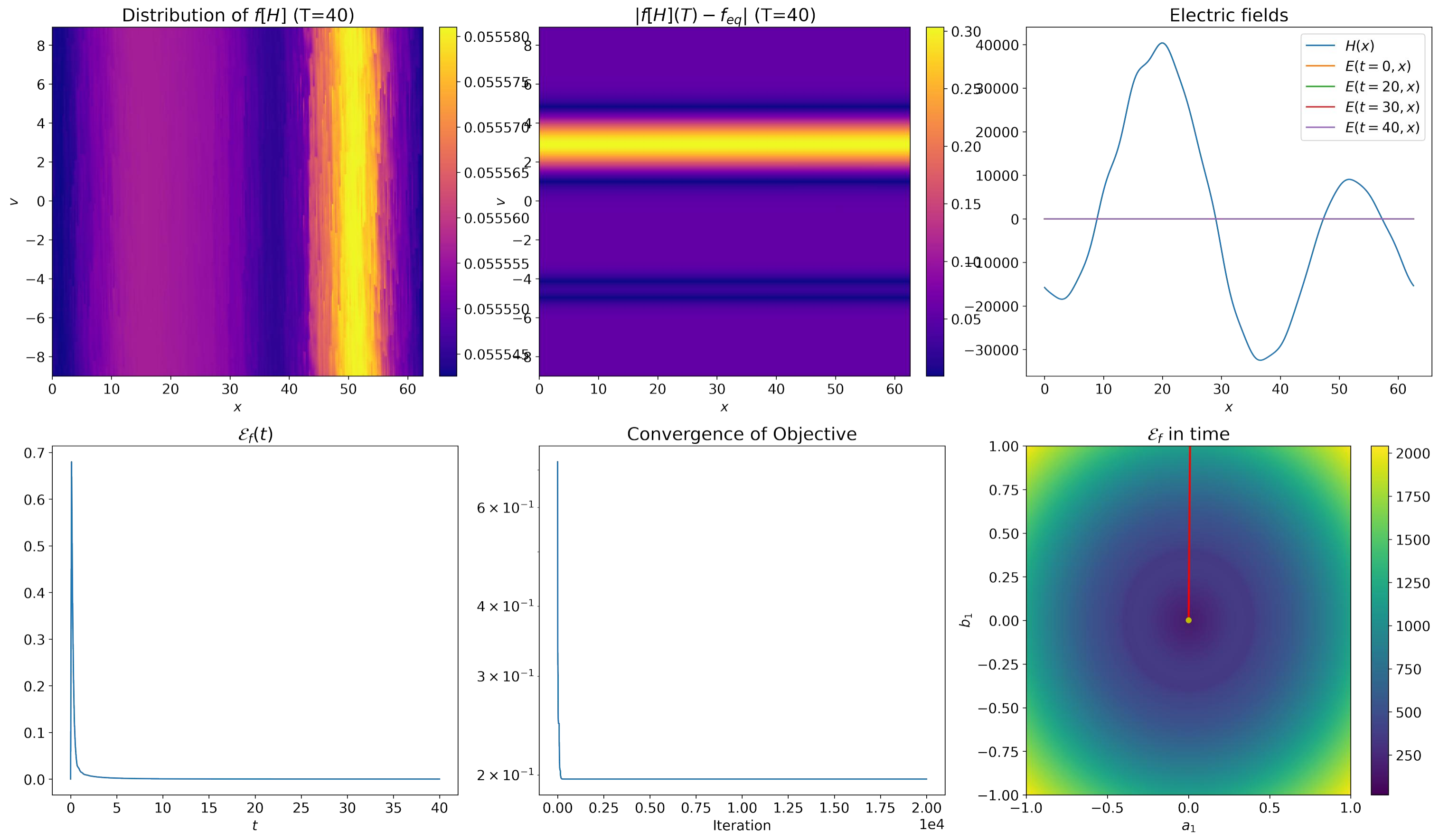}
    \caption{Simulation of~\eqref{eq:vlasov-poisoon_system_ext_1d} with over-parametrized $H$ obtained from~\eqref{eq:optimization_pb_simple} using~\eqref{eq:EET_obj} with local initialization using GD with line-search. From left to right and top to bottom: $f[H](T=30,x,v)$, $|f[H](T,x,v)-f_{\text{eq}}(v)|$, $H$ and $E_{f[H]}(t,x)$, $\mathcal{E}_{f[H]}(t)$, convergence of objective and, trajectory over the landscape of the objective (yellow dot is initial guess).}
    \label{fig:BoT_ee_GDL_local_over}
\end{figure}

\begin{figure}[H]
    \centering
    \includegraphics[width=0.85\linewidth]{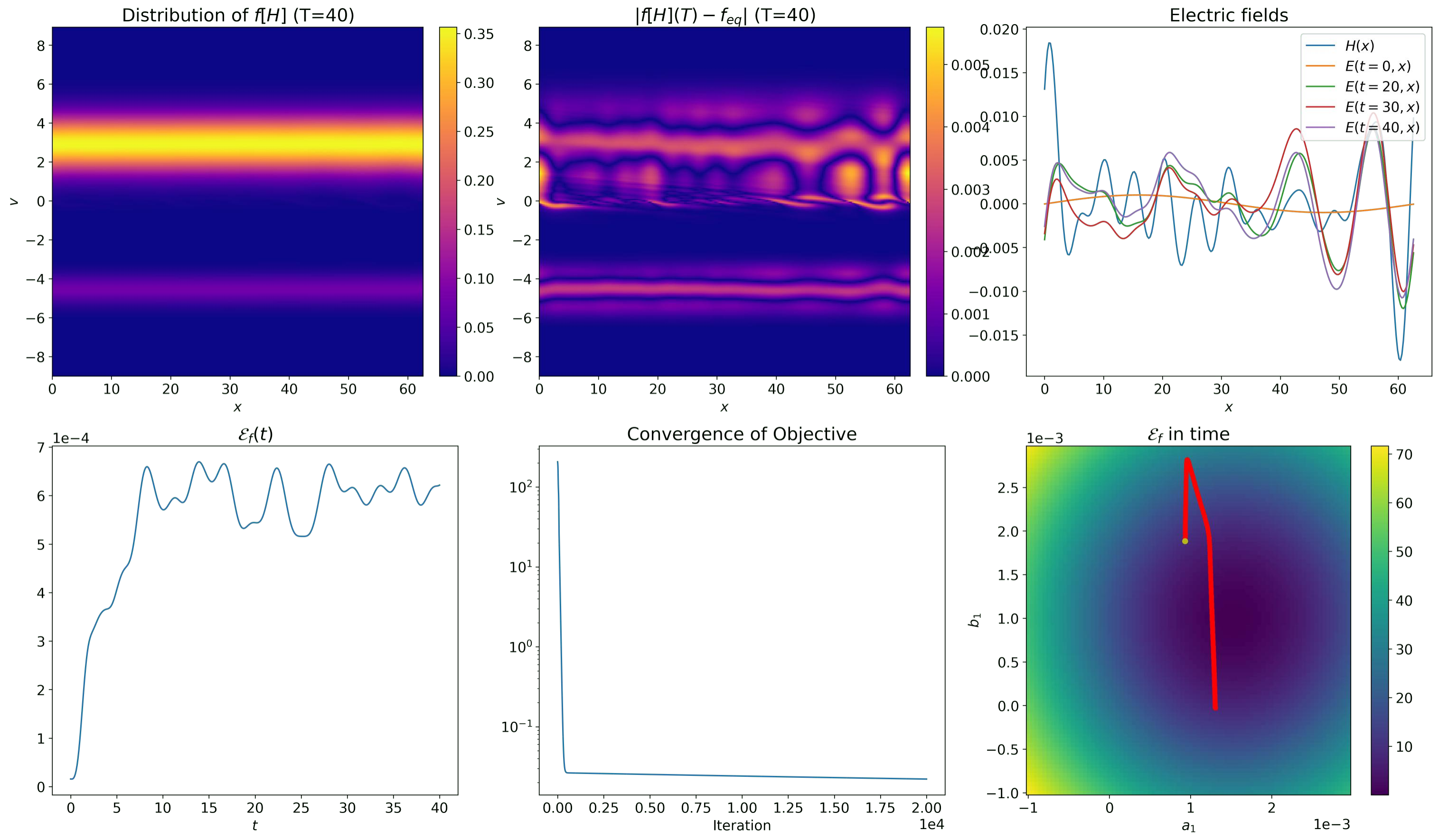}
    \caption{Simulation of~\eqref{eq:vlasov-poisoon_system_ext_1d} using~\eqref{eq:EET_obj} with over-parametrized $H$ obtained from~\eqref{eq:optimization_pb_simple} with local initialization using GD with constant stepsize. From left to right and top to bottom: $f[H](T=30,x,v)$, $|f[H](T,x,v)-f_{\text{eq}}(v)|$, $H$ and $E_{f[H]}(t,x)$, $\mathcal{E}_{f[H]}(t)$, convergence of objective and, trajectory over the landscape of the objective (yellow dot is initial guess).}
    \label{fig:BoT_ee_GD_local_over}
\end{figure}

\end{document}